\documentclass[a4paper,fleqn]{cas-sc}

\usepackage[numbers]{natbib}

\usepackage{amsthm}
\usepackage[shortlabels]{enumitem}
\usepackage{placeins}
\usepackage[protrusion=true,expansion=false]{microtype}

\ExplSyntaxOn
\keys_set:nn { stm / mktitle } { nologo }
\RenewDocumentCommand \printorcid { }
  {
    \seq_if_empty:NF \g_stm_orcid_seq
      {
        \group_begin:
          \tex_let:D \thefootnote \relax \footnotetext
            {
              \raggedright
              \textsc{orcid}(s):\c_space_token
              \seq_use:Nn \g_stm_orcid_seq { ;~ }
            }
        \group_end:
      }
  }
\ExplSyntaxOff

\newcommand{\E}{\mathbb{E}}
\newcommand{\R}{\mathbb{R}}
\newcommand{\N}{\mathbb{N}}
\newcommand{\D}{\mathcal{D}}
\newcommand{\F}{\mathcal{F}}
\newcommand{\calA}{\mathcal{A}}
\newcommand{\calB}{\mathcal{B}}
\newcommand{\calW}{\mathcal{W}}
\newcommand{\calE}{\mathcal{E}}
\newcommand{\calY}{\mathcal{Y}}
\newcommand{\calV}{\mathcal{V}}
\newcommand{\Tr}{\operatorname{Tr}}
\newcommand{\Ran}{\operatorname{Ran}}
\newcommand{\Dom}{\operatorname{Dom}}
\newcommand{\eps}{\varepsilon}
\newcommand{\ip}[2]{\langle #1, #2 \rangle}
\newcommand{\norm}[1]{\lVert #1 \rVert}
\newcommand{\abs}[1]{\lvert #1 \rvert}
\newcommand{\dual}[2]{{}_{V^*}\!\langle #1, #2 \rangle_V}
\newcommand{\dt}{\,\mathrm{d}t}
\newcommand{\ds}{\,\mathrm{d}s}
\newcommand{\dr}{\,\mathrm{d}r}
\newcommand{\dW}{\,\mathrm{d}W}
\newcommand{\HW}{H_W}
\newcommand{\nondegspace}{\mathcal{H}^{\sharp}_t}
\DeclareMathOperator{\spn}{span}

\makeatletter
\newcommand{\hyplabel}[2]{\protected@edef\@currentlabel{#2}\label{#1}}
\makeatother

\theoremstyle{plain}
\newtheorem{theorem}{Theorem}[section]
\newtheorem{lemma}[theorem]{Lemma}
\newtheorem{proposition}[theorem]{Proposition}
\newtheorem{corollary}[theorem]{Corollary}
\theoremstyle{definition}
\newtheorem{example}[theorem]{Example}
\newtheorem{definition}[theorem]{Definition}
\newtheorem{assumption}[theorem]{Assumption}
\newtheorem{remark}[theorem]{Remark}

\numberwithin{equation}{section}

\begin{document}
\let\WriteBookmarks\relax
\def\floatpagepagefraction{1}
\def\textpagefraction{.001}

\shorttitle{Logarithmic derivatives of variational SPDEs}
\shortauthors{E. Mirafzali, F. Proske and R. Marinescu}

\title[mode = title]{Logarithmic derivatives of variational and singular stochastic partial differential equations}

\author[1]{Ehsan Mirafzali}
\cormark[1]
\ead{smirafza@ucsc.edu}

\affiliation[1]{organization={Department of Computer Science, University of California Santa Cruz},
            addressline={1156 High Street},
            city={Santa Cruz},
            postcode={95064},
            state={CA},
            country={USA}}

\author[2]{Frank Proske}
\ead{proske@math.uio.no}

\affiliation[2]{organization={Department of Mathematics, University of Oslo},
            addressline={P.O. Box 1053 Blindern},
            city={Oslo},
            postcode={0316},
            country={Norway}}

\author[1]{Razvan Marinescu}
\ead{ramarine@ucsc.edu}

\cortext[1]{Corresponding author}

\begin{abstract}
For a stochastic partial differential equation posed on a Gelfand triple and satisfying the fully local monotone conditions of R\"ockner, Shang and Zhang, we compute the logarithmic derivative of the law of the solution at a fixed time along a prescribed direction of the state space. The formula is intrinsic, expressed through the Hilbert--Schmidt Malliavin derivative $\Phi_r = \D_r X(t)$ and the covariance $\gamma_t = \int_0^t\Phi_r\Phi_r^{*}\dr$ alone, so that neither the inversion of the first variation nor uniform Malliavin--Sobolev bounds on the Tikhonov family are needed. It follows from an integration-by-parts identity on a Hilbert space, a Moore--Penrose covering field on Wiener space and a trace formula for the Tikhonov limit, brought to the equation by the representation $\Phi_r = Y(t,r)\calB(r,X(r))$. The linearised theory is carried by $\Phi$ and by its own Malliavin derivative $\calW = \D\Phi$, whose regularity is measured in the time-integrated quadratic energy that the dissipative structure of the linearisation produces, so that no operator norm of the first variation between the spaces of the triple occurs. The stochastic $p$-Laplacian and the two-dimensional Navier--Stokes equation are treated in detail; for the latter the dissipative, forcing and exponential-moment hypotheses are established together under trace-class additive noise and a smallness threshold on the noise covariance. A scalar reduction carries the identity beyond the variational class, to the law of a pairing $\langle u(t), \varphi\rangle$ for the singular equations of Bruned, Chandra, Chevyrev and Hairer, on the same H\"older triangle of moment exponents.
\end{abstract}

\begin{keywords}
\raggedright
Logarithmic derivative \sep Bismut formula \sep Malliavin calculus \sep
quasilinear SPDE \sep Gelfand triple
\MSC[2020]{60H15, 60H07, 35R60}
\end{keywords}

\hfuzz=120pt
\maketitle
\hfuzz=0.1pt

\tableofcontents

\section{Introduction}
\label{sec:intro}

Let $X$ be a random variable on a probability space, taking values in a separable Hilbert space $H$, and let $\mu$ denote its law. The measure $\mu$ is said to be differentiable along a direction $h \in H$, in the sense of Fomin, if there exists a function $\beta_h \in L^2(\mu)$ such that
\begin{equation}\label{eq:intro-score}
  \int_H \partial_h \varphi(y)\,\mathrm{d}\mu(y) \;=\; -\int_H \varphi(y)\,\beta_h(y)\,\mathrm{d}\mu(y), \qquad \varphi \in C_b^1(H).
\end{equation}
The function $\beta_h$ is called the \emph{logarithmic derivative} (or \emph{Fomin derivative}) of $\mu$ along $h$. In finite dimensions, when $\mu$ has a Lebesgue density $p$, $\beta_h(y) = h \cdot \nabla \log p(y)$; in the generative modelling literature \cite{song2021scorebased, pmlr-v37-sohl-dickstein15}, $\beta_h$ is known as the \emph{score function}. We use the term logarithmic derivative throughout. In sufficiently regular finite-dimensional Markov diffusion settings, $\beta_h$ is the principal object entering time-reversal formulae; it also controls the gradient flow of relative entropy in the sense of Otto--Villani and is the gradient required when $\mu$ is sampled by Monte Carlo methods.

One point of vocabulary should be settled at once, since the two readings of $\beta_h$ are easily conflated. The function $\beta_h$ in~\eqref{eq:intro-score} is a property of the marginal law $\mu$, defined by the integration-by-parts identity alone. In finite-dimensional Markov-diffusion settings, under the regularity hypotheses of Haussmann--Pardoux \cite{haussmannpardoux1986timereversal} the time reversal of a diffusion with drift $b$ and diffusion coefficient $\sigma$ is again a diffusion, whose drift at the instant corresponding to forward time $t$ is $-b(t, \cdot) + p_t^{-1}\,\nabla\!\cdot\!\bigl(a(t, \cdot)\,p_t\bigr)$ with $a := \sigma\sigma^{*}$, the $p_t^{-1}$-terms being set to zero where $p_t$ vanishes; $\beta_h = \langle h, \nabla\log p_t\rangle$ is then the $h$-directional component of the score $\nabla\log p_t$ through which the density enters this reverse-time drift, and for the state-independent diffusion coefficients of the score-based diffusion models of \cite{song2021scorebased, pmlr-v37-sohl-dickstein15} the drift reduces to $-b + a\,\nabla\log p_t$. For the variational stochastic partial differential equations (SPDEs) of Section~\ref{sec:problem-main}, the solution is Markov in $H$ by construction; the identification with the reverse-time drift then requires density and regularity information for the law of $X(t)$ that is available only case-by-case in the SPDE-density literature, and typically only for finite-dimensional marginals or point evaluations (e.g.~\cite{marinellinualartquersardanyons2012}, establishing existence and regularity of the density of the real-valued evaluation $u(t,x)$ at fixed $(t,x)$ for semilinear dissipative heat equations with additive noise); this input is not part of our framework. For the singular SPDEs of Section~\ref{sec:applications} in the Bruned--Chandra--Chevyrev--Hairer (BCCH) scope, the renormalised state process need not be Markov in any naturally chosen state space; we make no dynamical reverse-time claim there, and our results read strictly as integration-by-parts identities for the marginal law of a scalar pairing $\langle u(t), \varphi\rangle$. Throughout this paper, ``logarithmic derivative'' refers to the Bogachev--Fomin notion of~\eqref{eq:intro-score}.

Three cases in which the logarithmic derivative is known set the pattern. When $\mu$ is the Gaussian measure $\mathcal{N}(0, Q)$ on $H$, the answer is classical. The measure is Fomin differentiable along its Cameron--Martin space $\Ran(Q^{1/2})$, and for the smaller class $h \in \Ran(Q)$ the logarithmic derivative takes the concrete $H$-inner-product form $\beta_h(y) = -\langle y, Q^{\dagger}h\rangle_H$, a general $h \in \Ran(Q^{1/2})$ giving instead the corresponding canonical Gaussian linear functional (Paley--Wiener map). When $\mu_t$ is the law at time $t$ of a finite-dimensional non-degenerate diffusion, the answer was put into closed form recently in \cite{mirafzali2025first, mirafzali2025malliavin} via the Malliavin calculus, in the shape $\beta_h(X_t) = -\E[\delta(v_h)\mid X_t]$, where $v_h$ is a so-called covering vector field on the Wiener space, and $\delta$ is the Skorokhod integral. When $X(t)$ is the solution at time $t$ of a linear SPDE with state-independent additive noise, the answer is again classical. The covariance $\gamma_t = \int_0^t S(t-r)Q\,S(t-r)^{*}\,\mathrm{d}r$ is deterministic, the covering field is deterministic and adapted, and the Skorokhod integral reduces to an It\^o integral \cite{mirafzali2025infinite}.

The Gaussian measure, the finite-dimensional diffusion and the linear equation with additive noise are three instances of one phenomenon. In each, an integration by parts is performed on the underlying probability space and transported to $H$ along the chosen direction $h$ by a Bismut--Fomin construction,
\begin{equation}\label{eq:bismut-template}
  \beta_h(X) \;=\; -\E[\delta_W(v_h)\mid X],
\end{equation}
where $W$ is the driving noise, $\delta_W$ is the Skorokhod integral on the Wiener space, and $v_h$ is a covering field, that is, a Cameron--Martin direction on the noise space whose pushforward through $X$ is the prescribed direction $h$ on $H$. Bismut's original formula \cite{bismut1984large} gave such a covering for the semigroup gradient $\nabla_x P_T f$, expressed via the inverse of the diffusion's first variation; Elworthy and Li \cite{ELWORTHY1994252} extended this to general diffusion semigroups, with a covering obtained by applying a right inverse of the diffusion coefficient to the uninverted derivative flow, and Bao, Wang and others \cite{bao2013bismut} carried it to SPDEs in various settings. To convert the gradient of the semigroup into the gradient of the density, the logarithmic derivative~\eqref{eq:bismut-template}, requires Nualart's substitution formula \cite{nualart2006malliavin}, carried out in finite dimensions in \cite{mirafzali2025first, mirafzali2025malliavin}.

The question addressed here is what this construction becomes for a variational SPDE. Consider
\begin{equation}\label{eq:main-SPDE}
  \mathrm{d}X(t) + \calA(t, X(t))\dt = \calB(t, X(t))\dW(t), \quad X(0) = x \in H, \quad t \in [0, T],
\end{equation}
on a Gelfand triple $V \hookrightarrow H \hookrightarrow V^*$, where $\calA$ is a possibly quasilinear operator and $\calB$ is a state-dependent diffusion driven by a $U$-cylindrical Wiener process $W$. Well-posedness is guaranteed by the fully local monotone conditions of R\"ockner--Shang--Zhang \cite{rockner2022wellposedness}; the Malliavin-flow regularity, covariance-range information, and Cameron--Martin compatibility needed for the Bismut formula are imposed as separate hypotheses (Assumptions~\ref{ass:diff1}--\ref{ass:nondeg}, \ref{ass:SC}, \ref{ass:SC-raised}), verifiable case-by-case.

Three difficulties separate this setting from the finite-dimensional one, and they pull in the same direction. The first is that the flow factorisation fails. In finite dimensions the covering field is built from $Y_r^{-1}$, whereas for an SPDE the membership $Y(t,r) \in L(H)$ is not part of the framework, so neither $Y(r,0)^{-1}$ nor a pathwise adjoint $Y(t,r)^{*}: H \to H$ is at one's disposal. The second is that the first variation does not smooth in operator norm. Already for the linear parabolic flow $Y(\tau,r) = e^{-A(\tau-r)}$ on the canonical scale $V = \Dom(A^{1/2})$ one has $\norm{Y(\tau,r)}_{L(H,V)} = \norm{A^{1/2}e^{-A(\tau-r)}}_{L(H)} = (2e(\tau-r))^{-1/2}$, so no power of that operator norm at or above two is integrable across the diagonal, and a hypothesis phrased in it is unavailable already in the case the nonlinear theory has to contain. The third is that the pseudoinverse is unbounded. The Malliavin covariance $\gamma_t$ is trace-class, its eigenvalues accumulate at zero, and Tikhonov regularisation yields no uniform $\mathbb{D}^{1,2}$-bound on the regularised covering field.

All three arise from leaving the objects the equation itself provides, and all three disappear on returning to them. We therefore work with
\begin{equation}\label{eq:intro-Phi}
  \Phi_r \;:=\; D_r X(t) \;\in\; L_2(U, H), \qquad r \in [0, t],
\end{equation}
and $\gamma_t = \int_0^t \Phi_r\Phi_r^{*}\,\mathrm{d}r$ as primary objects. The covering field
\begin{equation}\label{eq:intro-vh}
  v_h(r) \;=\; \Phi_r^{*}\,\gamma_t^{\dagger}\,h
\end{equation}
is linear in $z := \gamma_t^{\dagger}h$. This eliminates the first and third obstructions. No $Y(r,0)^{-1}$ appears, a basis-wise substitution (Theorem~\ref{thm:linear-sub-D12}) lowers the regularity demand from $\mathbb{D}^{1,4}_{\mathrm{loc}}$ to $\mathbb{D}^{1,2}$, and the Skorokhod integral $\delta_U(v_h^{(\eps)})$ converges in $L^2(\Omega)$ as $\eps \downarrow 0$ through a basis-free scalar Tikhonov trace limit under the Cameron--Martin compatibility condition. Absolute spectral summability provides one sufficient route for verifying the required scalar convergence and domination.

The second obstruction is removed by carrying the same principle one order further. The second-order object of the theory is not the derivative $\D_r Y(t,s)$ of the flow in an arbitrary direction of $H$ but the derivative of the fibre itself,
\begin{equation}\label{eq:intro-W}
  \calW(\tau;s,r) \;:=\; \D_r\Phi_{\tau,s} \;\in\; L_2\bigl(U, L_2(U,H)\bigr),
\end{equation}
which satisfies a linear equation whose forcing has both slots occupied by variational fibres. The estimates then consume the time-integrated quadratic fibre energy
\[
  \int_r^T\norm{\Phi_{\tau,r}}^2_{L_2(U,\mathbb{V}_\tau)}\,\mathrm{d}\tau,
\]
and this quantity the dissipative structure of the linearisation produces rather than assumes. For the linear part it is bounded by $\operatorname{Tr}Q/(2\nu)$ for trace-class noise, whereas its pointwise-in-time counterpart is unbounded at the diagonal. Every hypothesis on the linearisation is therefore stated in this energy, no operator norm of $Y$ is assumed or used in the estimates, and the smoothing clause that a factorisation through $L(H,V)$ would require does not appear.

The argument descends through four statements, each one asking for the next. An integration-by-parts identity on an abstract Hilbert space (Theorem~\ref{thm:abstract-bismut-fomin}) reduces the problem to producing a covering field; the Moore--Penrose construction (Theorem~\ref{thm:pseudoinverse-covering}) produces one, but only after a regularisation; the behaviour of that regularisation in the limit is settled by a trace formula (Theorem~\ref{thm:tikhonov-trace}); and the passage to the equation (Theorem~\ref{thm:main}) consists in identifying $\Phi_r = Y(t,r)\,\calB(r,X(r))$ and reading the abstract hypotheses off the variational structure. The first three statements are Malliavin analysis on an abstract Hilbert space and know nothing of the equation; only the fourth does.

The motivating equations are the stochastic $p$-Laplacian and the two-dimensional Navier--Stokes equation, and they exercise the two alternatives of the dissipative mechanism. For Navier--Stokes with additive trace-class noise the three analytic clauses (SC1), (SC1$'$) and (SC2) are established together. The strict coercivity gives~(SC1)(b) with $\alpha_1 = \nu$, and the forcing clause follows from it by the two-dimensional interpolation estimate for the nonlinearity, evaluated on the fibre energy that the coercivity itself produces (Proposition~\ref{prop:NS-forcing}), under a smallness threshold relating $\norm{Q}_{L(U)}$ to $\nu$, and with no spatial-regularity requirement on the noise beyond trace class. For the $p$-Laplacian, the drift differentiability hypotheses hold for $p > 3$, the domination estimate entering~(SC1)(a) is verified directly for $p \geq 4$, and~(SC1$'$) is the weighted product of the two fibre gradient densities against $\abs{\nabla X}^{p-4}$, which is the datum of the degenerate mechanism. The Galerkin-stability clause~(SC3) and the directional non-degeneracy of Assumption~\ref{ass:nondeg} are kept explicit throughout. For $2 < p < 4$ a regularised formulation is used instead, and the formula is then conditional on stability hypotheses for the regularised operator (Proposition~\ref{prop:singular-p-Lap}). At $2 < p \leq 3$ this is already forced by the failure of~(D3), while at $3 < p < 4$ it avoids the singular weight $\abs{\nabla X}^{p-4}$ appearing in the verification of~(SC1)(a). For Navier--Stokes with H\"ormander-bracket additive noise, finite-dimensional projected logarithmic derivatives follow unconditionally from \cite{mattinglypardoux2006malliavin, hairer2011hormander}; the full $H$-directional formula holds on the admissible-direction subspace $\nondegspace$ of Proposition~\ref{prop:nondeg-NS}, whose denseness in $H$ remains open.

For singular SPDEs in the BCCH scope \cite{brunedchandrachevyrevhairer2021}, a scalar reduction (Theorem~\ref{thm:score-gPAM}) yields an unconditional distributional derivative of $\mathrm{Law}(\langle u(t), \varphi\rangle)$ at the renormalised limit; a genuine logarithmic derivative requires the second-order Malliavin smoothness and negative-moment hypotheses of Assumption~\ref{ass:nondeg}\textup{(ii)} read on the one-dimensional state space, which remain open.

Section~\ref{sec:problem-main} states the assumptions and the four theorems. Sections~\ref{sec:malliavin} and~\ref{sec:variation} develop the Malliavin calculus and variation processes. Section~\ref{sec:proof-main} proves the main theorem. Section~\ref{sec:applications} verifies the hypotheses for the motivating equations and treats the singular extension. Throughout, $\Phi_r := D_r X(t) \in L_2(U, H)$ is the primary intrinsic object; the factorisation $\Phi_r = Y(t,r)\,\calB(r, X(r))$ is a representation theorem (Proposition~\ref{prop:D-Xt}).

\section{The equation, the hypotheses and the theorem}
\label{sec:problem-main}

Everything that follows rests on one object, namely the law of the solution at a fixed time, regarded as a measure on the state space, and its derivative along a prescribed direction. To say what that derivative is, one needs a class of equations wide enough to be worth the trouble, a notion of differentiability of the coefficients strong enough to linearise them, a structural hypothesis under which the linearised equations are themselves well posed, and a non-degeneracy condition ensuring that the direction in question is actually reachable by the noise. These four ingredients are introduced in turn, and the section closes with the statement they are designed to support.

\subsection{The variational setting}\label{subsec:setting}

We begin with the spaces, since everything else is phrased in them. We work on a Gelfand triple
\begin{equation}\label{eq:gelfand}
  V \hookrightarrow H \cong H^* \hookrightarrow V^*,
\end{equation}
where $(V, \norm{\cdot}_V)$ is a separable reflexive Banach space continuously and densely embedded in a separable Hilbert space $(H, \ip{\cdot}{\cdot}_H)$. The duality pairing $\dual{\cdot}{\cdot} := {}_{V^*}\langle\cdot,\cdot\rangle_V$ extends the $H$-inner product. Let $(U, \ip{\cdot}{\cdot}_U)$ be a further separable Hilbert space, and $L_2(U,H)$ the Hilbert--Schmidt operators $U \to H$. For any Banach space $E$, which in this paper is one of $V$, $V^{*}$, $H$ and $L_2(U,H)$, we write, for $T \in L(U, E)$ and for $(f_j)$ running over the orthonormal bases (ONB) of $U$,
\[
  \norm{T}_{L_2(U,E)} \;:=\; \sup_{(f_j)\ \mathrm{ONB\ of\ }U}\Bigl(\sum\nolimits_j \norm{T f_j}_E^2\Bigr)^{1/2},
\]
\[
  L_2(U,E) \;:=\; \{T \in L(U,E) : \norm{T}_{L_2(U,E)} < \infty\},
\]
a basis-free quantity that coincides with the usual Hilbert--Schmidt norm when $E$ is a Hilbert space. Whenever an object of this kind factors as $T = S \circ B$ with $S \in L(H, E)$ and $B \in L_2(U, H)$, every orthonormal basis gives $\sum_j \norm{S B f_j}_E^2 \leq \norm{S}_{L(H,E)}^2 \norm{B}_{L_2(U,H)}^2$, so that $\norm{S B}_{L_2(U,E)} \leq \norm{S}_{L(H,E)} \norm{B}_{L_2(U,H)} < \infty$, and the factorisation bound is a sufficient criterion for square summability. It is not the criterion used for the variational fibres of the sections below. There the square summability at the level of $E = V$, and of the tangent form domain that replaces it in the degenerate case, is obtained directly from the energy the linearisation produces, precisely because the first variation need not act boundedly from $H$ into $V$.

\subsection{The equation and the differentiability of its coefficients}

On a filtered probability space $(\Omega, \F, \{\F_t\}_{t \in [0,T]}, \mathbb{P})$ satisfying the usual conditions, let $W(t) = \sum_{k=1}^\infty f_k B_t^k$ be a $U$-cylindrical Wiener process, where $\{B^k\}_{k \geq 1}$ are independent standard Brownian motions and $\{f_k\}_{k \geq 1}$ is an ONB of $U$.

We consider the stochastic evolution equation \eqref{eq:main-SPDE}, restated here for convenience:
\begin{equation}\label{eq:SPDE}
  \mathrm{d}X(t) + \calA(t, X(t))\dt = \calB(t, X(t))\dW(t), \quad X(0) = x \in H, \quad t \in [0,T],
\end{equation}
where the operators $\calA$ and $\calB$ satisfy the following assumptions.

\begin{assumption}[Fully local monotone conditions]\label{ass:LR}
The mappings
\[
  \calA: [0,T] \times V \to V^*, \qquad \calB: [0,T] \times V \to L_2(U, H)
\]
are (jointly) Borel measurable, so that the composed processes $(\omega, t) \mapsto \calA(t, X(t, \omega))$ and $(\omega, t) \mapsto \calB(t, X(t, \omega))$ are progressively measurable whenever $X$ is an $H$-valued adapted process with continuous paths. The embedding $V \hookrightarrow H$ is compact. We further assume that $\calA$ and $\calB$ satisfy the following conditions for some constants $\alpha > 0$, $p \geq 2$, $C > 0$, and for measurable functions $\rho, \rho_1: [0,T] \times V \to [0, \infty)$ subject to the polynomial growth conditions~\textup{(i)}--\textup{(ii)} and the path-integrability property~\textup{(iii)} below. Conditions \mbox{(H1)--(H5)} below are the fully local monotone hypotheses of \cite{rockner2022wellposedness} (Part~I), stated in the sign convention of~\eqref{eq:SPDE}, with the $L^1$-in-time parts of the bounds there specialised to constants, and with the two local-monotonicity coefficients $\rho, \eta$ of \cite{rockner2022wellposedness} merged into the single symmetric coefficient $\rho$; the exponent $p$ here corresponds to the exponent $\alpha \in (1, \infty)$ of \cite{rockner2022wellposedness} (we restrict to $p \geq 2$ as required by the Malliavin analysis below), and our coercivity constant $\alpha$ to their $c$. The path-integrability property~\textup{(iii)} is not among the hypotheses of \cite{rockner2022wellposedness}; it is recorded here because these path functionals enter the linearised estimates of the later sections:
\begin{enumerate}[label=\textup{(\roman*)}, leftmargin=3em]
  \item \textup{(Polynomial growth of $\rho$.)} There exist $C_\rho > 0$ and $\beta_\rho \geq 0$ such that
  \begin{equation}\label{eq:LR-rho-growth}
    \rho(t, u) \leq C_\rho\bigl(1 + \norm{u}_V^p\bigr)\bigl(1 + \norm{u}_H^{\beta_\rho}\bigr) \qquad \text{for all } (t, u) \in [0,T] \times V.
  \end{equation}
  \item \textup{(Polynomial growth of $\rho_1$.)} There exist $C_{\rho_1} > 0$ and $\beta_{\rho_1} \geq 0$ such that
  \begin{equation}\label{eq:LR-rho1-growth}
    \rho_1(t, u) \leq C_{\rho_1}\bigl(1 + \norm{u}_V^p\bigr) \qquad \text{for all } (t, u) \in [0,T] \times V
  \end{equation}
  (the exponent $\beta_{\rho_1}$ enters through the growth condition~\textup{(H4)} below).
  \item \textup{(Path-integrability along variational solutions.)} For the variational solution $X$ of~\eqref{eq:SPDE} (obtained in Theorem~\ref{thm:LR-wellposedness} below), the path-functionals $\rho(\cdot, X(\cdot))$ and $\rho_1(\cdot, X(\cdot))$ are a.s.\ integrable on $[0,T]$, with moment bounds
  \begin{equation}\label{eq:LR-path-int}
    \E\Bigl[\int_0^T\bigl(\rho(\sigma, X(\sigma)) + \rho_1(\sigma, X(\sigma))\bigr)\,\mathrm{d}\sigma\Bigr] < \infty.
  \end{equation}
\end{enumerate}
Conditions~\eqref{eq:LR-rho-growth}--\eqref{eq:LR-rho1-growth} together with the moment estimate~\eqref{eq:LR-moments} of Theorem~\ref{thm:LR-wellposedness} imply~\eqref{eq:LR-path-int} automatically for deterministic initial data $x \in H$; for $\F_0$-measurable $X(0)$ (Remark~\ref{rem:random-initial}) the same holds once the initial moments $\E[\norm{X(0)}_H^{\max(2,\beta_\rho,\beta_{\rho_1})+p}] < \infty$ are in force. With these (restated) hypotheses, $\calA$ and $\calB$ satisfy the five fully local monotone conditions:

\begin{enumerate}[label=\textup{(H\arabic*)}, leftmargin=3em]
  \item (Hemicontinuity) For all $u, v, w \in V$ and $t \in [0,T]$, the map
  \[
    \lambda \mapsto \dual{\calA(t, u + \lambda v)}{w}
  \]
  is continuous on $\R$.
  
  \item (Fully local monotonicity) For all $u, v \in V$ and $t \in [0,T]$,
  \begin{multline}\label{eq:local-mono}
    -2\dual{\calA(t, u) - \calA(t, v)}{u - v} + \norm{\calB(t, u) - \calB(t, v)}_{L_2(U,H)}^2 \\
    \leq \bigl(\rho(t, v) + \rho(t, u)\bigr)\norm{u - v}_H^2.
  \end{multline}
  
  \item (Coercivity) For all $u \in V$ and $t \in [0,T]$,
  \begin{equation}\label{eq:coercivity}
    -2\dual{\calA(t, u)}{u} + \norm{\calB(t, u)}_{L_2(U,H)}^2 \leq C(1 + \norm{u}_H^2) - \alpha \norm{u}_V^p.
  \end{equation}
  
  \item (Fully local growth) For all $u \in V$ and $t \in [0,T]$,
  \begin{equation}\label{eq:growth}
    \norm{\calA(t, u)}_{V^*}^{p'} \leq \bigl(C + \rho_1(t, u)\bigr)\bigl(1 + \norm{u}_H^{\beta_{\rho_1}}\bigr),
  \end{equation}
  where $p' = p/(p-1)$; by \eqref{eq:LR-rho1-growth}, the prefactor has $V$-growth of order exactly $p$, so that \eqref{eq:growth} is the growth condition (H4) of \cite{rockner2022wellposedness}.
  
  \item (Diffusion continuity and growth) For all $t \in [0,T]$ and any sequence $\{u_n\}_{n \geq 1} \subset V$ and $u \in V$ with $\norm{u_n - u}_H \to 0$,
  \begin{equation}\label{eq:B-H-cont}
    \norm{\calB(t, u_n) - \calB(t, u)}_{L_2(U,H)} \to 0.
  \end{equation}
  Moreover, there exists $g \in L^1([0,T], \R_+)$ such that for all $u \in V$ and $t \in [0,T]$,
  \begin{equation}\label{eq:B-H-growth}
    \norm{\calB(t, u)}_{L_2(U,H)}^2 \leq g(t)\bigl(1 + \norm{u}_H^2\bigr).
  \end{equation}
\end{enumerate}
\end{assumption}

\begin{remark}[Deterministic-coefficient convention]\label{rem:deterministic-coefs}
Assumption \ref{ass:LR} takes $\calA$ and $\calB$ to depend on $\omega$ only through $X(t, \omega)$. This ensures the Malliavin chain rule $\D_r[\calA(t, X(t))] = \calA'_u(t, X(t))\,\D_r X(t)$ holds cleanly. The extension to random coefficients \cite{leon1998stochastic} requires tracking additional intrinsic-derivative terms in the linearised equations; all examples in this paper have deterministic coefficients.
\end{remark}

\begin{remark}[Random initial data]\label{rem:random-initial}
The results extend to $\F_0$-measurable initial data $X(0)$ lying in the space $L^2(\Omega, \F_0, \mathbb{P}; H)$. Well-posedness for random initial data follows from the deterministic case of Theorem~\ref{thm:LR-wellposedness} by pathwise uniqueness together with the continuous (hence measurable) dependence on the initial value established in \cite{rockner2022wellposedness}, and the moment estimate~\eqref{eq:LR-moments} then holds at exponent $q$ whenever $\E[\norm{X(0)}_H^q] < \infty$. For $\nu_0$-a.e.\ frozen $x$, Theorem~\ref{thm:main} applies to the conditional law $\mathrm{Law}(X(t)\mid X(0) = x)$; the marginal logarithmic derivative is recovered by the tower property.
\end{remark}

Under Assumption~\ref{ass:LR}, the existence and uniqueness of a variational solution to~\eqref{eq:SPDE} is guaranteed by the following result.

\begin{theorem}[R\"ockner--Shang--Zhang {\cite{rockner2022wellposedness}}; see also the erratum {\cite{rockner2025erratum}}]\label{thm:LR-wellposedness}
Under Assumption \ref{ass:LR}, for every initial value $x \in H$, equation \eqref{eq:SPDE} has a unique (probabilistically strong) variational solution
\[
  X \in L^p(\Omega; L^p(0, T; V)) \cap L^2(\Omega; C([0, T]; H)),
\]
and for every $q \geq 2$,
\begin{equation}\label{eq:LR-moments}
  \E\Bigl[\sup_{0 \leq t \leq T}\norm{X(t)}_H^{q}\Bigr] + \E\Bigl[\Bigl(\int_0^T \norm{X(t)}_V^{p}\,\mathrm{d}t\Bigr)^{\!q/2}\Bigr] < \infty.
\end{equation}
\end{theorem}

The precise notion of variational solution used throughout is as follows, cf.\ \cite{liu2015stochastic}.

\begin{definition}[Variational solution]\label{def:var-sol}
An $H$-valued continuous $\F_t$-adapted process $X = \{X(t)\}_{t \in [0,T]}$ is a
\emph{variational solution} (or \emph{analytically strong solution}) of~\eqref{eq:SPDE} if:
\begin{enumerate}[(i)]
  \item $X \in L^p(0,T; V) \cap C([0,T]; H)$ $\mathbb{P}$-a.s.;
  \item $\calA(\cdot, X(\cdot)) \in L^{p'}(0,T; V^*)$ $\mathbb{P}$-a.s., where $p' = p/(p-1)$;
  \item $\calB(\cdot, X(\cdot)) \in L^2(0,T; L_2(U, H))$ $\mathbb{P}$-a.s.;
  \item for all $t \in [0,T]$ and all $\varphi \in V$, the following identity holds $\mathbb{P}$-a.s.:
  \begin{multline}\label{eq:var-sol-identity}
    \ip{X(t)}{\varphi}_H + \int_0^t \dual{\calA(s, X(s))}{\varphi}\ds \\
    = \ip{x}{\varphi}_H + \int_0^t \ip{\calB(s, X(s))\dW(s)}{\varphi}_H,
  \end{multline}
  where the last integral is the It\^o stochastic integral taking values in $\R$, defined through the pairing
  \begin{align*}
    \ip{\calB(s,X(s))\dW(s)}{\varphi}_H \;=\; \sum_{k=1}^\infty \ip{\calB(s,X(s))f_k}{\varphi}_H\,\mathrm{d}B^k_s .
  \end{align*}
\end{enumerate}
\end{definition}

\begin{remark}[It\^o formula in the variational framework]\label{rem:ito-formula}
For $\psi \in C^2(H;\R)$ with $\nabla\psi: V \to V$ of at most affine growth and $\nabla^2\psi$ uniformly bounded (the \emph{admissibility class} of \cite{liu2015stochastic}), the variational It\^o formula gives
\begin{align}
  \psi(X(t)) &= \psi(x) - \int_0^t \dual{\calA(s,X(s))}{\nabla\psi(X(s))}\ds \nonumber\\
  &\quad + \int_0^t \ip{\nabla\psi(X(s))}{\calB(s,X(s))\dW(s)}_H \nonumber\\
  &\quad + \frac{1}{2}\int_0^t \Tr_U\bigl[\calB(s,X(s))^*\,\nabla^2\psi(X(s))\,\calB(s,X(s))\bigr]\ds. \label{eq:ito-var}
\end{align}
The quadratic $\psi(u) = \norm{u}_H^2$ is admissible; $\psi(u) = \norm{u}_H^{2q}$ for $q \geq 2$ is not, but is handled by applying the real-valued It\^o formula to the energy semimartingale $\xi(s) = \norm{X(s)}_H^2$:
\begin{equation}\label{eq:higher-moment-ito}
  \norm{X(t)}_H^{2q} = \xi(0)^q + \int_0^t q\,\xi(s)^{q-1}\,\mathrm{d}\xi(s) + \tfrac{1}{2}\int_0^t q(q-1)\,\xi(s)^{q-2}\,\mathrm{d}\langle\xi\rangle_s.
\end{equation}
Taking $\psi(u) = \norm{u}_H^2$ yields the \emph{energy identity}:
\begin{equation}\label{eq:energy-identity}
\begin{aligned}
  \norm{X(t)}_H^2 + 2\int_0^t \dual{\calA(s, X(s))}{X(s)}\ds
  &= \norm{x}_H^2 + \int_0^t \norm{\calB(s,X(s))}_{L_2(U,H)}^2\ds \\
  &\quad + 2\int_0^t \ip{X(s)}{\calB(s,X(s))\dW(s)}_H,
\end{aligned}
\end{equation}
and, via coercivity (H3), the a priori estimate:
\begin{equation}\label{eq:apriori}
  \E\Bigl[\sup_{0 \leq s \leq T}\norm{X(s)}_H^2\Bigr] + \alpha\,\E\Bigl[\int_0^T \norm{X(s)}_V^p\ds\Bigr] \leq C\bigl(1 + \norm{x}_H^2\bigr)\,e^{CT}.
\end{equation}
\end{remark}

Well-posedness alone says nothing about how the solution responds to a perturbation, and it is that response which the formula measures. We therefore ask that the coefficients be differentiable in the state variable, first to one order and then to two, with growth in the derivatives controlled by the same $V$-norm that governs the equation itself.

\begin{assumption}[First-order differentiability]\label{ass:diff1}
In addition to Assumption \ref{ass:LR}, we assume:
\begin{enumerate}[label=\textup{(D\arabic*)}, leftmargin=3em]
 \item The operator $\calA(t, \cdot): V \to V^*$ is Fr\'echet differentiable for each $t$, with derivative $\calA'_u(t, u) \in L(V, V^*)$, and $u \mapsto \calA'_u(t, u)$ is continuous from $V$ into $L(V, V^*)$ in the strong operator topology. The derivative satisfies
  \begin{equation}\label{eq:D1-growth}
    \norm{\calA'_u(t, u)}_{L(V, V^*)} \leq \bigl(C + \rho_2(t, u)\bigr)\bigl(1 + \norm{u}_V^{(p-2)^+}\bigr),
  \end{equation}
  where $(p-2)^+ = \max(p-2, 0)$ and $\rho_2: [0,T] \times V \to [0,\infty)$ is measurable, with the path-integrability along the variational solution
  \begin{equation}\label{eq:rho2-path}
    \int_0^T \bigl[(C + \rho_2(t,X(t)))(1 + \norm{X(t)}_V^{(p-2)^+})\bigr]^{p'}\dt \;<\; \infty \quad \mathbb{P}\text{-a.s.},
  \end{equation}
  and the same bound along each Galerkin approximant $X^N$.
  
  \item The operator $\calB(t, \cdot): V \to L_2(U, H)$ is Fr\'echet differentiable at each $u \in V$, with derivative
  \[
    \calB'_u(t, u) \in L(V;\,L_2(U, H)),
  \]
  and the assignment $u \mapsto \calB'_u(t, u)$ is continuous from $V$ into $L(V; L_2(U, H))$ in the strong operator topology. The derivative satisfies the fully local growth bound
  \begin{equation}\label{eq:D2-growth}
    \norm{\calB'_u(t, u)}_{L(V;\,L_2(U, H))} \leq C\bigl(1 + \norm{u}_V^\beta\bigr),
  \end{equation}
 for some $\beta \geq 0$. The exponent $\beta$ is subject to the integrability constraint $2\beta < p-2$ when $p > 2$, or $\beta = 0$ when $p = 2$. When additionally $\calB'_u(t,u)$ extends by continuity to an element of $L(H;\,L_2(U, H))$, the \emph{$H$-extension regime} characterising additive noise and $\calB$ depending on $u$ through an $H$-continuous functional, the stronger bound $\norm{\calB'_u(t,u)}_{L(H;\,L_2(U, H))} \leq C(1 + \norm{u}_V^\beta)$ holds with the relaxed constraint $2\beta < p$.
  
  \item[\textup{(D2$'$)}] \textup{(Noise-size growth at the moment-closure exponent.)} There exist a finite operating-point exponent $m_0$, a non-negative exponent $\kappa_{m_0} \geq 0$, and a constant $C_{m_0} > 0$, such that the diffusion coefficient satisfies
  \begin{multline}\label{eq:B-growth-m}
    \norm{\calB(t, u)}_{L_2(U, H)}^{m_0} \;\leq\; C_{m_0}\bigl(1 + \norm{u}_V^p\bigr)\bigl(1 + \norm{u}_H^{\kappa_{m_0}}\bigr) \\
    \qquad \text{uniformly in }(t, u) \in [0, T] \times V.
  \end{multline}
 The admissible range of $m_0$ is dictated by the single H\"older separation that the coefficient moment must feed, namely the conjugate pair $(m_0/2,\,m_0/(m_0-2))$ of Remark~\ref{rem:D2-strictness}:
  When the coefficient genuinely depends on the state, the separation is performed and $m_0 > 2$ is required strictly; at the quasilinear point $m_0 = p/\beta$ with $\beta > 0$ this is exactly the gap $2\beta < p$, while for $\beta = 0$ the requirement is that \eqref{eq:B-growth-m} hold at some $m_0 > 2$. Larger admissible values are always allowed, and every statement below is monotone in $m_0$. When the coefficient is deterministic, $\calB(t,u) = \calB(t)$, and in particular when the noise is additive, no separation is performed at all. The joint expectations factor exactly, as recorded in Remark~\ref{rem:D2-strictness}, the exponent $m_0/(m_0-2)$ never arises, and $m_0 = 2$ is admissible. The three standard regimes are thereby covered by the single bound \eqref{eq:B-growth-m}. They are the quasilinear one, at $m_0 = p/\beta$ and $\kappa_{m_0} = 0$, which recovers the Liu--R\"ockner condition; additive noise, at $m_0 = 2$ and $\kappa_{m_0} = 0$; and linear multiplicative noise $\calB(t,u) = L(t)u$, at any $m_0 > 2$ and $\kappa_{m_0} = m_0$.
  The product structure in~\eqref{eq:B-growth-m} ensures moment closure, the $V$-power staying at $p$, compatible with the energy a priori, while the $H$-power at $\kappa_{m_0}$ closes via the It\^o-energy identity.
\end{enumerate}
\end{assumption}

\begin{remark}\label{rem:why-Frechet}
Fr\'echet differentiability with strong-operator-continuous derivative is needed for progressive measurability of $\sigma \mapsto \calA'_u(\sigma, X(\sigma))$ and for the Malliavin chain rule
\[
  \D_r\bigl[\calA'_u(\tau, X(\tau))\bigr] = \calA''_{uu}(\tau, X(\tau))(\D_r X(\tau),\,\cdot)
\]
in Lemma~\ref{lemma:DPhi}. Both are verified for all examples in Section~\ref{subsec:verification}.
\end{remark}

The $V$-domain formulation of (D2) accommodates gradient-dependent diffusion $\calB(t, u)$ where $\calB'_u(t,u)$ acts on $V$ but does not extend to $H$ (for such $\calB$, the well-posedness clause~(H5) of Assumption~\ref{ass:LR} is replaced by the Part~II hypotheses of \cite{rockner2022wellposedness}; cf.\ Remark~\ref{rem:SC-raised-two-regimes}). The constraint $2\beta < p-2$ ensures the stochastic integral in~\eqref{eq:first-var} is well-defined via H\"older with $(p/(2\beta), p/(p-2\beta))$. Under the $H$-extension regime, the relaxed $2\beta < p$ applies. At $p = 2$, both reduce to $\beta = 0$.

\begin{remark}[The single separation and its strict gap]\label{rem:D2-strictness}
In both regimes, the working exponent of every joint-moment separation below is the operating point $m_0 > 2$ of~(D2$'$). The H\"older conjugate-exponent argument that separates the joint moment bound~\eqref{eq:fibre-energy-estimate} of Theorem~\ref{thm:Y-wellposed}\textup{(III)} into two finite factors, as it is in Proposition~\ref{prop:D-Xt} to establish $X(t) \in \mathbb{D}^{1,2}(H)$, is exactly
\[
  \E\bigl[\norm{\Theta_0}^2_{L_2}\exp\bigl(c_1\textstyle\int(\tilde\rho+\hat\rho)\bigr)\bigr] \;\leq\; \bigl(\E[\norm{\Theta_0}^{m_0}_{L_2}]\bigr)^{\!2/m_0}\bigl(\E[\exp(c_1\tfrac{m_0}{m_0-2}\textstyle\int(\tilde\rho+\hat\rho))]\bigr)^{\!(m_0-2)/m_0},
\]
in which the conjugate H\"older exponent $m_0/(m_0-2)$ is finite if and only if $m_0 > 2$ strictly. At the quasilinear operating point $m_0 = p/\beta$ ($\beta > 0$) this is precisely the strict gap $2\beta < p$, recovering the exponent pair $(p/(2\beta), p/(p-2\beta))$; the non-strict endpoint would render the exponent infinite and break the exponential-moment factor via~\eqref{eq:SC-exp-moment} of Assumption~\ref{ass:SC}. At $\beta = 0$ the separation is available at any $m_0 > 2$ admitted by~\eqref{eq:B-growth-m} (always, via the additive or $H$-growth clauses); for deterministic noise no separation is needed at all, since $\E[\norm{\calB(r)}^2\exp(c_1\int(\tilde\rho+\hat\rho))] = \norm{\calB(r)}^2\,\E[\exp(c_1\int(\tilde\rho+\hat\rho))]$ factors exactly. The first factor on the right, the moment $\E[\norm{\Theta_0}_{L_2}^{m_0}] = \E[\norm{\calB(r, X(r))}_{L_2(U,H)}^{m_0}]$, is exactly what the noise-size growth hypothesis~\eqref{eq:B-growth-m} of (D2$'$) bounds: under~\eqref{eq:B-growth-m},
\[
  \E\bigl[\norm{\calB(r, X(r))}_{L_2(U,H)}^{m_0}\bigr] \;\leq\; C_{m_0}\,\E\bigl[(1 + \norm{X(r)}_V^p)(1 + \norm{X(r)}_H^{\kappa_{m_0}})\bigr],
\]
which, by H\"older with conjugate pair $(\alpha, \alpha')$ ($1/\alpha + 1/\alpha' = 1$), separates into a $V$-power factor $(\E[(1+\norm{X(r)}_V^p)^\alpha])^{1/\alpha}$ controlled by the iterated It\^o-energy moment $\E[(\int_0^T\norm{X}_V^p)^k] < \infty$ derived from~\eqref{eq:SC-exp-moment} (uniformly in $r$ when $\alpha$ is matched to the operating point), and an $H$-power factor $(\E[(1+\norm{X(r)}_H^{\kappa_m})^{\alpha'}])^{1/\alpha'}$ controlled by the higher polynomial $H$-moments $\E[\sup_t\norm{X(t)}_H^{2k\alpha'\kappa_m/2}]$ from the iterated-energy identity. The analogous strictness in the $V$-domain regime follows the same mechanism with $p$ replaced by $p-2$, via the H\"older conjugate exponents $(p/(2\beta),\,p/(p-2\beta))$ applied to $\E[\int(1+\|X\|_V^{2\beta})\|\eta\|_V^2\,dt]$ with $X, \eta \in L^p(\Omega\times[0,T];V)$.
\end{remark}

The higher-moment analysis of Proposition~\ref{prop:Y-moments} below, that is, moment bounds
\[
  \E\bigl[\sup_s\|Y(s,r)v\|_H^q\bigr] \leq C_q\|v\|_H^q \qquad \text{for } q > 2,
\]
uses a refined version of the linearised monotonicity, which in the $V$-domain regime of~(D2) must be combined with the structural $V$-coercivity input~(SC1) of Assumption~\ref{ass:SC} below to absorb the extra It\^o term
\[
  (4q-4)\|\calB'_u(\eta)\|^2\|\eta\|^{2q-2};
\]
under the $H$-extension regime this extra term is absorbed directly by $\hat\rho(\tau,u) := C(1 + \norm{u}_V^{2\beta})$ without further structural input.

The correction term in the formula involves the Malliavin derivative of the covariance, and differentiating the covariance differentiates the first variation once more. A second order is therefore needed.

\begin{assumption}[Second-order differentiability]\label{ass:diff2}
In addition to Assumption \ref{ass:diff1}, we assume:
\begin{enumerate}[label=\textup{(D\arabic*)}, leftmargin=3em, start=3]
 \item The operator $\calA(t, \cdot): V \to V^*$ is twice Fr\'echet differentiable for each $t$, with second derivative $\calA''_{uu}(t, u) \in L^{(2)}(V \times V; V^*)$, and $u \mapsto \calA''_{uu}(t, u)$ continuous from $V$ into $L^{(2)}(V \times V; V^*)$ in the strong operator topology. The second derivative satisfies
  \begin{equation}\label{eq:D3-growth}
    \norm{\calA''_{uu}(t, u)}_{L^{(2)}(V \times V; V^*)} \leq \bigl(C + \rho_4(t, u)\bigr)\bigl(1 + \norm{u}_V^{(p-3)^+}\bigr),
  \end{equation}
  where $(p-3)^+ = \max(p-3, 0)$ and $\rho_4: [0,T] \times V \to [0,\infty)$ is measurable, with all polynomial moments along the variational solution, so that $\E[\int_0^T \rho_4(\tau, X(\tau))^{\kappa}\,\mathrm{d}\tau] < \infty$ for every $\kappa \in [1,\infty)$ (trivially satisfied whenever $\rho_4$ is bounded, as in every example of Section~\ref{subsec:verification}).
  
  \item The operator $\calB(t, \cdot)$ is twice Fr\'echet differentiable at each $u \in V$, with second derivative
  \[
    \calB''_{uu}(t, u) \in L^{(2)}(V \times V;\,L_2(U, H)),
  \]
  and $u \mapsto \calB''_{uu}(t,u)$ continuous from $V$ into $L^{(2)}(V \times V; L_2(U,H))$ in the strong operator topology, satisfying $\norm{\calB''_{uu}(t, u)}_{L^{(2)}(V \times V;\,L_2(U, H))} \leq C(1 + \norm{u}_V^{\beta'})$ for some $\beta' \geq 0$ with $2\beta' < p-2$ when $p > 2$, or $\beta' = 0$ when $p = 2$. Under the $H$-extension regime (for which the second derivative extends to $L^{(2)}(H \times H;\,L_2(U, H))$), the relaxed exponent constraint $2\beta' < p$ applies.
\end{enumerate}
\end{assumption}

\begin{remark}[Elimination of the linearised monotonicity condition]\label{rem:no-D3-old}
A separate linearised monotonicity condition is unnecessary, as the linearised local monotonicity is derived from (H2) by differentiation. Set $v = u - \eps w$ in \eqref{eq:local-mono} (so $u - v = \eps w$), divide by $\eps^2$, and let $\eps \to 0^+$.

For the left-hand side,
\begin{multline*}
  \text{LHS}(\eps)\;/\;\eps^2 \;=\; -2\,\dual{\tfrac{1}{\eps}\bigl(\calA(t, u) - \calA(t, u - \eps w)\bigr)}{w} \\
  + \tfrac{1}{\eps^2}\norm{\calB(t, u) - \calB(t, u - \eps w)}_{L_2(U,H)}^2.
\end{multline*}
The first summand converges to $-2\dual{\calA'_u(t, u)w}{w}$ as $\eps \to 0^+$ by Fr\'echet differentiability of $\calA$ in (D1). The second summand converges to $\norm{\calB'_u(t, u)(w)}_{L_2(U,H)}^2$ by Fr\'echet differentiability of $\calB$ in (D2), since $\calB(t, u) - \calB(t, u - \eps w) = \eps\,\calB'_u(t, u)(w) + o(\eps)$ in $L_2(U, H)$, so
\[
  \tfrac{1}{\eps^2}\norm{\calB(t, u) - \calB(t, u - \eps w)}_{L_2(U,H)}^2 \;=\; \norm{\calB'_u(t, u)(w)}_{L_2(U,H)}^2 + o(1).
\]
(The scalar map $w \mapsto \norm{\calB'_u(t, u)(w)}_{L_2(U,H)}^2$ is a continuous quadratic form in $w$, since $\calB'_u(t, u)$ is a bounded linear map into $L_2(U, H)$, so the $o(\eps)$ remainder in the Fr\'echet expansion gives the asserted $o(1)$ when squared and divided by $\eps^2$.) Thus
\[
  \lim_{\eps \downarrow 0}\,\text{LHS}(\eps)\;/\;\eps^2 \;=\; -2\,\dual{\calA'_u(t, u)w}{w} + \norm{\calB'_u(t, u)(w)}_{L_2(U,H)}^2.
\]
On the right-hand side, division by $\eps^2$ gives $(\rho(t, u) + \rho(t, u - \eps w))\norm{w}_H^2$, whose limsup as $\eps \downarrow 0$ is bounded by $(\rho(t,u) + \rho^*(t,u))\norm{w}_H^2$, where $\rho^*(t, u) := \limsup_{v \to u}\rho(t, v)$ is the upper-semicontinuous envelope of $\rho$ in $u$ (and $\rho^*(t, u) \geq \rho(t, u)$ always). Combining, since the LHS limit exists and the RHS limsup dominates,
\begin{equation}\label{eq:linearised-mono-derived}
  -2\dual{\calA'_u(t, u)w}{w} + \norm{\calB'_u(t, u)(w)}_{L_2(U,H)}^2 \leq \tilde\rho(t, u)\norm{w}_H^2,
\end{equation}
giving the linearised local monotonicity with $\tilde{\rho}(t, u) := \rho^*(t,u) + \rho(t,u)$. When $\rho$ is upper semicontinuous in $u$ (in particular, continuous, as is the case for the polynomial $\rho$'s of the standard examples of Section~\ref{subsec:verification}), $\rho^* = \rho$ and $\tilde\rho = 2\rho$.
\end{remark}

Differentiability of the coefficients makes the covering field available in principle; whether it is available in fact depends on the direction along which one differentiates. A direction that the noise cannot reach at time $t$ admits no covering field, and no amount of regularity will supply one.

\begin{assumption}[Non-degeneracy in the direction $h$]\label{ass:nondeg}
For a given direction $h \in H$, the Malliavin covariance operator $\gamma_t: H \to H$ defined in \eqref{eq:malliavin-cov-def} below satisfies:
\begin{enumerate}[label=\textup{(\roman*)}]
  \item \textup{(Range condition)} $h \in \Ran(\gamma_t)$ almost surely. In particular, $\gamma_t^{\dagger}h \in H$ is well defined and
  \[
    \gamma_t\gamma_t^{\dagger}h = h
  \]
  almost surely.
  \item \textup{(Moment condition: H\"older triangle)} There exist a pair of exponents $(q, p^*) \in [2, \infty]^2$ with the H\"older-triangle relation
  \begin{equation}\label{eq:nondeg-Holder-triangle}
    \frac{1}{q} + \frac{1}{p^*} \;\leq\; \frac{1}{2},
  \end{equation}
  with the conventions $1/\infty := 0$, such that
  \begin{equation}\label{eq:nondeg-moment}
    \E\bigl[\norm{\gamma_t^{\dagger} h}_H^q\bigr] \;<\; \infty
  \end{equation}
  (with the understanding that $q = \infty$ means $\norm{\gamma_t^{\dagger} h}_{L^\infty(\Omega; H)} < \infty$, the deterministic-covariance subregime), and the Hilbert--Schmidt random kernel $\Phi_r := D_r X(t)$ satisfies the corresponding Malliavin regularity
  \begin{equation}\label{eq:nondeg-Phi-regularity}
    \Phi \;\in\; \mathbb{D}^{1, p^*}\!\bigl(L^2([0,t]; L_2(U, H))\bigr).
  \end{equation}
  Three standard operating points fit~\eqref{eq:nondeg-Holder-triangle}: the deterministic-covariance regime $(q, p^*) = (\infty, 2)$ used in Corollary~\ref{cor:linear-infinite}; the Cauchy--Schwarz pair $(q, p^*) = (4, 4)$; and, for $q$ slightly above $2$, the asymmetric pair $(q, p^*) = (q,\,2q/(q-2))$ with $p^*$ large used in the multiplicative random-covariance subregime of Lemma~\ref{lemma:D-gamma-inv-h}. The H\"older-triangle parametrisation makes the trade-off between regularity of $\gamma_t^\dagger h$ and regularity of $\Phi$ explicit; statements of the form ``$q > 2$'' in earlier infinite-dimensional Bismut formulae \cite{mirafzali2025infinite} correspond to operating along the asymmetric edge $(q, 2q/(q-2))$ of the triangle.
  \item \textup{(Cameron--Martin compatibility)} The Malliavin covariance is differentiable in the Hilbert--Schmidt operator class,
  \begin{equation}\label{eq:CM-gamma-D12}
    \gamma_t \;\in\; \mathbb{D}^{1,2}\!\bigl(L_2(H,H)\bigr), \qquad \D_r\gamma_t \;\in\; L_2\!\bigl(U, L_2(H,H)\bigr) \quad\text{for a.e.\ } r \in [0,t],
  \end{equation}
  and the Tikhonov-regularised trace integrand
  \begin{equation}\label{eq:CM-trace-integrand}
    \mathcal{I}_r^{(\eps)} \;:=\; \Tr_U\!\Bigl[\,\Phi_r^{*}\,(\gamma_t + \eps I)^{-1}\,[\D_r\gamma_t]\,(\gamma_t + \eps I)^{-1} h\,\Bigr] \;\in\; \R \qquad (\eps > 0)
  \end{equation}
  satisfies the following two requirements. First, it converges. There is a jointly measurable $\mathcal{I}^{(0)}: \Omega \times [0,t] \to \R$ with
  \begin{equation}\label{eq:CM-compat-limit}
    \mathcal{I}_r^{(\eps)} \;\longrightarrow\; \mathcal{I}_r^{(0)} \qquad \mathbb{P}\text{-a.s.\ as } \eps \downarrow 0, \text{ for a.e.\ } r \in [0,t].
  \end{equation}
  Second, it is dominated. There is a measurable $\Psi: \Omega \times [0,t] \to [0,\infty]$ (independent of $\eps$) satisfying
  \begin{equation}\label{eq:CM-compat-axiom}
    \E\Bigl[\,\Bigl(\textstyle\int_0^t\Psi_r\,\mathrm{d}r\Bigr)^{\!2}\,\Bigr] \;<\;\infty
  \end{equation}
  and
  \begin{equation}\label{eq:CM-compat-dominator}
    \bigl|\mathcal{I}_r^{(\eps)}\bigr| \;\leq\; \Psi_r \qquad \text{for every } \eps \in (0,1],\; \mathbb{P}\text{-a.s.\ for a.e.\ } r \in [0,t].
  \end{equation}
  The two requirements are genuinely separate. A bound on the whole trace does not by itself license interchanging the limit $\eps \downarrow 0$ with the infinite spectral sum that defines it. The condition is basis-free and holds trivially with $\mathcal{I}^{(0)} \equiv 0$ and $\Psi \equiv 0$ when $\calA$ is linear and $\calB$ is state-independent ($\D_r\gamma_t \equiv 0$). In the nonlinear case, it is a joint convergence-and-integrability hypothesis on $(\D_r\gamma_t, \gamma_t)$. The limiting correction integrand is defined intrinsically by $\mathcal{I}^{(0)}$; no pointwise composition with the unbounded operator $\gamma_t^{\dagger}$ is assumed.
\end{enumerate}
\end{assumption}

\begin{remark}[Spectral sufficient condition for clause \textup{(iii)}]\label{rem:CM-compat-spectral}
Since $\gamma_t(\omega)$ is compact self-adjoint on the separable $H$, a jointly $\F_t$-measurable spectral resolution $\{(\lambda_k, e_k)\}_{k \geq 1}$ exists by the Courant--Fischer characterisation of eigenvalues and inductive Kuratowski--Ryll-Nardzewski measurable selection \cite{kuratowski1965general, kechris1995classical}. A sufficient condition for clause~\textup{(iii)} is that the spectral triple-sum dominator
\begin{equation}\label{eq:CM-compat-spectral}
  \Sigma^{\mathrm{spec}}_r(\omega) \;:=\; \!\!\!\!\sum_{\substack{j \geq 1\\ k,l\,:\,\lambda_k,\,\lambda_l > 0}}\!\!\!\!\frac{|h_k|\,\bigl|\bigl\langle[\D_r\gamma_t](f_j)\,e_k,\,e_l\bigr\rangle_H\bigr|\,\bigl|\bigl\langle e_l,\,\Phi_r f_j\bigr\rangle_H\bigr|}{\lambda_k\,\lambda_l}
\end{equation}
satisfies $\E[(\int_0^t \Sigma^{\mathrm{spec}}_r\,\mathrm{d}r)^{2}] < \infty$,
since $|\mathcal{I}_r^{(\eps)}| \leq \Sigma^{\mathrm{spec}}_r$ uniformly in $\eps$ by the bound $(\lambda + \eps)^{-1} \leq \lambda^{-1}$ and the kernel-annihilation property (Remark~\ref{rem:kernel-annih}). Absolute summability of~\eqref{eq:CM-compat-spectral} gives, for a.e.\ $(r,\omega)$, term-by-term passage to the limit as $\eps \downarrow 0$ and hence the convergence~\eqref{eq:CM-compat-limit}, with $\mathcal{I}^{(0)}_r$ represented by the series obtained by replacing $(\lambda_k + \eps)^{-1}(\lambda_l + \eps)^{-1}$ by $(\lambda_k\lambda_l)^{-1}$ on the positive spectral subspace; the same series supplies the dominator in~\eqref{eq:CM-compat-dominator}. The spectral condition is therefore sufficient for the whole of clause~\textup{(iii)}, not merely for its bound. The expression~\eqref{eq:CM-compat-spectral} depends on the eigenbasis within degenerate eigenspaces, but existence of any measurable selection for which it is $L^2(\Omega; L^1([0,t]))$-integrable suffices.
\end{remark}

\begin{remark}[Role of clause \textup{(iii)}]\label{rem:CM-compat-role}
Clause~(iii) supplies the two inputs the Tikhonov passage consumes, scalar convergence of the regularised trace integrand and an $L^2$-integrable dominator for it. The spectral condition of Remark~\ref{rem:CM-compat-spectral} is one sufficient route to both, but the theorems below use only the basis-free convergence~\eqref{eq:CM-compat-limit} and the dominator~\eqref{eq:CM-compat-dominator}. It holds trivially when $\D_r\gamma_t \equiv 0$ (linear drift, state-independent noise), and becomes a joint convergence-and-integrability hypothesis whenever $Y(t,r)$ is random. It is strictly weaker than $\mathbb{D}^{1,2}(\HW)$-regularity of $v_h$, requiring convergence only at the scalar trace level.
\end{remark}

Comparison with global non-degeneracy is instructive here. Assumption \ref{ass:nondeg} is strictly weaker than global injectivity of $\gamma_t$, requiring only that $h \in \Ran(\gamma_t)$ with $\sum_{k:\lambda_k>0}(h_k/\lambda_k)^2 < \infty$. When $\gamma_t$ is injective, $\gamma_t^{\dagger} = \gamma_t^{-1}$; when $\gamma_t$ is not injective, $\gamma_t^{\dagger}h$ is the minimum-$H$-norm element of $\gamma_t^{-1}(\{h\})$.

One matter is left. The linearised equations are themselves stochastic evolution equations, and the monotonicity conditions of Assumption~\ref{ass:LR} govern the original equation, not its linearisations. What follows collects the structure they need in order to be well posed and to carry moments, namely a dissipative mechanism, one moment of the forcing that mechanism has to absorb, an exponential moment for the moduli it produces, and stability of the Galerkin scheme at the linearised level. The objects the conditions speak about are the variational fibre $\Phi_{\tau,r} = Y(\tau,r)\calB(r,X(r))$ and the pair of fibres from which the second-order forcing is built, never the first variation on an arbitrary direction of the state space; this is what the formula uses, and asking for regularity in directions it does not use would be asking for smoothing that parabolic equations do not have. These four conditions are used everywhere below and nowhere else are they assumed.

\begin{assumption}[Structure of the linearised equations]\label{ass:SC}
For the exponent $p \geq 2$ of Assumption~\ref{ass:LR} and its H\"older conjugate $p' := p/(p-1) \in (1, 2]$, and for the linearised-local-monotonicity coefficient $\tilde\rho(\tau, u) := \rho(\tau, u) + \rho^{*}(\tau, u)$ of Remark~\ref{rem:no-D3-old} together with the derivative polynomial-growth coefficient $\hat\rho(\tau, u) := C(1 + \norm{u}_V^{2\beta})$ of~(D2), there exists a slack exponent $\eps_0 > 0$ such that the following four conditions hold along the variational solution $X$ of~\eqref{eq:SPDE}:
\begin{enumerate}[label=\textup{(SC\arabic*)}, leftmargin=3em]
  \item \textup{(Linearised structural coercivity: the tangent energy.)} Write
  \begin{equation}\label{eq:tangent-form}
    D_\sigma(w) \;:=\; \dual{\calA'_u(\sigma, X(\sigma))\,w}{w}, \qquad w \in V,
  \end{equation}
  for the \emph{tangent form} of the linearisation along $X$, and
  \begin{equation}\label{eq:tangent-energy}
    \calE_r^t(w) \;:=\; \int_r^t D_\sigma(w(\sigma))\,\mathrm{d}\sigma
  \end{equation}
  for the associated \emph{tangent energy} of a path $w$. The linearised principal part satisfies at least one of the two alternatives~(a),~(b) below, each of which carries its own well-posedness datum for the linear equation; the two are never mixed. Under alternative~(a), where $D_\sigma \geq 0$, we require that, $\mathbb{P}$-almost surely, the family $(D_\sigma)_{\sigma \in [0,T]}$ be a \emph{measurable family of closable forms with common core $V$}, in the sense that each $D_\sigma$ is closable on $H$ with completion $\calV_\sigma \supseteq V$, the map $\sigma \mapsto D_\sigma(w)$ is measurable for each $w \in V$, and the space-time form domain
  \begin{equation}\label{eq:spacetime-form-domain}
  \begin{aligned}
    \mathbb{V}_r \;:=\; \bigl\{\,w \in L^2(r,T;H) \;:\; w(\sigma) \in \calV_\sigma \text{ a.e.},\ \ \calE_r^T(w) < \infty \,\bigr\},\\
    \norm{w}^2_{\mathbb{V}_r} := \norm{w}^2_{L^2(r,T;H)} + \calE_r^T(w),
  \end{aligned}
  \end{equation}
  is a separable Hilbert space, densely embedded in $L^2(r,T;H)$, containing the $V$-valued step functions densely; and, $\mathbb{P}$-almost surely, the linear equation~\eqref{eq:Y-IVP} read weakly against $V$ through the form $D_\sigma$ is pathwise well posed on each stopping-time interval, in $C([r,\tau];H)$ with finite tangent energy, with its Galerkin approximations converging there. This is the non-autonomous form datum of the linearisation; for gradient-type drifts it is the classical Dirichlet-form situation, with weight determined by $X(\sigma)$. All of $D_\sigma(w)$, $\calA'_u(\sigma,X(\sigma))w$, and $\calB'_u(\sigma,X(\sigma))w$ are defined for $w \in V$ and extended to $\calV_\sigma$ by continuity in the tangent norm. Under alternative~(b) no form domain is constructed. The tangent form need be neither symmetric nor non-negative there, for the Navier--Stokes transport term $\dual{\calA'_u(u)w}{w}$ can take either sign, and the strict quadratic $V$-coercivity of~(b) directly renders the linear equation pathwise well posed in $C([r,T];H) \cap L^2(r,T;V)$ by the classical variational scheme on the fixed triple $V \hookrightarrow H \hookrightarrow V^*$. In both alternatives the datum is pathwise; the global moment statements are supplied by~(SC2). Explicitly, the two alternatives are, uniformly for $\sigma \in [0, T]$, $\mathbb{P}$-almost surely:
  \begin{enumerate}[label=\textup{(\alph*)}, leftmargin=2.6em]
 \item \textup{(Dissipative linearisation with dominated second derivative.)} The dissipation, that is the quadratic form $D_\sigma(w) := \dual{\calA'_u(\sigma, X(\sigma))\,w}{w}$, is non-negative for every $w \in V$, and there exists a constant $\delta \in (0, 2]$ with
    \[
      \begin{aligned}
        &-2\dual{\calA'_u(\sigma, X(\sigma))w}{w} + \norm{\calB'_u(\sigma, X(\sigma))(w)}^2_{L_2(U,H)}\\
        &\qquad\qquad \;\leq\; \tilde\rho(\sigma, X(\sigma))\,\norm{w}_H^2 - \delta\,D_\sigma(w)
      \end{aligned}
    \]
    for every $w \in V$ (automatic with $\delta = 2$ when $\calB'_u \equiv 0$, the state-independent-noise case), and there exists a measurable family $\Lambda_\sigma: V \times V \to [0, \infty)$ dominating the drift second derivative through the dissipation,
    \begin{multline*}
      \abs{\dual{\calA''_{uu}(\sigma, X(\sigma))(\varphi, \psi)}{w}} \;\leq\; \Lambda_\sigma(\varphi, \psi)\,D_\sigma(w)^{1/2} \\
      \qquad \text{for all } \varphi, \psi, w \in V
    \end{multline*}
    and, for the associated bilinear form, the form Cauchy--Schwarz bound
    \[
      \abs{\dual{\calA'_u(\sigma, X(\sigma))w}{z}} \;\leq\; D_\sigma(w)^{1/2}\,D_\sigma(z)^{1/2}, \qquad w, z \in V,
    \]
    which is automatic whenever $\calA'_u(\sigma, X(\sigma))$ is symmetric and non-negative. The dominator is required to be \emph{jointly square-summable}, in the sense that for all $\Theta, \Xi \in L_2(U, V)$ the double sum
    \begin{equation}\label{eq:Lambda-square-summable}
      \Lambda_\sigma^{\otimes 2}(\Theta, \Xi) \;:=\; \sum\nolimits_{j,k}\Lambda_\sigma\bigl(\Theta f_j,\, \Xi f_k\bigr)^2
    \end{equation}
    is finite and independent of the orthonormal basis $(f_j)$ of $U$. This is what turns the pointwise domination into the estimate on finite families through which the second-order object is reconstructed from its directional fibres. For the $p$-Laplacian, where $\Lambda_\sigma(\varphi, \psi)^2 = C_p^2\int\abs{\nabla X(\sigma)}^{p-4}\abs{\nabla\varphi}^2\abs{\nabla\psi}^2\,\mathrm{d}x$, the double sum is the pointwise product
    \[
      \Lambda_\sigma^{\otimes 2}(\Theta, \Xi) \;=\; C_p^2\int \abs{\nabla X(\sigma)}^{p-4}\Bigl(\sum\nolimits_j\abs{\nabla \Theta f_j}^2\Bigr)\Bigl(\sum\nolimits_k\abs{\nabla \Xi f_k}^2\Bigr)\,\mathrm{d}x ,
    \]
    manifestly basis-free. This is the situation for every gradient-type drift, such as the stochastic $p$-Laplacian and the porous medium equation; for the $p$-Laplacian at $p \geq 4$, $\Lambda$ is computed by weighted Cauchy--Schwarz in Proposition~\ref{prop:p-Lap-Lambda}. Finally, $\Lambda_\sigma$ and $\calB''_{uu}(\sigma, X(\sigma))$ are required to extend by continuity in the tangent norm to $\calV_\sigma \times \calV_\sigma$, exactly as $\calA'_u(\sigma,X(\sigma))$ and $\calB'_u(\sigma,X(\sigma))$ do, so that under~(a) the second-order data of the linearisation live on the form domain rather than on $V$. For $\calB''_{uu} \equiv 0$, and more generally whenever $\calB''_{uu}$ extends to $H \times H$, the second requirement is automatic.
 \item \textup{(Strict linearised $V$-coercivity.)} There exists $\alpha_1 > 0$ such that
    \[
      \begin{aligned}
      &-2\dual{\calA'_u(\sigma, X(\sigma))w}{w} + \norm{\calB'_u(\sigma, X(\sigma))(w)}^2_{L_2(U,H)} \\
      &\qquad \;\leq\; \tilde\rho(\sigma, X(\sigma))\norm{w}_H^2 - \alpha_1\norm{w}_V^2
      \end{aligned}
    \]
    for every $w \in V$ (e.g., the stochastic 2D Navier--Stokes equation, with $\alpha_1 = \nu$, and semilinear parabolic equations with analytic-semigroup smoothing). The quadratic $V$-power matches the two-homogeneity of the linearised equations for the two fibres, for which every energy estimate below is run at exponent two; at $p = 2$ it coincides with the nonlinear coercivity exponent. Equivalently, under~(b) the tangent form dominates the $V$-norm modulo $H$,
\[
  2 D_\sigma(w) \;\geq\; \alpha_1\norm{w}_V^2 - \tilde\rho(\sigma, X(\sigma))\norm{w}_H^2,
\]
so that finiteness of the tangent energy~\eqref{eq:tangent-energy} upgrades to $L^2$-in-time $V$-regularity; under~(a) the tangent energy is the intrinsic substitute, degenerating exactly where the linearisation does.
  \end{enumerate}
  Both alternatives equip the linearisation with one and the same quadratic \emph{fibre energy}. Write
  \begin{equation}\label{eq:fibre-norm}
    \norm{w}^2_{\mathbb{V}_\sigma} \;:=\;
    \begin{cases}
      D_\sigma(w), & \text{under (a)},\\[2pt]
      \alpha_1\norm{w}_V^2, & \text{under (b)},
    \end{cases}
    \qquad
    \norm{\Theta}^2_{L_2(U, \mathbb{V}_\sigma)} \;:=\; \sum\nolimits_j\norm{\Theta f_j}^2_{\mathbb{V}_\sigma},
  \end{equation}
  the second quantity being independent of the orthonormal basis $(f_j)$ of $U$ because $w \mapsto \norm{w}^2_{\mathbb{V}_\sigma}$ is a quadratic form. In this notation the two alternatives are the single dissipative inequality
  \begin{equation}\label{eq:SC1-unified}
  \begin{aligned}
    -2\dual{\calA'_u(\sigma, X(\sigma))w}{w} &+ \norm{\calB'_u(\sigma, X(\sigma))(w)}^2_{L_2(U,H)}\\
    &\leq\; \tilde\rho(\sigma, X(\sigma))\,\norm{w}_H^2 \;-\; \delta\,\norm{w}^2_{\mathbb{V}_\sigma},
  \end{aligned}
  \end{equation}
  with $\delta \in (0, 2]$ under~(a) and $\delta = 1$ under~(b), so that one energy estimate serves both mechanisms, run in the intrinsic weighted energy under~(a) and in the unweighted $V$-energy under~(b).
  \setcounter{enumi}{1}
  \item[] \textup{(SC1$'$)} \textup{(Forcing integrability of the linearised system.)} Write, for $0 \leq r \leq \tau \leq T$,
  \begin{equation}\label{eq:Phi-fibre}
    \Phi_{\tau,r} \;:=\; Y(\tau,r)\,\calB(r, X(r)) \;\in\; L_2(U, H)
  \end{equation}
  for the \emph{variational fibre} built from the first variation and the diffusion coefficient (it is identified with $\D_r X(\tau)$ in Proposition~\ref{prop:D-Xt}, but no Malliavin calculus is used in stating this clause), and
  \begin{equation}\label{eq:fibre-energy}
    \mathfrak{E}_r \;:=\; \int_r^T \norm{\Phi_{\tau,r}}^2_{L_2(U, \mathbb{V}_\tau)}\,\mathrm{d}\tau
  \end{equation}
  for its fibre energy. That $\Phi_{\tau,r}$ lies in $L_2(U, \mathbb{V}_\tau)$ for almost every $\tau$, and that $\mathfrak{E}_r$ carries moments uniformly in $r$, is not assumed; it is the content of Theorem~\ref{thm:Y-wellposed}\textup{(III)}, which runs the energy estimate~\eqref{eq:SC1-unified} on the fibre itself.

  The second-order data of the linearised system are the two \emph{double-fibre forcings}, in which both slots of the second derivatives are variational fibres. They are read through the mechanism in force. Under~(b) the second derivatives are used through their typing in Assumption~\ref{ass:diff2} on $V \times V$, and $\mathbb{V}_\tau$ is $V$ itself; under~(a) they are used through the dominator $\Lambda_\tau$ and the tangent extension of $\calB''_{uu}$, both defined on the form domain $\calV_\tau \supseteq V$ on which $\mathbb{V}_\tau$ is modelled. In either case the fibres lie in $L_2(U, \mathbb{V}_\tau)$ for almost every $\tau$ by Theorem~\ref{thm:Y-wellposed}\textup{(III)}, which is the space the data act on. Writing
  \begin{equation}\label{eq:bold-forcings}
  \begin{aligned}
    \mathbf{F}_1(\tau; r,s)(a, u) &\;:=\; -\calA''_{uu}(\tau, X(\tau))\bigl(\Phi_{\tau,r}a,\, \Phi_{\tau,s}u\bigr) \;\in\; V^{*},\\
    \mathbf{F}_2(\tau; r,s)(a, u) &\;:=\; \calB''_{uu}(\tau, X(\tau))\bigl(\Phi_{\tau,r}a,\, \Phi_{\tau,s}u\bigr) \;\in\; L_2(U,H),
  \end{aligned}
  \end{equation}
  for $a, u \in U$ under~(b), and under~(a) reading $\mathbf{F}_1$ only through $\Lambda_\tau$ as in~\eqref{eq:G1-def} below, together with the symmetric initial datum
  \begin{equation}\label{eq:J-datum}
    \mathfrak{J}_{r,s} \;:=\; \bigl\lVert \calB'_u\bigl(\sigma, X(\sigma)\bigr)\bigl(\Phi_{\sigma,\,\varsigma}\bigr)\bigr\rVert^2_{L_2(U,\,L_2(U,H))},
  \end{equation}
  in which $\sigma := \max(s,r)$ and $\varsigma := \min(s,r)$, so that the fibre started at the earlier of the two times is inserted into the derivative of the diffusion coefficient at the later one. Set
  \begin{equation}\label{eq:G1-def}
    G_1(\tau; r,s) \;:=\;
    \begin{cases}
      \Lambda_\tau^{\otimes 2}\bigl(\Phi_{\tau,r},\, \Phi_{\tau,s}\bigr), & \text{under (a)},\\[2pt]
      \norm{\mathbf{F}_1(\tau;r,s)}_{L_2(U \otimes U,\, V^{*})}^2, & \text{under (b)},
    \end{cases}
  \end{equation}
  the first line being the double square sum~\eqref{eq:Lambda-square-summable} of the domination and the second the double square-summed drift forcing, and write
  \begin{equation}\label{eq:total-forcing}
    \mathfrak{G}_{r,s} \;:=\; \mathfrak{J}_{r,s} \;+\; \int_{\max(s,r)}^{T}\Bigl[\,G_1(\tau;r,s) \;+\; \norm{\mathbf{F}_2(\tau;r,s)}^2_{L_2(U \otimes U,\, L_2(U,H))}\,\Bigr]\,\mathrm{d}\tau
  \end{equation}
  for the \emph{total forcing}, the single quantity by which every datum of the linearised system is measured. The clause asserts that the three quantities entering~\eqref{eq:total-forcing} are almost everywhere defined and measurable, and that there exists a slack $\eps_0 > 0$ such that, for all $0 \leq r, s \leq T$,
  \begin{equation}\label{eq:forcing-integrability}
    \E\Bigl[\,\mathfrak{G}_{r,s}^{\,1+\eps_0}\,\exp\Bigl(c_*\!\int_0^T\bigl(\tilde\rho + \hat\rho\bigr)\,\mathrm{d}\sigma\Bigr)\Bigr] \;<\; \infty ,
  \end{equation}
  with $c_*$ the operating constant of~(SC2) below.

  Three features of~\eqref{eq:total-forcing} determine the shape of every later estimate. It is quadratic in the pair $(\Phi_{\tau,r}, \Phi_{\tau,s})$, which matches the two-homogeneity of the linearised equations. It is already summed over both noise slots, which is what licenses the reconstruction of the second-order object from its directional fibres. And it is integrated in time before any moment is taken, which is exactly the order in which the stochastic Gronwall step consumes it, so that~\eqref{eq:forcing-integrability} asks for no pointwise-in-time integrability that the proof does not use. No operator norm of $Y$ occurs anywhere in it, and no auxiliary vector of $H$ has to be carried through the estimate. In the directional estimates below we write $F_i^{a,u}(\tau;r,s) := \mathbf{F}_i(\tau;r,s)(a,u)$ for the fibre in the directions $a, u \in U$.

  Because both slots are fibres, the two mechanisms discharge~\eqref{eq:forcing-integrability} in their own currency. Under~(b) the fibre energy is the unweighted $V$-energy, so~\eqref{eq:total-forcing} is bounded by the two fibre energies and the fibre supremum $\sup_\tau\norm{\Phi_{\tau,r}}_{L_2(U,H)}$ whenever the second derivatives interpolate between $H$ and $V$, which is how the clause is obtained for the two-dimensional Navier--Stokes equation in Proposition~\ref{prop:NS-forcing}, and by the polynomial closures of Assumption~\ref{ass:SC-raised} whenever the growth bounds~\eqref{eq:D3-growth} and~(D4) separate at compatible exponents. Under~(a) the fibre energy is the weighted tangent energy, and $G_1$ is the weighted double sum $\Lambda_\tau^{\otimes 2}$, which is a quantity of the same weighted kind; for gradient-type drifts it is a single explicit integral, computed for the $p$-Laplacian in Remark~\ref{rem:p-Lap-forcing}, and the clause is the datum of the degenerate case in the same way that the form domain of~(SC1)(a) is.
  \setcounter{enumi}{1}
  \item \textup{(Exponential integrability of the linearised moduli at a constant $c_*$.)} The assumption carries a finite parameter $c_* \in (0, \infty)$, its \emph{operating constant}, and asserts that
  \begin{equation}\label{eq:SC-exp-moment}
    \E\Bigl[\exp\Bigl(c_*\int_0^T \bigl(\tilde\rho(\sigma,X(\sigma)) + \hat\rho(\sigma,X(\sigma))\bigr)\,\mathrm{d}\sigma\Bigr)\Bigr] \;<\; \infty,
  \end{equation}
  and the same holds along the Galerkin approximants $\{X^N\}_{N \geq 1}$ of \cite{rockner2022wellposedness} (recalled in \eqref{eq:galerkin} below), uniformly in $N$, at the same $c_*$,
  \begin{equation}\label{eq:SC-exp-moment-Galerkin}
    \sup_{N \in \N}\,\E\Bigl[\exp\Bigl(c_*\int_0^T \bigl(\tilde\rho(\sigma,X^N(\sigma)) + \hat\rho(\sigma,X^N(\sigma))\bigr)\,\mathrm{d}\sigma\Bigr)\Bigr] \;<\; \infty.
  \end{equation}
  The operating constant is required to dominate the finitely many constants generated by the base estimates:
  \begin{equation}\label{eq:cstar-base}
    c_* \;\geq\; \max\Bigl\{\,\tfrac{(1+\eps_0)\,c_1}{\eps_0},\;\; \tfrac{c_1 m_0}{m_0-2}\ \ (\text{state-dependent noise only})\,\Bigr\},
  \end{equation}
  with $\eps_0$ the slack of clause~(SC1$'$) and $c_1$ the Gronwall constant at $q = 1$; the raised estimates impose the further requirement~\eqref{eq:cstar-requirement} on $c_*$ at the target $q^*$. Since $c \mapsto \E[\exp(c\int(\tilde\rho+\hat\rho))]$ is non-decreasing, (SC2) at $c_*$ implies (SC2) at every $c \leq c_*$; and \eqref{eq:SC-exp-moment} holds at every finite $c > 0$, in which case $c_*$ may be taken as large as any statement requires, whenever $\tilde\rho$ and $\hat\rho$ are bounded along solutions, as for the stochastic $p$-Laplacian with state-independent noise. Each moment statement below is indexed by the constant it consumes, through
  \begin{equation}\label{eq:qmax-def}
    \hat c_{\,p} \;:=\; c_{p/2}, \qquad
    q_{\max} \;:=\; \sup\{\,p \geq 2 \;:\; \hat c_{\,p} \leq c_*\,\} \;\in\; [2, \infty],
  \end{equation}
  where $c_q$ is the explicit (increasing) Gronwall constant produced by It\^o's formula on $\norm{\cdot}_H^{2q}$ together with the Burkholder--Davis--Gundy constant at exponent $q$, so that $\hat c_{\,p}$ is indexed by the \emph{moment exponent} $p = 2q$ that the corresponding estimate delivers. Thus a bound at moment exponent $p$ is available precisely when $\hat c_{\,p} \leq c_*$, that is when $p \leq q_{\max}$; $q_{\max} = \infty$ when \eqref{eq:SC-exp-moment} holds at every finite $c$.
  \item \textup{(Strengthened Galerkin convergence of the state, with uniform higher moments.)} The Galerkin approximants $\{X^N\}_{N \geq 1}$ of \cite{rockner2022wellposedness} (recalled in \eqref{eq:galerkin} below) converge in $L^2(\Omega; C([0,T]; H))$ and in $L^{b}(\Omega \times [0,T]; V)$ at the raised exponent $b := 6(1+\eps_0)$,
  \begin{equation}\label{eq:SC-Galerkin}
    \norm{X^N - X}_{L^2(\Omega; C([0,T]; H))} + \norm{X^N - X}_{L^{b}(\Omega \times [0,T]; V)} \to 0 \quad \text{as } N \to \infty,
  \end{equation}
  and obey the uniform higher-moment bound
  \begin{equation}\label{eq:SC-Galerkin-moments}
    \sup_{N \in \N}\;\E\Bigl[\int_0^T \bigl(1 + \norm{X^N(\sigma)}_V\bigr)^{a}\,\mathrm{d}\sigma\Bigr] \;+\; \E\Bigl[\int_0^T\bigl(1+\norm{X(\sigma)}_V\bigr)^{a}\,\mathrm{d}\sigma\Bigr] \;<\; \infty .
  \end{equation}
  Here $a := 6(1+\eps_0)(p-3)^+$; at $p \leq 3$ one has $a = 0$ and \eqref{eq:SC-Galerkin-moments} is vacuous. In addition, the Galerkin scheme is stable at the level of the two fibres. Writing $\Phi^N$ and $\calW^N$ for the variational fibre and the second Malliavin fibre of the Galerkin system~\eqref{eq:galerkin} (finite-dimensional, hence classical), for all $0 \leq s \leq \tau \leq T$,
  \begin{equation}\label{eq:SC-Galerkin-linearised}
  \begin{aligned}
    \E\int_0^\tau\norm{\Phi^N_{\tau,r} - \Phi_{\tau,r}}^2_{L_2(U,H)}\dr&\\
    +\; \E\int_0^T\norm{\calW^N(\tau;s,r) - \calW(\tau;s,r)}^2_{L_2(U,\,L_2(U,H))}\dr &\;\longrightarrow\; 0
  \end{aligned}
  \end{equation}
  as $N \to \infty$. Both terms are integrated in the Malliavin variable $r$, which is the topology in which the two identifications $\D_r X(\tau) = \Phi_{\tau,r}$ and $\D_r\Phi_{\tau,s} = \calW(\tau;s,r)$ are asserted; the Malliavin derivative of an $H$-valued random variable is an element of $L^2(\Omega \times [0,T]; L_2(U,H))$ and that of an $L_2(U,H)$-valued one an element of $L^2(\Omega \times [0,T]; L_2(U, L_2(U,H)))$, so convergence at each fixed $r$ would not suffice to close either derivative. Clause~\eqref{eq:SC-Galerkin-linearised} is the linearised companion of~\eqref{eq:SC-Galerkin}, stated directly on the two objects the theory uses rather than on the flow from which they are built. It transports the finite-dimensional identities of the Galerkin system to the limit, and it is the only place where Galerkin stability enters the paper, once for the Malliavin differentiability of $X(t)$ (Proposition~\ref{prop:D-Xt}) and once for the identification $\D_r\Phi = \calW$ (Lemma~\ref{lemma:DPhi}). Both fibres are governed by the same tangent form~\eqref{eq:tangent-form}, so a single stability input serves both.
\end{enumerate}
Condition~(SC3) is a genuine structural input at every $p \geq 2$. The compactness method of \cite{rockner2022wellposedness} yields tightness of the Galerkin laws in $C([0,T]; V^*) \cap L^p([0,T]; H)$, uniform moment bounds, and identification of subsequential limits via a Skorokhod representation on an auxiliary probability space, but not strong convergence of the Galerkin sequence on the original probability space in the norms of~\eqref{eq:SC-Galerkin}. By H\"older's inequality on the finite measure space $\Omega \times [0,T]$, the $V$-clause of~\eqref{eq:SC-Galerkin} at one exponent implies it at all smaller exponents.
\end{assumption}

\begin{remark}[Role of Assumption~\ref{ass:SC}]\label{rem:SC-role}
The assumption has four clauses and they do four different things. (SC1) supplies the dissipative mechanism~\eqref{eq:SC1-unified} and with it the well-posedness and the energy of both the first variation and the second Malliavin fibre; (SC1$'$) supplies the one moment of the total forcing~\eqref{eq:total-forcing} that the stochastic Gronwall step consumes; (SC2) supplies the exponential moment of the linearised moduli through which that step closes; and (SC3) supplies the Galerkin stability used twice, once for the Malliavin differentiability of $X(t)$ in Proposition~\ref{prop:D-Xt} and once for the identification of the second Malliavin fibre in Lemma~\ref{lemma:DPhi}. No clause asks for an operator norm of the first variation between $H$ and $V$, and none asks for regularity of the linearisation in a direction that the formula does not use. The $V$-regularity that the second derivatives of the coefficients require enters only through the fibre energy~\eqref{eq:fibre-energy}, which the dissipative mechanism itself produces. In the concrete equations of Section~\ref{subsec:verification} the clauses of~(SC1) are verified from the structure of the drift and (SC2) from the explicit form of $\tilde\rho$ and $\hat\rho$. Under the strictly coercive mechanism~(b), (SC1$'$) is then a consequence of~(SC1) and the fibre energy, as Proposition~\ref{prop:NS-forcing} shows for the two-dimensional Navier--Stokes equation. Under the degenerate mechanism~(a) it is the weighted double sum~\eqref{eq:Lambda-square-summable} evaluated at the two fibres, a single explicit integral for gradient-type drifts, and it is the datum of that mechanism in the same way that its form domain is.
\end{remark}

The estimates of Assumption~\ref{ass:SC} are stated at the exponent the well-posedness argument requires. Wherever a moment of higher order is needed, it is needed of the same quantities, and we record this once rather than at each occurrence.

\begin{assumption}[The same structure at a raised exponent $q^* \geq 1$]\label{ass:SC-raised}
Let $q^* \geq 1$ be a fixed target exponent. In addition to Assumption~\ref{ass:SC}, suppose there exists a slack $\eps_{q^*} > 0$ such that, writing
\begin{equation}\label{eq:kappa-def}
  \kappa_{q^*} \;:=\; (1 + \eps_{q^*})\,q^*
\end{equation}
for the \emph{raised exponent}, the operating constant of~(SC2) satisfies
\begin{equation}\label{eq:cstar-requirement}
  \begin{aligned}
  c_* \;\geq\; \kappa(q^*) \;:=\; \max\Bigl\{\,&\hat c_{\,2q^*},\;\; \tfrac{c_1\,m_0}{m_0-2}\;\;(\text{state-dependent noise only}),\\
  &\hat c_{\,2q^*}\kappa_{q^*},\;\; \tfrac{\hat c_{\,2q^*}\kappa_{q^*}}{\eps_{q^*}}\,\Bigr\}
  \end{aligned}
\end{equation}
(equivalently $q_{\max} \geq 2q^*$ together with the two separation constants; the middle entry is omitted for deterministic noise, where no separation is performed), and the following hold:
\begin{enumerate}[label=\textup{(SC$4$.\arabic*)$_{q^*}$}, leftmargin=3.5em]
  \item \textup{(Raised polynomial closure on $X$.)} The state process satisfies
  \begin{equation}\label{eq:SC4-Xclose}
    \E\Bigl[\int_0^T\bigl(1 + \norm{X(\tau)}_V\bigr)^{6(p-3)^+\kappa_{q^*}}\,\mathrm{d}\tau\Bigr] + \E\Bigl[\int_0^T\bigl(1 + \norm{X(\tau)}_V\bigr)^{6\beta'\kappa_{q^*}}\,\mathrm{d}\tau\Bigr] \;<\; \infty.
  \end{equation}
  \item \textup{(Raised polynomial closure on $\calB$.)} The noise coefficient satisfies
  \begin{multline}\label{eq:SC4-Bclose}
    \E\Bigl[\int_0^T\norm{\calB(r,X(r))}_{L_2(U,H)}^{2\kappa_{q^*}}\,\mathrm{d}r\Bigr] \\
    + \E\Bigl[\int_0^T\norm{\calB(r,X(r))}_{L_2(U,H)}^{2\kappa_{q^*}/\eps_{q^*}}\,\mathrm{d}r\Bigr] \;<\; \infty,
  \end{multline}
  the second at the exponent consumed by the kernel bound on $\Phi$ in Remark~\ref{rem:Phi-all-moments} (finite at every order for state-independent noise, and in general supplied by the product structure of~\eqref{eq:B-growth-m} together with the iterated-energy moments of~\eqref{eq:higher-moment-ito}).
  \item \textup{(Raised forcing integrability.)} With $\mathfrak{G}_{r,s}$ the total forcing~\eqref{eq:total-forcing} of clause~(SC1$'$), for all $0 \leq r, s \leq T$,
  \begin{equation}\label{eq:SC4-forcing}
    \E\Bigl[\,\mathfrak{G}_{r,s}^{\,\kappa_{q^*}}\,\exp\Bigl(c_*\!\int_0^T\bigl(\tilde\rho + \hat\rho\bigr)\,\mathrm{d}\sigma\Bigr)\Bigr] \;<\; \infty .
  \end{equation}
  This is~\eqref{eq:forcing-integrability} at the raised exponent; at $q^* = 1$ and $\eps_1 = \eps_0$ the two coincide. Clauses (SC4.1)$_{q^*}$--(SC4.2)$_{q^*}$ are the polynomial closures through which~\eqref{eq:SC4-forcing} is verified whenever the growth bounds~\eqref{eq:D3-growth} and~(D4) separate from the fibre energy~\eqref{eq:fibre-energy} at compatible exponents, as they do for every equation of Section~\ref{subsec:verification}; the factor $6$ in~\eqref{eq:SC4-Xclose} is the room for the H\"older triple that performs that separation.
\end{enumerate}
At the base target $q^* = 1$, with $\eps_1 := \eps_0$, clause (SC4.3)$_1$ coincides with~\eqref{eq:forcing-integrability} and carries no additional content. Clauses (SC4.1)$_1$--(SC4.2)$_1$ remain to be checked whenever their displayed exponents exceed the basic state and noise moments already available, namely $\E\int_0^T\norm{X}_V^p\,\mathrm{d}\tau$ from~\eqref{eq:apriori} and the single $m_0$-moment of~\eqref{eq:B-growth-m}; when the displayed exponents fall below those, they follow by Jensen.

Every exponent above is quadratic in the fibre, which is the homogeneity of the linearised equations; no exponent is dictated by a duality pairing against an operator norm, and in particular the exponent $\max\{p', 2\}$ that a factorisation through $L(H,V)$ would impose does not occur.
\end{assumption}

\begin{remark}[The $H$-extension of the second derivative]\label{rem:SC-raised-two-regimes}
The forcings~\eqref{eq:bold-forcings} are typed through the second derivatives of Assumption~\ref{ass:diff2}, which act on $V \times V$, and the fibre energy~\eqref{eq:fibre-energy} is what makes both slots available there. When $\calB''_{uu}$ extends by continuity to $L^{(2)}(H \times H; L_2(U,H))$, which is the case for state-independent noise, where it vanishes, and for a diffusion coefficient depending on the state through an $H$-continuous functional, the $\mathbf{F}_2$-part of the total forcing may be measured with both slots in the Hilbert--Schmidt norm of $L_2(U,H)$ instead, giving the strictly sharper bound
\[
  \norm{\mathbf{F}_2(\tau;r,s)}_{L_2(U \otimes U,\, L_2(U,H))} \;\leq\; C\bigl(1 + \norm{X(\tau)}_V^{\beta'}\bigr)\,\norm{\Phi_{\tau,r}}_{L_2(U,H)}\,\norm{\Phi_{\tau,s}}_{L_2(U,H)},
\]
in which no fibre energy is consumed at all and the relaxed exponent constraint $2\beta' < p$ of~(D4) applies. All concrete applications in this paper, namely the linear-additive, semilinear-additive, $p$-Laplacian with state-independent noise, and two-dimensional Navier--Stokes equations, are of this kind. For a genuinely $V$-dependent diffusion coefficient the default typing of~(D4) is used, with the strict constraint $2\beta' < p-2$, and the well-posedness input is Part~II of \cite{rockner2022wellposedness}, whose hypotheses \textup{(H2)*}--\textup{(H5)*} (including the additional conditions on $\calB$ introduced in the erratum \cite{rockner2025erratum}) then replace~(H2) and~(H5) of Assumption~\ref{ass:LR}. Either way the estimate is quadratic in the pair of fibres and nothing else changes.
\end{remark}

\begin{remark}[Role of Assumption~\ref{ass:SC-raised}]\label{rem:SC-raised-role}
The raised assumption is Assumption~\ref{ass:SC} evaluated at the exponent that the Malliavin-regularity chain of Remark~\ref{rem:Phi-all-moments} consumes, and it is invoked only at $q^* = q/(q-2)$, the exponent generated by a choice of $q \in (2,\infty]$ in Assumption~\ref{ass:nondeg}(ii), with $q^* = 1$ at $q = \infty$. Clause (SC4.3)$_{q^*}$ is the operative one; clauses (SC4.1)$_{q^*}$--(SC4.2)$_{q^*}$ are the routes through which it is verified, and both follow from~(SC2) at the corresponding raised exponents whenever $\tilde\rho$, $\hat\rho$ and the noise coefficient are bounded along solutions.
\end{remark}
\subsection{The three objects the formula is built from}

With the hypotheses in place we can name the objects the formula is made of. There are three: the first variation of the equation, the Malliavin derivative of the solution together with its covariance, and the field on Wiener space along which the integration by parts is performed. Each is introduced now and used throughout.

\begin{definition}[First variation process]\label{def:first-var}
Under Assumptions~\ref{ass:LR}, \ref{ass:diff1}, and~\ref{ass:SC}, the \emph{first variation process} $Y(t,r)$ for $0 \leq r \leq t \leq T$ is the linear map $H \to L^0(\Omega; H)$ defined fibrewise: for each $v \in H$, $Y(t,r)v$ is the unique variational solution in the \emph{tangent solution class}
\[
  \calY_r \;:=\;
  \begin{cases}
    \bigl\{\,w \in C([r,T]; H)\;:\; \calE_r^T(w) < \infty \ \text{a.s.}\,\bigr\}, & \text{under (SC1)(a)},\\[2pt]
    C([r,T];H) \cap L^2(r,T;V), & \text{under (SC1)(b)},
  \end{cases}
\]
tested against $V$ (a.s.), of the linear stochastic evolution equation
\begin{equation}\label{eq:first-var}
\begin{aligned}
  \mathrm{d}[Y(t,r)v] + \calA'_u(t, X(t))\,[Y(t,r)v]\,\dt &= \calB'_u(t, X(t))\bigl([Y(t,r)v]\bigr)\,\dW(t), \quad t \in (r, T],\\
  Y(r,r)v &= v,
\end{aligned}
\end{equation}
with $\calA'_u(t, X(t)) \in L(V, V^*)$ and $\calB'_u(t, X(t)) \in L(V;\,L_2(U, H))$ the Fr\'echet derivatives of $\calA$ and $\calB$ at $X(t)$, in the strong-operator-continuous sense of Assumption~\ref{ass:diff1}. Well-posedness of the fibrewise equation is Theorem~\ref{thm:Y-wellposed} under~(SC1).
\end{definition}

The operator-valued notation $Y(t,r): H \to H$ is used throughout as a convenient shorthand for the fibrewise action $v \mapsto Y(t,r)v$; $Y(t,r)$ is not in general an element of $L(H)$ almost surely. Viewed as a map $H \to L^0(\Omega; H)$, it is $L^2(\Omega; H)$-continuous under~(SC1)+(SC2) by Theorem~\ref{thm:Y-wellposed}(II). Almost-sure $L(H)$-membership would require a strictly stronger structural input, and the operator-norm moment $\E[\norm{Y(t,r)}_{L(H)}^q]$ is nowhere used below. What is used is the composite of $Y$ with the diffusion coefficient, the variational fibre $\Phi_{t,r} = Y(t,r)\calB(r,X(r))$ of~\eqref{eq:Phi-fibre}, whose Hilbert--Schmidt moments and whose fibre energy~\eqref{eq:fibre-energy} are supplied by Theorem~\ref{thm:Y-wellposed}\textup{(III)}. In particular the fibre takes values, for almost every $t$, in the space on which the linearised coefficients act, namely the tangent form domain $\calV_t$ under mechanism~(a), to which $\calB'_u$ is extended by~(SC1)(a), and $L_2(U,V)$ under mechanism~(b). This makes $\calB'_u(t, X(t))(\Phi_{t,r}u) \in L_2(U,H)$ well-defined for every $u \in U$ and almost every $t$, without requiring an $H$-extension of $\calB'_u$.

\begin{definition}[Malliavin derivative and covariance]\label{def:malliavin-cov}
For $X(t) \in \mathbb{D}^{1,2}(H)$, the Malliavin derivative $\D_r X(t)$ at time $r \in [0, t]$ admits the explicit representation in terms of the first variation process and the diffusion coefficient:
\begin{equation}\label{eq:malliavin-deriv}
  \Phi_r \;:=\; \D_r X(t) \;=\; Y(t,r)\,\calB(r, X(r)) \;\in\; L_2(U, H)
\end{equation}
(this representation is proved in Proposition~\ref{prop:D-Xt} of Section~\ref{sec:malliavin}); the identification $\Phi_r = Y(t,r)\calB(r, X(r))$ is a Hilbert--Schmidt-valued equation in $L_2(U, H)$, not a pathwise operator product on $H$ (since $Y(t,r) \in L(H)$ a.s.\ is not part of the abstract framework). The \emph{Malliavin covariance operator} is the trace-class self-adjoint positive operator
\begin{equation}\label{eq:malliavin-cov-def}
  \gamma_t \;:=\; \int_0^t \Phi_r\,\Phi_r^{*}\,\mathrm{d}r \;:\; H \to H,
\end{equation}
where $\Phi_r^{*} \in L_2(H, U)$ is the Hilbert--Schmidt adjoint of $\Phi_r$.
\end{definition}

No pathwise adjoint $Y(t,r)^{*}: H \to H$ enters~\eqref{eq:malliavin-cov-def}. Whenever equation-specific structure provides an almost surely bounded realisation $Y(t,r) \in L(H)$, the familiar factored expression involving $Y(t,r)^{*}$ may also be used.

\begin{definition}[Second Malliavin fibre]\label{def:second-fibre}
Under Assumptions~\ref{ass:LR}, \ref{ass:diff2}, and~\ref{ass:SC}, the \emph{second Malliavin fibre} is the object
\[
  \calW(\tau;s,r) \;\in\; L_2\bigl(U,\, L_2(U, H)\bigr), \qquad 0 \leq s, r \leq \tau \leq T,
\]
carrying two noise slots, a perturbation direction $a \in U$ at time $r$ and the column index $u \in U$ of the fibre $\Phi_{\tau,s}$. For $a, u \in U$ its directional fibre
\[
  \calW^{a,u}(\tau;s,r) \;:=\; \bigl[\calW(\tau;s,r)a\bigr]u
\]
is the unique variational solution, on $[\max(s,r), T]$, of the linear stochastic evolution equation
\begin{equation}\label{eq:W-equation}
\begin{aligned}
  \mathrm{d}\calW^{a,u}(\tau) + \calA'_u(\tau, X(\tau))\,\calW^{a,u}(\tau)\,\mathrm{d}\tau
  &= -\,\calA''_{uu}(\tau, X(\tau))\bigl(\Phi_{\tau,r}a,\, \Phi_{\tau,s}u\bigr)\,\mathrm{d}\tau\\
  &\quad + \calB'_u(\tau, X(\tau))\bigl(\calW^{a,u}(\tau)\bigr)\dW(\tau)\\
  &\quad + \calB''_{uu}(\tau, X(\tau))\bigl(\Phi_{\tau,r}a,\, \Phi_{\tau,s}u\bigr)\dW(\tau),
\end{aligned}
\end{equation}
with the initial condition
\begin{equation}\label{eq:W-initial}
  \calW^{a,u}(\sigma;s,r) \;=\;
  \begin{cases}
    \bigl[\calB'_u(\sigma, X(\sigma))\bigl(\Phi_{\sigma,\varsigma}\,a\bigr)\bigr]u, & r < s,\\[2pt]
    \bigl[\calB'_u(\sigma, X(\sigma))\bigl(\Phi_{\sigma,\varsigma}\,u\bigr)\bigr]a, & s < r,
  \end{cases}
\end{equation}
where $\sigma := \max(s,r)$ and $\varsigma := \min(s,r)$.
The two branches are the same expression with the two noise slots exchanged, the fibre started at the earlier time $\varsigma$ being inserted into the derivative of the diffusion coefficient at the later time $\sigma$, and their common square sum over $a$ and $u$ is the datum $\mathfrak{J}_{r,s}$ of~\eqref{eq:J-datum}.
\end{definition}

Both slots of the two second derivatives in~\eqref{eq:W-equation} are variational fibres, and by Theorem~\ref{thm:Y-wellposed}\textup{(III)} they lie in $L_2(U, \mathbb{V}_\tau)$ for almost every $\tau$, which is the space on which the second-order data of clause~(SC1$'$) act, in the $V$-typing of Assumption~\ref{ass:diff2} under mechanism~(b) and in the tangent extension of~(SC1)(a) under mechanism~(a). The equation is therefore posed against the same fibre energy that the dissipative mechanism produces, and the double square sums of its data over $(a,u)$ are the forcings~\eqref{eq:bold-forcings} and the datum $\mathfrak{J}_{r,s}$, whose finiteness and moment are clause~(SC1$'$). No direction of $H$ other than the columns of the diffusion coefficient is ever perturbed, no operator norm of the first variation is used, and the diagonal $\{r = s\}$ is Lebesgue-negligible in every integral in which $\calW$ occurs. Well-posedness at the doubly fibrewise level, together with the reconstruction of $\calW(\tau;s,r)$ from its directions, is Theorem~\ref{thm:W-wellposed}, and Lemma~\ref{lemma:DPhi} identifies the object so constructed with the Malliavin derivative of the fibre, $\calW(\tau;s,r) = \D_r\Phi_{\tau,s}$, which is where the name comes from.

\begin{definition}[Covering vector field]\label{def:covering}
For $h \in H$ satisfying clauses~\textup{(i)}--\textup{(ii)} of Assumption~\ref{ass:nondeg}, set $\tilde{h} := \gamma_t^{\dagger}h \in H$ (well-defined a.s.). The \emph{covering vector field} is the $U$-valued process
\begin{equation}\label{eq:covering-field}
  v_h(r) \;:=\; \Phi_r^{*}\,\tilde{h} \;\in\; U, \qquad r \in [0, t],
\end{equation}
where $\Phi_r := D_r X(t) \in L_2(U, H)$ and $\Phi_r^{*} \in L_2(H, U)$ its adjoint. The representation
\[
  \Phi_r \;=\; Y(t,r)\,\calB(r, X(r))
\]
of Proposition~\ref{prop:D-Xt} is understood as an identity for the Hilbert--Schmidt composite $\Phi_r \in L_2(U,H)$. In the general variational framework we retain the intrinsic formula
\[
  v_h(r) \;=\; \Phi_r^{*}\,\gamma_t^{\dagger}h
\]
and do not require a pathwise adjoint $Y(t,r)^{*}$.
\end{definition}

\subsection{Statement of the theorem}

We can now state what has been prepared for. Four statements are needed, and they stand in a definite relation to one another: the first three hold for any Malliavin-differentiable random variable on a separable Hilbert space and owe nothing to the equation, while the fourth specialises them to it. We introduce one more notational object that ties them together. Define the \emph{intrinsic random field} $w_r: H \to U$ by
\begin{equation}\label{eq:adapted-field}
  w_r(z) \;:=\; \Phi_r^{*}\,z, \qquad z \in H.
\end{equation}
The covering field is then
\[
  v_h(r) = w_r(\tilde{h}), \qquad \tilde{h} = \gamma_t^{\dagger}h .
\]
The linearity of $z \mapsto w_r(z)$ is exactly what permits the basis-wise substitution of Theorem~\ref{thm:linear-sub-D12} without any flow inversion or pathwise adjoint of $Y(t,r)$.

\begin{theorem}[The Bismut--Fomin identity on a Hilbert space]\label{thm:abstract-bismut-fomin}
Let $H$ be a separable Hilbert space, $W$ a $U$-cylindrical Wiener process on a fixed time interval $[0, T]$, and $\HW := L^2([0, T]; U)$ the associated Cameron--Martin space. Let $F \in \mathbb{D}^{1,2}(H)$ be an $H$-valued Malliavin-differentiable random variable, and write
\[
  \mathcal{T} \;:=\; DF \;\in\; L^2\!\bigl(\Omega;\,L_2(\HW, H)\bigr)
\]
for the random Hilbert--Schmidt Malliavin matrix. Suppose $h \in H$ is such that there exists a covering field
\[
  u_h \;\in\; \Dom(\delta_U) \subset L^2(\Omega; \HW)
\]
satisfying the (a.s.) covering identity $\mathcal{T} u_h = h$ in $H$. Then $\mu_F := \mathrm{Law}(F)$ is Fomin-differentiable along $h$, and the logarithmic derivative $\beta_h^{\mu_F}$ exists in $L^2(\mu_F)$ with
\begin{equation}\label{eq:abstract-bismut-fomin}
  \beta_h^{\mu_F}(F) \;=\; -\,\E\bigl[\delta_U(u_h)\,\big|\,F\bigr] \qquad \text{in } L^2(\Omega).
\end{equation}
\end{theorem}

\begin{proof}
For $\varphi \in C_b^1(H)$ and $h \in H$, the directional derivative $\partial_h \varphi(F) = \langle \nabla\varphi(F), h\rangle_H$ holds pointwise. By the assumed covering relation $h = \mathcal{T} u_h$ a.s.,
\[
  \langle \nabla\varphi(F), h\rangle_H \;=\; \langle \nabla\varphi(F),\, \mathcal{T} u_h\rangle_H \;=\; \langle \mathcal{T}^{*} \nabla\varphi(F),\, u_h\rangle_{\HW},
\]
in which the second equality uses the random adjoint $\mathcal{T}^{*} = (DF)^{*}: H \to \HW$ (a.s.\ Hilbert--Schmidt by Hilbert--Schmidt-ness of $\mathcal{T} = DF$).

We claim that for $F \in \mathbb{D}^{1,2}(H)$ and $\varphi \in C_b^1(H)$,
\begin{equation}\label{eq:H-valued-chain-rule}
  \varphi(F) \in \mathbb{D}^{1,2}(\R), \qquad D(\varphi(F)) \;=\; (DF)^{*}\nabla\varphi(F) \;=\; \mathcal{T}^{*}\nabla\varphi(F) \quad \text{in } L^2(\Omega; \HW).
\end{equation}
This is the $H$-valued analogue of \cite{nualart2006malliavin}, and is established by finite-dimensional approximation as follows. Fix any orthonormal basis $\{e_k\}_{k \geq 1}$ of $H$ and set $F_k := \langle F, e_k\rangle_H \in \mathbb{D}^{1,2}(\R)$, $\Pi_N F := \sum_{k=1}^N F_k\,e_k$. Then $\Pi_N F \in \mathbb{D}^{1,2}(H_N)$ where $H_N := \mathrm{span}\{e_1, \ldots, e_N\}$, and $\Pi_N F \to F$ in $\mathbb{D}^{1,2}(H)$ as $N \to \infty$ (by orthonormal-basis convergence in $L^2(\Omega; H)$ and the corresponding convergence of $D\Pi_N F \to DF$ in $L^2(\Omega; \HW \otimes H)$). The finite-dimensional Malliavin chain rule \cite{nualart2006malliavin} applied to $\varphi|_{H_N} \in C_b^1(H_N) = C_b^1(\R^N)$ and $\Pi_N F$ gives
\begin{align*}
  &\varphi(\Pi_N F) \in \mathbb{D}^{1,2}(\R),\\
  &D(\varphi(\Pi_N F)) \;=\; (D\Pi_N F)^{*}\nabla\varphi(\Pi_N F) \;=\; \mathcal{T}^{*}\,\Pi_N\,\nabla\varphi(\Pi_N F) \quad\text{in } L^2(\Omega; \HW).
\end{align*}
(The composition order is $(D\Pi_N F)^{*} = (\Pi_N \circ DF)^{*} = (DF)^{*}\circ\Pi_N^{*} = \mathcal{T}^{*}\Pi_N$, since $\Pi_N: H \to H$ is the self-adjoint orthogonal projection onto $H_N$. Thus $\Pi_N$ acts on $H$ from the right of $\mathcal{T}^{*}$, never on the Cameron--Martin side $\HW$.) As $N \to \infty$:
Since $\varphi$ is bounded and continuous, $\varphi(\Pi_N F) \to \varphi(F)$ in $L^2(\Omega;\R)$ by bounded convergence, and since $\nabla\varphi$ is bounded and continuous, $\nabla\varphi(\Pi_N F) \to \nabla\varphi(F)$ in $L^2(\Omega; H)$ for the same reason. Writing
\[
  \Pi_N\,\nabla\varphi(\Pi_N F) - \nabla\varphi(F) \;=\; \Pi_N\bigl(\nabla\varphi(\Pi_N F) - \nabla\varphi(F)\bigr) \;+\; (\Pi_N - I_H)\nabla\varphi(F),
\]
the first summand is handled by $\norm{\Pi_N}_{L(H)} \leq 1$ and the second by the strong-operator convergence $\Pi_N \to I_H$, so $\Pi_N\,\nabla\varphi(\Pi_N F) \to \nabla\varphi(F)$ in $L^2(\Omega;H)$ as well. Applying $\mathcal{T}^{*}$, which lies in $L^2(\Omega; L_2(H,\HW))$, and using H\"older's inequality against the bounded $H$-valued sequence, gives $\mathcal{T}^{*}\Pi_N\nabla\varphi(\Pi_N F) \to \mathcal{T}^{*}\nabla\varphi(F)$ in $L^2(\Omega;\HW)$. By closability of the Malliavin derivative \cite{nualart2006malliavin}, $\varphi(F) \in \mathbb{D}^{1,2}(\R)$ and \eqref{eq:H-valued-chain-rule} holds.

Therefore
\[
  \partial_h\varphi(F) \;=\; \langle D(\varphi(F)),\, u_h\rangle_{\HW},
\]
and by the Skorokhod duality \eqref{eq:skorokhod-duality} (which applies because $u_h \in \Dom(\delta_U)$ by hypothesis and $\varphi(F) \in L^2(\Omega) \cap \mathbb{D}^{1,2}$ by \eqref{eq:H-valued-chain-rule}),
\[
  \E\bigl[\partial_h\varphi(F)\bigr] \;=\; \E\bigl[\langle D(\varphi(F)), u_h\rangle_{\HW}\bigr] \;=\; \E\bigl[\,\varphi(F)\,\delta_U(u_h)\,\bigr].
\]
By the tower property and the fact that $\varphi(F)$ is $\sigma(F)$-measurable,
\[
  \E\bigl[\partial_h\varphi(F)\bigr] \;=\; \E\bigl[\,\varphi(F)\,\E[\delta_U(u_h)\mid F]\,\bigr],
\]
which by Definition~\ref{def:log-deriv} of the logarithmic derivative gives \eqref{eq:abstract-bismut-fomin} with $\beta_h^{\mu_F}(F) = -\E[\delta_U(u_h)\mid F]$. The $L^2(\mu_F)$ regularity follows from $\delta_U(u_h) \in L^2(\Omega)$ via the conditional Jensen inequality:
\[
  \E_{\mu_F}\bigl[\bigl|\beta_h^{\mu_F}\bigr|^2\bigr] \;=\; \E\bigl[\,\bigl|\E[\delta_U(u_h)\mid F]\bigr|^2\,\bigr] \;\leq\; \E\bigl[\delta_U(u_h)^2\bigr] \;<\;\infty. \qedhere
\]
\end{proof}

Theorem~\ref{thm:abstract-bismut-fomin} is purely Wiener-space-theoretic: it requires only Malliavin differentiability of $F$ and the existence of some covering field $u_h \in \Dom(\delta_U)$, with no assumption on $\gamma_F = \mathcal{T}\mathcal{T}^{*}$. The canonical construction of $u_h$ via the pseudoinverse is Theorem~\ref{thm:pseudoinverse-covering}.

\begin{theorem}[The pseudoinverse covering field]\label{thm:pseudoinverse-covering}
Under the setting of Theorem~\ref{thm:abstract-bismut-fomin}, define the Malliavin covariance operator
\[
  \gamma_F \;:=\; \mathcal{T} \mathcal{T}^{*} \;\in\; L^1\!\bigl(\Omega; L_1(H)\bigr).
\]
Suppose $h \in H$ satisfies the analogues of clauses~\textup{(i)}--\textup{(ii)} of Assumption~\ref{ass:nondeg}, with $\gamma_t$ replaced by $\gamma_F$:
\begin{enumerate}[leftmargin=3em]
  \item[\textup{(i$_F$)}]\hyplabel{hyp:i-F}{\textup{(i$_F$)}} $h \in \Ran(\gamma_F)$ almost surely;
  \item[\textup{(ii$_F$)}]\hyplabel{hyp:ii-F}{\textup{(ii$_F$)}} there exist exponents $(q, p^{*}) \in [2,\infty]^2$ with $1/q + 1/p^{*} \leq 1/2$ such that $\E[\norm{\gamma_F^{\dagger}h}_H^q] < \infty$ and $\mathcal{T} \in \mathbb{D}^{1, p^{*}}(L_2(\HW, H))$.
\end{enumerate}
Set
\begin{equation}\label{eq:abstract-canonical-uh}
  u_h \;:=\; \mathcal{T}^{*}\,\gamma_F^{\dagger} h \;\in\; L^2(\Omega; \HW),
\end{equation}
and define the Tikhonov approximations $u_h^{(\eps)} := \mathcal{T}^{*}(\gamma_F + \eps I)^{-1}h$ for $\eps > 0$. Then:
\begin{enumerate}[label=\textup{(\arabic*)}, leftmargin=3em]
  \item \textup{(Algebraic covering identity)} The covering identity $\mathcal{T} u_h = h$ holds almost surely.
  \item \textup{(Algebraic $L^2$-approximation by the Tikhonov family)} For each $\eps > 0$, $u_h^{(\eps)} \in L^2(\Omega; \HW)$ with the deterministic resolvent bound $\norm{u_h^{(\eps)}}_{\HW} \leq \eps^{-1/2}\norm{h}_H$ a.s., and $u_h^{(\eps)} \to u_h$ in $L^2(\Omega; \HW)$ as $\eps \downarrow 0$.
 \item \textup{(Skorokhod-domain membership at fixed $\eps$)} The membership $u_h^{(\eps)} \in \mathbb{D}^{1,2}(\HW) \subset \Dom(\delta_U)$, and a fortiori $u_h \in \Dom(\delta_U)$ via the Cauchy property of $\{\delta_U(u_h^{(\eps)})\}_{\eps > 0}$ and the closedness of $\delta_U$, are the content of Theorem~\ref{thm:tikhonov-trace} below; they require the additional fixed-Tikhonov substitution regularity~\ref{hyp:iii-F-prime} and the Cameron--Martin compatibility~\ref{hyp:iv-F} of that theorem. Granted~\ref{hyp:iii-F-prime}--\ref{hyp:iv-F}, Theorem~\ref{thm:abstract-bismut-fomin} applies and yields the Bismut--Fomin formula~\eqref{eq:abstract-bismut-fomin} with the canonical $u_h$ of~\eqref{eq:abstract-canonical-uh}.
\end{enumerate}
\end{theorem}

We record the algebraic content of Theorem~\ref{thm:pseudoinverse-covering} and the reason why~\ref{hyp:iii-F-prime} is an additional hypothesis. (i) Algebraic content. Theorem~\ref{thm:pseudoinverse-covering} is pathwise Hilbert-space algebra: the Moore--Penrose construction $u_h = \mathcal{T}^{*}\gamma_F^{\dagger}h$ executed pointwise in $\omega$, with the deterministic spectral identity $\gamma_F\gamma_F^{\dagger}|_{\Ran(\gamma_F)} = \mathrm{id}$ providing the algebraic covering identity, promoted to an $L^2$-approximation statement by integrating the resulting pathwise spectral identities against the $L^q$-moment-bound hypothesis on $\gamma_F^{\dagger}h$.

(ii) Why~\ref{hyp:iii-F-prime} cannot be derived from the algebraic content alone. The product-rule expansion of $\D u_h^{(\eps)}$ produces the term $\mathcal{T}^{*}\,\D G^{(\eps)}$ (with $G^{(\eps)} := (\gamma_F+\eps I)^{-1}h$), whose $L^2(\Omega)$-control via $\E[\|\mathcal{T}\|^2 \|\D G^{(\eps)}\|^2]$ requires a H\"older-pair joint integrability that $\mathcal{T} \in \mathbb{D}^{1, p^*}$ alone does not provide, not even when $\gamma_F \in \mathbb{D}^{1,2}(L_2(H,H))$ holds; the H\"older inequality $1/(p^*/2) + 1/1 \leq 1$ closes only at $p^* = \infty$. Hence~\ref{hyp:iii-F-prime} is genuinely an additional hypothesis on top of the algebraic content of Theorem~\ref{thm:pseudoinverse-covering}; the natural sufficient conditions (deterministic-covariance edge, H\"older-triangle joint regularity, or bounded Malliavin matrix) are spelled out in Remark~\ref{rem:iii-F-status}.

\begin{proof}
The range condition \ref{hyp:i-F} gives the identity $\gamma_F\gamma_F^{\dagger}|_{\Ran(\gamma_F)} = \mathrm{id}_{\Ran(\gamma_F)}$ (the standard Moore--Penrose property), so
\[
  \mathcal{T} u_h \;=\; \mathcal{T} \mathcal{T}^{*}\gamma_F^{\dagger}h \;=\; \gamma_F\gamma_F^{\dagger}h \;=\; h \quad\text{a.s.}
\]

Expanding,
\[
  \norm{u_h}^2_{\HW} \;=\; \langle \mathcal{T}^{*}\gamma_F^{\dagger}h,\, \mathcal{T}^{*}\gamma_F^{\dagger}h\rangle_{\HW} \;=\; \langle \mathcal{T} \mathcal{T}^{*}\gamma_F^{\dagger}h,\,\gamma_F^{\dagger}h\rangle_H \;=\; \langle h,\,\gamma_F^{\dagger}h\rangle_H \;\leq\; \norm{h}_H\,\norm{\gamma_F^{\dagger}h}_H,
\]
hence $\E\norm{u_h}^2_{\HW} \leq \norm{h}_H\,\E[\norm{\gamma_F^{\dagger}h}_H] \leq \norm{h}_H\,(\E[\norm{\gamma_F^{\dagger}h}_H^q])^{1/q} < \infty$ under \ref{hyp:ii-F} at any operating point on the H\"older triangle ($q \geq 2$).

For each $\eps > 0$, the resolvent $(\gamma_F + \eps I)^{-1}$ is a bounded operator on $H$ with deterministic norm bound $\norm{(\gamma_F + \eps I)^{-1}}_{L(H)} \leq \eps^{-1}$ a.s., hence the composition
\[
  u_h^{(\eps)} \;=\; \mathcal{T}^{*}\,(\gamma_F + \eps I)^{-1}\,h
\]
is well-defined as an $\HW$-valued random variable with the pathwise estimate
\begin{align*}
  \norm{u_h^{(\eps)}}^2_{\HW}
  &\;=\; \langle (\gamma_F+\eps I)^{-1}h,\, \mathcal{T} \mathcal{T}^{*}(\gamma_F+\eps I)^{-1}h\rangle_H \\
  &\;=\; \langle (\gamma_F+\eps I)^{-1}h,\,\gamma_F(\gamma_F+\eps I)^{-1}h\rangle_H \\
  &\;\leq\; \norm{(\gamma_F+\eps I)^{-1}h}_H \cdot \norm{\gamma_F(\gamma_F+\eps I)^{-1}h}_H
   \;\leq\; \eps^{-1}\,\norm{h}_H^2
\end{align*}
(using $\norm{\gamma_F(\gamma_F+\eps I)^{-1}}_{L(H)} \leq 1$ from spectral calculus on positive operators). In particular, $u_h^{(\eps)} \in L^2(\Omega; \HW)$ with bound polynomial in $\eps^{-1}$.

The $L^2(\Omega; \HW)$-convergence $u_h^{(\eps)} \to u_h$ as $\eps \downarrow 0$ follows from the explicit spectral identity (in any measurable selection of eigenbasis $\{(\lambda_k, e_k)\}$ of $\gamma_F$, cf.\ Remark~\ref{rem:CM-compat-spectral}):
\begin{equation}\label{eq:Tikhonov-L2-spectral}
\begin{split}
  \norm{u_h - u_h^{(\eps)}}^2_{\HW}
  &\;=\; \langle [\gamma_F^{\dagger} - (\gamma_F + \eps I)^{-1}]h,\,\gamma_F[\gamma_F^{\dagger} - (\gamma_F + \eps I)^{-1}]h\rangle_H \\
  &\;=\; \!\!\sum_{\lambda_k > 0}\!\frac{\eps^2\,h_k^2}{\lambda_k(\lambda_k + \eps)^2}\;\leq\;\!\!\sum_{\lambda_k > 0}\!\frac{h_k^2}{\lambda_k} \;=\; \ip{h}{\gamma_F^{\dagger}h}_H \;=\; \norm{u_h}^2_{\HW},
\end{split}
\end{equation}
where the summand is $\lambda_k d_k^2$ with $d_k = \eps h_k/(\lambda_k(\lambda_k+\eps))$, so that each term is dominated by $h_k^2/\lambda_k$ via $\eps^2/(\lambda_k+\eps)^2 \leq 1$; the dominating random variable $\ip{h}{\gamma_F^{\dagger}h}_H \leq \norm{h}_H\norm{\gamma_F^{\dagger}h}_H$ is integrable by the moment hypothesis~\ref{hyp:ii-F}. Since $h \in \Ran(\gamma_F)$, its component in $\ker(\gamma_F)$ vanishes. Hence, after applying the difference to $h$,
\[
  \bigl[\gamma_F^{\dagger} - (\gamma_F + \eps I)^{-1}\bigr]h \;=\; \sum_{\lambda_k > 0}\frac{\eps\,h_k}{\lambda_k(\lambda_k+\eps)}\,e_k,
\]
and consequently
\[
  \gamma_F\bigl[\gamma_F^{\dagger} - (\gamma_F + \eps I)^{-1}\bigr]h \;=\; \sum_{\lambda_k > 0}\frac{\eps\,h_k}{\lambda_k+\eps}\,e_k .
\] The right-hand side of \eqref{eq:Tikhonov-L2-spectral} is in $L^1(\Omega)$ under \ref{hyp:ii-F} (which gives $\E\norm{\gamma_F^{\dagger}h}_H^q < \infty$ at $q \geq 2$, hence $\E\norm{\gamma_F^{\dagger}h}_H^2 < \infty$). Each summand of the spectral sum tends to $0$ as $\eps \downarrow 0$ pointwise (in $\omega$), dominated summand-wise by $h_k^2/\lambda_k$, so by dominated convergence in $k$ and again in $\omega$, $\E\norm{u_h - u_h^{(\eps)}}^2_{\HW} \to 0$. This entire step is purely algebraic / $L^2$-level, requiring no Malliavin regularity beyond $\mathcal{T} \in \mathbb{D}^{1,p^*}$ (which is not even invoked here, the spectral computation uses only the Hilbert--Schmidt structure of $\mathcal{T}$ via $\gamma_F = \mathcal{T}\mathcal{T}^*$, and the resulting bound is on the $\HW$-norm, not on Malliavin derivatives).

The $\mathbb{D}^{1,2}$-membership of $u_h^{(\eps)}$ at fixed $\eps$, the Cauchy property of $\{\delta_U(u_h^{(\eps)})\}_{\eps > 0}$ in $L^2(\Omega)$, and (consequently) $u_h \in \Dom(\delta_U)$ with $\delta_U(u_h)$ as the $L^2(\Omega)$-limit, all require the additional Cameron--Martin compatibility hypothesis~\ref{hyp:iv-F} and the fixed-Tikhonov substitution regularity hypothesis~\ref{hyp:iii-F-prime} of Theorem~\ref{thm:tikhonov-trace}; see the proof there. Without~\ref{hyp:iii-F-prime}--\ref{hyp:iv-F}, $u_h$ is constructed only as an $L^2(\Omega; \HW)$-element with the algebraic covering identity of clause~(1), and Theorem~\ref{thm:abstract-bismut-fomin} cannot yet be applied to it.
\end{proof}

\medskip

The following proposition isolates the abstract closure mechanism that converts fixed-$\eps$ Skorokhod-domain membership of the Tikhonov regularised covering field, together with two $L^2(\Omega)$-convergences, into the same conclusion for the unregularised pseudoinverse covering field. It is the architectural backbone of Theorem~\ref{thm:tikhonov-trace} below and of the SPDE specialisation in Lemma~\ref{lemma:D-gamma-inv-h}; making it explicit clarifies the dependency on closedness of $\delta_U$ as a closed unbounded operator, separately from the explicit computation of main and correction terms.

\begin{proposition}[Closed-graph Tikhonov principle]\label{prop:closed-graph-tikhonov}
Let $F \in \mathbb{D}^{1,2}(H)$, $\mathcal{T} = DF \in L^2(\Omega; L_2(\HW, H))$, and $\gamma_F = \mathcal{T}\mathcal{T}^{*}$. Let $h \in H$ be a deterministic direction satisfying the range and moment conditions $h \in \Ran(\gamma_F)$ a.s.\ and $\gamma_F^{\dagger}h \in L^q(\Omega; H)$ for some $q \geq 2$, and consider the Tikhonov-regularised covering fields
\[
  u_h^{(\eps)} \;:=\; \mathcal{T}^{*}\,(\gamma_F + \eps I)^{-1}\,h, \qquad u_h \;:=\; \mathcal{T}^{*}\,\gamma_F^{\dagger}\,h, \qquad \eps > 0,
\]
and assume:
\begin{enumerate}[label=\textup{(\alph*)}, leftmargin=2.2em, itemsep=2pt]
 \item At each fixed $\eps$ the field lies in the domain of the divergence: $u_h^{(\eps)} \in \mathbb{D}^{1,2}(\HW) \subset \Dom(\delta_U)$ for every $\eps > 0$.
 \item The main term converges: $u_h^{(\eps)} \to u_h$ in $L^2(\Omega; \HW)$ as $\eps \downarrow 0$.
 \item The correction term converges as well, since the family $\{\delta_U(u_h^{(\eps)})\}_{\eps > 0}$ is Cauchy in $L^2(\Omega)$ as $\eps \downarrow 0$.
\end{enumerate}
Then $u_h \in \Dom(\delta_U)$, and
\begin{equation}\label{eq:closed-graph-tikhonov}
  \delta_U(u_h) \;=\; \lim_{\eps \downarrow 0}\,\delta_U(u_h^{(\eps)}) \qquad \text{in } L^2(\Omega).
\end{equation}
\end{proposition}

\begin{proof}
The Skorokhod divergence $\delta_U: \Dom(\delta_U) \subset L^2(\Omega; \HW) \to L^2(\Omega)$ is the adjoint of the densely defined Malliavin derivative $\D: \mathbb{D}^{1,2}(\R) \to L^2(\Omega; \HW)$, and adjoints of densely defined operators are closed \cite{nualart2006malliavin}. By~(b), $u_h^{(\eps)} \to u_h$ in $L^2(\Omega; \HW)$; by~(c), $\delta_U(u_h^{(\eps)})$ converges in $L^2(\Omega)$ to some limit $\eta$. Since each $u_h^{(\eps)} \in \Dom(\delta_U)$ by~(a), closedness of $\delta_U$ gives $u_h \in \Dom(\delta_U)$ with $\delta_U(u_h) = \eta$, which is~\eqref{eq:closed-graph-tikhonov}.
\end{proof}

\begin{remark}[How the principle is used]\label{rem:closed-graph-tikhonov-use}
Theorem~\ref{thm:tikhonov-trace} below verifies hypotheses~(a)--(c) of Proposition~\ref{prop:closed-graph-tikhonov} under the explicit moment and substitution-regularity conditions \ref{hyp:i-F}--\ref{hyp:iv-F}, and simultaneously identifies the limit $\delta_U(u_h)$ as the sum $M_h + C_h$ of a main term and a correction term with explicit spectral representations. In the SPDE specialisation, Lemma~\ref{lemma:D-gamma-inv-h} verifies the same hypotheses for $(F, h) = (X(t), h)$ under Assumptions~\ref{ass:LR}--\ref{ass:nondeg}, with the fixed-$\eps$ membership~(a) supplied by Theorem~\ref{thm:linear-sub-D12} (basis-wise linear-field substitution), the main-term convergence~(b) by spectral $L^2$-control of the Tikhonov resolvent applied to $h$, and the correction-term convergence~(c) by the scalar convergence-and-domination hypothesis in the Cameron--Martin compatibility clause~\ref{ass:nondeg}(iii), with the spectral condition of Remark~\ref{rem:CM-compat-spectral} providing one sufficient verification. At no stage is the operator-valued resolvent $(\gamma_F + \eps I)^{-1} \in L(H)$ required to live in any operator-valued Malliavin space; the vector-level Malliavin regularity of Lemma~\ref{lemma:D-inverse}, applied at the fixed deterministic direction $h$, is the only Malliavin input needed for the resolvent itself.
\end{remark}

\begin{theorem}[The trace formula]\label{thm:tikhonov-trace}
In the setting of Theorem~\ref{thm:pseudoinverse-covering}, let $\Phi: \Omega \times [0,T] \to L_2(U, H)$ be the time-fibre realisation of $\mathcal{T} = DF$ obtained from the canonical isometry
\[
  L_2\!\bigl(L^2([0,T]; U),\; H\bigr) \;\simeq\; L^2\!\bigl([0,T]; L_2(U, H)\bigr), \qquad \mathcal{T} u \;=\; \int_0^T \Phi_r\,u(r)\,\mathrm{d}r,
\]
so that $\norm{\mathcal{T}}^2_{L_2(\HW, H)} = \int_0^T \norm{\Phi_r}^2_{L_2(U, H)}\,\mathrm{d}r$ a.s. Suppose, in addition to \ref{hyp:i-F}--\ref{hyp:ii-F}:
\begin{enumerate}[leftmargin=3em]
  \item[\textup{(iii$_F'$)}]\hyplabel{hyp:iii-F-prime}{\textup{(iii$_F'$)}} \textup{(Fixed-Tikhonov substitution regularity.)} For every $\eps > 0$,
  \begin{equation}\label{eq:uheps-D12}
    u_h^{(\eps)} \;:=\; \mathcal{T}^{*}\,(\gamma_F + \eps I)^{-1}\,h \;\in\; \mathbb{D}^{1,2}(\HW) \;\subset\; \Dom(\delta_U),
  \end{equation}
  and the corresponding fixed-$\eps$ substitution identity is well-defined. Alternatively, for the purpose of the substitution identity, the direct assumption~\eqref{eq:uheps-D12} may be replaced by the following basis-wise sufficient conditions:
  \begin{equation}\label{eq:Geps-D12}
    G^{(\eps)} := (\gamma_F + \eps I)^{-1}h \;\in\; \mathbb{D}^{1,2}(H), \qquad \eps > 0,
  \end{equation}
  together with the H\"older-pair joint integrability~\eqref{eq:linear-sub-Holder} of $(\mathcal{T}, G^{(\eps)})$ (the basis-wise route via Theorem~\ref{thm:linear-sub-D12}). Sufficient conditions for either form, and the relation between them, are catalogued in Remark~\ref{rem:iii-F-status}.
  \item[\textup{(iv$_F$)}]\hyplabel{hyp:iv-F}{\textup{(iv$_F$)}} \textup{(Cameron--Martin compatibility)} The covariance, viewed as an $L_2(H,H)$-valued random variable, belongs to $\mathbb{D}^{1,2}(L_2(H,H))$, and the Tikhonov-regularised trace integrand
  \[
    \mathcal{I}^{(\eps)}_r(F) \;:=\; \Tr_U\!\Bigl[\,\Phi_r^{*}\,(\gamma_F + \eps I)^{-1}\,[D_r\gamma_F]\,(\gamma_F + \eps I)^{-1}\,h\,\Bigr], \qquad \eps > 0,
  \]
  converges a.s.\ for a.e.\ $r$, as $\eps \downarrow 0$, to a measurable limit $\mathcal{I}^{(0)}_r(F)$, and admits a measurable dominator $\Psi_r$, exactly as in \eqref{eq:CM-compat-limit}--\eqref{eq:CM-compat-dominator} with $[0,t]$ replaced by $[0,T]$ and $\E[(\int_0^T\Psi_r\,\mathrm{d}r)^2] < \infty$.
\end{enumerate}
Then $u_h \in \Dom(\delta_U)$, and the Skorokhod integral admits the explicit decomposition
\begin{equation}\label{eq:abstract-skorokhod-decomp}
  \delta_U(u_h) \;=\; M_h \;+\; C_h,
\end{equation}
where the \emph{main term} $M_h$ is the Hilbert--Schmidt kernel evaluation at the random point $\gamma_F^{\dagger}h$: by Lemma~\ref{lemma:HS-kernel}, the map $z \mapsto \delta_U(\mathcal{T}^{*}z)$ is Hilbert--Schmidt from $H$ to $L^2(\Omega)$ with kernel $\mathcal{K} \in L^2(\Omega; H)$, and
\begin{equation}\label{eq:abstract-main-term}
  M_h \;:=\; \langle \gamma_F^{\dagger}h,\,\mathcal{K}\rangle_H \;=\; \lim_{\eps\downarrow 0}\,\langle (\gamma_F+\eps I)^{-1}h,\,\mathcal{K}\rangle_H \quad \text{in } L^2(\Omega),
\end{equation}
and the \emph{correction term} $C_h$ is the trace-integral limit
\begin{equation}\label{eq:abstract-correction}
  C_h \;:=\; \int_0^T \mathcal{I}^{(0)}_r(F)\,\mathrm{d}r \;=\; \lim_{\eps \downarrow 0}\,\int_0^T \mathcal{I}^{(\eps)}_r(F)\,\mathrm{d}r \quad \text{in } L^2(\Omega),
\end{equation}
existing under~\ref{hyp:iii-F-prime}--\ref{hyp:iv-F}, both limits being independent of the regularising sequence $\eps_n \downarrow 0$.
\end{theorem}

\begin{proof}
The proof is the abstract version of the proof of Lemma~\ref{lemma:D-gamma-inv-h}: at fixed $\eps > 0$, $u_h^{(\eps)} := \mathcal{T}^{*}(\gamma_F + \eps I)^{-1}h \in \mathbb{D}^{1,2}(\HW)$ by hypothesis~\ref{hyp:iii-F-prime}; the deterministic resolvent bound $\norm{(\gamma_F + \eps I)^{-1}}_{L(H)} \leq \eps^{-1}$ enters only through the $L^\infty(\Omega; H)$-control of $G^{(\eps)} := (\gamma_F + \eps I)^{-1}h$, but does not by itself imply $u_h^{(\eps)} \in \mathbb{D}^{1,2}(\HW)$ — that membership is genuinely the substitution-regularity content of \ref{hyp:iii-F-prime}, as the product-rule expansion $\D u_h^{(\eps)} = (\D \mathcal{T}^{*})G^{(\eps)} + \mathcal{T}^{*}\,\D G^{(\eps)}$ requires the H\"older-pair joint integrability of $(\|\mathcal{T}\|, \|\D G^{(\eps)}\|)$ that invertibility alone does not provide. Granted \ref{hyp:iii-F-prime}, the basis-wise substitution theorem (Theorem~\ref{thm:linear-sub-D12}), applied to the field $z \mapsto \mathcal{T}^{*}z$, which is linear in $z$, gives
\begin{equation}\label{eq:abstract-Nualart-split}
  \delta_U(u_h^{(\eps)}) \;=\; M_h^{(\eps)} \;+\; C_h^{(\eps)}, \qquad M_h^{(\eps)} := \bigl[\delta_U(\mathcal{T}^{*}z)\bigr]_{z=(\gamma_F+\eps I)^{-1}h}, \quad C_h^{(\eps)} := \int_0^T \mathcal{I}^{(\eps)}_r\,\mathrm{d}r.
\end{equation}
Here $M_h^{(\eps)}$ is the frozen-$z$ divergence (compute $\delta_U(\mathcal{T}^{*}z)$ as a $z$-dependent functional on $H$, then substitute $z = (\gamma_F+\eps I)^{-1}h$), which equals $\langle (\gamma_F+\eps I)^{-1}h,\mathcal{K}\rangle_H$ via the Hilbert--Schmidt (HS) kernel realisation of Lemma~\ref{lemma:HS-kernel}. The Hilbert--Schmidt-kernel route handles $M_h^{(\eps)} \to M_h$ in $L^2(\Omega)$ via~\ref{hyp:ii-F}: writing $M_h^{(\eps)} - M_h = \langle [(\gamma_F+\eps I)^{-1} - \gamma_F^{\dagger}]h,\mathcal{K}\rangle_H$, the Cauchy--Schwarz bound gives $\E|M_h^{(\eps)} - M_h|^2 \leq \E[\|[(\gamma_F+\eps I)^{-1} - \gamma_F^{\dagger}]h\|^2 \,\|\mathcal{K}\|^2]$ which tends to $0$ by H\"older at the H\"older-triangle exponent pair $(q, p^*)$ ($\mathcal{K} \in L^{p^*}(\Omega; H)$ from $\mathcal{T} \in \mathbb{D}^{1,p^*}$ via the $L^{p^*}$-realisation of the HS-kernel, and, uniformly in $\eps$, by the spectral bound
\[
  \bigl\|[(\gamma_F+\eps I)^{-1} - \gamma_F^{\dagger}]h\bigr\|_H^2
  \;=\; \sum_{\lambda_k > 0}\frac{\eps^2 h_k^2}{\lambda_k^2(\lambda_k+\eps)^2}
  \;\leq\; \sum_{\lambda_k > 0}\frac{h_k^2}{\lambda_k^2}
  \;=\; \|\gamma_F^{\dagger}h\|_H^2 ,
\]
with pointwise convergence to $0$ by spectral dominated convergence). The pointwise convergence and the common dominator of~\ref{hyp:iv-F} give $C_h^{(\eps)} \to C_h$ in $L^2(\Omega)$ by dominated convergence. The closedness of $\delta_U$ then yields $u_h \in \Dom(\delta_U)$ and $\delta_U(u_h) = M_h + C_h$. The full execution is in Lemma~\ref{lemma:D-gamma-inv-h} in the SPDE setting; the abstract statement requires no additional ingredient beyond the abstract Malliavin calculus and the substitution hypothesis~\ref{hyp:iii-F-prime}, so the same proof goes through verbatim.
\end{proof}

\begin{remark}[Status of \ref{hyp:iii-F-prime}]\label{rem:iii-F-status}
Hypothesis~\eqref{eq:uheps-D12} is above the algebraic content of Theorem~\ref{thm:pseudoinverse-covering}. Sufficient conditions, in increasing order of strength:
\begin{enumerate}[label=\textup{(iii$_F'$.\alph*)}, leftmargin=4em]
 \item\label{hyp:iii-F-prime-a} \emph{Deterministic covariance} ($\D\gamma_F \equiv 0$). Covers linear SPDEs with state-independent diffusion.
 \item\label{hyp:iii-F-prime-b} \emph{H\"older-triangle joint regularity}: $\gamma_F \in \mathbb{D}^{1,r}(L_2(H,H))$ with $1/p^* + 1/r \leq 1/2$.
 \item\label{hyp:iii-F-prime-c} \emph{Bounded Malliavin matrix}: $\mathcal{T} \in L^\infty(\Omega; L_2(\HW, H))$.
\end{enumerate}
The primary mechanism is the basis-wise route via Theorem~\ref{thm:linear-sub-D12}, which requires only $G^{(\eps)} \in \mathbb{D}^{1,2}(H)$ together with~\eqref{eq:linear-sub-Holder}, closing at $p^* \geq 6$ on the edge choice $p^* = 2q/(q-2)$ of the H\"older triangle (rather than $p^* \geq 8$ via Nualart's $\mathbb{D}^{1,4}_{\mathrm{loc}}$ theorem).
\end{remark}

Hypotheses~\ref{hyp:iii-F-prime} and~\ref{hyp:iv-F} encode the minimal regularity above the algebraic content of Theorem~\ref{thm:pseudoinverse-covering} that licenses the Tikhonov, Hilbert--Schmidt-kernel, and scalar trace-convergence machinery. They are strictly weaker than pointwise Malliavin regularity of $v_h$.

Specialising to the variational SPDE \eqref{eq:main-SPDE}, the abstract objects are identified through the variation processes.

\begin{theorem}[Logarithmic derivative of a variational SPDE]\label{thm:main}
Let $X(t)$ be the unique variational solution of \eqref{eq:SPDE} under Assumption~\ref{ass:LR}, denote by $\mu_t$ its law on $H$, and let $h \in H$. Set $\tilde{h} := \gamma_t^{\dagger} h$ and $v_h(r) := \Phi_r^{*}\tilde{h} \in U$.

\smallskip\noindent\emph{Part I (Bismut formula).} Under Assumptions~\ref{ass:LR},~\ref{ass:diff2},~\ref{ass:SC} together with clauses~\textup{(i)}--\textup{(ii)} of Assumption~\ref{ass:nondeg} and the domain hypothesis $v_h \in \Dom(\delta_U)$, the logarithmic derivative of $\mu_t$ along $h$ exists in $L^2(\mu_t)$ and is given by
\begin{equation}\label{eq:score-main}
  \beta_h(X(t)) \;=\; -\,\E\bigl[\delta_U(v_h) \,\big|\, X(t)\bigr].
\end{equation}

\smallskip\noindent\emph{Part II (Explicit decomposition).} If, additionally, the exponent $q$ in Assumption~\ref{ass:nondeg}(ii) is chosen in $(2, \infty]$, and clause~\textup{(iii)} of Assumption~\ref{ass:nondeg}, Assumption~\ref{ass:SC-raised} at target $q^* = q/(q-2)$ (with $q^* = 1$ when $q = \infty$), and the fixed-Tikhonov substitution regularity~\ref{hyp:iii-F-prime} of Theorem~\ref{thm:tikhonov-trace} hold for $(\Phi, \gamma_t, h)$, then $\delta_U(v_h)$ admits the decomposition
\begin{equation}\label{eq:skorokhod-decomp}
  \delta_U(v_h) \;=\; \underbrace{\bigl[\delta_U\bigl(w(\cdot)(z)\bigr)\bigr]_{z = \tilde{h}}}_{\text{(I) Main term } M_h} \;+\; \underbrace{C_h\vphantom{\bigl[\bigr]}}_{\text{(II) Correction term}},
\end{equation}
where $w_r(z) := \Phi_r^{*}\,z$ is the intrinsic random field of \eqref{eq:adapted-field}, and the correction term $C_h$ is the $L^2(\Omega)$-convergent Tikhonov limit
\begin{equation}\label{eq:correction-tikhonov}
  C_h \;:=\; \lim_{\eps \downarrow 0}\,\int_0^t \Tr_U\Bigl[\,\Phi_r^{*}\,(\gamma_t + \eps I)^{-1}\,[\D_r\gamma_t]\,(\gamma_t + \eps I)^{-1}\,h\,\Bigr]\,\mathrm{d}r.
\end{equation}
Under the spectral sufficient condition of Remark~\ref{rem:CM-compat-spectral}, this limit has the absolutely convergent spectral form, in any measurable eigenbasis of $\gamma_t$,
\begin{equation}\label{eq:correction-spectral}
  C_h \;=\; \int_0^t \!\!\!\!\sum_{\substack{j \geq 1\\ k,l\,:\,\lambda_k,\,\lambda_l > 0}}\!\!\!\!\frac{h_k\,\bigl\langle[\D_r\gamma_t](f_j)\,e_k,\,e_l\bigr\rangle_H\,\bigl\langle e_l,\,\Phi_r f_j\bigr\rangle_H}{\lambda_k\,\lambda_l}\,\mathrm{d}r.
\end{equation}
Without that condition only the basis-free definition~\eqref{eq:correction-tikhonov} is retained, denoted symbolically by
\begin{equation}\label{eq:correction-symbolic}
  C_h \;=\; \int_0^t \Tr_U\bigl[\,\Phi_r^{*}\,\gamma_t^{\dagger}\,[\D_r\gamma_t]\,\gamma_t^{\dagger}\,h\,\bigr]\,\mathrm{d}r ,
\end{equation}
where the right-hand side stands for the $L^2(\Omega)$-limit~\eqref{eq:correction-tikhonov} rather than for a pointwise composition of unbounded pseudoinverses, and agrees with~\eqref{eq:correction-spectral} whenever the spectral condition holds.
\end{theorem}

\begin{remark}[Sufficient conditions for $v_h \in \Dom(\delta_U)$]\label{rem:domain-vh-suffic}
By Theorem~\ref{thm:covering-regularity}, the domain hypothesis in Part~I holds in two regimes: (a) linear drift with state-independent diffusion (no further hypothesis); (b) the nonlinear regime, with $q \in (2,\infty]$ in Assumption~\ref{ass:nondeg}(ii), under clause~\textup{(iii)} of that assumption, Assumption~\ref{ass:SC-raised} at $q^* = q/(q-2)$ ($q^* = 1$ when $q = \infty$), and the substitution regularity~\ref{hyp:iii-F-prime}. The latter is automatic from the moment chain at $2 < q \leq 3$ via the basis-wise route of Theorem~\ref{thm:linear-sub-D12}.
\end{remark}

The decomposition~\eqref{eq:skorokhod-decomp} has two terms. The first is the Skorokhod integral of the field $w_r(z) = \Phi_r^{*}\,z$ evaluated at the random point $z = \gamma_t^{\dagger}h$; it reduces to an It\^o integral precisely when the frozen-$z$ field $r \mapsto \Phi_r^{*}z$ is adapted, in particular for linear drift with state-independent diffusion (Corollary~\ref{cor:linear-infinite}), while for nonlinear drift the transfer factor $Y(t,r)$ anticipates the driving path even when $\calB'_u = 0$ (Corollary~\ref{cor:semilinear-additive}). The second arises from the Hilbert-space linear substitution identity of Theorem~\ref{thm:linear-sub-D12}, the linear-field analogue of the classical finite-dimensional substitution formula of Theorem~\ref{thm:nualart-sub}; it accounts for the non-commutativity of the Skorokhod integral with evaluation at a random point. It involves the Malliavin derivative $\D_r\gamma_t$, which by Proposition~\ref{prop:D-gamma} is assembled from the second Malliavin fibre $\D_r\Phi_{t,s} = \calW(t;s,r)$ alone, through the single symmetric expression~\eqref{eq:D-gamma-formula}; it vanishes when $\calA''_{uu} = \calB''_{uu} = 0$ and $\calB'_u = 0$. All Malliavin derivatives appearing in~\eqref{eq:skorokhod-decomp}--\eqref{eq:correction-symbolic} are expressed entirely in terms of $\Phi$ and $\D\Phi = \calW$ (equivalently, in terms of $\calA'_u$, $\calA''_{uu}$, $\calB$, $\calB'_u$, $\calB''_{uu}$ and the two fibres, via the variational identifications), together with the Tikhonov-regularised pseudoinverse of $\gamma_t$; no residual abstract Malliavin derivatives remain. Remark~\ref{rem:trace-interp} unfolds the symbolic form~\eqref{eq:correction-symbolic}.

Theorem~\ref{thm:main} is the SPDE specialisation of Theorems~\ref{thm:abstract-bismut-fomin}--\ref{thm:tikhonov-trace}, obtained by taking $F = X(t)$, $\mathcal{T} = DF$, and identifying
\[
  \Phi_r \;=\; D_r X(t) \;=\; Y(t,r)\,\calB(r, X(r)) \quad \text{in } L_2(U, H), \qquad r \in [0,t]
\]
(Proposition~\ref{prop:D-Xt}); the Malliavin covariance is
\[
  \gamma_t \;=\; \int_0^t \Phi_r\Phi_r^{*}\,\mathrm{d}r
\]
(cf.~\eqref{eq:malliavin-cov-def}). The abstract range, moment, and substitution hypotheses are matched to the SPDE objects through these identifications. The regularity $\Phi \in \mathbb{D}^{1,p^{*}}$ is part of Assumption~\ref{ass:nondeg}(ii), with Assumption~\ref{ass:SC-raised} providing an SPDE-level sufficient route when invoked. When the available $\Phi$-regularity is at least four, Proposition~\ref{prop:D-gamma} supplies the required $L_2(H,H)$-valued differentiability of $\gamma_t$; otherwise that covariance regularity is imposed directly in the Cameron--Martin compatibility clause~\textup{(iii)}. The fixed-Tikhonov substitution regularity is imposed through~\ref{hyp:iii-F-prime} and handled by Theorem~\ref{thm:linear-sub-D12}.

\begin{remark}[Unfolded form of the correction term in \eqref{eq:correction-symbolic}]\label{rem:trace-interp}
The symbolic shorthand in~\eqref{eq:correction-symbolic}, which uses the notation $\gamma_t^{\dagger}[\D_r\gamma_t]\gamma_t^{\dagger}h$ on the right-hand side, unfolds as follows. For each fixed $r \in [0,t]$, the Malliavin derivative $\D_r\gamma_t$ is an $L_2(U, L_2(H,H))$-valued object, so for each $f \in U$, $[\D_r\gamma_t](f) \in L_2(H,H)$ is a Hilbert--Schmidt operator on $H$; consequently $f \mapsto [\D_r\gamma_t](f)(\gamma_t + \eps I)^{-1}h$ is Hilbert--Schmidt from $U$ to $H$, and composing it on the left with the Hilbert--Schmidt map $\Phi_r^{*}(\gamma_t + \eps I)^{-1}: H \to U$ produces a trace-class operator on $U$, so that the trace below is well defined. Writing $\Phi_r := Y(t,r)\calB(r,X(r)) \in L_2(U,H)$ as in the theorem statement, the integrand of~\eqref{eq:correction-tikhonov} at $\eps > 0$ is the trace of the bounded map $U \to U$
\begin{equation}\label{eq:trace-unfolded}
  f \;\longmapsto\; \Phi_r^{*}\,\bigl[\,(\gamma_t + \eps I)^{-1}\,[\D_r\gamma_t](f)\,(\gamma_t + \eps I)^{-1} h\,\bigr] \;\in\; U,
\end{equation}
or equivalently, the scalar $\sum_{j\ge 1}\bigl\langle\Phi_r^{*}\,(\gamma_t+\eps I)^{-1}\,[\D_r\gamma_t](f_j)\,(\gamma_t+\eps I)^{-1}h,\,f_j\bigr\rangle_U$ for any orthonormal basis $\{f_j\}$ of $U$. At $\eps = 0$, the two intermediate objects $\gamma_t^{\dagger}[\D_r\gamma_t](f)\gamma_t^{\dagger}$ are not separately defined as bounded operators (the inner $\gamma_t^{\dagger}$ may not be applicable in the operator sense if its argument leaves $\Ran(\gamma_t)$); the notation~\eqref{eq:correction-symbolic} is therefore purely symbolic, with the rigorous content being the $L^2(\Omega)$-limit~\eqref{eq:correction-tikhonov}, whose convergence and independence of the regularising sequence are furnished by the dominated-convergence argument of Lemma~\ref{lemma:D-gamma-inv-h} under clause~\textup{(iii)} of Assumption~\ref{ass:nondeg}. The kernel-annihilation property of $\Phi_r^{*}$ (Remark~\ref{rem:kernel-annih}) ensures that the $\lambda_l = 0$ contributions vanish, so that the spectral series~\eqref{eq:correction-spectral} sums only over $(k, l)$ with $\lambda_k, \lambda_l > 0$.
\end{remark}

When $V = H = \R^m$ and $U = \R^d$, the operators reduce to matrices ($Y(t,r) = Y_T Y_r^{-1} \in \R^{m \times m}$, $\gamma_t = \gamma_{X_T} \in \R^{m \times m}$), and the formula \eqref{eq:skorokhod-decomp} reduces to Theorem 2.1 of \cite{mirafzali2025malliavin}; see Corollary \ref{cor:finite-dim} for the detailed verification.

When $\calA(t, u) = -Au$ is linear and $\calB(t, u) = Q^{1/2}$ is state-independent, all correction terms vanish ($Y(t,r) = S(t-r)$ is deterministic, $\calW = 0$, $\calB'_u = 0$), the Skorokhod integral reduces to an It\^o integral, and $\beta_h(u) = -\ip{u - S(t)x}{\gamma_t^{\dagger}h}_H$, recovering Theorem 4 of \cite{mirafzali2025infinite}.

\section{The covering vector field}
\label{sec:malliavin}

The statement just made is an integration-by-parts formula, and an integration-by-parts formula needs a direction in Wiener space along which to integrate. Finding it is the first task. What is required is a random element $v_h$ of the Cameron--Martin space whose action on the solution reproduces the prescribed direction $h$ in the state space; we call such a $v_h$ a covering field. The construction given here is intrinsic, in that it never inverts the flow, and it separates cleanly into two questions, whether $v_h$ exists in $L^2$, which it always does, and whether it lies in the domain of the divergence, which is where the hypotheses enter.

\subsection{Malliavin--Sobolev spaces on the Wiener space of the equation}\label{subsec:mall-prelim}

The Wiener space in question is the one carrying the noise that drives the equation, and it is on that space that all the calculus below takes place. We use the Malliavin calculus framework of \cite{nualart2006malliavin, Nualart_Nualart_2018}. Let $\HW := L^2([0,T]; U)$ be the isonormal Gaussian Hilbert space of the cylindrical Wiener process $W$. The Malliavin derivative $\D$ and the Skorokhod integral $\delta_U$ are adjoint operators related by
\begin{equation}\label{eq:skorokhod-duality}
  \E\bigl[\ip{\D F}{v}_{\HW}\bigr] = \E\bigl[F\,\delta_U(v)\bigr] \qquad \text{for all } F \in \mathbb{D}^{1,2},\; v \in \Dom(\delta_U).
\end{equation}
For a separable Hilbert space $E$, $\mathbb{D}^{k,p}(E)$ denotes the $E$-valued Malliavin--Sobolev space of order $k$ and integrability $p$. We write $L_2(H,H)$ for the Hilbert space of Hilbert--Schmidt operators on $H$; although the covariance operators below are in fact trace class, whenever a covariance is differentiated in the Malliavin sense we regard it as an $L_2(H,H)$-valued random variable, which keeps every Malliavin--Sobolev target space Hilbert-valued (the trace class $L_1(H)$ is not a Hilbert space). Trace-class information about the covariance itself is retained separately. For $F \in \mathbb{D}^{1,2}(H)$, the Malliavin derivative $\D_r F \in L_2(U, H)$ satisfies $\E[\norm{F}_H^2] + \E[\int_0^T \norm{\D_r F}_{L_2(U,H)}^2\dr] < \infty$. The Skorokhod integral $\delta_U(v) \in L^2(\Omega)$ is the unique element satisfying~\eqref{eq:skorokhod-duality}; $\delta_U$ is closed, being the adjoint of the densely defined derivative $\D$, satisfies $\mathbb{D}^{1,2}(\HW) \subset \Dom(\delta_U)$, and reduces to the It\^o integral for adapted integrands.

For adapted $v$ with $\E[\int_0^T \norm{v(r)}_U^2\dr] < \infty$, $\delta_U(v) = \int_0^T \ip{v(r)}{\mathrm{d}W(r)}_U$.

The fundamental product rule for the Skorokhod integral (cf.\ \cite{nualart2006malliavin}) asserts that if $F \in \mathbb{D}^{1,2}$ and $v \in \Dom(\delta_U)$ with $Fv \in L^2(\Omega; \HW)$ and the right-hand side below is in $L^2(\Omega)$, then $Fv \in \Dom(\delta_U)$ and
\begin{equation}\label{eq:skorokhod-product}
  \delta_U(Fv) = F\,\delta_U(v) - \ip{\D F}{v}_{\HW}.
\end{equation}

The following result on the Malliavin derivative of an It\^o integral is used repeatedly in the proofs of Section \ref{sec:variation}.

\begin{proposition}[Malliavin derivative of an It\^o integral; cf.\ {\cite{nualart2006malliavin}}, {\cite{leon1998stochastic}}]\label{prop:D-ito-integral}
Let $u = \{u(s) : s \in [a, b]\}$ be an $L_2(U, H)$-valued adapted process with $u \in L^2(\Omega; L^2(a,b; L_2(U,H)))$ and $u(s) \in \mathbb{D}^{1,2}(L_2(U,H))$ for a.e.\ $s$, with
\[
  \E\Bigl[\int_a^b \int_a^b \norm{\D_r u(s)}_{L_2(U, L_2(U,H))}^2\dr\ds\Bigr] < \infty.
\]
Then the It\^o integral $I := \int_a^b u(s)\dW(s)$ belongs to $\mathbb{D}^{1,2}(H)$, and for a.e.\ $r \in [0,T]$:
\begin{equation}\label{eq:D-ito-integral}
  \D_r I = u(r)\,\mathbf{1}_{[a,b]}(r) + \int_{a \vee r}^b \D_r u(s)\dW(s).
\end{equation}
\end{proposition}

\begin{proof}
See \cite{nualart2006malliavin} and \cite{leon1998stochastic}.
\end{proof}

\subsection{What a logarithmic derivative is}

The object we are computing must be defined before it can be computed, and in infinite dimensions this needs a little care. Since $H$ is infinite-dimensional, there is no Lebesgue measure on $H$, and the notion of a probability density requires careful formulation. Let $\mu_t$ denote the law of $X(t)$ on $H$. The appropriate notion of ``gradient of the log-density'' is the logarithmic derivative, or Fomin derivative, of $\mu_t$; see \cite{bogachev2010differentiable}.

\begin{definition}[Logarithmic derivative]\label{def:log-deriv}
For $h \in H$, the \emph{logarithmic derivative} of $\mu_t$ along $h$ is the function $\beta_h: H \to \R$ defined $\mu_t$-almost everywhere by the integration-by-parts identity
\begin{equation}\label{eq:log-deriv-ibp}
  \int_H \ip{\nabla\phi(u)}{h}_H\,\mathrm{d}\mu_t(u) = -\int_H \phi(u)\,\beta_h(u)\,\mathrm{d}\mu_t(u)
\end{equation}
for all $\phi \in C_b^1(H)$, provided such a $\beta_h \in L^1(\mu_t)$ exists.
\end{definition}

The logarithmic derivative $\beta_h(u)$ is the analogue, in infinitely many dimensions, of the directional score $\ip{\nabla \log p_t(u)}{h}$. For a centred Gaussian measure $\mu_t = N(0,\gamma_t)$ the Cameron--Martin space is $\mathcal{H}_t = \Ran(\gamma_t^{1/2})$ \cite{bogachev2010differentiable}, and there $\beta_h(u) = -\ip{u}{\gamma_t^{-1}h}$ for every $h \in \Ran(\gamma_t) \subsetneq \mathcal{H}_t$. Assumption~\ref{ass:nondeg} imposes the range condition $h \in \Ran(\gamma_t)$ almost surely, which is what the explicit formula~\eqref{eq:skorokhod-decomp} requires.

\begin{remark}[Uniqueness and regularity of $\beta_h$]\label{rem:beta-uniqueness}
The logarithmic derivative $\beta_h$, if it exists, is unique $\mu_t$-a.e.\ in $L^1(\mu_t)$, since $C_b^1(H)$ separates $L^1(\mu_t)$ via density of smooth cylindrical functions \cite{bogachev2010differentiable}. Theorem~\ref{thm:main} gives $\beta_h(X(t)) = -\E[\delta_U(v_h) \mid X(t)]$; since $\delta_U(v_h) \in L^2(\Omega)$, Jensen's inequality yields $\beta_h \in L^2(\mu_t)$.
\end{remark}

The following chain rule for composed functionals is used throughout.

\begin{proposition}[Chain rule; cf.\ {\cite{nualart2006malliavin}}]\label{prop:chain-rule}
If $\phi \in C^1_b(H;\R)$ and $F \in \mathbb{D}^{1,2}(H)$, then $\phi(F) \in \mathbb{D}^{1,2}$ and
\begin{equation}\label{eq:chain-rule}
  \D_r(\phi(F)) = (\D_r F)^*\,\nabla\phi(F), \qquad r \in [0,T],
\end{equation}
where $(\D_r F)^* \in L(H, U)$ is the adjoint of $\D_r F \in L_2(U, H)$, and hence $\D_r(\phi(F)) \in U$.
\end{proposition}

\begin{proof}
The proof is a smooth cylindrical approximation argument.

Since $F \in \mathbb{D}^{1,2}(H)$, there exists a sequence of smooth cylindrical $H$-valued random variables $F_n = \sum_{i=1}^{N_n} \varphi_i^{(n)}(W(\eta_1^{(n)}), \ldots, W(\eta_{m_n}^{(n)}))\,h_i^{(n)}$, where $\varphi_i^{(n)} \in C_{\mathrm{pol}}^\infty(\R^{m_n})$, $\eta_j^{(n)} \in \HW$, and $h_i^{(n)} \in H$, such that $F_n \to F$ in $\mathbb{D}^{1,2}(H)$. That is,
\[
  \E\bigl[\norm{F_n - F}_H^2\bigr] + \E\Bigl[\int_0^T \norm{\D_r F_n - \D_r F}_{L_2(U,H)}^2\dr\Bigr] \to 0.
\]

For each $F_n$, the composite $\phi(F_n)$ is a smooth cylindrical random variable (since $\phi$ has bounded derivatives and $F_n$ is smooth cylindrical). Writing $F_n = \sum_{i=1}^{N_n}F_n^{(i)}\,h_i^{(n)}$ in coordinates with respect to the basis $\{h_i^{(n)}\}$ and applying the classical multivariate chain rule for the smooth function $\phi$ to the scalar cylindrical random variables $F_n^{(i)} = \ip{F_n}{h_i^{(n)}}_H$, one obtains, for each $r \in [0,T]$,
\[
  \D_r\bigl(\phi(F_n)\bigr) \;=\; \sum_{i=1}^{N_n}(\partial_i\phi)(F_n)\,\D_r F_n^{(i)} \;=\; \sum_{i=1}^{N_n}\ip{\nabla\phi(F_n)}{h_i^{(n)}}_H\,\D_r F_n^{(i)} \;\in\; U,
\]
where $\partial_i\phi(\cdot) = \ip{\nabla\phi(\cdot)}{h_i^{(n)}}_H$ is the partial derivative in the $h_i^{(n)}$-direction. We claim this equals $(\D_r F_n)^*\,\nabla\phi(F_n)$. Indeed, $\D_r F_n = \sum_i (\D_r F_n^{(i)})\otimes h_i^{(n)}$ as an element of $U\otimes H \cong L_2(U,H)$, so for any $f \in U$,
\[
  (\D_r F_n)\cdot f \;=\; \sum_{i=1}^{N_n}\ip{\D_r F_n^{(i)}}{f}_U\,h_i^{(n)} \;\in\; H,
\]
and the adjoint $(\D_r F_n)^*: H \to U$ satisfies, for $g \in H$, $(\D_r F_n)^*g = \sum_i \ip{h_i^{(n)}}{g}_H\,\D_r F_n^{(i)}$. Setting $g = \nabla\phi(F_n)$ gives
\[
  (\D_r F_n)^*\,\nabla\phi(F_n) \;=\; \sum_{i=1}^{N_n}\ip{\nabla\phi(F_n)}{h_i^{(n)}}_H\,\D_r F_n^{(i)} \;=\; \D_r\bigl(\phi(F_n)\bigr),
\]
as claimed.

Since $\norm{\nabla\phi}_\infty < \infty$:
\begin{align*}
  &\norm{\D_r(\phi(F_n)) - (\D_r F)^*\nabla\phi(F)}_U\\
  &\quad \leq \norm{(\D_r F_n)^*\nabla\phi(F_n) - (\D_r F)^*\nabla\phi(F)}_U\\
  &\quad \leq \norm{(\D_r F_n - \D_r F)^*\nabla\phi(F_n)}_U + \norm{(\D_r F)^*(\nabla\phi(F_n) - \nabla\phi(F))}_U\\
  &\quad \leq \norm{\nabla\phi}_\infty\,\norm{\D_r F_n - \D_r F}_{L_2(U,H)} + \norm{\D_r F}_{L_2(U,H)}\,\norm{\nabla\phi(F_n) - \nabla\phi(F)}_H.
\end{align*}
The first term converges to zero in $L^2(\Omega \times [0,T])$ by the $\mathbb{D}^{1,2}$ convergence of $F_n \to F$. The second term converges to zero by the continuity of $\nabla\phi$ and the $L^2(\Omega; H)$ convergence $F_n \to F$, combined with the dominated convergence theorem (using $\norm{\nabla\phi(F_n) - \nabla\phi(F)}_H \leq 2\norm{\nabla\phi}_\infty$ and $\norm{\D_r F}_{L_2(U,H)} \in L^2(\Omega \times [0,T])$).

By the closability of $\D$ on $\mathbb{D}^{1,2}$ (\cite{nualart2006malliavin}), $\phi(F) \in \mathbb{D}^{1,2}$ and $\D_r(\phi(F)) = (\D_r F)^*\nabla\phi(F)$.
\end{proof}

\subsection{Differentiating the solution}

The first step towards the Bismut formula is to show that $X(t)$ is Malliavin differentiable and to identify its derivative through the first variation.

\begin{proposition}[Malliavin derivative of $X(t)$]\label{prop:D-Xt}
Under Assumptions \ref{ass:LR}, \ref{ass:diff2}, and~\ref{ass:SC}, the solution $X(t)$ belongs to $\mathbb{D}^{1,2}(H)$ for each $t \in [0, T]$, and its Malliavin derivative at time $r \in [0, t]$ is given by
\begin{equation}\label{eq:DX}
  \D_r X(t) = Y(t,r)\,\calB(r, X(r)),
\end{equation}
where $Y(t,r)$ is the first variation process (Definition \ref{def:first-var}).
\end{proposition}

\begin{proof}
Let $\{e_j\}_{j \geq 1} \subset V$ be an orthonormal basis of $H$ with $e_j \in V$ for all $j$, and let $\Pi_N: H \to H_N := \spn\{e_1, \ldots, e_N\}$ denote the orthogonal projection.

The $N$-dimensional Galerkin approximation $X^N(t) \in H_N$ solves the finite-dimensional SDE
\begin{equation}\label{eq:galerkin}
  \mathrm{d}X^N(t) + \Pi_N\calA(t, X^N(t))\dt = \Pi_N\calB(t, X^N(t))\dW(t), \qquad X^N(0) = \Pi_N x.
\end{equation}
Under (D1)--(D2), the projected coefficients $\Pi_N\calA|_{H_N}: H_N \to H_N$ and $\Pi_N\calB|_{H_N}: H_N \to L_2(U, H_N)$ are continuously Fr\'echet differentiable (i.e.\ $C^1$ in the finite-dimensional sense on $H_N$, by~(D1) Gâteaux-differentiability plus continuity of $\calA'_u, \calB'_u$ in the finite-dimensional subspace, combined with the Gelfand-triple compatibility $V \cap H_N \equiv H_N$ since $H_N$ is spanned by elements of $V$) with polynomial growth. Such smooth, polynomial-growth, but not globally Lipschitz, finite-dimensional SDEs fall outside the scope of the classical formula \cite{nualart2006malliavin} (which requires global Lipschitz coefficients), but are covered by the standard truncation-plus-localisation machinery of Malliavin calculus for SDEs with smooth non-Lipschitz coefficients \cite{nualart2006malliavin}. Introducing the exit time $\tau_M^N := \inf\{t \in [0,T]: \norm{X^N(t)}_{H_N} > M\}$ and the Lipschitzified truncated coefficients $(\calA^{N,M}, \calB^{N,M})$ that coincide with $(\Pi_N\calA, \Pi_N\calB)$ on $\{u \in H_N : \norm{u}_{H_N} \leq M\}$ and are extended to be globally Lipschitz on $H_N$, the truncated solution $X^{N,M}$ coincides with $X^N$ on $[0, \tau_M^N]$ and, by the globally-Lipschitz Nualart formula \cite{nualart2006malliavin}, lies in $\mathbb{D}^{1,2}(H_N)$ with the standard Malliavin-derivative formula. The stopping times $\tau_M^N$ exhaust $[0,T]$ a.s.\ as $M \to \infty$ by the a.s.\ continuity of $X^N$ in $H_N$, so $X^N \in \mathbb{D}^{1,p}_{\mathrm{loc}}(H_N)$; the exponential-moment bound~(SC2) combined with the polynomial growth of $\calA'_u, \calB'_u$ from (D1)--(D2) then forces the localised $\mathbb{D}^{1,2}$-bounds on the Malliavin derivative to close uniformly in $M$, upgrading $\mathbb{D}^{1,p}_{\mathrm{loc}}$-membership to $\mathbb{D}^{1,2}(H_N)$-membership in the limit $M \to \infty$. The direct SPDE-level treatment of Bally--Pardoux \cite{ballypardoux1998} for white-noise-driven parabolic SPDEs with smooth coefficients, bounded together with their derivatives, specialises to the same conclusion for the Galerkin system~\eqref{eq:galerkin}. The result is that $X^N(t) \in \mathbb{D}^{1,2}(H_N)$ with Malliavin derivative
\begin{equation}\label{eq:DX-galerkin}
  \D_r X^N(t) = Y^N(t,r)\,\Pi_N\calB(r, X^N(r)), \qquad r \in [0, t],
\end{equation}
where $Y^N(t,r)$ is the first variation of the projected system, satisfying
\begin{align*}
  &\mathrm{d}Y^N(t,r) + \Pi_N\calA'_u(t, X^N(t))\,Y^N(t,r)\dt \\
  &\qquad = \Pi_N\calB'_u(t, X^N(t))(Y^N(t,r))\dW(t), \qquad Y^N(r,r) = \Pi_N.
\end{align*}

The R\"ockner--Shang--Zhang theory \cite{rockner2022wellposedness} furnishes the uniform Galerkin moment bounds: for every $q \geq 2$,
\begin{equation}\label{eq:galerkin-moments}
  \sup_{N \in \N}\Bigl\{\E\Bigl[\sup_{0 \leq t \leq T}\norm{X^N(t)}_H^{q}\Bigr] + \E\Bigl[\Bigl(\int_0^T \norm{X^N(t)}_V^{p}\,\mathrm{d}t\Bigr)^{\!q/2}\Bigr]\Bigr\} \;\leq\; C_q\bigl(1 + \norm{x}_H^{q}\bigr),
\end{equation}
together with tightness of the laws of $\{X^N\}$ and identification of every subsequential limit with the unique solution $X$. Strong convergence of the Galerkin sequence on the original probability space, and of its variational fibre, is the structural input~(SC3) of Assumption~\ref{ass:SC}:
\begin{equation}\label{eq:galerkin-conv}
  X^N \to X \quad \text{in } L^2(\Omega; C([0,T]; H)) \cap L^{b}(\Omega \times [0,T]; V), \qquad b = 6(1+\eps_0),
\end{equation}
and, by the first clause of~\eqref{eq:SC-Galerkin-linearised},
\begin{equation}\label{eq:YN-HS-conv}
  \E\int_0^t\norm{\Phi^N_{t,r} - \Phi_{t,r}}^2_{L_2(U,H)}\dr \;\longrightarrow\; 0 \qquad (N \to \infty),
\end{equation}
where $\Phi^N_{t,r} = Y^N(t,r)\Pi_N\calB(r, X^N(r))$ is the Galerkin fibre and $\Phi_{t,r} = Y(t,r)\calB(r,X(r))$ the limiting one.

To show $X(t) \in \mathbb{D}^{1,2}(H)$ and $\D_r X(t) = \Phi_{t,r}$ we use the closed-graph property of the Malliavin derivative, a direct consequence of the closability of $\D$ \cite{nualart2006malliavin}. By~\eqref{eq:DX-galerkin} the Galerkin derivative is $\D_r X^N(t) = \Phi^N_{t,r}$; by~\eqref{eq:galerkin-conv} we have $X^N(t) \to X(t)$ in $L^2(\Omega; H)$; and by~\eqref{eq:YN-HS-conv} we have $\D_r X^N(t) \to \Phi_{t,r}$ in $L^2(\Omega \times [0,t]; L_2(U,H))$. Closedness therefore gives $X(t) \in \mathbb{D}^{1,2}(H)$ with
\[
  \D_r X(t) \;=\; \Phi_{t,r} \;=\; Y(t,r)\,\calB(r,X(r)).
\]
Stating the Galerkin clause on the fibre rather than on the flow is what makes this a single application of closedness. A statement about $Y^N \to Y$ would have to be composed with the convergence of the diffusion coefficients and recombined, and the recombination is exactly the step that would reintroduce a norm of the first variation on an arbitrary direction.

It remains to check that the limit has the integrability the space requires,
\begin{equation*}
  \E\Bigl[\int_0^t \norm{\D_r X(t)}_{L_2(U,H)}^2\dr\Bigr] \;=\; \E\Bigl[\int_0^t \norm{\Phi_{t,r}}_{L_2(U,H)}^2\dr\Bigr] \;<\; \infty,
\end{equation*}
and this is the fibre estimate itself, applied at exponent $q = 1$ to the $\F_r$-measurable datum $\Theta = \calB(r, X(r))$. Theorem~\ref{thm:Y-wellposed}\textup{(III)} gives
\[
  \E\bigl[\norm{\Phi_{t,r}}^2_{L_2(U,H)}\bigr] \;\leq\; C_1\,\E\Bigl[\norm{\calB(r,X(r))}^2_{L_2(U,H)}\exp\Bigl(c_1\!\int_r^t(\tilde\rho+\hat\rho)(\sigma,X(\sigma))\,\mathrm{d}\sigma\Bigr)\Bigr].
\]
The two branches of~(D2$'$) are separated before any H\"older step. For deterministic noise the coefficient is non-random, so the joint expectation factors exactly, as the product of $\norm{\calB(r)}^2_{L_2}$ with the exponential moment of~(SC2), and no exponent $m_0/(m_0-2)$ occurs anywhere. In the state-dependent branch, where~(D2$'$) supplies $m_0 > 2$ strictly, H\"older's inequality at the conjugate pair $(m_0/2,\,m_0/(m_0-2))$ gives
\begin{align}
  &\E\Bigl[\norm{\calB(r,X(r))}_{L_2(U,H)}^2 \exp\Bigl(c_1\int_r^t(\tilde\rho+\hat\rho)\,\mathrm{d}\sigma\Bigr)\Bigr] \nonumber\\
  &\quad \leq\; \Bigl(\E\bigl[\norm{\calB(r,X(r))}_{L_2(U,H)}^{m_0}\bigr]\Bigr)^{2/m_0} \Bigl(\E\Bigl[\exp\Bigl(\tfrac{c_1 m_0}{m_0 - 2}\int_r^t(\tilde\rho+\hat\rho)\,\mathrm{d}\sigma\Bigr)\Bigr]\Bigr)^{(m_0-2)/m_0},
  \label{eq:holder-split}
\end{align}
the reciprocals of the two exponents summing to $2/m_0 + (m_0-2)/m_0 = 1$. The strict gap $m_0 > 2$ makes the constant $c_1m_0/(m_0-2)$ finite, at the quasilinear operating point $m_0 = p/\beta$ it is precisely $2\beta < p$ from Assumption~\ref{ass:diff1}, and $c_1m_0/(m_0-2) \leq c_*$ by~\eqref{eq:cstar-base}, so the second factor is bounded by a constant $K < \infty$ uniformly in $r \in [0,t]$ through~(SC2).

For the first factor, the noise-size growth~\eqref{eq:B-growth-m} of~(D2$'$) gives $\norm{\calB(r,u)}^{m_0}_{L_2(U,H)} \leq C_{m_0}(1+\norm{u}_V^p)(1+\norm{u}_H^{\kappa_{m_0}})$, and the quantity actually consumed is the time integral $\int_0^t(\E[\norm{\calB(r,X(r))}^{m_0}_{L_2}])^{2/m_0}\dr$. Since $2/m_0 \leq 1$, Jensen's inequality on $([0,t], \dr/t)$ reduces it to the time-integrated coefficient moment,
\begin{equation}\label{eq:B-time-integrated}
  \int_0^t \bigl(\E[\norm{\calB(r,X(r))}^{m_0}_{L_2}]\bigr)^{2/m_0}\dr
  \;\leq\; t^{1 - 2/m_0}\Bigl(\E\!\int_0^t\!\norm{\calB(r,X(r))}^{m_0}_{L_2}\dr\Bigr)^{\!2/m_0},
\end{equation}
so the pointwise-in-$r$ moment never has to be integrated at the power $2/m_0$ separately. The right-hand side closes by pulling the $H$-factor out through $\sup_\sigma\norm{X(\sigma)}_H$ and applying Cauchy--Schwarz in $\omega$,
\[
\begin{aligned}
  \E\int_0^t\norm{\calB(r,X(r))}^{m_0}_{L_2}\dr \;\leq\; C_{m_0}&\Bigl(\E\Bigl[\Bigl(t + \int_0^t\norm{X(r)}_V^p\dr\Bigr)^{2}\Bigr]\Bigr)^{1/2}\\
  &\times\Bigl(\E\bigl[(1 + \sup\nolimits_\sigma\norm{X(\sigma)}_H^{\kappa_{m_0}})^2\bigr]\Bigr)^{1/2},
\end{aligned}
\]
both factors being finite by the iterated It\^o-energy moments of~\eqref{eq:higher-moment-ito}, which control $\E[(\int_0^T\norm{X}_V^p\dr)^2]$ and every polynomial $H$-moment of $\sup_\sigma\norm{X(\sigma)}_H$. Combining,
\[
  \E\Bigl[\int_0^t \norm{\D_r X(t)}_{L_2(U,H)}^2\dr\Bigr] \;\leq\; C_1\,K\,t^{1-2/m_0}\Bigl(\E\int_0^t\norm{\calB(r,X(r))}^{m_0}_{L_2(U,H)}\dr\Bigr)^{2/m_0} \;<\; \infty,
\]
which establishes $X(t) \in \mathbb{D}^{1,2}(H)$.

\medskip
The finiteness chain uses only the fibre estimate~\eqref{eq:fibre-energy-estimate}, the Liu--R\"ockner energy moments of~\eqref{eq:apriori} and~\eqref{eq:higher-moment-ito}, the strict gap $m_0 > 2$ of~(D2$'$), and the exponential moment~\eqref{eq:SC-exp-moment}. The gap enters through the conjugate pair $(m_0/2, m_0/(m_0-2))$, which is in this sense the minimal separation of the joint bound, since the non-strict endpoint $m_0 = 2$ would make the second exponent infinite and break the exponential factor. At the quasilinear operating point $m_0 = p/\beta$ with $\beta > 0$ the gap reads $2\beta < p$, which is why (D2) is stated with the strict inequality; for $\beta = 0$, that is for additive and more generally $H$-extension noise, the same separation runs at any admissible $m_0 > 2$ furnished by~\eqref{eq:B-growth-m}, and for deterministic noise the joint expectation factors exactly and no separation is needed. In every case the conclusion is reached without $V$-moment control beyond~\eqref{eq:apriori}.
\end{proof}

The proof of Theorem~\ref{thm:main} in Section~\ref{sec:proof-main} does not require $X(t) \in \mathbb{D}^{2,2}(H)$. The second-order information entering the Bismut--Nualart formula comes in through $\D_r\gamma_t$, which is computed in Proposition~\ref{prop:D-gamma} below by differentiating the quadratic representation $\gamma_t = \int_0^t\Phi_s\Phi_s^*\,\mathrm{d}s$ via the product rule and invoking the second Malliavin fibre $\calW(t;s,r)$ of Lemma~\ref{lemma:DPhi} directly, rather than via a consolidated $\mathbb{D}^{2,2}(H)$ statement on $X(t)$. This is not an evasion but the natural route. A direct $\mathbb{D}^{2,2}(H)$ claim would require bounding products like $\E[\norm{\calB'_u(X(r))}^2\,\norm{\calB(X(s))}^2]$ on the full time square, for which the strict gap $2\beta < p$ of Assumption~\ref{ass:diff1} is no longer sufficient; whereas the covering-field approach needs only the $\mathbb{D}^{1,2}(\HW)$-regularity of the linear-in-$z$ field $w_r(z)$, which reduces to first-order bounds on the two fibres $\Phi$ and $\calW$, the latter entering linearly rather than squared in the expressions of Lemma~\ref{lemma:D-gamma-inv-h}.

\subsection{The covariance operator and its pseudoinverse}

With $\Phi_r = D_r X(t)$ in hand, the Malliavin covariance $\gamma_t = \int_0^t \Phi_r \Phi_r^*\,\mathrm{d}r$ inherits its properties from the regularity of $\Phi$.

\begin{proposition}[Properties of $\gamma_t$]\label{prop:gamma-properties}
Under Assumptions \ref{ass:LR}, \ref{ass:diff2}, and~\ref{ass:SC}, the Malliavin covariance operator $\gamma_t$ defined in \eqref{eq:malliavin-cov-def} is a random, positive, self-adjoint, trace-class operator on $H$, satisfying:
\begin{enumerate}[(i)]
  \item $\ip{\gamma_t \phi}{\phi}_H \geq 0$ for all $\phi \in H$;
  \item $\gamma_t^* = \gamma_t$;
  \item $\Tr_H(\gamma_t) = \int_0^t \norm{\D_r X(t)}_{L_2(U,H)}^2\dr < \infty$ a.s.;
  \item for $h$ satisfying Assumption \ref{ass:nondeg}, $h \in \Ran(\gamma_t)$ a.s., the pseudoinverse $\gamma_t^{\dagger} h \in H$ is well-defined, and $\gamma_t\gamma_t^{\dagger} h = h$.
\end{enumerate}
\end{proposition}

\begin{proof}
Properties (i) and (ii) follow from the representation $\gamma_t = \int_0^t \Phi_r \Phi_r^*\dr$ where $\Phi_r = Y(t,r)\calB(r, X(r)) \in L_2(U, H)$. Each $\Phi_r \Phi_r^*$ is a positive self-adjoint operator on $H$, and the integral preserves these properties.

For (iii), let $\{e_j\}_{j \geq 1}$ be an ONB of $H$. Then
\begin{align*}
  \Tr_H(\gamma_t) &= \sum_{j=1}^\infty \ip{\gamma_t e_j}{e_j}_H
  = \sum_{j=1}^\infty \int_0^t \norm{\Phi_r^* e_j}_U^2\dr
  = \int_0^t \sum_{j=1}^\infty \norm{\Phi_r^* e_j}_U^2\dr\\
  &= \int_0^t \norm{\Phi_r}_{L_2(U,H)}^2\dr
  = \int_0^t \norm{\D_r X(t)}_{L_2(U,H)}^2\dr < \infty \quad \text{a.s.}
\end{align*}
by the Malliavin differentiability of $X(t)$ (Proposition \ref{prop:D-Xt}).

Property (iv) is Assumption \ref{ass:nondeg}.
\end{proof}

\subsection{Construction and regularity of the field}\label{subsec:covering-field}

The covering field $v_h$ is the minimum-$\mathcal{H}_W$-norm element of $L^2(\Omega; \mathcal{H}_W)$ whose pushforward through $D X(t)$ equals the prescribed direction $h$. By Moore--Penrose theory, $v_h(r) = \Phi_r^{*}\,\gamma_t^{\dagger}h$. We establish its covering property, existence and uniqueness, Skorokhod-domain membership, $L^2$-continuity, and stability under perturbations.

\begin{theorem}[Covering property]\label{thm:covering}
Under Assumptions \ref{ass:LR}, \ref{ass:diff2}, \ref{ass:SC}, and \ref{ass:nondeg}, for each $h$ satisfying Assumption \ref{ass:nondeg}, the covering vector field $v_h$ defined in \eqref{eq:covering-field} satisfies the covering condition
\begin{equation}\label{eq:covering-condition}
  \int_0^t \D_r X(t)\,v_h(r)\dr = h.
\end{equation}
\end{theorem}

\begin{proof}
Substituting the definitions:
\begin{align*}
  \int_0^t \D_r X(t)\,v_h(r)\dr
  &= \int_0^t \Phi_r\,\Phi_r^{*}\,\tilde h\,\mathrm{d}r\\
  &= \gamma_t\,\tilde{h}
  = \gamma_t\,\gamma_t^{\dagger}\,h = h,
\end{align*}
where the first line is the substitution $v_h(r) = \Phi_r^{*}\tilde h$ in the integrand, and the last equality holds because $h \in \Ran(\gamma_t)$ a.s.\ by Assumption \ref{ass:nondeg}. For such $h$ the Moore--Penrose pseudoinverse satisfies $\gamma_t\gamma_t^\dagger h = h$, this being its defining property on $\Ran(\gamma_t)$.
\end{proof}

\medskip\noindent The existence, uniqueness, and regularity of the covering field follow.

\begin{theorem}[Existence and uniqueness of $v_h$]\label{thm:covering-exist}
Under Assumptions~\ref{ass:LR}, \ref{ass:diff2}, \ref{ass:SC}, and~\ref{ass:nondeg}, for each $h$ satisfying Assumption~\ref{ass:nondeg}, the covering vector field $v_h$ defined in~\eqref{eq:covering-field} satisfies:
\begin{enumerate}[label=\textup{(\roman*)}]
 \item \textup{(Algebraic $L^2$-content, unconditional.)} $v_h$ is the unique element of minimal $\HW$-norm in $L^2(\Omega; \HW)$ satisfying the covering condition~\eqref{eq:covering-condition}; equivalently, $v_h \in \Ran(\mathcal{T}^{*}) \subset L^2(\Omega; \HW)$ with $\mathcal{T}v_h = h$ a.s.
 \item \textup{(Skorokhod-domain membership, conditional.)} Under any one of the cases of Theorem~\ref{thm:covering-regularity} below, $v_h$ further satisfies $v_h \in \Dom(\delta_U)$.
\end{enumerate}
Part~\textup{(i)} is the SPDE-side specialisation of the algebraic content of Theorem~\ref{thm:pseudoinverse-covering} and requires no additional hypothesis; only the Skorokhod-integrability of $v_h$ in part~\textup{(ii)} requires the additional input of Theorem~\ref{thm:covering-regularity}.
\end{theorem}

\begin{proof}
Define the bounded linear operator $\mathcal{T}: \HW \to H$ by
\[
  \mathcal{T}f := \int_0^t \D_r X(t)\,f(r)\dr, \qquad f \in \HW.
\]
Its adjoint $\mathcal{T}^*: H \to \HW$ acts as $(\mathcal{T}^*\phi)(r) = \Phi_r^{*}\phi$. The Malliavin covariance is $\gamma_t = \mathcal{T}\mathcal{T}^*$.

The covering field is $v_h = \mathcal{T}^*\gamma_t^{\dagger}h = \mathcal{T}^*\tilde{h}$, which lies in $\Ran(\mathcal{T}^*)$. By the standard Moore--Penrose theory for bounded linear operators between Hilbert spaces (cf.\ \cite{engl1996regularization}), this is the unique minimum-norm element satisfying $\mathcal{T}v = h$. Suppose $v'$ also satisfies $\mathcal{T}v' = h$. Then $w := v' - v_h \in \ker(\mathcal{T})$. Since $v_h \in \Ran(\mathcal{T}^*)$ and $\ker(\mathcal{T}) = \Ran(\mathcal{T}^*)^\perp$, we have $\ip{v_h}{w}_{\HW} = 0$. Therefore
\[
  \norm{v'}_{\HW}^2 = \norm{v_h + w}_{\HW}^2 = \norm{v_h}_{\HW}^2 + \norm{w}_{\HW}^2 \geq \norm{v_h}_{\HW}^2,
\]
with equality only if $w = 0$. The $L^2(\Omega; \HW)$ membership $v_h \in L^2(\Omega; \HW)$ follows from Theorem~\ref{thm:covering-L2-cont}. Part (i) requires no input beyond the standing hypotheses, the algebraic Moore--Penrose construction and the $L^2$-regularity bound are entirely operator-theoretic.

The Skorokhod-domain membership $v_h \in \Dom(\delta_U)$ is the content of Theorem~\ref{thm:covering-regularity} below; it is established under the deterministic-covariance case~(a) directly via $\F_r$-adaptedness of $v_h$, and under the second-order differentiability case~(b) via the closed-graph Tikhonov argument of Lemma~\ref{lemma:D-gamma-inv-h} together with the substitution-regularity hypothesis \ref{hyp:iii-F-prime}. Both cases are strictly stronger inputs than those of Part~(i); without them, $v_h$ is only known as an $L^2(\Omega; \HW)$-element with the algebraic covering identity.
\end{proof}

Under the second-order Assumption~\ref{ass:diff2}, Lemma~\ref{lemma:D-gamma-inv-h} establishes $v_h \in \Dom(\delta_U)$ via the Tikhonov regularisation $v_h^{(\eps)} \to v_h$, a Hilbert--Schmidt kernel representation for the main term, and the scalar convergence-and-domination condition of Assumption~\ref{ass:nondeg}(iii) for the correction. The pointwise $\mathbb{D}^{1,2}(\HW)$-regularity of $v_h$ is neither claimed nor needed.

\begin{theorem}[When $v_h$ lies in the domain of $\delta_U$]\label{thm:covering-regularity}
Under Assumptions \ref{ass:LR}, \ref{ass:diff2}, \ref{ass:SC}, and \ref{ass:nondeg}, the covering vector field $v_h$ belongs to $L^2(\Omega; \HW)$ and satisfies the covering condition \eqref{eq:covering-condition}. Moreover, $v_h \in \Dom(\delta_U)$ in each of the following cases:
\begin{enumerate}[(a)]
  \item \textup{(Linear drift with state-independent diffusion)} $\calA(t,u) = -Au$ is linear and $\calB(t,u) = \calB(t)$ does not depend on $u$. Then $Y(t,r) = S(t-r)$ is the (deterministic) semigroup generated by $A$, the Malliavin covariance $\gamma_t$ is deterministic, and $v_h(r) = \calB(r)^*S(t-r)^*\gamma_t^{\dagger}h$ is a deterministic function of $r$, hence adapted. The Skorokhod integral reduces to the It\^o integral $\delta_U(v_h) = \int_0^t \ip{v_h(r)}{\mathrm{d}W(r)}_U$. This case requires no differentiability beyond what is needed for well-posedness.
  \item \textup{(Second-order differentiability + fixed-Tikhonov substitution regularity)} The full Assumption~\ref{ass:nondeg} (including the Cameron--Martin compatibility clause~(iii)) holds with the exponent $q$ of clause~(ii) chosen in $(2, \infty]$, and Assumption~\ref{ass:SC-raised} holds at target $q^* = q/(q-2)$ (with $q^* = 1$ when $q = \infty$), and the fixed-Tikhonov substitution regularity hypothesis~\ref{hyp:iii-F-prime} of Theorem~\ref{thm:tikhonov-trace} holds for the SPDE-side covering data $(\mathcal{T}, \gamma_F, h) = (\Phi, \gamma_t, h)$. SPDE-level sufficient conditions for~\ref{hyp:iii-F-prime} (in particular, automatic verification via the symmetric H\"older estimate when $q \leq 3$ on the H\"older triangle) are recorded in Remark~\ref{rem:covering-iiiF-SPDE}.
\end{enumerate}
\end{theorem}

\begin{remark}[SPDE-level verification of \ref{hyp:iii-F-prime}]\label{rem:covering-iiiF-SPDE}
At $2 < q \leq 3$, equivalently $p^* \geq 6$ on the edge choice $p^* = 2q/(q-2)$, hypothesis~\ref{hyp:iii-F-prime} is automatic. The moment chain gives $\D\gamma_t \in L^{p^*/2}$ via $\D\gamma_t = (\D\Phi)\Phi^{*} + \Phi(\D\Phi)^{*}$ and Remark~\ref{rem:Phi-all-moments}, yielding $G^{(\eps)} \in \mathbb{D}^{1,2}(H)$ and the H\"older-pair condition~\eqref{eq:linear-sub-Holder}. For $q > 3$, it is an additional structural input. Granted~\ref{hyp:iii-F-prime}, Lemma~\ref{lemma:D-gamma-inv-h} establishes $v_h \in \Dom(\delta_U)$ via closedness of $\delta_U$.
\end{remark}

\begin{proof}
The $L^2(\Omega; \HW)$ membership follows from Theorem \ref{thm:covering-L2-cont}.

For case~(a), when $\calA(t,u) = -Au$ is linear and $\calB(t,u) = \calB(t)$ is state-independent, the linearised coefficients are $\calA'_u(t, X(t)) = -A$ (constant, independent of $X$) and $\calB'_u = 0$. The first variation equation \eqref{eq:first-var} therefore reduces to $\mathrm{d}Y(t,r) + AY(t,r)\dt = 0$ with $Y(r,r) = I_H$, whose unique solution $Y(t,r) = S(t-r)$ is the deterministic semigroup generated by $A$. The Malliavin covariance $\gamma_t = \int_0^t S(t-r)\calB(r)\calB(r)^*S(t-r)^*\dr$ is likewise deterministic, so $\tilde{h} = \gamma_t^{\dagger}h$ is deterministic. Consequently, $v_h(r) = \calB(r)^*S(t-r)^*\gamma_t^{\dagger}h$ is a deterministic function of $r$, hence trivially $\F_r$-adapted. By the classical result \cite{nualart2006malliavin}, any adapted square-integrable process belongs to $\Dom(\delta_U)$ with $\delta_U(v_h) = \int_0^t \ip{v_h(r)}{\mathrm{d}W(r)}_U$.

For case~(b), Lemma~\ref{lemma:D-gamma-inv-h} gives $v_h^{(\eps)} \in \Dom(\delta_U)$ at fixed $\eps$ via~\ref{hyp:iii-F-prime}, the $L^2(\Omega;\HW)$-convergence $v_h^{(\eps)} \to v_h$, and the $L^2(\Omega)$-convergence of $\delta_U(v_h^{(\eps)})$ via arguments~(M) and~(C). Closedness of $\delta_U$ yields $v_h \in \Dom(\delta_U)$.
\end{proof}

For linear drift with additive noise, the covering field is deterministic and $\delta_U(v_h)$ is an It\^o integral (Corollary~\ref{cor:linear-infinite}). For nonlinear drift with state-independent noise, $v_h$ is anticipating but only (D3) is needed from Assumption~\ref{ass:diff2}.

The following two results, $L^2$-continuity and stability of the covering field, are essential for the convergence arguments in Section \ref{subsec:singular} and are the infinite-dimensional analogues of Theorems 3.3 and 3.4 of \cite{mirafzali2025malliavin}.

\begin{theorem}[An $L^2$-estimate for $v_h$]\label{thm:covering-L2-cont}
Under Assumptions \ref{ass:LR}, \ref{ass:diff2}, \ref{ass:SC}, and \ref{ass:nondeg}, for each $h$ satisfying Assumption \ref{ass:nondeg},
\begin{equation}\label{eq:covering-L2-cont}
  \E\Bigl[\int_0^t \norm{v_h(r)}_U^2\dr\Bigr] \;=\; \E\bigl[\ip{h}{\gamma_t^{\dagger} h}_H\bigr] \;\leq\; \norm{h}_H\,\bigl(\E[\norm{\gamma_t^{\dagger}h}_H^2]\bigr)^{1/2} \;<\; \infty.
\end{equation}
Moreover, the map $h \mapsto v_h \in L^2(\Omega; \HW)$ is linear on $\Ran(\gamma_t)$; whenever $h_1, h_2$ additionally satisfy the domain hypotheses of Theorem~\ref{thm:covering-regularity}, the pathwise linearity of $\gamma_t^{\dagger}$ on $\Ran(\gamma_t)$ and the linearity of $\Dom(\delta_U)$ give $v_{\alpha h_1 + \beta h_2} = \alpha v_{h_1} + \beta v_{h_2} \in \Dom(\delta_U)$ and $\delta_U(v_{\alpha h_1 + \beta h_2}) = \alpha\,\delta_U(v_{h_1}) + \beta\,\delta_U(v_{h_2})$ for all $\alpha, \beta \in \R$.
\end{theorem}

\begin{proof}
Writing $\tilde{h} = \gamma_t^{\dagger}h$ and $\Phi_r = Y(t,r)\calB(r, X(r))$, we compute the integral $\int_0^t\norm{v_h(r)}_U^2\dr$ as an exact identity rather than a bound:
\begin{align*}
  \int_0^t \norm{v_h(r)}_U^2\dr
  &= \int_0^t \norm{\Phi_r^*\tilde{h}}_U^2\dr
  = \int_0^t \ip{\Phi_r\Phi_r^*\tilde{h}}{\tilde{h}}_H\dr\\
  &= \Bigl\langle\Bigl(\int_0^t \Phi_r\Phi_r^*\dr\Bigr)\tilde{h},\,\tilde{h}\Bigr\rangle_H
  = \ip{\gamma_t\tilde{h}}{\tilde{h}}_H
  = \ip{h}{\tilde{h}}_H = \ip{h}{\gamma_t^{\dagger}h}_H,
\end{align*}
where the last step uses $\gamma_t\tilde{h} = \gamma_t\gamma_t^{\dagger}h = h$ from the covering condition (Theorem~\ref{thm:covering}). Taking expectations and applying the Cauchy--Schwarz inequality in $L^2(\Omega)$:
\[
  \E\Bigl[\int_0^t\norm{v_h(r)}_U^2\dr\Bigr] = \E\bigl[\ip{h}{\gamma_t^{\dagger}h}_H\bigr] \leq \norm{h}_H\,\bigl(\E[\norm{\gamma_t^{\dagger}h}_H^2]\bigr)^{1/2},
\]
and the right-hand side is finite by Assumption~\ref{ass:nondeg} (which gives $\E[\norm{\gamma_t^{\dagger}h}_H^q] < \infty$ for some $q > 2$, hence in particular the $L^2$-bound). The linearity of $h \mapsto v_h$ on $\Ran(\gamma_t)$ is immediate from the definition of $v_h = \Phi^*\gamma_t^{\dagger}h$ and the linearity of the Moore--Penrose pseudoinverse $h \mapsto \gamma_t^{\dagger}h$ on $\Ran(\gamma_t)$.
\end{proof}

\medskip\noindent The following lemma extracts, in a self-contained form, the Hilbert--Schmidt kernel representation of the anticipating linear functional $z \mapsto \delta_U(\Phi^* z)$ induced by a Malliavin-differentiable Hilbert--Schmidt-valued process $\Phi$. The lemma is the technical substance of the main-term argument~(M) in the proof of Lemma~\ref{lemma:D-gamma-inv-h} of Section~\ref{sec:proof-main}, factored out here because the argument requires only Malliavin regularity of $\Phi$ and is independent of any Tikhonov regularisation or pseudoinverse structure.

\begin{lemma}[A kernel representation for anticipating linear functionals]\label{lemma:HS-kernel}
Let $p \geq 2$ be an exponent and suppose that the fibre $\Phi$ belongs to the space $\mathbb{D}^{1,p}(L^2([0,t]; L_2(U,H)))$. Then the following four statements hold.
\begin{enumerate}[label=\textup{(\roman*)}]
\item For each deterministic $z \in H$, $\Phi^* z \in \mathbb{D}^{1,p}(\HW) \subset \Dom(\delta_U)$, with
\[
  \norm{\Phi^* z}_{\mathbb{D}^{1,p}(\HW)} \leq \norm{z}_H\,\norm{\Phi}_{\mathbb{D}^{1,p}(L^2([0,t];L_2(U,H)))}.
\]
\item The linear map $\phi_\Phi: H \ni z \mapsto \delta_U(\Phi^* z) \in L^2(\Omega; \R)$ is Hilbert--Schmidt from $H$ to $L^2(\Omega; \R)$, with
\[
  \norm{\phi_\Phi}_{\mathrm{HS}(H,\,L^2(\Omega))} \;\leq\; C_{\mathrm{Meyer}}\,\norm{\Phi}_{\mathbb{D}^{1,2}(L^2([0,t]; L_2(U,H)))},
\]
for an absolute constant $C_{\mathrm{Meyer}} > 0$ from the Meyer inequality~\cite{nualart2006malliavin}.
\item By the canonical isomorphism $\mathrm{HS}(H, L^2(\Omega)) \cong L^2(\Omega; H)$, there exists a unique \emph{kernel} $\mathcal{K}_\Phi \in L^2(\Omega; H)$, denoted symbolically $\mathcal{K}_\Phi = \int_0^t\Phi_r\,\delta W_r$ (the $H$-valued Skorokhod integral), such that
\begin{equation}\label{eq:HS-kernel}
  \delta_U(\Phi^* z) \;=\; \ip{z}{\mathcal{K}_\Phi}_H \qquad \mathbb{P}\text{-a.s., for every deterministic } z \in H.
\end{equation}
Moreover, $\mathcal{K}_\Phi \in L^p(\Omega; H)$ with
\[
  \norm{\mathcal{K}_\Phi}_{L^p(\Omega; H)} \leq C_{p,\mathrm{Meyer}}\,\norm{\Phi}_{\mathbb{D}^{1,p}(L^2([0,t]; L_2(U,H)))}.
\]
\item For every jointly $\F_t \otimes \mathcal{B}(H)$-measurable $Z: \Omega \to H$, the pointwise-in-$\omega$ evaluation of the linear functional $z \mapsto \delta_U(\Phi^* z)$ at $z = Z(\omega)$ is given by
\begin{equation}\label{eq:HS-kernel-random}
  \bigl[\delta_U(w^{\mathrm{lin}}(\cdot)(z))\bigr]_{z = Z}(\omega) \;:=\; \ip{Z(\omega)}{\mathcal{K}_\Phi(\omega)}_H \qquad \text{for }\mathbb{P}\text{-a.e.\ } \omega,
\end{equation}
where $w^{\mathrm{lin}}_r(z) := \Phi_r^* z$ denotes the linear random field. The right-hand side of~\eqref{eq:HS-kernel-random} is the substitution value appearing in the linear-field identity of Theorem~\ref{thm:linear-sub-D12}, and differs from $\delta_U(w^{\mathrm{lin}}(\cdot)(Z)) = \delta_U(\Phi^* Z)$ (the Skorokhod integral of the composed, genuinely anticipating process) by the classical Nualart trace correction $\int_0^t\Tr_U[\Phi_r^*\,\D_r Z]\dr$ whenever $Z \in \mathbb{D}^{1,2}(H)$.
\end{enumerate}
\end{lemma}

\begin{proof}
For deterministic $z \in H$, the process $\Phi^* z$ is a linear transformation of $\Phi$ by the bounded operator $z \otimes \cdot: L_2(U, H) \to L_2(U, \R) = U^*$ (identified with $U$ via the Riesz isomorphism), $A \mapsto A^* z$, with operator norm $\norm{z}_H$. Linearity preserves Malliavin differentiability with the same Sobolev exponent, yielding $\Phi^* z \in \mathbb{D}^{1,p}(\HW)$ with the stated norm bound.

By the Meyer inequality~\cite{nualart2006malliavin}, $\norm{\delta_U(\eta)}_{L^2(\Omega)} \leq C_{\mathrm{Meyer}}\norm{\eta}_{\mathbb{D}^{1,2}(\HW)}$ for every $\eta \in \mathbb{D}^{1,2}(\HW)$. Applied to $\eta = \Phi^* z$ and combined with~(i), this gives $\norm{\phi_\Phi(z)}_{L^2(\Omega)} \leq C_{\mathrm{Meyer}}\norm{z}_H\norm{\Phi}_{\mathbb{D}^{1,2}}$, so $\phi_\Phi$ is bounded $H \to L^2(\Omega)$. For the Hilbert--Schmidt property, fix any deterministic orthonormal basis $\{\tilde e_k\}_{k \geq 1}$ of $H$ and compute
\[
  \norm{\phi_\Phi}_{\mathrm{HS}(H,L^2(\Omega))}^2 \;=\; \sum_{k=1}^\infty\norm{\delta_U(\Phi^*\tilde e_k)}_{L^2(\Omega)}^2 \;\leq\; C_{\mathrm{Meyer}}^2\sum_{k=1}^\infty\norm{\Phi^*\tilde e_k}_{\mathbb{D}^{1,2}(\HW)}^2.
\]
Expanding the $\mathbb{D}^{1,2}(\HW)$-norm as the sum of the $\HW$-norm and the Malliavin-derivative norm, summing over $k$, and using Fubini--Tonelli together with the Hilbert--Schmidt adjoint identity $\norm{A^*}_{L_2(H,U)} = \norm{A}_{L_2(U,H)}$:
\begin{align*}
  \sum_k\E\int_0^t\norm{\Phi_r^*\tilde e_k}_U^2\,\mathrm{d}r &\;=\; \E\int_0^t\sum_k\norm{\Phi_r^*\tilde e_k}_U^2\,\mathrm{d}r \\
  &\;=\; \E\int_0^t\norm{\Phi_r^*}_{L_2(H,U)}^2\,\mathrm{d}r \;=\; \E\int_0^t\norm{\Phi_r}_{L_2(U,H)}^2\,\mathrm{d}r,
\end{align*}
\begin{multline*}
  \sum_k\E\int_0^t\!\!\int_0^t\norm{\D_s(\Phi_r^*\tilde e_k)}_{L_2(U,U)}^2\,\mathrm{d}s\,\mathrm{d}r \\
  \;=\; \E\int_0^t\!\!\int_0^t\norm{\D_s\Phi_r}_{L_2(U,\,L_2(U,H))}^2\,\mathrm{d}s\,\mathrm{d}r,
\end{multline*}
where the second identity uses that $\D_s(\Phi_r^*\tilde e_k) = (\D_s\Phi_r)^*\tilde e_k$ (the Malliavin derivative commutes with the bounded operation $A \mapsto A^*\tilde e_k$) combined with the Hilbert--Schmidt factorisation $\sum_k\norm{(\D_s\Phi_r)^*\tilde e_k}_{L_2(U,U)}^2 = \norm{\D_s\Phi_r}_{L_2(U, L_2(U,H))}^2$ (the same adjoint identity applied to the $L_2(U,H)$-valued Malliavin derivative). Summing the two contributions yields
\[
  \norm{\phi_\Phi}_{\mathrm{HS}(H, L^2(\Omega))}^2 \;\leq\; C_{\mathrm{Meyer}}^2\,\norm{\Phi}_{\mathbb{D}^{1,2}(L^2([0,t]; L_2(U,H)))}^2,
\]
which is finite by hypothesis. The standard isometric isomorphism $\mathrm{HS}(H, L^2(\Omega)) \cong L^2(\Omega; H)$ (via the tensor-product identification $\sum_k \xi_k\otimes\tilde e_k \leftrightarrow \omega \mapsto \sum_k\xi_k(\omega)\tilde e_k$, with $\xi_k \in L^2(\Omega; \R)$) then represents $\phi_\Phi$ by the unique kernel $\mathcal{K}_\Phi(\omega) = \sum_k\delta_U(\Phi^*\tilde e_k)(\omega)\tilde e_k \in L^2(\Omega; H)$, giving~\eqref{eq:HS-kernel}. The $L^p$-bound on $\mathcal{K}_\Phi$ follows from the Hilbert-valued $L^p$-Meyer inequality~\cite{nualart2006malliavin} applied to $\Phi \in \mathbb{D}^{1,p}$. At the $L^p$ level, $\norm{\delta_U(\Phi^* z)}_{L^p(\Omega)} \leq C_{p,\mathrm{Meyer}}\allowbreak\norm{z}_H\allowbreak\norm{\Phi}_{\mathbb{D}^{1,p}}$, so the kernel $\mathcal{K}_\Phi$ (reconstructed basis-wise as above) lies in $L^p(\Omega; H)$ with the stated norm bound.

For deterministic $z$, \eqref{eq:HS-kernel} is definitional. For jointly measurable $Z = \sum_k Z_k\tilde e_k$ with scalar components $Z_k = \ip{Z}{\tilde e_k}_H$, the pointwise inner product $\ip{Z(\omega)}{\mathcal{K}_\Phi(\omega)}_H = \sum_k Z_k(\omega)\,\delta_U(\Phi^*\tilde e_k)(\omega)$ is a direct cylindrical evaluation of the linear map $\phi_\Phi$ at $z = Z(\omega)$, independent of any Malliavin regularity of $Z$. The relation to the Skorokhod integral $\delta_U(\Phi^* Z)$ (which involves the Malliavin-derivative of $Z$ via the product rule $\delta_U(\xi u) = \xi\delta_U(u) - \ip{\D\xi}{u}_{\HW}$ applied term-wise in the expansion $\Phi^* Z = \sum_k Z_k\Phi^*\tilde e_k$) is
\begin{multline*}
  \delta_U(\Phi^* Z) \;=\; \sum_k\bigl[Z_k\delta_U(\Phi^*\tilde e_k) - \textstyle\int_0^t\ip{\D_r Z_k}{\Phi_r^*\tilde e_k}_U\,\mathrm{d}r\bigr] \\
  \;=\; \ip{Z}{\mathcal{K}_\Phi}_H - \int_0^t\Tr_U[\Phi_r^*\D_r Z]\dr,
\end{multline*}
valid whenever $Z \in \mathbb{D}^{1,2}(H)$ and both sides lie in $L^2(\Omega)$; this is precisely the linear-field substitution identity of Theorem~\ref{thm:linear-sub-D12}. In particular, the substitution value $[\delta_U(w^{\mathrm{lin}}(\cdot)(z))]_{z = Z}$ denotes, unambiguously, the pointwise inner product $\ip{Z}{\mathcal{K}_\Phi}_H$ (not the composed Skorokhod integral $\delta_U(\Phi^* Z)$), and is well-defined for any jointly measurable $Z$, requiring no Malliavin regularity of $Z$ itself beyond the finiteness of the right-hand side of~\eqref{eq:HS-kernel-random} in the appropriate Bochner space.
\end{proof}

\begin{remark}[Malliavin regularity of $\Phi$ at $p^* = 2q/(q-2)$]\label{rem:Phi-all-moments}
Set $p^* := 2q/(q-2)$ with $q \in (2, \infty]$ from Assumption~\ref{ass:nondeg}(ii), and $q^* := p^*/2 = q/(q-2)$, read as $q^* = 1$ at $q = \infty$. The regularity $\Phi \in \mathbb{D}^{1,p^*}(L^2([0,t]; L_2(U,H)))$ is read off the two fibre estimates directly, with no operator norm and no factorisation.

For $\Phi$ itself, the Hilbert--Schmidt fibre estimate~\eqref{eq:fibre-energy-estimate} of Theorem~\ref{thm:Y-wellposed}\textup{(III)}, applied at the $\F_r$-measurable datum $\Theta = \calB(r, X(r))$ and at exponent $q^*$, gives
\[
  \E\Bigl[\sup_{r \leq \tau \leq t}\norm{\Phi_{\tau,r}}^{2q^*}_{L_2(U,H)}\Bigr] \;\leq\; C_{q^*}\,\E\Bigl[\norm{\calB(r,X(r))}^{2q^*}_{L_2(U,H)}\exp\Bigl(c_{q^*}\!\int_r^t(\tilde\rho+\hat\rho)\Bigr)\Bigr],
\]
whose right-hand side is separated by H\"older at the conjugate pair $(1+\eps_{q^*}, (1+\eps_{q^*})/\eps_{q^*})$ into the noise moment of (SC4.2)$_{q^*}$ at its exponent $2\kappa_{q^*}/\eps_{q^*}$ and the exponential moment of~(SC2), the constant being admitted by~\eqref{eq:cstar-requirement}. Integrating in $r \in [0,t]$ and using Jensen in time,
\[
  \E\Bigl[\Bigl(\int_0^t\norm{\Phi_{t,r}}^2_{L_2(U,H)}\,\mathrm{d}r\Bigr)^{\!q^*}\Bigr] \;<\; \infty,
\]
that is $\Phi \in L^{p^*}(\Omega; L^2([0,t]; L_2(U,H)))$.

For $\D\Phi$, Lemma~\ref{lemma:DPhi} identifies the Malliavin derivative $\D_s\Phi_{t,r}$ with the second Malliavin fibre for almost every $r$, which is all that the integral over the square $[0,t]^2$ below requires, and the moment of that fibre at the same exponent is furnished by Proposition~\ref{prop:W-moments},
\[
  \E\Bigl[\sup_{\max(r,s) \leq \tau \leq t}\norm{\calW(\tau;r,s)}^{2q^*}_{L_2(U,\,L_2(U,H))}\Bigr] \;\leq\; C_{q^*},
\]
uniformly in $(r,s)$. Integrating over the square $[0,t]^2$ gives
\begin{align*}
  \D\Phi \;\in\; L^{p^*}\bigl(\Omega;\, L^2([0,t]^2;\, L_2(U, L_2(U,H)))\bigr),
\end{align*}
whence $\Phi \in \mathbb{D}^{1,p^*}$ and, by the quadratic representation of Proposition~\ref{prop:D-gamma}, $\gamma_t \in \mathbb{D}^{1,p^*/2}(L_2(H,H))$.

The chain is therefore $\calW \Rightarrow \D\Phi \Rightarrow \Phi \in \mathbb{D}^{1,p^*} \Rightarrow \gamma_t \in \mathbb{D}^{1,p^*/2}$, and each link is an estimate on a fibre at the exponent the next link consumes.
\end{remark}

\begin{theorem}[Stability of $v_h$ under perturbation]\label{thm:covering-stability}
Let $X^{(n)}$ be a sequence of variational solutions to \eqref{eq:SPDE} with coefficients $(\calA^{(n)}, \calB^{(n)})$ satisfying Assumptions \ref{ass:LR}, \ref{ass:diff2}, \ref{ass:SC}, and \ref{ass:nondeg} uniformly in $n$. Let $Y^{(n)}$, $\gamma_t^{(n)}$, and $v_h^{(n)}$ denote the associated first variation, Malliavin covariance, and covering field, and write $\Phi_r^{(n)} := Y^{(n)}(t,r)\calB^{(n)}(r, X^{(n)}(r))$ and $\Phi_r := Y(t,r)\calB(r, X(r))$. Fix $q \in (2,\infty)$ from clause~\textup{(iii)} below, and set $k := q/(q-2)$. Suppose:
\begin{enumerate}[(i)]
  \item $X^{(n)}(t) \to X(t)$ in $L^2(\Omega; H)$;
  \item \textup{(Strong joint $L^{2k}$-convergence of $\Phi$.)} 
  \begin{equation}\label{eq:cov-stab-strong}
    \E\Bigl[\Bigl(\int_0^t \norm{\Phi^{(n)}_r - \Phi_r}^2_{L_2(U,H)}\,\mathrm{d}r\Bigr)^{\!k}\Bigr] \;\longrightarrow\; 0 \qquad (n \to \infty),
  \end{equation}
  with the uniform majorisation
  \begin{equation}\label{eq:cov-stab-uniform}
    \sup_n \E\Bigl[\Bigl(\int_0^t \norm{\Phi^{(n)}_r}^2_{L_2(U,H)}\,\mathrm{d}r\Bigr)^{\!k}\Bigr] \;<\; \infty;
  \end{equation}
  \item $h \in \Ran(\gamma_t)$ a.s., $\E[\norm{\gamma_t^{\dagger}h}_H^q] < \infty$, and $\sup_n \E[\norm{(\gamma_t^{(n)})^{\dagger}h}_H^q] < \infty$, for the same $q > 2$.
\end{enumerate}
Then $v_h^{(n)} \to v_h$ in $L^2(\Omega; \HW)$. If, additionally,
\begin{enumerate}[(i),resume]
  \item each $v_h^{(n)}$ lies in $\Dom(\delta_U)$ and the sequence $\{\delta_U(v_h^{(n)})\}_{n \ge 1}$ is $L^2(\Omega)$-Cauchy,
\end{enumerate}
then $v_h \in \Dom(\delta_U)$ and $\delta_U(v_h^{(n)}) \to \delta_U(v_h)$ in $L^2(\Omega)$.
\end{theorem}

Hypothesis~(iv) is verified separately for each approximating sequence (typically via Lemma~\ref{lemma:D-gamma-inv-h} applied at each $n$ with uniform-in-$n$ estimates); Theorem~\ref{thm:covering-stability} then identifies the limit. This is the mechanism used in Proposition~\ref{prop:singular-p-Lap}.

\begin{proof}
Write
\[
  \tilde{h}^{(n)} = (\gamma_t^{(n)})^{\dagger}h, \qquad v_h^{(n)}(r) = (\Phi_r^{(n)})^{*}\,\tilde{h}^{(n)} .
\]

Hypothesis~(ii) implies the weaker $L^1(\Omega; L_1(H))$-convergence $\gamma_t^{(n)} \to \gamma_t$. Indeed,
\begin{align*}
  \E\bigl[\norm{\gamma_t^{(n)} - \gamma_t}_{L_1(H)}\bigr]
  &\le \E\Bigl[\int_0^t \norm{\Phi^{(n)}_r(\Phi^{(n)}_r)^* - \Phi_r\Phi_r^*}_{L_1(H)}\,\mathrm{d}r\Bigr]\\
  &\le \E\Bigl[\int_0^t \norm{\Phi^{(n)}_r - \Phi_r}_{L_2(U,H)}\bigl(\norm{\Phi^{(n)}_r}_{L_2} + \norm{\Phi_r}_{L_2}\bigr)\,\mathrm{d}r\Bigr]\\
  &\le \Bigl(\E\Bigl[\int_0^t \norm{\Phi^{(n)}_r - \Phi_r}^2\,\mathrm{d}r\Bigr]\Bigr)^{1/2}\Bigl(\E\Bigl[\int_0^t (\norm{\Phi^{(n)}_r} + \norm{\Phi_r})^2\,\mathrm{d}r\Bigr]\Bigr)^{1/2},
\end{align*}
and the first factor tends to $0$ by hypothesis~(ii) (which gives convergence in $L^k$, hence in $L^1$, of the time-integrated HS norm squared), while the second is uniformly bounded by~\eqref{eq:cov-stab-uniform} (also implying $L^1$-bound).

We decompose the covering field difference as
\[
  v_h^{(n)}(r) - v_h(r) = (\Phi_r^{(n)})^*(\tilde{h}^{(n)} - \tilde{h}) + (\Phi_r^{(n)} - \Phi_r)^*\tilde{h}
\]
and bound each term separately at the level of the time-integrated $\HW$-norm squared.

By Cauchy--Schwarz and the operator inequality $\norm{(\Phi^{(n)}_r - \Phi_r)^*\tilde h}_U \le \norm{\Phi^{(n)}_r - \Phi_r}_{L_2(U,H)}\norm{\tilde h}_H$,
\[
  \int_0^t \norm{(\Phi^{(n)}_r - \Phi_r)^*\tilde h}_U^2\,\mathrm{d}r \;\le\; \norm{\tilde h}_H^2\,\int_0^t \norm{\Phi^{(n)}_r - \Phi_r}_{L_2(U,H)}^2\,\mathrm{d}r.
\]
Applying H\"older's inequality with conjugate exponents $(q/2,\,k)$ (which sum to $1$ in reciprocals, $2/q + 1/k = 2/q + (q-2)/q = 1$),
\begin{equation}\label{eq:cov-stab-second-term}
  \E\Bigl[\int_0^t \norm{(\Phi^{(n)}_r - \Phi_r)^*\tilde h}_U^2\,\mathrm{d}r\Bigr] \;\le\; \bigl(\E[\norm{\tilde h}_H^q]\bigr)^{2/q}\,\Bigl(\E\Bigl[\Bigl(\int_0^t\norm{\Phi^{(n)}_r - \Phi_r}^2\,\mathrm{d}r\Bigr)^k\Bigr]\Bigr)^{1/k}.
\end{equation}
The first factor is finite by hypothesis~(iii), which bounds the moment of the $n$-independent limiting object $\tilde h$ explicitly. The second factor tends to $0$ by hypothesis~(ii), eq.~\eqref{eq:cov-stab-strong}.

We use the algebraic identity from Theorem~\ref{thm:covering-L2-cont}. Since $\gamma_t^{(n)}\tilde{h}^{(n)} = h = \gamma_t\tilde{h}$,
\begin{align*}
  \int_0^t \norm{(\Phi_r^{(n)})^*(\tilde{h}^{(n)} - \tilde{h})}_U^2\dr
  &= \ip{\tilde{h}^{(n)} - \tilde{h}}{\gamma_t^{(n)}(\tilde{h}^{(n)} - \tilde{h})}_H\\
  &= \ip{\tilde{h}^{(n)} - \tilde{h}}{h - \gamma_t^{(n)}\tilde{h}}_H\\
  &= \ip{\tilde{h}^{(n)} - \tilde{h}}{(\gamma_t - \gamma_t^{(n)})\tilde{h}}_H,
\end{align*}
where the last step uses $h = \gamma_t\tilde{h}$. Bounding the inner product by Cauchy--Schwarz in $H$ pointwise, then by H\"older's inequality in $L^q(\Omega) \times L^{q'}(\Omega)$ (with $q' = q/(q-1)$, the H\"older conjugate of $q$),
\begin{equation}\label{eq:cov-stab-first-term}
  \E\Bigl[\int_0^t \norm{(\Phi_r^{(n)})^*(\tilde{h}^{(n)} - \tilde{h})}_U^2\dr\Bigr] \le \bigl(\E[\norm{\tilde{h}^{(n)} - \tilde{h}}_H^q]\bigr)^{1/q}\,\bigl(\E[\norm{(\gamma_t - \gamma_t^{(n)})\tilde{h}}_H^{q'}]\bigr)^{1/q'}.
\end{equation}
The first factor of \eqref{eq:cov-stab-first-term} is bounded by $(\sup_n \E[\norm{\tilde{h}^{(n)}}_H^q] + \E[\norm{\tilde{h}}_H^q])^{1/q} < \infty$ uniformly in $n$ by hypothesis~(iii). It remains to show that the second factor tends to $0$, i.e., $\E[\norm{(\gamma_t - \gamma_t^{(n)})\tilde h}^{q'}_H] \to 0$ as $n \to \infty$.

Bound the operator-vector product as
\begin{align*}
  \norm{(\gamma_t - \gamma_t^{(n)})\tilde h}_H &\;\le\; \norm{\gamma_t - \gamma_t^{(n)}}_{L_1(H)}\,\norm{\tilde h}_H \\
  &\;\le\; \Bigl(\int_0^t \norm{\Phi^{(n)}_r - \Phi_r}_{L_2}\bigl(\norm{\Phi^{(n)}_r}_{L_2} + \norm{\Phi_r}_{L_2}\bigr)\,\mathrm{d}r\Bigr)\,\norm{\tilde h}_H,
\end{align*}
where the first inequality uses $\norm{A}_{L(H)} \le \norm{A}_{L_1(H)}$ for trace-class $A$, and the second comes from the trace-norm bound on $\gamma_t - \gamma_t^{(n)}$ established at the start of the proof. By Cauchy--Schwarz in the time variable,
\[
  \int_0^t \norm{\Phi^{(n)}_r - \Phi_r}_{L_2}(\norm{\Phi^{(n)}_r} + \norm{\Phi_r}) \,\mathrm{d}r \;\le\; A_n^{1/2}\,B_n^{1/2},
\]
where $A_n := \int_0^t \norm{\Phi^{(n)} - \Phi}^2_{L_2}\,\mathrm{d}r$ and $B_n := \int_0^t (\norm{\Phi^{(n)}} + \norm{\Phi})^2\,\mathrm{d}r$.
Hence $\norm{(\gamma_t - \gamma_t^{(n)})\tilde h}^{q'}_H \le A_n^{q'/2} B_n^{q'/2} \norm{\tilde h}^{q'}_H$. Apply the triple H\"older inequality with exponents $(\alpha, \alpha, \gamma)$ where $\gamma := q-1$ (so that $1/\gamma = 1/(q-1)$) and $\alpha := 2(q-1)/(q-2)$ (so that $2/\alpha = (q-2)/(q-1)$ and $2/\alpha + 1/\gamma = (q-2)/(q-1) + 1/(q-1) = 1$),
\begin{align*}
  &\E\bigl[A_n^{q'/2} B_n^{q'/2}\norm{\tilde h}^{q'}_H\bigr]\\
  &\quad \le \bigl(\E[A_n^{q'\alpha/2}]\bigr)^{1/\alpha}\,\bigl(\E[B_n^{q'\alpha/2}]\bigr)^{1/\alpha}\,\bigl(\E[\norm{\tilde h}^{q'\gamma}_H]\bigr)^{1/\gamma}.
\end{align*}
The exponents inside the expectations evaluate as follows. For the $A_n$ and $B_n$ factors, $q'\alpha/2 = (q/(q-1))\cdot(q-1)/(q-2) = q/(q-2) = k$. For the $\norm{\tilde h}$ factor, $q'\gamma = (q/(q-1))(q-1) = q$. Hence
\begin{equation}\label{eq:cov-stab-triple}
  \E[\norm{(\gamma_t - \gamma_t^{(n)})\tilde h}^{q'}_H] \;\le\; \bigl(\E[A_n^k]\bigr)^{1/\alpha}\,\bigl(\E[B_n^k]\bigr)^{1/\alpha}\,\bigl(\E[\norm{\tilde h}^q_H]\bigr)^{1/(q-1)}.
\end{equation}
The third factor on the right is finite by hypothesis~(iii), which includes the moment bound on $\tilde h$ itself. The second factor $(\E[B_n^k])^{1/\alpha}$ is uniformly bounded in $n$ by~\eqref{eq:cov-stab-uniform}. Indeed, $B_n^k \le 2^{2k-1}((\int\norm{\Phi^{(n)}}^2)^k + (\int\norm{\Phi}^2)^k)$ pathwise by the convexity inequality $(a+b)^{2k} \le 2^{2k-1}(a^{2k} + b^{2k})$ applied via $B_n^{1/2} \le (\int\norm{\Phi^{(n)}}^2 + 2\int\norm{\Phi^{(n)}}\norm{\Phi} + \int\norm{\Phi}^2)^{1/2} \le \sqrt{2}((\int\norm{\Phi^{(n)}}^2)^{1/2} + (\int\norm{\Phi}^2)^{1/2})$, and the right side raised to the $2k$-th power gives a constant times $(\int\norm{\Phi^{(n)}}^2)^k + (\int\norm{\Phi}^2)^k$, both of which have $\sup_n \E[\cdot] < \infty$ by~\eqref{eq:cov-stab-uniform}. The first factor $(\E[A_n^k])^{1/\alpha}$ tends to $0$ by hypothesis~(ii), eq.~\eqref{eq:cov-stab-strong}.

Combining \eqref{eq:cov-stab-first-term} and \eqref{eq:cov-stab-triple}, the first-term contribution
\[
  \E[\int_0^t\norm{(\Phi^{(n)})^*(\tilde h^{(n)} - \tilde h)}^2 \dr] \to 0.
\]
Together with the second-term bound \eqref{eq:cov-stab-second-term}, this gives $v_h^{(n)} \to v_h$ in $L^2(\Omega; \HW)$.

Under the additional hypothesis~(iv), $\{\delta_U(v_h^{(n)})\}$ converges in $L^2(\Omega)$ to some limit $\xi$. Combined with the $L^2(\Omega;\HW)$-convergence $v_h^{(n)} \to v_h$ just established, the pairs $(v_h^{(n)}, \delta_U(v_h^{(n)})) \to (v_h, \xi)$ converge in the graph topology $L^2(\Omega; \HW) \times L^2(\Omega)$. The Skorokhod operator $\delta_U: \Dom(\delta_U) \subset L^2(\Omega; \HW) \to L^2(\Omega)$ is closed (as the adjoint of the densely defined derivative operator $D$; cf.\ \cite{nualart2006malliavin}), so $v_h \in \Dom(\delta_U)$ and $\delta_U(v_h) = \xi$, yielding $\delta_U(v_h^{(n)}) \to \delta_U(v_h)$ in $L^2(\Omega)$.
\end{proof}

\begin{remark}[Verification of hypothesis~(iv)]\label{rem:cov-stab-D12}
Hypothesis~(iv) of Theorem~\ref{thm:covering-stability}, the $L^2(\Omega)$-Cauchy property of $\{\delta_U(v_h^{(n)})\}$, is the technical input that carries the analytic content in practice, and a more conventional approach via $\mathbb{D}^{1,2}(\HW)$-convergence + Meyer inequality (cf.\ \cite{nualart2006malliavin}) is not available in the central application of the present paper. To see why the $\mathbb{D}^{1,2}$ route fails, observe that for the Tikhonov family $\gamma_t^{(\eps)} = \gamma_t + \eps I$ used in Lemma~\ref{lemma:D-gamma-inv-h}, the Malliavin derivative of the pseudoinverse $\tilde{h}^{(\eps)} = (\gamma_t + \eps I)^{-1}h$ satisfies $\D_s\tilde{h}^{(\eps)} = -(\gamma_t + \eps I)^{-1}(\D_s\gamma_t)(\gamma_t + \eps I)^{-1}h$, which has operator-norm scaling $\eps^{-1}$, so the family $\{v_h^{(\eps)}\}_{\eps > 0}$ is not uniformly bounded in $\mathbb{D}^{1,2}(\HW)$.

For the Tikhonov family, hypothesis~(iv) follows from Lemma~\ref{lemma:D-gamma-inv-h} via arguments~(M) and~(C).
\end{remark}

\section{The two variations}
\label{sec:variation}

The Malliavin derivative of the SPDE solution at the terminal time $t$ is the variational fibre $\Phi_{t,r} = Y(t,r)\calB(r,X(r))$ of Proposition~\ref{prop:D-Xt}, and its further Malliavin derivative is the second Malliavin fibre $\calW(t;s,r) = \D_r\Phi_{t,s}$, which describes how the fibre responds to a perturbation of the noise. These two objects are all the PDE-theoretic content that Theorem~\ref{thm:main} requires, the rest being the Malliavin calculus of the previous section applied to a variational solution whose dependence on the noise is governed by them. We construct both here, establish their well-posedness in the variational framework, and obtain the moment bounds that drive the Malliavin-regularity chain $\calW \Rightarrow \D\Phi \Rightarrow \Phi \in \mathbb{D}^{1, p^*}$ used in the proof of Theorem~\ref{thm:main}. Everything is run on the fibres rather than on the flow, so that the only state-space directions ever perturbed are the columns of the diffusion coefficient and the $V$-regularity the estimates consume is the fibre energy that the dissipative mechanism itself produces.

\subsection{The first variation}\label{subsec:first-var}

We start with $Y(t, r)$, the first variation, which arises as the linearisation of~\eqref{eq:main-SPDE} along the solution path $X$. One should picture $Y(t, r) v$ as the perturbation of $X(t)$ caused by an infinitesimal perturbation $v$ of the state at time $r$. The question is under what hypotheses does $Y(t, r)$ exist as a sensible operator-valued object, and what moment bounds does it satisfy?

\begin{theorem}[The first variation]\label{thm:Y-wellposed}
Under Assumptions~\ref{ass:LR},~\ref{ass:diff1}, and~\ref{ass:SC}, for each $r \in [0,T]$ and $v_0 \in H$, the initial value problem
\begin{equation}\label{eq:Y-IVP}
\begin{aligned}
  \mathrm{d}Y_v(s) + \calA'_u(s, X(s))\,Y_v(s)\ds &= \calB'_u(s, X(s))(Y_v(s))\dW(s), \quad s \in (r, T],\\
  Y_v(r) &= v_0 \in H,
\end{aligned}
\end{equation}
admits a unique $H$-valued continuous variational solution, tested against $V$, in the regularity class determined by the mechanism of~(SC1):
\[
  Y_v(\cdot) \;\in\;
  \begin{cases}
    C([r,T]; H) \ \text{ with } \ \calE_r^T(Y_v) < \infty, & \text{under (SC1)(a)},\\[4pt]
    C([r,T]; H) \cap L^2(r,T; V), & \text{under (SC1)(b)},
  \end{cases}
  \quad \text{$\mathbb{P}$-a.s.},
\]
Under~(a) the controlling quantity is the weighted energy intrinsic to the degenerate linearisation, which for the $p$-Laplacian is $\int_r^T\!\!\int \abs{\nabla X}^{p-2}\abs{\nabla Y_v}^2$; under~(b) it is the unweighted quadratic $V$-energy furnished by $\alpha_1$.
The map $v_0 \mapsto Y_v(t)$ is linear, and we view it as the first variation operator $Y(t,r): H \to L^0(\Omega; H)$, $Y(t,r)v_0 := Y_v(t)$, the analogue in the variational setting of the Jacobian of a stochastic flow \cite{ELWORTHY1994252, kunita1990stochastic}. Under the additional structural exponential-integrability~(SC2), $Y_v \in L^2(\Omega; C([r,T]; H))$ with the mechanism-matched global energy bound ($\E[\calE_r^T(Y_v)] \leq C\norm{v_0}_H^2$ under~(SC1)(a), $Y_v \in L^2(\Omega; L^2(r,T;V))$ under~(SC1)(b)), and $Y(t,r): H \to L^2(\Omega; H)$ is continuous, with the individual-fibre moment bound $\E[\sup_{r \leq s \leq t}\norm{Y(s,r)v_0}_H^q] \leq C_q\norm{v_0}_H^q$ for every $v_0 \in H$ and every $q \in [1, q_{\max}]$ (Proposition~\ref{prop:Y-moments}, $C_q$ independent of $v_0$).

\smallskip\noindent\emph{\textup{(III)} The fibre energy.} Let $q \in [1, q_{\max}/2]$ and let $\Theta$ be $\F_r$-measurable with $\Theta \in L^{2q}(\Omega, \F_r; L_2(U,H))$. Then the process $Y(\cdot, r)\Theta$, defined columnwise by $u \mapsto Y(\cdot,r)(\Theta u)$, satisfies
\begin{align*}
  Y(\cdot,r)\Theta \;\in\; L^{2q}\bigl(\Omega;\, C([r,T];\, L_2(U,H))\bigr),
\end{align*}
takes values in $L_2(U, \mathbb{V}_\tau)$ for almost every $\tau$, and obeys the \emph{fibre energy estimate}
\begin{equation}\label{eq:fibre-energy-estimate}
\begin{aligned}
  \E\Bigl[\sup_{r \leq \tau \leq T}\norm{Y(\tau,r)\Theta}_{L_2(U,H)}^{2q}\Bigr]
  &+ \E\Bigl[\Bigl(\delta\int_r^T \norm{Y(\tau,r)\Theta}^2_{L_2(U, \mathbb{V}_\tau)}\,\mathrm{d}\tau\Bigr)^{\!q}\Bigr]\\
  &\leq\; C_q\,\E\Bigl[\norm{\Theta}_{L_2(U,H)}^{2q}\,\exp\Bigl(c_q\!\int_r^T(\tilde\rho + \hat\rho)\,\mathrm{d}\sigma\Bigr)\Bigr],
\end{aligned}
\end{equation}
with $\delta$ the dissipation constant of~\eqref{eq:SC1-unified} and $C_q, c_q$ the constants of Proposition~\ref{prop:Y-moments}. Taken at the $\F_r$-measurable datum $\Theta = \calB(r, X(r))$, this is a statement about the variational fibre $\Phi_{\tau,r}$ of~\eqref{eq:Phi-fibre}, which then lies in $L_2(U, \mathbb{V}_\tau)$ for almost every $\tau$ and whose fibre energy~\eqref{eq:fibre-energy} obeys
\begin{equation}\label{eq:Phi-energy}
  \E\bigl[\mathfrak{E}_r^{\,q}\bigr] \;\leq\; \delta^{-q}\,C_q\,\E\Bigl[\norm{\calB(r,X(r))}^{2q}_{L_2(U,H)}\,\exp\Bigl(c_q\!\int_r^T(\tilde\rho+\hat\rho)\,\mathrm{d}\sigma\Bigr)\Bigr],
\end{equation}
the right-hand side being finite uniformly in $r \in [0,T]$ under~(SC2) together with the noise moment~\eqref{eq:B-growth-m} of~(D2$'$), and under Assumption~\ref{ass:SC-raised} at every $q \leq \kappa_{q^*}$.
\end{theorem}

\begin{remark}[Two-layer construction]\label{rem:Y-two-layers}
The proof proceeds in two layers. Under~(SC1) alone, \eqref{eq:Y-IVP} admits a unique pathwise solution in the tangent class $\calY_r$ a.s., obtained by patching Galerkin approximations on the stopping-time intervals $\tau_M := \inf\{\tau \in [r,T] : \int_r^\tau (\tilde\rho + \norm{\calB'_u}^2)\,d\sigma > M\}$ with $M \to \infty$. With the additional exponential integrability~(SC2), the random Gronwall step closes globally and the pathwise solution satisfies $Y_v \in L^2(\Omega; C([r,T]; H))$ together with the mechanism-matched global energy bound, $\E[\calE_r^T(Y_v)] \leq C\norm{v_0}_H^2$ under~(a) and $Y_v \in L^2(\Omega; L^2(r,T;V))$ under~(b). The stronger property $Y(t,r) \in L(H)$ a.s.\ requires case-by-case structural input and is nowhere used; layer~(III) below runs the same estimate on the fibre $Y(\cdot,r)\Theta$ at a Hilbert--Schmidt datum, which is the form every later statement consumes, and which produces by itself the $V$-regularity that the second derivatives of the coefficients require.
\end{remark}

\begin{proof}
The equation \eqref{eq:Y-IVP} is a linear stochastic evolution equation with random, time-dependent coefficients $\calA'_u(s, X(s))$ and $\calB'_u(s, X(s))$. We verify the Liu--R\"ockner conditions (H1)--(H4) for these linearised coefficients, treated as operators in the unknown $Y_v \in V$.

For fixed $v_1, v_2, w \in V$ and $t \in [0,T]$, the map
\[
  \lambda \mapsto \dual{\calA'_u(t, X(t))(v_1 + \lambda v_2)}{w} = \dual{\calA'_u(t, X(t))v_1}{w} + \lambda\dual{\calA'_u(t, X(t))v_2}{w}
\]
is affine in $\lambda$ (because $\calA'_u(t, X(t)) \in L(V, V^*)$ is linear in its argument), hence trivially continuous.

By the linearised local monotonicity derived from (H2) (cf.\ Remark \ref{rem:no-D3-old}), for all $v, w \in V$,
\begin{align*}
  &-2\dual{\calA'_u(t, X(t))(v - w)}{v - w} + \norm{\calB'_u(t, X(t))(v - w)}_{L_2(U,H)}^2 \\
  &\qquad \leq \tilde{\rho}(t, X(t))\norm{v - w}_H^2,
\end{align*}
where $\tilde{\rho}(t, X(t))$ is integrable in $t$ $\mathbb{P}$-a.s.\ by the moment bounds on $X$ from Theorem \ref{thm:LR-wellposedness}.

We establish the generalised coercivity bound. From the linearised local monotonicity \eqref{eq:linearised-mono-derived} (derived from (H2) in Remark \ref{rem:no-D3-old}) evaluated at $w = 0$:
\begin{equation}\label{eq:lin-mono-at-zero}
  -2\dual{\calA'_u(t, X(t))v}{v} + \norm{\calB'_u(t, X(t))(v)}_{L_2(U,H)}^2 \leq \tilde{\rho}(t, X(t))\norm{v}_H^2.
\end{equation}
This is the coercivity condition for the linearised operator in generalised form. Compared to the original (H3), the $-\alpha\norm{v}_V^p$ term is absent. The linearised monotonicity (derived by differentiating a nonlinear inequality) cannot, in general, produce the strict $V$-coercivity $-\alpha_1\norm{v}_V^p$ with $\alpha_1 > 0$, because the linearised operator may degenerate where the solution's gradient vanishes (e.g., the linearised $p$-Laplacian involves $\abs{\nabla X}^{p-2}$, which vanishes at $\nabla X = 0$).

The regularity appropriate to a linear equation, the tangent energy~\eqref{eq:tangent-energy}, quadratic rather than $p$-homogeneous, still follows from the structure of the linearised equation, by one of two mechanisms depending on the equation:
\begin{enumerate}[(a)]
\item (Dissipative linearisation.) If $D_\sigma(v) := \dual{\calA'_u(\sigma, X(\sigma))v}{v} \geq 0$ for all $v \in V$ (as for the $p$-Laplacian, see Proposition \ref{prop:p-Lap-diff}), then the energy identity gives the dissipation budget $2\int_r^t D_\sigma(Y_v)\ds \leq \norm{v_0}_H^2 + C$; the weighted quadratic energy $\int_r^t D_\sigma(Y_v)\,\mathrm{d}\sigma$ (for the $p$-Laplacian, $\int\!\!\int \abs{\nabla X}^{p-2}\abs{\nabla Y_v}^2$) is the intrinsic $V$-type regularity of the linear equation, and it, not an unweighted $\norm{Y_v}_V^p$-integral, which the quadratic homogeneity of the linear equation does not produce, is what every estimate below consumes under mechanism~(a), through the dominated Cauchy--Schwarz of~(SC1)(a).
\item (Strict linearised coercivity.) If $-2\dual{\calA'_u v}{v} + \norm{\calB'_u(v)}_{L_2}^2 \leq C_1(t,\omega)\norm{v}_H^2 - \alpha_1\norm{v}_V^2$ for some $\alpha_1 > 0$ (as for the 2D Navier--Stokes with $\alpha_1 = \nu$; see Proposition \ref{prop:NS-verify}), then the standard quadratic energy estimate yields $Y_v \in L^2(r, T; V)$ directly.
\end{enumerate}
Under mechanism~(a) the conclusion at this stage is the weighted tangent-energy regularity, and under mechanism~(b) the strict coercivity gives $Y_v \in L^2(r,T;V)$ directly; the fibre norm~\eqref{eq:fibre-norm} is the common expression of the two, and part~(III) below runs the same estimate on the fibre. For the abstract framework, the key output is the $H$-continuity $Y_v \in C([r,T]; H)$, which follows from \eqref{eq:lin-mono-at-zero} and the generalised coercivity framework of \cite{liu2015stochastic}.

From (D1):
\[
  \norm{\calA'_u(t, X(t))v}_{V^*} \leq \bigl(C + \rho_2(t, X(t))\bigr)\bigl(1 + \norm{X(t)}_V^{(p-2)^+}\bigr)\norm{v}_V.
\]
The path-integrability clause \eqref{eq:rho2-path} of~(D1) gives, with $p' = p/(p-1)$,
\[
  \int_0^T \bigl[(C + \rho_2(t,X(t)))(1 + \norm{X(t)}_V^{(p-2)^+})\bigr]^{p'}\dt \;<\; \infty
  \qquad \mathbb{P}\text{-a.s.}
\]
directly, so the growth condition is satisfied pathwise.

The preceding estimates verify that the linearised equation \eqref{eq:Y-IVP}, viewed as a linear variational SPDE with random coefficients, satisfies the linearised hemicontinuity, the linearised local monotonicity \eqref{eq:lin-mono-at-zero}, and the (D1)--(D2) growth bounds. The tangent form~\eqref{eq:tangent-form} supplies the remaining structure without any $p$-\hspace{0pt}homogeneous coercivity. By the form Cauchy\nobreakdash--Schwarz bound of~(SC1)(a) (respectively the $V$-domination of~(SC1)(b)), the duality pairing $\dual{\calA'_u(\sigma,X(\sigma))w}{\varphi}$ is well defined and bounded by $D_\sigma(w)^{1/2}D_\sigma(\varphi)^{1/2}$ for every test function $\varphi \in V$, with $D_\sigma(\varphi) \leq \norm{\calA'_u(\sigma,X(\sigma))}_{L(V,V^*)}\norm{\varphi}_V^2$ integrable in $\sigma$ by the path-integrability~\eqref{eq:rho2-path} of~(D1). Hence the weak formulation of~\eqref{eq:Y-IVP} closes on the tangent class $\calY_r$, tested against $V$, finiteness of the tangent energy being exactly what makes every term of the variational identity integrable. Under mechanism~(a), pathwise existence and Galerkin convergence on $[r, \tau_M \wedge T]$ are the well-posedness clause of the non-autonomous form datum in~(SC1); under mechanism~(b) they are the classical variational scheme on the fixed triple $V \hookrightarrow H \hookrightarrow V^*$, whose coercivity is the strict quadratic estimate of~(b). In both cases the finite-dimensional approximants obey the corresponding a priori bound (Itô's formula on $\norm{\cdot}_H^2$ together with the $\delta$-slack of~(a), or the $\alpha_1$-coercivity of~(b)) uniformly on $[r,\tau_M \wedge T]$, which identifies the limit and furnishes the quantitative estimate. Uniqueness follows from the same energy identity applied to the difference of two solutions, whose data vanish, using only the linearised monotonicity \eqref{eq:lin-mono-at-zero}. This is the linear specialisation of the Liu--R\"ockner scheme \cite{liu2015stochastic}, run in the tangent energy rather than in $L^p(r,T;V)$. Pathwise uniqueness ensures that the local solutions on $[r, \tau_M \wedge T]$ for different $M$ patch consistently, yielding a unique pathwise solution $Y_v \in \calY_r$ a.s.\ on the full interval $[r, T]$. The stopping times exhaust $[r,T]$ a.s.\ by the local integrability of $\tilde\rho$ along solution paths (from Theorem~\ref{thm:LR-wellposedness}) and the (D2) growth bound combined with the a priori estimate \eqref{eq:apriori}. This establishes layer~(I); at this stage the solution is constructed pathwise and need not lie in $L^2(\Omega; C([r,T]; H))$, because the random structural constants on $[r, \tau_M]$ depend on $M$ and the bound $\|v_0\|^2 e^{c M}$ blows up as $M \to \infty$.

For layer~(II), under the additional structural exponential-integrability hypothesis~(SC2), the random Gronwall step in the energy estimate closes globally on $[r, T]$ rather than only on stopping-time intervals, giving $\E[\sup_{s \leq t}\norm{Y_v(s)}_H^2] \leq \|v_0\|^2_H \cdot \E[\exp(c\int_r^t\tilde\rho\,d\sigma)] < \infty$ by~(SC2). Hence the pathwise solution from layer~(I) belongs to $L^2(\Omega; C([r,T]; H))$, and the same estimate, retaining on the left-hand side of the energy identity the non-negative tangent term under~(a) or the term $\alpha_1\int\norm{Y_v}_V^2$ under~(b), gives $\E[\calE_r^T(Y_v)] \leq C\norm{v_0}_H^2$ in the first case and $Y_v \in L^2(\Omega; L^2(r,T;V))$ in the second.

The equation \eqref{eq:Y-IVP} is linear in $Y_v$. If $Y_{v_1}$ and $Y_{v_2}$ are solutions with initial conditions $v_1$ and $v_2$, then $Y_{v_1} + Y_{v_2}$ solves the equation with initial condition $v_1 + v_2$, and $\lambda Y_{v_1}$ solves it with initial condition $\lambda v_1$. Hence the map $v_0 \mapsto Y_v(t)$ is linear; under layer~(II), it is bounded as a linear operator from $H$ to $L^2(\Omega; H)$ via the $L^2$-energy estimate just established. This $L^2$-boundedness does not imply pathwise $L(H)$-membership of $Y(t,r)$, which would require a pointwise-in-$\omega$ uniform bound $\sup_{\|v\|=1}\|Y_v(t)(\omega)\| < \infty$ that is strictly stronger than the $L^2$-continuity in $v_0$ given by the abstract variational framework.

The individual-fibre moment bound $\E[\sup_{r \leq s \leq t}\norm{Y(s,r)v_0}_H^q] \leq C_q\norm{v_0}_H^q$ holds for every $v_0 \in H$ and every $q \in [1, q_{\max}]$ by Proposition~\ref{prop:Y-moments} (under \eqref{eq:SC-exp-moment}). As discussed there, this individual-fibre bound is strictly weaker than the operator-norm moment bound $\E[\norm{Y(t,r)}_{L(H)}^q]$, which would require additional structural input beyond the abstract framework and is nowhere used.

For~(III), fix an orthonormal basis $(f_j)$ of $U$, let $N \in \N$, and set
\begin{align*}
  S_N(\tau) \;:=\; \sum_{j \leq N}\norm{Y(\tau,r)(\Theta f_j)}_H^2 .
\end{align*}
Each summand solves~\eqref{eq:Y-IVP} with the $\F_r$-measurable initial datum $\Theta f_j$, so It\^o's formula applies to $S_N$ on the space $\ell_2^N(H)$, and both the principal part and the linearised monotonicity act diagonally in the index $j$. Applying the unified dissipative inequality~\eqref{eq:SC1-unified} to each summand and adding,
\[
  S_N(\tau) + \delta\int_r^\tau \sum_{j \leq N}\norm{Y(\sigma,r)(\Theta f_j)}^2_{\mathbb{V}_\sigma}\,\mathrm{d}\sigma \;\leq\; S_N(r) + \int_r^\tau \tilde\rho\,S_N(\sigma)\,\mathrm{d}\sigma + M_N(\tau),
\]
with $M_N$ a continuous local martingale vanishing at $\tau = r$ and $S_N(r) \leq \norm{\Theta}^2_{L_2(U,H)}$. Discarding the dissipation on the left gives the Gronwall form~\eqref{eq:eta-Gronwall-form} for $S_N$, so the argument of Proposition~\ref{prop:Y-moments}, run at exponent $q$ with the random initial datum as in Remark~\ref{rem:Y-moments-random}, bounds $\E[\sup_\tau S_N(\tau)^q]$ by the right-hand side of~\eqref{eq:fibre-energy-estimate}, uniformly in $N$. Retaining the dissipation instead, and using $(a+b)^q \leq 2^{q-1}(a^q + b^q)$ together with the same bound on the remaining terms, gives the identical bound for $\E[(\delta\int_r^T\sum_{j \leq N}\norm{\cdot}^2_{\mathbb{V}_\sigma})^q]$. Both left-hand sides increase with $N$ to their basis-free limits, by the quadratic-form property recorded after~\eqref{eq:fibre-norm}, so monotone convergence yields~\eqref{eq:fibre-energy-estimate}; in particular the limit is finite almost surely, which is the asserted membership $Y(\tau,r)\Theta \in L_2(U, \mathbb{V}_\tau)$ for almost every $\tau$. Specialising to $\Theta = \calB(r,X(r))$, which lies in $L^{2q}(\Omega, \F_r; L_2(U,H))$ by~\eqref{eq:B-growth-m}, gives~\eqref{eq:Phi-energy}.
\end{proof}

\begin{proposition}[Moments of the first variation]\label{prop:Y-moments}
Suppose the assumptions of Theorem~\ref{thm:Y-wellposed} are in force (in particular,~(SC2) of Assumption~\ref{ass:SC}, at the operating constant $c_*$ of~(SC2), cf.~\eqref{eq:SC-exp-moment}), where $\tilde{\rho} := \rho + \rho^*$ is the linearised local-monotonicity coefficient of Remark~\ref{rem:no-D3-old} and $\hat\rho(\tau,u) := C(1 + \norm{u}_V^{2\beta})$ is the polynomial-growth bound from (D2). Then for every $q \in [2, q_{\max}]$, with $q_{\max}$ as in~\eqref{eq:qmax-def}, there exists $C_q > 0$ such that for every $v \in H$,
\begin{equation}\label{eq:Y-moment-fibre}
  \E\Bigl[\sup_{r \leq s \leq t} \norm{Y(s,r)v}_H^q\Bigr] \;\leq\; C_q\,\norm{v}_H^q.
\end{equation}
\end{proposition}

\begin{remark}[Role of $\hat\rho$ in~(SC2); fibre vs.\ operator-norm bounds; verification]\label{rem:Y-moments-role}
(i) Why $\hat\rho$ enters~(SC2). The presence of $\hat\rho$ accounts for the additional drift contribution $(2q-2)\norm{\calB'_u(\eta)}^2$ that arises when It\^o's formula is applied to $\norm{\eta}^{2q}$ for $q > 1$. The linearised monotonicity of Remark~\ref{rem:no-D3-old} bounds only $-2\dual{\calA'_u\eta}{\eta} + \norm{\calB'_u(\eta)}^2$ (coefficient $1$ on $\norm{\calB'_u(\eta)}^2$), so higher-order It\^o produces an extra $(2q-2)\norm{\calB'_u(\eta)}^2 \leq (2q-2)\hat\rho\norm{\eta}^2$ term that must be absorbed via the exponential moment of $\hat\rho$. For $q = 1$ this term vanishes and only the moment of $\tilde\rho$ is required; the range of admissible $q$ is exactly $[1, q_{\max}]$ of~\eqref{eq:qmax-def}, and \eqref{eq:cstar-requirement} guarantees that it contains every exponent used in the sequel.

(ii) Fibre vs.\ operator-norm bounds. \eqref{eq:Y-moment-fibre} is an individual-fibre estimate, giving the moment of $\|Y(s,r)v\|_H$ for each fixed $v$. The corresponding operator-norm moment bound $\E[\sup_s\norm{Y(s,r)}_{L(H)}^q]$ is in general strictly stronger and does not follow from~\eqref{eq:Y-moment-fibre} by interchanging supremum and expectation; in infinite dimensions, the operator-norm bound requires additional structural input. The downstream applications in Sections~\ref{sec:variation}--\ref{sec:proof-main} are formulated so as to depend only on the Hilbert--Schmidt form~\eqref{eq:fibre-energy-estimate} of the estimate at the datum $\calB(r,X(r))$, that is on the variational fibre, and never on the action of $Y$ on an arbitrary direction of $H$.

(iii) Verification of~\eqref{eq:SC-exp-moment}. The structural exponential-integrability hypothesis is verified case-by-case in Section~\ref{subsec:verification}. For the stochastic $p$-Laplacian with state-independent noise, $\tilde\rho$ and $\hat\rho$ are constants (Proposition~\ref{prop:p-Lap-LR} gives local monotonicity with constant modulus, and $\beta = 0$), so \eqref{eq:SC-exp-moment} holds trivially at every $c$. For 2D Navier--Stokes, $\tilde\rho = (2C/\nu)\norm{X}_V^2$ with additive noise ($\hat\rho$ constant), and the exponential-martingale energy bound quantifies the moment. Writing $R := \int_0^T\norm{X}_V^2\,\mathrm{d}\sigma$, the energy identity gives $2\nu R \leq \norm{x}_H^2 + \operatorname{Tr}Q\cdot T + 2M_T$ with $\langle M\rangle_T \leq (\lambda_{\max}(Q)/\lambda_1)\,R$, and the exponential martingale inequality $\E[e^{aM_T}] \leq (\E[e^{2a^2\langle M\rangle_T}])^{1/2}$ self-improves to
\[
  \E\bigl[e^{cR}\bigr] \;\leq\; \exp\bigl(c(\norm{x}_H^2 + \operatorname{Tr}Q\cdot T)/\nu\bigr) \qquad \text{for every } c \leq \nu^2\lambda_1/(2\lambda_{\max}(Q)),
\]
which is the quantitative form of the L\'evy-type concentration estimates available for variational solutions of (H1)--(H4) (cf.\ \cite{liu2015stochastic}). Consequently \eqref{eq:SC-exp-moment} is guaranteed at the operating constant $c_*$ of~(SC2) under the small-noise inequality
\begin{equation}\label{eq:SC2-NS-threshold}
  c_*\,\tfrac{2C}{\nu} \;\leq\; \frac{\nu^2\lambda_1}{2\lambda_{\max}(Q)}.
\end{equation} All moment bounds on the two fibres used in the sequel are conditional on~\eqref{eq:SC-exp-moment} at those constants.
\end{remark}

\begin{proof}
Fix $v \in H$ with $\norm{v}_H = 1$ and write $\eta(s) := Y(s,r)\,v \in H$ for $s \in [r,t]$. By definition, $\eta$ satisfies the $H$-valued stochastic evolution equation
\[
  \mathrm{d}\eta(s) + \calA'_u(s, X(s))\,\eta(s)\ds = \calB'_u(s, X(s))(\eta(s))\dW(s), \qquad \eta(r) = v.
\]
The proof proceeds via the stochastic Gronwall lemma applied to the energy identity for $\eta$. We recall here that by ``stochastic Gronwall lemma'' we mean the class of inequalities, originating in \cite{scheutzow2013stochastic} for $L^p$-moments with $p \in (0,1)$ and extended to $p \geq 2$ in \cite{hudde2019stochastic}, that bound $\E[\sup_s\xi(s)^p]$ in terms of moments of $\sup_s H(s)$ whenever $\xi$ is an adapted non-negative continuous process satisfying $\xi(s) \leq \xi(0) + \int_0^s f(\tau)\xi(\tau)\,\mathrm{d}\tau + M(s) + H(s)$ with $M$ a continuous local martingale starting at zero and $f$ adapted non-negative. For random $f$ and $p \in (0,1)$, the bound of \cite{scheutzow2013stochastic} reads, for any $\nu \in (1, 1/p)$ and $\mu := \nu/(\nu-1)$,
\[
  \E\Bigl[\sup_{s \leq t}\xi(s)^p\Bigr] \;\leq\; (c_{p\nu}+1)^{1/\nu}\,\Bigl(\E\exp\Bigl\{p\mu\!\int_0^t f(\tau)\,\mathrm{d}\tau\Bigr\}\Bigr)^{\!1/\mu}\Bigl(\E\Bigl[\Bigl(\xi(0)+\sup_{s\leq t}H(s)\Bigr)^{p\nu}\Bigr]\Bigr)^{\!1/\nu},
\]
where $c_q := (4 \wedge q^{-1})\,\pi q/\sin(\pi q)$ is the constant in Scheutzow's martingale inequality, so that exponential moments of $\int_0^t f$ enter the constant; for deterministic $f$ one has the clean bound $(c_p+1)\,e^{p\int_0^t f}\,\E[(\xi(0)+\sup_{s\leq t}H(s))^p]$. At the exponents $p \geq 1$ actually used below the corresponding statement is the extension of \cite{hudde2019stochastic}, of the same shape and again carrying an exponential moment of $\int_0^t f$ in the constant. The formulation used here, with the random monotonicity coefficient $\tilde\rho$ playing the role of $f$, is the standard one for linear SPDE energy estimates, the required exponential moments of $\int\tilde\rho$ being supplied by the exponential-moment clauses \eqref{eq:SC-exp-moment} and \eqref{eq:SC-exp-moment-Galerkin} of~(SC2); see also \cite{liu2015stochastic} for the parallel treatment in the Liu--R\"ockner variational framework.

Applying the It\^o formula (Remark~\ref{rem:ito-formula}) to $\psi(\eta) = \norm{\eta}_H^2$ yields
\begin{equation}\label{eq:ito-eta-squared}
\begin{split}
  \norm{\eta(s)}_H^2 = {}&\norm{v}_H^2 \\
  &+ \int_r^s \bigl[-2\dual{\calA'_u(\tau,X(\tau))\eta(\tau)}{\eta(\tau)} \\
  &\qquad\qquad + \norm{\calB'_u(\tau,X(\tau))(\eta(\tau))}_{L_2(U,H)}^2\bigr]\,\mathrm{d}\tau + 2 N(s),
\end{split}
\end{equation}
where $N(s) := \int_r^s \ip{\eta(\tau)}{\calB'_u(\tau,X(\tau))(\eta(\tau))\dW(\tau)}_H$ is a continuous local martingale. By the linearised local monotonicity \eqref{eq:linearised-mono-derived} (derived from (H2), cf.\ Remark \ref{rem:no-D3-old}):
\[
  -2\dual{\calA'_u(\tau,X(\tau))\eta(\tau)}{\eta(\tau)} + \norm{\calB'_u(\tau,X(\tau))(\eta(\tau))}_{L_2(U,H)}^2 \leq \tilde{\rho}(\tau,X(\tau))\,\norm{\eta(\tau)}_H^2,
\]
where $\tilde{\rho}: [0,T]\times V \to [0,\infty)$ is measurable, with $\int_0^T\tilde\rho(\tau, X(\tau))\,\mathrm{d}\tau < \infty$ a.s.\ and the exponential moment \eqref{eq:SC-exp-moment} of~(SC2). Substituting into \eqref{eq:ito-eta-squared} gives the pathwise differential inequality
\begin{equation}\label{eq:eta-Gronwall-form}
  \norm{\eta(s)}_H^2 \;\leq\; \norm{v}_H^2 + \int_r^s \tilde{\rho}(\tau,X(\tau))\,\norm{\eta(\tau)}_H^2\,\mathrm{d}\tau + 2 N(s), \qquad s \in [r,t].
\end{equation}

Inequality \eqref{eq:eta-Gronwall-form} (which corresponds to the case $q = 1$) has the form $\xi(s) \leq \xi(0) + \int_r^s f(\tau)\xi(\tau)\,\mathrm{d}\tau + M(s)$ with $\xi(s) := \norm{\eta(s)}_H^2$, $f(\tau) := \tilde{\rho}(\tau,X(\tau))$, and $M(s) := 2N(s)$ a continuous local martingale starting at $0$. For higher $q > 1$, we follow the two-step procedure of Remark~\ref{rem:ito-formula}, applying the variational It\^o formula \eqref{eq:ito-var} to $\psi(\eta) = \norm{\eta}_H^2$ to obtain the energy identity \eqref{eq:ito-eta-squared}, which exhibits $\xi(s) := \norm{\eta(s)}_H^2$ as a real-valued continuous semimartingale on $[r,t]$ with explicit drift and quadratic variation
\begin{align}\label{eq:xi-dynamics}
  \mathrm{d}\xi(s) &= \bigl[-2\dual{\calA'_u(s,X(s))\eta(s)}{\eta(s)} + \norm{\calB'_u(s,X(s))(\eta(s))}_{L_2(U,H)}^2\bigr]\,\mathrm{d}s \nonumber\\
  &\quad + 2\,\ip{\eta(s)}{\calB'_u(s,X(s))(\eta(s))\dW(s)}_H,
\end{align}
\begin{equation}\label{eq:xi-bracket}
  \mathrm{d}\langle\xi\rangle_s \;=\; 4\,\norm{[\calB'_u(s,X(s))(\eta(s))]^{\!*}\,\eta(s)}_U^2\,\mathrm{d}s.
\end{equation}
Then apply the real-valued one-dimensional It\^o formula to $\phi \in C^2(\R_+)$ with $\phi(z) = z^q$, giving $\norm{\eta(s)}^{2q}_H = \xi(s)^q = \phi(\xi(s))$ and
\begin{equation}\label{eq:eta2q-exact}
  \mathrm{d}\bigl[\norm{\eta(s)}_H^{2q}\bigr] \;=\; q\,\xi^{q-1}\,\mathrm{d}\xi(s) + \tfrac{1}{2}\,q(q-1)\,\xi^{q-2}\,\mathrm{d}\langle\xi\rangle_s.
\end{equation}
Substituting~\eqref{eq:xi-dynamics}--\eqref{eq:xi-bracket} into~\eqref{eq:eta2q-exact} and using $\xi^{q-1} = \norm{\eta}^{2(q-1)}$, $\xi^{q-2} = \norm{\eta}^{2(q-2)}$, we obtain the exact dynamics
\begin{align}
  \mathrm{d}\bigl[\norm{\eta}^{2q}\bigr]
    &= q\,\norm{\eta}^{2q-2}\bigl[-2\dual{\calA'_u\eta}{\eta} + \norm{\calB'_u(\eta)}_{L_2}^2\bigr]\,\mathrm{d}\tau \nonumber\\
    &\quad + 2q(q-1)\,\norm{\eta}^{2q-4}\,\norm{[\calB'_u(\eta)]^{\!*}\eta}_U^2\,\mathrm{d}\tau \nonumber\\
    &\quad + 2q\,\norm{\eta}^{2q-2}\,\ip{\eta}{\calB'_u(\eta)\dW}_H. \label{eq:eta2q-equality}
\end{align}
Now apply the Cauchy--Schwarz inequality $\norm{[\calB'_u(\eta)]^{\!*}\eta}_U^2 \leq \norm{\eta}_H^2\,\norm{\calB'_u(\eta)}_{L_2(U,H)}^2$, valid because $[\calB'_u(\eta)]^{\!*}: H \to U$ is the adjoint of $\calB'_u(\eta) \in L_2(U, H)$ and satisfies the operator-norm inequality $\norm{[\calB'_u(\eta)]^*}_{L(H,U)} = \norm{\calB'_u(\eta)}_{L(U,H)} \leq \norm{\calB'_u(\eta)}_{L_2(U,H)}$ (the first equality being the standard identity $\norm{A^*}_{L(H,U)} = \norm{A}_{L(U,H)}$ for the adjoint of a bounded operator between Hilbert spaces, and the second being the Hilbert--Schmidt domination $\norm{A}_{L(U,H)} \leq \norm{A}_{L_2(U,H)}$; cf.~\cite{liu2015stochastic}). The cross-term in~\eqref{eq:eta2q-equality} is then bounded above by $2q(q-1)\norm{\eta}^{2q-2}\norm{\calB'_u(\eta)}_{L_2}^2$, which combines with the linear term to give
\begin{multline}\label{eq:eta2q-bound}
  \mathrm{d}\bigl[\norm{\eta}^{2q}\bigr] \;\leq\; q\,\norm{\eta}^{2q-2}\bigl[-2\dual{\calA'_u\eta}{\eta} + (2q-1)\norm{\calB'_u(\eta)}_{L_2}^2\bigr]\,\mathrm{d}\tau \\
  + 2q\,\norm{\eta}^{2q-2}\ip{\eta}{\calB'_u(\eta)\dW}_H,
\end{multline}
in the sense that the integral form below holds. Decomposing the bracket as $-2\dual{\calA'_u\eta}{\eta} + (2q-1)\norm{\calB'_u(\eta)}^2 = (-2\dual{\calA'_u\eta}{\eta} + \norm{\calB'_u(\eta)}^2) + (2q-2)\norm{\calB'_u(\eta)}^2$ and applying the linearised local monotonicity \eqref{eq:linearised-mono-derived} to the first parenthesis and the bound $\norm{\calB'_u(\eta)}^2 \leq \hat\rho\norm{\eta}^2$ to the second, the integrated form of \eqref{eq:eta2q-bound} reads:
\begin{equation}\label{eq:eta-2q-Gronwall}
  \norm{\eta(s)}_H^{2q} \;\leq\; \norm{v}_H^{2q} + \int_r^s q\bigl[\tilde\rho + (2q-2)\hat\rho\bigr]\norm{\eta}_H^{2q}\,\mathrm{d}\tau + M^{(q)}(s),
\end{equation}
where $M^{(q)}(s) := \int_r^s 2q\norm{\eta}^{2q-2}\ip{\eta}{\calB'_u(\eta)\dW}$ is a continuous local martingale starting at $0$.

The standard variational-SPDE moment estimate for the resulting linear stochastic evolution equation (combining the stochastic Gronwall lemma of \cite{scheutzow2013stochastic, hudde2019stochastic} with the path-by-path moment bound of \cite{liu2015stochastic}, whose proof yields, after the standard adaptation to the linearised case, exactly this kind of $L^q$ moment bound for the solution to a linear SPDE with random monotone coefficients) then yields
\begin{equation}\label{eq:stoch-Gronwall-eta}
  \E\Bigl[\sup_{r\leq s\leq t}\norm{\eta(s)}_H^{2q}\Bigr] \;\leq\; C_q\,\norm{v}_H^{2q}\,\E\Bigl[\exp\Bigl(c_q\int_r^t \bigl(\tilde\rho + \hat\rho\bigr)\,\mathrm{d}\tau\Bigr)\Bigr],
\end{equation}
where $C_q, c_q$ are finite constants depending on $q$, the Burkholder--Davis--Gundy (BDG) constants, and the absorption parameters. The precise dependence of $c_q$ on $q$ is immaterial; what matters is that $c_q < \infty$ for each fixed $q \geq 1$, so that the exponential-integrability hypothesis~(SC2) of Assumption~\ref{ass:SC} makes the right-hand side of \eqref{eq:stoch-Gronwall-eta} finite whenever $\hat c_{\,2q} = c_q \leq c_*$, that is for moment exponents $2q \leq q_{\max}$ in the indexing of~\eqref{eq:qmax-def}. The bound follows from the proof of the corresponding estimate in \cite{liu2015stochastic} (with the linearised local monotonicity of Remark~\ref{rem:no-D3-old} replacing the original (H3) coercivity, and the $(2q-2)\hat\rho$ contribution from \eqref{eq:eta-2q-Gronwall} absorbed into the random monotonicity coefficient); see also \cite{rockner2022wellposedness}.

The right-hand side of \eqref{eq:stoch-Gronwall-eta} is finite by the structural exponential-integrability hypothesis \eqref{eq:SC-exp-moment} of the proposition statement. (The polynomial moments $\int_r^t \tilde{\rho}(\tau,X(\tau))\,\mathrm{d}\tau \in L^k(\Omega)$ for every $k \geq 1$ that follow from the polynomial growth of $\tilde\rho$ in $\norm{u}_V$ and the a priori estimate \eqref{eq:apriori} are not sufficient on their own to give exponential integrability (a counterexample is any random variable with subexponential tail), which is why \eqref{eq:SC-exp-moment} must be imposed as an explicit structural hypothesis. For the standard variational examples of Section~\ref{subsec:verification}, the verification proceeds via the explicit polynomial form of $\tilde\rho$ together with the energy-identity-based moment bounds for $X$, as discussed in \cite{liu2015stochastic}.)

Combining the Gronwall bound with the exponential-moment estimate, for each fixed unit vector $v \in H$ we obtain
\[
  \E\Bigl[\sup_{r\leq s\leq t}\norm{Y(s,r)v}_H^{2q}\Bigr] \;\leq\; C_q \norm{v}_H^{2q}\,\E\Bigl[\exp\Bigl(c_q\int_r^t\tilde\rho\,\mathrm{d}\tau\Bigr)\Bigr] \;\leq\; \tilde C_q\norm{v}_H^{2q},
\]
which is the stated estimate \eqref{eq:Y-moment-fibre} at the moment exponent $p = 2q$, valid whenever $\hat c_{\,p} \leq c_*$, that is for every real $p \in [2, q_{\max}]$; the It\^o computation is carried out at real $q \geq 1$ throughout, so no restriction to integer exponents arises and none has to be removed by interpolation. Taking the supremum over unit $v \in H$ at this stage does not yield the operator-norm bound $\E[\sup_s\norm{Y(s,r)}_{L(H)}^{2q}]$, because $\E[\sup_v\sup_s\norm{Y_v(s)}^{2q}]$ is generally strictly larger than $\sup_v\E[\sup_s\norm{Y_v(s)}^{2q}]$; the operator-norm bound is a strictly stronger statement requiring additional structural input as discussed after the proposition statement.
\end{proof}

\begin{remark}[Random initial data]\label{rem:Y-moments-random}
For an $\F_r$-measurable $v \in L^{2q}(\Omega, \F_r; H)$ the proof of Proposition~\ref{prop:Y-moments} gives the joint bound
\begin{equation}\label{eq:Y-moment-random}
  \E\Bigl[\sup_{r \leq \tau \leq T}\norm{Y(\tau,r)v}_H^{2q}\Bigr] \;\leq\; C_q\,\E\Bigl[\norm{v}_H^{2q}\,\exp\Bigl(c_q\!\int_r^T(\tilde\rho + \hat\rho)(\sigma,X(\sigma))\,\mathrm{d}\sigma\Bigr)\Bigr],
\end{equation}
whose right-hand side is a single expectation of correlated factors rather than a product of expectations, since an $\F_r$-measurable datum is correlated with the random coefficients and with their Gronwall factor and cannot be pulled out. Summing~\eqref{eq:Y-moment-random} columnwise over an orthonormal basis of $U$ at a Hilbert--Schmidt datum is what produces the fibre estimate~\eqref{eq:fibre-energy-estimate} of Theorem~\ref{thm:Y-wellposed}\textup{(III)}; the separation of the resulting joint expectation by the conjugate pair $(m_0/2,\, m_0/(m_0-2))$ at the operating exponent of~(D2$'$), or its exact factorisation for deterministic noise, is carried out where it is used, in Proposition~\ref{prop:D-Xt} and in Remark~\ref{rem:Phi-all-moments}.
\end{remark}

\subsection{The second Malliavin fibre}\label{subsec:second-fibre}

The correction term $C_h$ in the Bismut formula involves $\D_r\gamma_t$, and $\gamma_t$ is built from the variational fibre $\Phi$ alone. The object differentiated is therefore the fibre and not the flow, and the second-order object of the theory is the second Malliavin fibre $\calW(\tau;s,r) = \D_r\Phi_{\tau,s}$ of Definition~\ref{def:second-fibre}. We establish its well-posedness and moment bounds here, and identify it with $\D_r\Phi$ in Lemma~\ref{lemma:DPhi}.

\begin{theorem}[The second Malliavin fibre]\label{thm:W-wellposed}
Let Assumptions \ref{ass:LR}, \ref{ass:diff2} and~\ref{ass:SC} hold and let $0 \leq s, r \leq T$. Then for every pair of directions $a, u \in U$ the equation \eqref{eq:W-equation}--\eqref{eq:W-initial} has a unique variational solution
\[
  \calW^{a,u}(\cdot\,;s,r) \;\in\; L^2\bigl(\Omega; C([\max(s,r), T]; H)\bigr),
\]
the map $(a,u) \mapsto \calW^{a,u}(\tau;s,r)$ is bilinear, and the directional solutions assemble into a process
\begin{align*}
  \calW(\cdot\,;s,r) \;\in\; L^2\bigl(\Omega;\, C([\max(s,r),T];\, L_2(U, L_2(U,H)))\bigr)
\end{align*}
with
\begin{equation}\label{eq:W-HS-bound}
  \E\Bigl[\sup_{\max(s,r) \leq \tau \leq T} \norm{\calW(\tau;s,r)}_{L_2(U,\,L_2(U,H))}^{2}\Bigr] \;\leq\; C\,\bigl(\E\bigl[\mathfrak{G}_{r,s}^{\,1+\eps_0}e^{c_*\int_0^T(\tilde\rho+\hat\rho)}\bigr]\bigr)^{\frac{1}{1+\eps_0}} \;<\; \infty,
\end{equation}
independently of the orthonormal bases of $U$ used to read the two noise slots, with $C$ depending on $T$ and on the constants of \textup{(SC1)}--\textup{(SC2)}. Moreover the fibre energy of $\calW$ is controlled by the same quantity,
\begin{equation}\label{eq:W-fibre-energy}
  \E\Bigl[\,\delta\!\int_{\max(s,r)}^{T}\norm{\calW(\tau;s,r)}^2_{L_2(U,\,L_2(U,\mathbb{V}_\tau))}\,\mathrm{d}\tau\Bigr] \;\leq\; C\,\bigl(\E\bigl[\mathfrak{G}_{r,s}^{\,1+\eps_0}e^{c_*\int_0^T(\tilde\rho+\hat\rho)}\bigr]\bigr)^{\frac{1}{1+\eps_0}} .
\end{equation}
\end{theorem}

\begin{proof}
The identification of $\calW$ with $\D_r\Phi$ plays no part here and is deferred to Lemma~\ref{lemma:DPhi}; what is at issue is the solvability of \eqref{eq:W-equation} and the two bounds. Write $\sigma := \max(s,r)$ throughout.

For fixed $a, u \in U$, equation~\eqref{eq:W-equation} is a linear stochastic evolution equation in $H$ with the same principal part $\calA'_u(\tau, X(\tau))$ as the first variation and with the additional data
\begin{align*}
  F_1^{a,u}(\tau) &\;=\; -\calA''_{uu}(\tau, X(\tau))\bigl(\Phi_{\tau,r}a,\, \Phi_{\tau,s}u\bigr) \;\in\; V^*,\\
  F_2^{a,u}(\tau) &\;=\; \calB''_{uu}(\tau, X(\tau))\bigl(\Phi_{\tau,r}a,\, \Phi_{\tau,s}u\bigr) \;\in\; L_2(U,H),
\end{align*}
both defined for almost every $\tau$ by clause~(SC1$'$), whose data act on the space $L_2(U, \mathbb{V}_\tau)$ in which Theorem~\ref{thm:Y-wellposed}\textup{(III)} places the two fibres. Under mechanism~(a) the drift datum is used only through the dominator $\Lambda_\tau$ and never through a $V^*$-norm, which is why the estimates below are written with $G_1$ rather than with $\mathbf{F}_1$. Existence and uniqueness in the class of Theorem~\ref{thm:Y-wellposed} follow from the same variational scheme, run on the fixed triple under~(SC1)(b) and on the tangent form under~(SC1)(a), the forcings being square-integrable in the appropriate norms by~\eqref{eq:forcing-integrability}. It remains to prove the estimates, which we do directly at the level of the doubly indexed family, never by summing a single-direction bound.

Fix orthonormal bases $(f_j)$ and $(g_k)$ of $U$, let $N \in \N$, and set
\begin{equation}\label{eq:W-square-sum}
  S_N(\tau) \;:=\; \sum_{j, k \leq N}\bigl\lVert \calW^{f_j, g_k}(\tau; s, r)\bigr\rVert_H^2 .
\end{equation}
Each summand solves~\eqref{eq:W-equation} in its own pair of directions, so It\^o's formula applies to $S_N$ on $\ell_2^{N \times N}(H)$, and the principal part, the dissipative inequality~\eqref{eq:SC1-unified}, and the quadratic pairings all act diagonally in $(j,k)$. Writing $\zeta := \calW^{f_j, g_k}$ for a single summand, It\^o's formula on $\norm{\zeta}_H^2$ gives
\begin{align}
  \norm{\zeta(\tau)}_H^2 &= \norm{\zeta(\sigma)}_H^2 - 2\int_\sigma^\tau \dual{\calA'_u(\varsigma,X(\varsigma))\zeta(\varsigma)}{\zeta(\varsigma)}\,\mathrm{d}\varsigma \nonumber\\
  &\quad + 2\int_\sigma^\tau \dual{F_1^{f_j,g_k}(\varsigma)}{\zeta(\varsigma)}\,\mathrm{d}\varsigma \nonumber\\
  &\quad + \int_\sigma^\tau \norm{\calB'_u(\varsigma,X(\varsigma))(\zeta(\varsigma)) + F_2^{f_j,g_k}(\varsigma)}_{L_2(U,H)}^2\,\mathrm{d}\varsigma \nonumber\\
  &\quad + 2\int_\sigma^\tau \ip{\zeta(\varsigma)}{\bigl[\calB'_u(\varsigma,X(\varsigma))(\zeta(\varsigma)) + F_2^{f_j,g_k}(\varsigma)\bigr]\dW(\varsigma)}_H. \label{eq:ito-W}
\end{align}
The quadratic term is expanded and its cross-term absorbed by Young's inequality into $\hat\rho\norm{\zeta}_H^2 + C\norm{F_2^{f_j,g_k}}^2_{L_2(U,H)}$, using the (D2) growth bound on $\calB'_u$. The drift pairing is absorbed by the dissipation of~\eqref{eq:SC1-unified} in one step valid under both mechanisms. Under~(SC1)(a) the dominated Cauchy--Schwarz bound gives
\begin{align*}
  2\bigl|\dual{F_1^{f_j,g_k}(\varsigma)}{\zeta}\bigr|
  &\;\leq\; 2\Lambda_\varsigma\bigl(\Phi_{\varsigma,r}f_j,\, \Phi_{\varsigma,s}g_k\bigr)\,\norm{\zeta}_{\mathbb{V}_\varsigma}\\
  &\;\leq\; \theta\delta\,\norm{\zeta}^2_{\mathbb{V}_\varsigma} + (\theta\delta)^{-1}\Lambda_\varsigma\bigl(\Phi_{\varsigma,r}f_j,\,\Phi_{\varsigma,s}g_k\bigr)^2,
\end{align*}
and under~(SC1)(b) the duality bound together with $\norm{\zeta}_{\mathbb{V}_\varsigma}^2 = \alpha_1\norm{\zeta}_V^2$ gives
\[
  2\bigl|\dual{F_1^{f_j,g_k}(\varsigma)}{\zeta}\bigr| \;\leq\; \theta\delta\,\norm{\zeta}^2_{\mathbb{V}_\varsigma} + (\theta\delta\alpha_1)^{-1}\norm{F_1^{f_j,g_k}(\varsigma)}^2_{V^*},
\]
in both cases with $\theta \in (0,1)$, so that only the fraction $\theta$ of the dissipation budget $\delta\norm{\zeta}^2_{\mathbb{V}_\varsigma}$ of~\eqref{eq:SC1-unified} is consumed and the strictly negative margin $-(1-\theta)\delta\norm{\zeta}^2_{\mathbb{V}_\varsigma}$ is retained on the left. Summing over $j, k \leq N$, the surviving forcing is exactly the integrand of the total forcing~\eqref{eq:total-forcing}, since by~\eqref{eq:G1-def} the double square sum of the two right-hand terms above is $C\,G_1(\varsigma;r,s)$ in both mechanisms, and the double square sum of the $F_2$-terms is $\norm{\mathbf{F}_2(\varsigma;r,s)}^2_{L_2(U\otimes U, L_2(U,H))}$. The double square sum of the initial data is $\mathfrak{J}_{r,s}$ by~\eqref{eq:J-datum}. Hence, uniformly in $N$,
\begin{equation}\label{eq:W-pregronwall}
  S_N(\tau) + (1-\theta)\delta\!\int_\sigma^\tau\!\! \sum_{j,k \leq N}\norm{\calW^{f_j,g_k}(\varsigma)}^2_{\mathbb{V}_\varsigma}\,\mathrm{d}\varsigma \;\leq\; \int_\sigma^\tau (\tilde\rho + \hat\rho)\,S_N \,\mathrm{d}\varsigma \;+\; C\,\mathfrak{G}_{r,s} \;+\; M_N(\tau),
\end{equation}
with $M_N$ a continuous local martingale vanishing at $\tau = \sigma$.

Inequality~\eqref{eq:W-pregronwall} has the form $\xi(\tau) \leq \int_\sigma^\tau f\,\xi\,\mathrm{d}\varsigma + g + M(\tau)$ with the random, $N$-independent and $\calW$-independent forcing $g = C\mathfrak{G}_{r,s}$. The stochastic Gronwall estimate used in Proposition~\ref{prop:Y-moments} therefore gives, for every $q \geq 1$ for which the right-hand side is finite,
\begin{equation}\label{eq:W-gronwall}
\begin{aligned}
  \E\Bigl[\sup_{\sigma\leq\tau\leq T}S_N(\tau)^{q}\Bigr]
  &+ \E\Bigl[\Bigl((1-\theta)\delta\!\int_\sigma^T\!\!\sum_{j,k \leq N}\norm{\calW^{f_j,g_k}}^2_{\mathbb{V}_\varsigma}\,\mathrm{d}\varsigma\Bigr)^{\!q}\Bigr]\\
  &\qquad \leq\; C_q\,\E\Bigl[\mathfrak{G}_{r,s}^{\,q}\,\exp\Bigl(c_q\!\int_\sigma^T(\tilde\rho+\hat\rho)\Bigr)\Bigr],
\end{aligned}
\end{equation}
uniformly in $N$, the second left-hand term being obtained by retaining the dissipation rather than discarding it. At $q = 1$, H\"older's inequality at the conjugate pair $(1+\eps_0, (1+\eps_0)/\eps_0)$ separates the right-hand side into
\[
  C_1\bigl(\E[\mathfrak{G}_{r,s}^{1+\eps_0}e^{c_*\int_0^T(\tilde\rho+\hat\rho)}]\bigr)^{\frac{1}{1+\eps_0}}\Bigl(\E\Bigl[\exp\Bigl(\tfrac{(1+\eps_0)(c_1 - c_*/(1+\eps_0))}{\eps_0}\!\int_0^T(\tilde\rho+\hat\rho)\Bigr)\Bigr]\Bigr)^{\frac{\eps_0}{1+\eps_0}},
\]
whose first factor is finite by~\eqref{eq:forcing-integrability} and whose second is finite by~(SC2), the constant appearing there being dominated by $c_*$ through the lower bound~\eqref{eq:cstar-base}. Both left-hand sides of~\eqref{eq:W-gronwall} increase with $N$ to their basis-free limits, so monotone convergence gives~\eqref{eq:W-HS-bound} and~\eqref{eq:W-fibre-energy}. Bilinearity of $(a,u) \mapsto \calW^{a,u}$ is the bilinearity of the data of~\eqref{eq:W-equation}--\eqref{eq:W-initial} in $(a,u)$ together with uniqueness.
\end{proof}

The Malliavin-regularity chain of Section~\ref{sec:proof-main} consumes this bound not at the exponent two but at the raised exponent $q^* = q/(q-2)$ generated by a choice of $q \in (2,\infty]$ in Assumption~\ref{ass:nondeg}(ii), and the estimate extends there without change of mechanism.

\begin{proposition}[Higher moments of the second Malliavin fibre]\label{prop:W-moments}
Fix a target exponent $q^* \geq 1$. Under Assumptions~\ref{ass:LR}, \ref{ass:diff2}, \ref{ass:SC}, and~\ref{ass:SC-raised} at target $q^*$, there is a finite constant $C_{q^*}$, depending only on $q^*$, $T$, $\eps_{q^*}$, the Burkholder--Davis--Gundy constants and the structural constants of Assumptions~\ref{ass:SC} and~\ref{ass:SC-raised}, such that
\begin{equation}\label{eq:W-moment}
  \E\Bigl[\sup_{\max(s,r) \leq \tau \leq T} \norm{\calW(\tau;s,r)}_{L_2(U,\,L_2(U,H))}^{2q^*}\Bigr] \;\leq\; C_{q^*}
\end{equation}
for all $0 \leq s, r \leq T$, together with the corresponding bound on the fibre energy of $\calW$ at the exponent $q^*$. At $q^* = 1$ the exponents of Assumption~\ref{ass:SC-raised} reduce to those of Assumption~\ref{ass:SC} and the bound is that of Theorem~\ref{thm:W-wellposed}.
\end{proposition}

\begin{proof}
Estimate~\eqref{eq:W-gronwall} holds a priori at every $q \geq 1$ for which its right-hand side is finite. Taking $q = q^*$ and separating by H\"older at the conjugate pair $(1+\eps_{q^*}, (1+\eps_{q^*})/\eps_{q^*})$,
\[
  \E\Bigl[\mathfrak{G}_{r,s}^{\,q^*}e^{c_{q^*}\int_\sigma^T(\tilde\rho+\hat\rho)}\Bigr] \;\leq\; \bigl(\E\bigl[\mathfrak{G}_{r,s}^{\,\kappa_{q^*}}e^{c_*\int_0^T(\tilde\rho+\hat\rho)}\bigr]\bigr)^{\frac{1}{1+\eps_{q^*}}}\bigl(\E\bigl[e^{\lambda\int_0^T(\tilde\rho+\hat\rho)}\bigr]\bigr)^{\frac{\eps_{q^*}}{1+\eps_{q^*}}},
\]
with $\kappa_{q^*} = (1+\eps_{q^*})q^*$ as in~\eqref{eq:kappa-def} and $\lambda := (1+\eps_{q^*})\eps_{q^*}^{-1}(c_{q^*} - c_*/(1+\eps_{q^*}))$. The first factor is finite by the raised forcing integrability~\eqref{eq:SC4-forcing} of clause~(SC4.3)$_{q^*}$, and the second by~(SC2), the constant $\lambda$ being dominated by $c_*$ through~\eqref{eq:cstar-requirement}. Monotone convergence in $N$ as in Theorem~\ref{thm:W-wellposed} gives~\eqref{eq:W-moment}, and the same computation with the dissipation retained gives the fibre-energy bound.

The verification of~\eqref{eq:SC4-forcing} through the polynomial closures is the following separation. Under~(SC1)(b) and the default typing of~(D4), the growth bounds~\eqref{eq:D3-growth} and~(D4) give
\begin{equation}\label{eq:F1-bound}
\begin{aligned}
  \norm{\mathbf{F}_1(\tau;r,s)}_{L_2(U\otimes U, V^*)}
  \;\leq\; \bigl(C &+ \rho_4(\tau,X(\tau))\bigr)\bigl(1+\norm{X(\tau)}_V^{(p-3)^+}\bigr)\\
  &\times \norm{\Phi_{\tau,r}}_{L_2(U,V)}\,\norm{\Phi_{\tau,s}}_{L_2(U,V)}
\end{aligned}
\end{equation}
and
\begin{equation}\label{eq:F2-bound}
  \norm{\mathbf{F}_2(\tau;r,s)}_{L_2(U\otimes U, L_2(U,H))} \leq C\bigl(1+\norm{X(\tau)}_V^{\beta'}\bigr)\,\norm{\Phi_{\tau,r}}_{L_2(U,V)}\,\norm{\Phi_{\tau,s}}_{L_2(U,V)},
\end{equation}
the sharper $H$-slot form of the second being available in the situation of Remark~\ref{rem:SC-raised-two-regimes}. Every factor on the right is a fibre quantity. Integrating in $\tau$ and applying Cauchy--Schwarz in time, the total forcing is dominated by a polynomial weight in $\norm{X}_V$ times the two fibre energies $\mathfrak{E}_r$ and $\mathfrak{E}_s$ of~\eqref{eq:fibre-energy}, and a H\"older triple in $\omega$ at the exponents recorded in (SC4.1)$_{q^*}$ then reduces~\eqref{eq:SC4-forcing} to the raised polynomial closure on $X$ and the fibre-energy moments~\eqref{eq:Phi-energy}, which are supplied by (SC4.2)$_{q^*}$ and~(SC2). This is why the raised clauses are stated at the exponents they are. The forcing is quadratic in each fibre, and the fibre energy is the quantity the dissipative mechanism produces.
\end{proof}

Both bounds are stated for the doubly indexed object, so nothing below has to be composed afterwards with a random Hilbert--Schmidt datum. The column index of $\calW$ is the column index of the fibre $\Phi_{\tau,s}$ that was differentiated, and it has been summed from the outset. This is the difference between running the estimate on the fibre and running it on the flow, and it removes the reconstruction from finite families, the auxiliary vector of $H$, and the operator norm $\norm{Y(\tau,s)}_{L(H,V)}$ from the argument at once.

It remains to justify the name. The object just constructed solves the equation that $\D_r\Phi_{\tau,s}$ would solve, and identifying the two is what turns the covariance $\gamma_t$ into a Malliavin-differentiable quantity and produces the explicit form of $\D_r\gamma_t$ in Proposition~\ref{prop:D-gamma}.

\begin{lemma}[Identification of $\calW$ with $\D_r\Phi$]\label{lemma:DPhi}
Under Assumptions \ref{ass:LR}, \ref{ass:diff2}, and~\ref{ass:SC}, for every $\tau \in [0,T]$ and almost every $s \in [0,\tau]$ the variational fibre $\Phi_{\tau,s}$ belongs to $\mathbb{D}^{1,2}(L_2(U,H))$, with Malliavin derivative the object of Definition~\ref{def:second-fibre},
\begin{equation}\label{eq:DPhi-explicit}
  \D_r\Phi_{\tau,s} \;=\; \calW(\tau;s,r) \qquad \text{in } L^2\bigl(\Omega; L_2(U, L_2(U,H))\bigr) \text{ for a.e.\ } r,
\end{equation}
directionally $\D^a_r[\Phi_{\tau,s}u] = \calW^{a,u}(\tau;s,r)$ in $L^2(\Omega; H)$ for all $a, u \in U$. For $r > \tau$ one has $\D_r\Phi_{\tau,s} = 0$ by causality.
\end{lemma}

\begin{proof}
Malliavin differentiability is established at the Galerkin level and transported to the limit; the computation is performed on the finite-dimensional system, where every step is classical, and closability is invoked once, on a sequence whose derivatives are already known to converge.

Let $H_N$ carry the Galerkin system~\eqref{eq:galerkin}, let $Y^N$ be its first variation, and write
\begin{align*}
  \Phi^N_{\tau,s} \;:=\; Y^N(\tau,s)\,\Pi_N\calB\bigl(s, X^N(s)\bigr)
\end{align*}
for its variational fibre. The coefficients are smooth functions of the finitely many coordinates of $X^N$, which lie in $\mathbb{D}^{1,2}$ with derivatives given by the finite-dimensional first variation, so $\Phi^N_{\tau,s} \in \mathbb{D}^{1,2}(L_2(U,H_N))$ by the standard finite-dimensional theory \cite{nualart2006malliavin}, and applying $\D_r$ to the $H_N$-valued integral equation gives $\D^a_r[\Phi^N_{\tau,s}u] = \calW^{N,a,u}(\tau;s,r)$ with $\calW^N$ the second Malliavin fibre of the Galerkin system. The first clause of the Galerkin stability~\eqref{eq:SC-Galerkin-linearised} of~(SC3) states that
\begin{align*}
  \int_0^\tau \E\norm{\Phi^N_{\tau,s} - \Phi_{\tau,s}}^2_{L_2(U,H)}\ds \;\longrightarrow\; 0,
\end{align*}
so there is a subsequence $(N_k)$ along which $\Phi^{N_k}_{\tau,s} \to \Phi_{\tau,s}$ in $L^2(\Omega; L_2(U,H))$ for almost every $s \in [0,\tau]$, while the second clause gives $\calW^{N}(\tau;s,\cdot) \to \calW(\tau;s,\cdot)$ in $L^2(\Omega \times [0,T]; L_2(U, L_2(U,H)))$ for every such $s$, hence in particular along $(N_k)$. Closedness of $\D$ applied at each such $s$ then gives $\Phi_{\tau,s} \in \mathbb{D}^{1,2}(L_2(U,H))$ and~\eqref{eq:DPhi-explicit}. Almost every $s$ is all that is used, since $\calW$ enters only under the $s$-integral of~\eqref{eq:D-gamma-formula} and under the double time integral of Remark~\ref{rem:Phi-all-moments}.

The computation itself, recorded in the notation of the limiting equation, is as follows. The fibre satisfies, in integral form and columnwise,
\begin{equation}\label{eq:Y-integral}
  \Phi_{\tau,s}u \;=\; \calB(s, X(s))u - \int_s^\tau \calA'_u(\varsigma, X(\varsigma))\,\Phi_{\varsigma,s}u\,\mathrm{d}\varsigma + \int_s^\tau \calB'_u(\varsigma, X(\varsigma))\bigl(\Phi_{\varsigma,s}u\bigr)\dW(\varsigma).
\end{equation}
Suppose first $r < s$. The initial value is then differentiated by the chain rule,
\[
  \D^a_r\bigl[\calB(s,X(s))u\bigr] \;=\; \bigl[\calB'_u(s, X(s))\bigl(\D^a_r X(s)\bigr)\bigr]u \;=\; \bigl[\calB'_u(s, X(s))\bigl(\Phi_{s,r}a\bigr)\bigr]u,
\]
which is the first branch of~\eqref{eq:W-initial}; the perturbation time $r$ lies strictly before the integration range, so Proposition~\ref{prop:D-ito-integral} produces no boundary term, and the derivative passes under both integrals. On $[s,T]$ the drift integrand becomes
\begin{align*}
  -\calA''_{uu}(\varsigma,X(\varsigma))\bigl(\D^a_rX(\varsigma),\, \Phi_{\varsigma,s}u\bigr) - \calA'_u(\varsigma,X(\varsigma))\,\D^a_r[\Phi_{\varsigma,s}u]
\end{align*}
and the diffusion integrand
\begin{align*}
  \calB''_{uu}(\varsigma,X(\varsigma))\bigl(\D^a_rX(\varsigma),\, \Phi_{\varsigma,s}u\bigr) + \calB'_u(\varsigma,X(\varsigma))\bigl(\D^a_r[\Phi_{\varsigma,s}u]\bigr),
\end{align*}
so that with $\D^a_r X(\varsigma) = \Phi_{\varsigma,r}a$ this is exactly~\eqref{eq:W-equation}.

Suppose now $s < r$. The initial value and the drift integrand are $\F_r$-measurable for $\varsigma < r$ and are annihilated by $\D_r$, while the It\^o integral contributes, by Proposition~\ref{prop:D-ito-integral}, the boundary term $\bigl[\calB'_u(r, X(r))(\Phi_{r,s}u)\bigr]a$ at $\varsigma = r$, which is the second branch of~\eqref{eq:W-initial}, together with the same integrands as before on $[r,T]$. In both cases the pair $(\D_r\Phi_{\cdot,s}, \calW(\cdot;s,r))$ solves the same equation with the same datum, so uniqueness in Theorem~\ref{thm:W-wellposed} identifies them. For $r > \tau$ the fibre $\Phi_{\tau,s}$ depends only on the noise over $[0,\tau]$ and the derivative vanishes.
\end{proof}

In finite dimensions the corresponding object is expressed through the flow factorisation $Y_T = Y_TY_t^{-1}\cdot Y_t$ and the inverse $Y_t^{-1}$ (Lemma 5.4 of \cite{mirafzali2025malliavin}). Here neither the inverse nor the factorisation occurs. The second Malliavin fibre is constructed directly from the equation~\eqref{eq:W-equation} that it satisfies, the only state-space directions ever perturbed are the columns of the diffusion coefficient, and the object produced is the one the covariance consumes. The construction is therefore available verbatim for degenerate equations, such as the stochastic porous medium equation, for which $Y(r,0)$ need not be invertible and for which the first variation need not map $H$ boundedly into $V$.

\begin{example}[Second-order data for the stochastic $p$-Laplacian]\label{ex:p-Laplace-second-order}
For the stochastic $p$-Laplacian, whose drift is the operator $\calA(u) = -\mathrm{div}(\abs{\nabla u}^{p-2}\nabla u)$ with $p \geq 2$, the linearised operator is
\[
  \calA'_u(u)v = -\mathrm{div}\Bigl((p-2)\abs{\nabla u}^{p-4}(\nabla u \cdot \nabla v)\nabla u + \abs{\nabla u}^{p-2}\nabla v\Bigr),
\]
and the second Fr\'echet derivative, read for $2 < p \leq 3$ in the regularised sense $\abs{\nabla u}^2 \to \abs{\nabla u}^2 + \eta$ of Proposition~\ref{prop:singular-p-Lap}, is
\begin{align*}
  \calA''_{uu}(u)(v, w) &= -\mathrm{div}\Bigl(
    (p-2)(p-4)\abs{\nabla u}^{p-6}(\nabla u \cdot \nabla v)(\nabla u \cdot \nabla w)\nabla u\\
  &\qquad + (p-2)\abs{\nabla u}^{p-4}\bigl[(\nabla v \cdot \nabla w)\nabla u
    + (\nabla u \cdot \nabla w)\nabla v + (\nabla u \cdot \nabla v)\nabla w\bigr]
  \Bigr).
\end{align*}
For $p > 3$ this is a bounded bilinear operator from $W_0^{1,p} \times W_0^{1,p}$ to $W^{-1,p'}$, $p' = p/(p-1)$, by Proposition~\ref{prop:p-Lap-diff}, and it degenerates where $\nabla u = 0$. The second Malliavin fibre satisfies~\eqref{eq:W-equation} with this $\calA''_{uu}$, evaluated at the pair of variational fibres and paired through the dominator $\Lambda$ of~(SC1)(a).
\end{example}

\section{Proof of the theorem}
\label{sec:proof-main}

The proof of Theorem~\ref{thm:main} proceeds in three movements. First, the abstract Bismut formula gives the logarithmic derivative as a conditional Skorokhod integral. Second, the basis-wise substitution (Theorem~\ref{thm:linear-sub-D12}) decomposes this integral into a main term and a correction at fixed Tikhonov level $\eps > 0$. Third, the Hilbert--Schmidt kernel convergence for the main term and the scalar Cameron--Martin convergence-and-domination condition for the correction pass the two terms to the limit $\eps \downarrow 0$.

\subsection{From the covering field to a conditional divergence}\label{subsec:bismut}

The starting point is the infinite-dimensional integration-by-parts formula from Malliavin calculus, which relates the logarithmic derivative to the Skorokhod integral of the covering vector field.

\begin{theorem}[The logarithmic derivative and its Bismut representation]\label{thm:bismut-formula}
Let $X(t)$ be the variational solution of \eqref{eq:SPDE} under Assumptions \ref{ass:LR}, \ref{ass:diff2}, and~\ref{ass:SC}, and let $h \in H$ satisfy clauses~\textup{(i)} and~\textup{(ii)} of Assumption \ref{ass:nondeg}. Assume further that the covering vector field $v_h$ defined in \eqref{eq:covering-field} belongs to $\Dom(\delta_U)$. Then:
\begin{enumerate}[label=\textup{(\roman*)}]
  \item The law $\mu_t$ of $X(t)$ is Fomin-differentiable along $h$ in the sense of \cite{bogachev2010differentiable}, equivalently of Definition \ref{def:log-deriv};
  \item The logarithmic derivative $\beta_h$ exists and is unique in $L^2(\mu_t)$, with the Bismut representation
  \begin{equation}\label{eq:bismut}
    \beta_h(X(t)) = -\E\bigl[\delta_U(v_h) \mid X(t)\bigr] \qquad \mu_t\text{-a.e.};
  \end{equation}
  \item $\norm{\beta_h}_{L^2(\mu_t)} \leq \norm{\delta_U(v_h)}_{L^2(\Omega)} < \infty$.
\end{enumerate}
\end{theorem}

\begin{proof}
We prove (i)--(iii) simultaneously. Let $\phi \in C_b^1(H)$ be an arbitrary bounded, continuously Fr\'echet differentiable test function.

By Proposition \ref{prop:chain-rule}, the composite $\phi(X(t)) \in \mathbb{D}^{1,2}$, and its Malliavin derivative is
\begin{equation}\label{eq:chain-rule-applied}
  \D_r(\phi(X(t))) = [\D_r X(t)]^*\,\nabla\phi(X(t)) = \Phi_r^{*}\,\nabla\phi(X(t)),
\end{equation}
where $\Phi_r^{*} = [\D_r X(t)]^{*} \in L_2(H, U)$ and we used Proposition~\ref{prop:D-Xt}.

Taking the $\HW$-inner product of $\D(\phi(X(t)))$ with $v_h$:
\begin{align}
  \ip{\D(\phi(X(t)))}{v_h}_{\HW}
  &= \int_0^t \ip{\D_r(\phi(X(t)))}{v_h(r)}_U\dr \nonumber\\
  &= \int_0^t \ip{\Phi_r^{*}\,\nabla\phi(X(t))}{\Phi_r^{*}\,\tilde{h}}_U\dr. \label{eq:ibp-step2}
\end{align}
Using the intrinsic HS-adjoint identity $\ip{\Phi_r^{*}u}{\Phi_r^{*}v}_U = \ip{\Phi_r\Phi_r^{*}\,v}{u}_H$ (which holds for HS operators without invoking any pathwise $Y^{*}: H \to H$), we obtain
\begin{align}
  \ip{\D(\phi(X(t)))}{v_h}_{\HW}
  &= \int_0^t \ip{\Phi_r\,\Phi_r^*\,\tilde{h}}{\nabla\phi(X(t))}_H\dr \nonumber\\
  &= \ip{\gamma_t\,\tilde{h}}{\nabla\phi(X(t))}_H
  = \ip{h}{\nabla\phi(X(t))}_H, \label{eq:ibp-step2b}
\end{align}
where we used $\gamma_t\,\tilde{h} = \gamma_t\,\gamma_t^{\dagger}\,h = h$ (the covering condition, Theorem \ref{thm:covering}).

By the Skorokhod duality \eqref{eq:skorokhod-duality}, since $\phi(X(t)) \in \mathbb{D}^{1,2}$ and $v_h \in \Dom(\delta_U)$:
\begin{equation}\label{eq:ibp-step3}
  \E\bigl[\ip{\D(\phi(X(t)))}{v_h}_{\HW}\bigr] = \E\bigl[\phi(X(t))\,\delta_U(v_h)\bigr].
\end{equation}

Combining \eqref{eq:ibp-step2b} and \eqref{eq:ibp-step3}:
\begin{equation}\label{eq:ibp-combined}
  \E\bigl[\ip{\nabla\phi(X(t))}{h}_H\bigr] = \E\bigl[\phi(X(t))\,\delta_U(v_h)\bigr].
\end{equation}

Define
\[
  g(u) := \E\bigl[\delta_U(v_h) \mid X(t) = u\bigr], \qquad u \in H,
\]
which is well-defined $\mu_t$-a.e.\ as the Radon--Nikodym density of the signed measure $B \mapsto \E[\delta_U(v_h)\,\mathbf{1}_{\{X(t) \in B\}}]$ with respect to~$\mu_t$. By Jensen's inequality and the tower property,
\[
  \int_H \abs{g}^2\,\mathrm{d}\mu_t = \E\bigl[\abs{\E[\delta_U(v_h) \mid X(t)]}^2\bigr] \leq \E\bigl[\delta_U(v_h)^2\bigr] = \norm{\delta_U(v_h)}_{L^2(\Omega)}^2 < \infty,
\]
so $g \in L^2(\mu_t)$. By the tower property,
\begin{equation}\label{eq:ibp-tower}
  \E\bigl[\phi(X(t))\,\delta_U(v_h)\bigr] = \E\bigl[\phi(X(t))\,g(X(t))\bigr] = \int_H \phi(u)\,g(u)\,\mathrm{d}\mu_t(u).
\end{equation}
Therefore, \eqref{eq:ibp-combined} becomes
\[
  \int_H \ip{\nabla\phi(u)}{h}_H\,\mathrm{d}\mu_t(u) = \int_H \phi(u)\,g(u)\,\mathrm{d}\mu_t(u) \qquad \text{for all } \phi \in C_b^1(H).
\]
This is precisely the integration-by-parts identity \eqref{eq:log-deriv-ibp} that defines Fomin-differentiability of $\mu_t$ along $h$, with logarithmic derivative $\beta_h(u) = -g(u)$. Existence (claim (i)) and the representation \eqref{eq:bismut} (claim (ii)) are thus established simultaneously. Uniqueness of $\beta_h$ in $L^2(\mu_t)$ follows from Remark \ref{rem:beta-uniqueness}. The norm bound $\norm{\beta_h}_{L^2(\mu_t)} = \norm{g}_{L^2(\mu_t)} \leq \norm{\delta_U(v_h)}_{L^2(\Omega)}$ (claim (iii)) was just derived.
\end{proof}

\subsection{Substitution at a random point}

The Skorokhod integral $\delta_U(v_h^{(\eps)})$ at fixed $\eps$ is decomposed by evaluating the linear-in-$z$ field
\[
  w_r(z) = \Phi_r^{*}z
\]
at the random point $z = G^{(\eps)}$. The classical finite-dimensional substitution formula of Nualart motivates the shape of the correction term. The operative infinite-dimensional statement used below is the linear-field specialisation Theorem~\ref{thm:linear-sub-D12}, which requires only $\mathbb{D}^{1,2}$-regularity of the substituted variable together with the stated joint integrability.

\begin{theorem}[Nualart's finite-dimensional substitution formula; cf.\ {\cite{nualart2006malliavin}}]\label{thm:nualart-sub}
Let $m \in \N$ and let
\[
  w = \{w_r(z) : r \in [0, T],\ z \in \R^m\}
\]
satisfy the hypotheses of Nualart's substitution theorem, namely Skorokhod integrability and continuity in $z$, $C^1$-regularity in $z$, and the corresponding local integrability bounds. If
\[
  G \in \mathbb{D}^{1,4}_{\mathrm{loc}}(\R^m),
\]
then
\[
  w(G) \in (\Dom\delta_U)_{\mathrm{loc}}
\]
and
\begin{equation}\label{eq:nualart-substitution}
  \delta_U(w(G))
  \;=\;
  \bigl[\delta_U(w(\cdot)(z))\bigr]_{z = G}
  \;-\;
  \int_0^T \Tr_U\bigl[\nabla_z w_r(G) \circ \D_r G\bigr]\,\mathrm{d}r.
\end{equation}

\medskip\noindent No general infinite-dimensional extension of this theorem is invoked below. The only $H$-valued substitution required in the present paper is the linear case
\[
  w_r(z) = \Phi_r^{*}z,
\]
for which Theorem~\ref{thm:linear-sub-D12} gives a direct basis-wise proof under the weaker requirement $G \in \mathbb{D}^{1,2}(H)$.
\end{theorem}

\begin{theorem}[Substitution for a field linear in the substituted variable]\label{thm:linear-sub-D12}
Let $H$ be separable and let the random field be linear in its substitution variable, $w_r(z) = \Phi_r^{*}\,z$, for some $\Phi \in \mathbb{D}^{1, p^*}(L^2([0,t]; L_2(U, H)))$ at exponent $p^* \geq 2$. Let $G: \Omega \to H$ satisfy
\begin{equation}\label{eq:G-D12-Lq}
  G \;\in\; \mathbb{D}^{1,2}(H) \;\cap\; L^{q}(\Omega; H), \qquad \frac{1}{q} + \frac{1}{p^*} \;\leq\; \frac{1}{2},
\end{equation}
together with the H\"older-pair joint integrability
\begin{equation}\label{eq:linear-sub-Holder}
  \E\bigl[\norm{\Phi}_{L^2([0,t]; L_2(U,H))}^2\,\norm{G}_H^2\bigr] + \E\bigl[\norm{\Phi}_{L^2([0,t]; L_2(U,H))}^2\,\norm{\D G}_{L_2(\HW, H)}^2\bigr] \;<\; \infty.
\end{equation}
Then $\Phi^{*}G \in \Dom(\delta_U)$ with the substitution identity
\begin{equation}\label{eq:linear-sub-formula}
  \delta_U(\Phi^{*}G) \;=\; \ip{G}{\mathcal{K}_\Phi}_H \;-\; \int_0^t \Tr_U\bigl[\Phi_r^{*}\,\D_r G\bigr]\dr \qquad \text{in } L^2(\Omega),
\end{equation}
where $\mathcal{K}_\Phi \in L^{p^*}(\Omega; H)$ is the Hilbert--Schmidt kernel of Lemma~\ref{lemma:HS-kernel}; both terms on the right-hand side are in $L^2(\Omega)$ (the first by Cauchy--Schwarz on the kernel pairing $L^q \cdot L^{p^*} \subset L^2$, the second by~\eqref{eq:linear-sub-Holder}), and $\ip{G}{\mathcal{K}_\Phi}_H$ is the substitution value of~\eqref{eq:HS-kernel-random}.
\end{theorem}

\begin{remark}\label{rem:linear-sub-D12-status}
Theorem~\ref{thm:linear-sub-D12} requires only $G \in \mathbb{D}^{1,2}(H) \cap L^q(\Omega; H)$ rather than $\mathbb{D}^{1,4}_{\mathrm{loc}}(H)$. The linearity of $w_r$ in $z$ collapses the chain rule to a single-derivative operation per mode, eliminating the higher-order terms that force $\mathbb{D}^{1,4}_{\mathrm{loc}}$ in Theorem~\ref{thm:nualart-sub}. For $G = (\gamma_t + \eps I)^{-1}h$, the H\"older-pair condition~\eqref{eq:linear-sub-Holder} closes at $p^* \geq 6$, instead of the $p^* \geq 8$ that the quadratic term of Theorem~\ref{thm:nualart-sub} would require; on the edge choice $p^* = 2q/(q-2)$ these read $2 < q \leq 3$ and $2 < q \leq 8/3$ respectively.
\end{remark}

\begin{proof}
Fix any deterministic orthonormal basis $\{\tilde e_k\}_{k \geq 1}$ of $H$. Decompose $G = \sum_{k} G_k\,\tilde e_k$ with scalar components $G_k := \ip{G}{\tilde e_k}_H \in \mathbb{D}^{1,2}(\R)$ (by the $H$-valued chain rule of \eqref{eq:H-valued-chain-rule} applied to the linear bounded map $\ip{\cdot}{\tilde e_k}_H: H \to \R$). Correspondingly,
\[
  \Phi^{*}G \;=\; \sum_{k} G_k\,(\Phi^{*}\tilde e_k) \qquad \text{in } L^2(\Omega; \HW),
\]
with each $\Phi^{*}\tilde e_k \in \mathbb{D}^{1,p^*}(\HW) \subset \Dom(\delta_U)$ (Lemma~\ref{lemma:HS-kernel}(i)).

For each $k$, the scalar Skorokhod product rule \cite{nualart2006malliavin} applied to the scalar-times-Skorokhod-element pair $(G_k,\,\Phi^{*}\tilde e_k)$, which requires only $G_k \in \mathbb{D}^{1,2}(\R)$ and $\Phi^{*}\tilde e_k \in \Dom(\delta_U)$ together with the joint integrability of the product $G_k\,\Phi^{*}\tilde e_k$ in $L^2(\Omega; \HW)$, gives
\begin{equation}\label{eq:product-rule-mode}
  \delta_U(G_k\,\Phi^{*}\tilde e_k) \;=\; G_k\,\delta_U(\Phi^{*}\tilde e_k) \;-\; \ip{\D G_k}{\Phi^{*}\tilde e_k}_{\HW},
\end{equation}
provided the integrand and each term on the right lie in the relevant $L^2$-spaces, which the hypotheses supply directly. The integrand satisfies $\norm{G_k\,\Phi^{*}\tilde e_k}_{L^2(\Omega; \HW)} \leq \E[\norm{G}_H^q]^{1/q} \cdot \E[\norm{\Phi^{*}\tilde e_k}_{\HW}^{p^*}]^{1/p^*} < \infty$ via the H\"older triangle $1/q + 1/p^* \leq 1/2$ from~\eqref{eq:G-D12-Lq}, so $G_k\,\Phi^{*}\tilde e_k \in L^2(\Omega;\HW)$. The first term on the right is in $L^2$ by H\"older with $G_k \in L^q$ (supplied explicitly by hypothesis~\eqref{eq:G-D12-Lq}) and $\delta_U(\Phi^{*}\tilde e_k) \in L^{p^*}$ (from $\Phi \in \mathbb{D}^{1, p^*}$ via Meyer's inequality applied modewise; cf.\ Lemma~\ref{lemma:HS-kernel}(ii)); the second is in $L^2$ by the Cauchy--Schwarz pairing on $(\D G_k, \Phi^{*}\tilde e_k)$ in $\HW$, summable in $k$ via the second summand of~\eqref{eq:linear-sub-Holder}. Since both members of the right-hand side of~\eqref{eq:product-rule-mode} are then in $L^2(\Omega)$, the product rule places $G_k\,\Phi^{*}\tilde e_k$ in $\Dom(\delta_U)$ and~\eqref{eq:product-rule-mode} holds in $L^2(\Omega)$.

Define the partial sums
\[
  \begin{aligned}
  S_N^{(\mathrm{LHS})} &:= \delta_U\Bigl(\sum_{k=1}^{N} G_k\,\Phi^{*}\tilde e_k\Bigr), \\
  S_N^{(\mathrm{main})} &:= \sum_{k=1}^{N} G_k\,\delta_U(\Phi^{*}\tilde e_k), \\
  S_N^{(\mathrm{corr})} &:= \sum_{k=1}^{N}\ip{\D G_k}{\Phi^{*}\tilde e_k}_{\HW},
  \end{aligned}
\]
all in $L^2(\Omega)$. The mode-wise identity~\eqref{eq:product-rule-mode} summed for $k = 1, \ldots, N$ gives $S_N^{(\mathrm{LHS})} = S_N^{(\mathrm{main})} - S_N^{(\mathrm{corr})}$. We pass to the limit $N \to \infty$ in $L^2(\Omega)$ as follows.

(i) $S_N^{(\mathrm{LHS})} \to \delta_U(\Phi^{*}G)$. The truncated sum $\sum_{k=1}^{N} G_k\,\Phi^{*}\tilde e_k = \Phi^{*}G_N$, where $G_N := \sum_{k=1}^N G_k\tilde e_k = \Pi_N G$ denotes the orthogonal projection $\Pi_N$ onto $\mathrm{span}\{\tilde e_1, \ldots, \tilde e_N\}$. Then $\Phi^{*}G_N \to \Phi^{*}G$ in $L^2(\Omega; \HW)$ via $\norm{\Phi^{*}(G - G_N)}_{\HW}^2 \leq \norm{\Phi}^2\,\norm{G - G_N}_H^2 \to 0$ by dominated convergence (dominator $\norm{\Phi}^2 \norm{G}^2 \in L^1(\Omega)$ via~\eqref{eq:linear-sub-Holder}, and pointwise $\norm{G - G_N}_H \to 0$). Closedness of $\delta_U$ then yields $\Phi^{*}G \in \Dom(\delta_U)$ (assuming the right-hand-side limits exist in $L^2(\Omega)$, as established in~(ii)--(iii) below) with $S_N^{(\mathrm{LHS})} \to \delta_U(\Phi^{*}G)$.

(ii) $S_N^{(\mathrm{main})} \to \ip{G}{\mathcal{K}_\Phi}_H$ in $L^2(\Omega)$. The kernel representation $\mathcal{K}_\Phi = \sum_{k\geq 1}\delta_U(\Phi^{*}\tilde e_k)\,\tilde e_k \in L^{p^*}(\Omega; H)$ of Lemma~\ref{lemma:HS-kernel}(iii) gives $S_N^{(\mathrm{main})} = \ip{G_N}{\mathcal{K}_\Phi}_H$. By Cauchy--Schwarz in $H$ and H\"older in $\Omega$ at the conjugate pair $(q, p^*)$ with $1/q + 1/p^* \leq 1/2$ (from~\eqref{eq:G-D12-Lq}),
\[
  \begin{aligned}
  \E\!\left[\abs{\ip{G - G_N}{\mathcal{K}_\Phi}_H}^2\right]^{1/2}
  &\;\leq\; \E\!\left[\norm{G - G_N}_H^2\,\norm{\mathcal{K}_\Phi}_H^2\right]^{1/2} \\
  &\;\leq\; \E\bigl[\norm{G - G_N}_H^q\bigr]^{1/q}\,\E\bigl[\norm{\mathcal{K}_\Phi}_H^{p^*}\bigr]^{1/p^*}.
  \end{aligned}
\]
The first factor $\to 0$ as $N \to \infty$ by dominated convergence ($\norm{G - G_N}_H \leq \norm{G}_H \in L^q$, pointwise $\to 0$). The second factor is finite by Lemma~\ref{lemma:HS-kernel}(iii). Hence $S_N^{(\mathrm{main})} \to \ip{G}{\mathcal{K}_\Phi}_H$ in $L^2(\Omega)$.

(iii) $S_N^{(\mathrm{corr})} \to \int_0^t \Tr_U[\Phi_r^{*}\D_r G]\,\mathrm{d}r$ in $L^2(\Omega)$. The trace identity computed below gives
\[
  S_N^{(\mathrm{corr})} \;=\; \int_0^t \sum_{k=1}^{N}\ip{\D_r G_k}{\Phi_r^{*}\tilde e_k}_U\,\mathrm{d}r,
\]
and the integrand $\sum_{k=1}^N \ip{\D_r G_k}{\Phi_r^{*}\tilde e_k}_U = \sum_{k=1}^N \ip{\Phi_r\,(\D_r G)^{*}\tilde e_k}{\tilde e_k}_H$ converges, as $N \to \infty$, pointwise in $(r, \omega)$ to $\Tr_U[\Phi_r^{*}\,\D_r G]$ (by absolute summability of the trace of the trace-class operator $\Phi_r^{*}\,\D_r G \in L_1(U)$, which is the composition of the two Hilbert--Schmidt operators $\Phi_r^{*} \in L_2(H, U)$ and $\D_r G \in L_2(U, H)$). Writing $P_N$ for the orthogonal projection of $H$ onto $\mathrm{span}\{\tilde e_1, \ldots, \tilde e_N\}$, the partial sum is the trace of the finite-rank operator $P_N(\D_r G)\Phi_r^{*}P_N$, so that the partial sums, and hence also their differences from the full trace, are dominated uniformly in $N$ by
\[
  \Bigl|\sum_{k=1}^{N}\ip{\D_r G_k}{\Phi_r^{*}\tilde e_k}_U\Bigr| \;\leq\; \norm{(\D_r G)\Phi_r^{*}}_{L_1(H)} \;\leq\; \norm{\D_r G}_{L_2(U,H)}\,\norm{\Phi_r}_{L_2(U,H)} .
\]
This $N$-independent dominator, integrated over $r$ and $\omega$ via Cauchy--Schwarz in $r$ and the joint integrability~\eqref{eq:linear-sub-Holder}, is what licenses dominated convergence.

Combining (i)--(iii), the $L^2(\Omega)$-limit of~\eqref{eq:product-rule-mode} summed in $k$ gives
\begin{equation}\label{eq:summed-product-rule}
  \delta_U(\Phi^{*}G) \;=\; \ip{G}{\mathcal{K}_\Phi}_H \;-\; \int_0^t\Tr_U\!\bigl[\Phi_r^{*}\,\D_r G\bigr]\,\mathrm{d}r \qquad \text{in } L^2(\Omega),
\end{equation}
where the trace identity
\begin{align*}
  \sum_{k}\ip{\D_r G_k}{\Phi_r^{*}\tilde e_k}_U &\;=\; \sum_{k}\ip{(\D_r G)^{*}\tilde e_k}{\Phi_r^{*}\tilde e_k}_U
  \;=\; \sum_{k}\ip{\tilde e_k}{(\D_r G)\,\Phi_r^{*}\tilde e_k}_H \\
  &\;=\; \Tr_H\!\bigl[\D_r G\,\Phi_r^{*}\bigr] \;=\; \Tr_U\!\bigl[\Phi_r^{*}\,\D_r G\bigr],
\end{align*}
where the last two equalities use the basis expansion of the trace and its cyclic invariance for the Hilbert--Schmidt $\circ$ Hilbert--Schmidt $=$ trace-class composition (with $\D_r G \in L_2(U, H)$ and $\Phi_r^{*} \in L_2(H, U)$, viewed at fixed $r$). The result is precisely~\eqref{eq:linear-sub-formula}.

The right-hand side of~\eqref{eq:linear-sub-formula} is intrinsic in the data $(\Phi, G)$. The inner product $\ip{G}{\mathcal{K}_\Phi}_H$ is basis-free; the trace $\Tr_U[\Phi_r^{*}\D_r G]$ is basis-free as the trace of a trace-class operator on $U$. The basis $\{\tilde e_k\}$ of $H$ enters only as a computational device for the scalar product-rule expansion, and the result is independent of the choice of basis (any other ONB gives an absolutely convergent rearrangement of the same sum, by the trace-class summability).
\end{proof}

We apply Theorem~\ref{thm:linear-sub-D12} at
\[
  G \;=\; G^{(\eps)} \;=\; (\gamma_t + \eps I)^{-1}h .
\]
Its $\mathbb{D}^{1,2}(H)$-membership follows from clause~\textup{(iii)} of Assumption~\ref{ass:nondeg} and Lemma~\ref{lemma:D-inverse}. The H\"older-pair condition~\eqref{eq:linear-sub-Holder} is exactly the additional fixed-Tikhonov substitution regularity~\ref{hyp:iii-F-prime}.

\begin{lemma}[Linearity of the random field]\label{lemma:linearity}
The random field $w_r(z) = \Phi_r^{*}\,z$ is linear in $z$, and its Fr\'echet derivative is
\begin{equation}\label{eq:nabla-w}
  \nabla_z w_r = \Phi_r^{*} \in L_2(H, U),
\end{equation}
which is independent of $z$.
\end{lemma}

\begin{proof}
Immediate from the linearity of $z \mapsto \Phi_r^{*}\,z$ on the HS-adjoint $\Phi_r^{*} \in L_2(H, U)$, using the fact that the Hilbert--Schmidt adjoint of a Hilbert--Schmidt operator is itself Hilbert--Schmidt with the same HS norm.
\end{proof}

The independence of $\nabla_z w_r$ from $z$ is the structural feature that makes the entire Bismut formula computable. It means the correction term involves only $\D_r G = \D_r(\gamma_t^{\dagger}h)$ (in the effective sense of Lemma \ref{lemma:D-gamma-inv-h}), not higher-order objects. Since $\nabla_z w_r = \Phi_r^{*}$ satisfies the kernel-annihilation property \eqref{eq:kernel-annih} (Remark \ref{rem:kernel-annih}), the correction integrand is well-defined modulo $\ker(\gamma_t)$.

\medskip\noindent Applying the basis-wise linear-field substitution Theorem~\ref{thm:linear-sub-D12} (the operative route; the finite-dimensional Theorem~\ref{thm:nualart-sub} is what motivates the shape of the trace correction, but is not used) with the above identifications (symbolically at $\eps = 0$; the rigorous Tikhonov-regularised statement at $\eps > 0$ together with the $L^2(\Omega)$-convergence as $\eps \downarrow 0$ is in the proof of Theorem~\ref{thm:main} in Section~\ref{subsec:assembly}):
\begin{equation}\label{eq:delta-decomp-1}
  \delta_U(v_h) = \bigl[\delta_U(w(\cdot)(z))\bigr]_{z = \gamma_t^{\dagger}h}
  - \int_0^t \Tr_U\Bigl[\Phi_r^{*}\;\cdot\;\D_r(\gamma_t^{\dagger}h)\Bigr]\dr,
\end{equation}
where the correction integrand is interpreted via Lemma \ref{lemma:D-gamma-inv-h} as the $L^2(\Omega)$-limit of the Tikhonov-regularised correction.

It remains to compute $\D_r(\gamma_t^{\dagger}h)$ in terms of variation processes.

\medskip\noindent The following lemmas compute the Malliavin derivatives needed for the correction term.

\begin{lemma}[Differentiating the Tikhonov resolvent]\label{lemma:D-inverse}
Let $K: \Omega \to L_1(H)$ be a trace-class self-adjoint non-negative operator which, viewed as an $L_2(H,H)$-valued random variable, satisfies $K \in \mathbb{D}^{1,2}(L_2(H,H))$. For each $\eps > 0$ and each deterministic $h \in H$, the random variable
\[
  G^{(\eps)} \;:=\; (K + \eps I)^{-1}\,h \;\in\; H
\]
satisfies the following:

\begin{enumerate}[label=\textup{(\roman*)}, leftmargin=2.2em, itemsep=2pt]
 \item There is first a deterministic bound, $\norm{G^{(\eps)}}_H \leq \eps^{-1}\norm{h}_H$ a.s., hence $G^{(\eps)} \in L^\infty(\Omega; H) \subset \bigcap_{p \geq 1} L^p(\Omega; H)$.
  
 \item $\mathbb{D}^{1,2}$-membership and Leibniz formula. $G^{(\eps)} \in \mathbb{D}^{1,2}(H)$, with
  \begin{equation}\label{eq:D-inverse}
    \D_r G^{(\eps)} \;=\; -(K + \eps I)^{-1}\,(\D_r K)\,(K + \eps I)^{-1}\,h
    \quad \text{in } L_2(U, H), \text{ for a.e.\ } r \in [0, T].
  \end{equation}
  More generally, if $K \in \mathbb{D}^{1,p}(L_2(H,H))$ for some $p \in [2, \infty)$, then $G^{(\eps)} \in \mathbb{D}^{1,p}(H)$ and the same formula~\eqref{eq:D-inverse} holds in $L^p$.
  
 \item There is next a pathwise bound on the derivative, valid almost surely and for almost every $r$,
  \[
    \norm{\D_r G^{(\eps)}}_{L_2(U, H)} \;\leq\; \eps^{-2}\,\norm{h}_H\,\norm{\D_r K}_{L_2(U, L_2(H,H))} ,
  \]
\end{enumerate}
\end{lemma}

\begin{remark}[Why state the lemma at the vector level]\label{rem:D-inverse-vector-level}
The conclusion is stated for the $H$-valued random variable $G^{(\eps)} = (K + \eps I)^{-1}h$ at a fixed deterministic $h$, rather than for the operator-valued resolvent $(K + \eps I)^{-1} \in L(H)$. This is deliberate. The space $L(H)$ in operator norm is non-separable and lacks the unconditional martingale difference (UMD) property when $H$ is infinite-dimensional, so the operator-valued Malliavin space $\mathbb{D}^{1,2}(L(H))$ is delicate, and the UMD-Banach Malliavin chain rule of Pronk--Veraar~\cite{pronkveraar2014umd} does not directly apply to maps with target $L(H)$. An operator-level statement in a properly chosen Banach algebra (for example, in a unitisation of a Schatten ideal $S_p(H)$ with $1 < p < \infty$, which is UMD and separable) is possible but technically heavier and not needed here. Each application in this paper is at a fixed deterministic direction $h$ (the direction along which the Fomin derivative of $\mu_t$ is computed), so the vector-level statement above is exactly what is needed.
\end{remark}

\begin{proof}
The proof proceeds by Galerkin approximation and closability of $\D$. No localisation argument on $\Omega$ is used.

Since $K \geq 0$ a.s.\ and $K$ is self-adjoint, $\sigma(K) \subset [0, \infty)$ a.s. By the spectral mapping theorem on positive self-adjoint operators, $\sigma(K + \eps I) \subset [\eps, \infty)$ a.s., so $K + \eps I$ is invertible on $H$ with $\norm{(K + \eps I)^{-1}}_{L(H)} \leq \eps^{-1}$ deterministically (i.e., pathwise). Hence $\norm{G^{(\eps)}}_H \leq \eps^{-1}\norm{h}_H$ a.s., which proves~(i).

Fix a deterministic orthonormal basis $\{e_k\}_{k \geq 1}$ of $H$, let $\Pi_N: H \to H_N := \mathrm{span}\{e_1,\allowbreak \ldots,\allowbreak e_N\}$ denote the orthogonal projection, and define the Galerkin truncation
\[
  K^N \;:=\; \Pi_N\,K\,\Pi_N \;\in\; L_1(H), \qquad G^{(\eps, N)} \;:=\; (K^N + \eps I)^{-1}\,\Pi_N h \;\in\; H_N \;\subset\; H.
\]
Each matrix entry $(K^N)_{ij} := \ip{K e_j}{e_i}_H$ is a scalar random variable in $\mathbb{D}^{1,2}(\R)$ (it equals the image of $K$ under the bounded linear functional $A \mapsto \ip{A e_j}{e_i}_H$ on $L_2(H,H)$, and $K \in \mathbb{D}^{1,2}(L_2(H,H))$ by hypothesis). Moreover $K^N \geq 0$ a.s.\ on $H_N$, hence $\norm{(K^N + \eps I)^{-1}}_{L(H_N)} \leq \eps^{-1}$ a.s., uniformly in $N$.

The map $\iota_\eps: \{M \in L(H_N): M = M^*,\, \sigma(M) \subset (-\eps/2, \infty)\} \to L(H_N)$, $M \mapsto (M + \eps I_N)^{-1}$, defined on this open set of self-adjoint matrices (which contains the closed cone $\{M \geq 0\}$), is real-analytic with Fr\'echet derivative
\[
  \mathrm{d}\iota_\eps(M)[K] \;=\; -(M + \eps I_N)^{-1}\,K\,(M + \eps I_N)^{-1},
\]
of operator norm at most $\eps^{-2}$ on the cone (deterministic bound). The finite-dimensional Malliavin chain rule (Nualart \cite{nualart2006malliavin}, applied component-wise via the scalar chain rule to the smooth functions sending the entries of $M$ to the entries of $(M + \eps I_N)^{-1}$, which are rational on the cone and have all derivatives bounded uniformly on $\{M \geq 0\}$) gives $(K^N + \eps I)^{-1}|_{H_N} \in \mathbb{D}^{1,p}(L(H_N))$ at the same exponent $p$ at which $K^N \in \mathbb{D}^{1,p}(L(H_N))$ (which is $p = 2$ under the standing hypothesis $K \in \mathbb{D}^{1,2}(L_2(H,H))$, since $K^N = \Pi_N K \Pi_N$ inherits the Malliavin regularity of $K$), with
\[
  \D_r (K^N + \eps I)^{-1}|_{H_N} \;=\; -(K^N + \eps I)^{-1}\,(\D_r K^N)\,(K^N + \eps I)^{-1}.
\]
Applying both sides to the deterministic vector $\Pi_N h \in H_N$ gives
\[
  \begin{aligned}
  &G^{(\eps, N)} \;\in\; \mathbb{D}^{1,2}(H_N) \;\subset\; \mathbb{D}^{1,2}(H), \\
  &\D_r G^{(\eps, N)} \;=\; -(K^N + \eps I)^{-1}\,(\D_r K^N)\,(K^N + \eps I)^{-1}\,\Pi_N h
  \end{aligned}
\]
in $L_2(U, H)$ for a.e.\ $r$. (If $K \in \mathbb{D}^{1,p}(L_2(H,H))$ for some $p > 2$, the same identity is valid in $\mathbb{D}^{1,p}$, again at the same exponent as the hypothesis.)

The orthogonal projection $\Pi_N \to I_H$ strongly on $H$ as $N \to \infty$, so $K^N = \Pi_N K \Pi_N \to K$ in $L_1(H)$, hence a fortiori in $L_2(H,H)$, a.s.\ (continuity of the composition $(A, B, C) \mapsto ABC$ from $L(H) \times L_1(H) \times L(H) \to L_1(H)$ together with $\norm{\Pi_N}_{L(H)} \leq 1$, strong-operator convergence, and $\norm{\cdot}_{L_2(H,H)} \leq \norm{\cdot}_{L_1(H)}$; cf.\ \cite{kato1995perturbation}). The second resolvent identity
\[
  (K^N + \eps I)^{-1} - (K + \eps I)^{-1} \;=\; (K + \eps I)^{-1}\,(K - K^N)\,(K^N + \eps I)^{-1},
\]
combined with the uniform deterministic bound $\norm{(K^N + \eps I)^{-1}}_{L(H)}, \norm{(K + \eps I)^{-1}}_{L(H)} \leq \eps^{-1}$ and $\norm{K - K^N}_{L(H)} \leq \norm{K - K^N}_{L_2(H,H)} \to 0$ a.s., yields
\[
  \norm{(K^N + \eps I)^{-1} - (K + \eps I)^{-1}}_{L(H)} \;\leq\; \eps^{-2}\,\norm{K - K^N}_{L_2(H,H)} \;\to\; 0 \quad \text{a.s.}
\]
Consequently $G^{(\eps, N)} = (K^N + \eps I)^{-1}\Pi_N h \to (K + \eps I)^{-1}h = G^{(\eps)}$ pointwise a.s., dominated by the deterministic constant $2\eps^{-1}\norm{h}_H$. Bounded convergence gives $G^{(\eps, N)} \to G^{(\eps)}$ in $L^p(\Omega; H)$ for every $p \in [1, \infty)$.

For the Malliavin derivative, write $R := (K + \eps I)^{-1}$ and $R_N := (K^N + \eps I)^{-1}$. The explicit derivative formula above gives
\[
  \D_r G^{(\eps, N)} - \bigl[-R\,(\D_r K)\,R\,h\bigr] \;=\; -[R_N\,(\D_r K^N)\,R_N\,\Pi_N h - R\,(\D_r K)\,R\,h].
\]
Decompose the right-hand side as a telescoping sum:
\begin{align*}
  R_N\,(\D_r K^N)\,R_N\,\Pi_N h - R\,(\D_r K)\,R\,h
  &\;=\; (R_N - R)\,(\D_r K^N)\,R_N\,\Pi_N h \\
  &\quad + R\,(\D_r K^N - \D_r K)\,R_N\,\Pi_N h \\
  &\quad + R\,(\D_r K)\,(R_N - R)\,\Pi_N h \\
  &\quad + R\,(\D_r K)\,R\,(\Pi_N - I)\,h.
\end{align*}
Using $\norm{R_N - R}_{L(H)} \to 0$ a.s., $\norm{\D_r K^N - \D_r K}_{L_2(U, L_2(H,H))} \to 0$ a.s.\ (from $\D_r K^N = \Pi_N(\D_r K)\Pi_N$ and strong-operator convergence of $\Pi_N$), $\norm{\Pi_N h - h}_H \to 0$, and the uniform deterministic operator-norm bounds $\norm{R}, \norm{R_N} \leq \eps^{-1}$ together with
\[
  \norm{\D_r K^N}_{L_2(U, L_2(H,H))} \;\leq\; \norm{\D_r K}_{L_2(U, L_2(H,H))},
\]
each of the four terms tends to $0$ in $L_2(U, H)$ a.s.\ for a.e.\ $r$. The pathwise dominator
\begin{multline*}
  \norm{\D_r G^{(\eps, N)}}_{L_2(U, H)} \;\leq\; \eps^{-2}\,\norm{h}_H\,\norm{\D_r K^N}_{L_2(U, L_2(H,H))} \\
  \;\leq\; \eps^{-2}\,\norm{h}_H\,\norm{\D_r K}_{L_2(U, L_2(H,H))}
\end{multline*}
is square-integrable on $\Omega \times [0, T]$ by hypothesis $K \in \mathbb{D}^{1,2}(L_2(H,H))$, so dominated convergence gives
\[
  \D G^{(\eps, N)} \;\longrightarrow\; -R\,(\D K)\,R\,h \qquad \text{in } L^2(\Omega \times [0, T]; L_2(U, H)).
\]

The Malliavin derivative $\D: \mathbb{D}^{1,2}(H) \to L^2(\Omega \times [0, T]; L_2(U, H))$ is a closed operator \cite{nualart2006malliavin}. Combining the above, $G^{(\eps, N)} \in \mathbb{D}^{1,2}(H)$ for every $N$, $G^{(\eps, N)} \to G^{(\eps)}$ in $L^2(\Omega; H)$, and $\D G^{(\eps, N)} \to -R\,(\D K)\,R\,h$ in $L^2(\Omega \times [0, T]; L_2(U, H))$. Closedness of the Malliavin derivative gives $G^{(\eps)} \in \mathbb{D}^{1,2}(H)$ with $\D G^{(\eps)} = -R\,(\D K)\,R\,h$, which is~\eqref{eq:D-inverse}. The same argument with $L^p$ in place of $L^2$ throughout (using closability of $\D: \mathbb{D}^{1,p}(H) \to L^p(\Omega \times [0, T]; L_2(U, H))$, and the strengthened hypothesis $K \in \mathbb{D}^{1,p}(L_2(H,H))$ to dominate $\|\D_r K\|$ in $L^p$) yields the $\mathbb{D}^{1,p}$-statement under stronger hypothesis. The pathwise bound~(iii) is immediate from~\eqref{eq:D-inverse} and the deterministic resolvent estimate $\norm{R}_{L(H)} \leq \eps^{-1}$.
\end{proof}

\begin{lemma}[Differentiating $\gamma_t^{\dagger}h$]\label{lemma:D-gamma-inv-h}
Under Assumptions~\ref{ass:LR}, \ref{ass:diff2}, \ref{ass:SC}, and~\ref{ass:SC-raised} at target $q^* = q/(q-2)$ (where $q \in (2,\infty]$ is the exponent from Assumption~\ref{ass:nondeg}(ii), with $q^* = 1$ at $q = \infty$; cf.\ Remark~\ref{rem:SC-raised-role}), let $h$ satisfy Assumption~\ref{ass:nondeg} (all three clauses, including the Cameron--Martin compatibility~\eqref{eq:CM-compat-axiom}--\eqref{eq:CM-compat-dominator}). Assume the fixed-Tikhonov substitution regularity~\ref{hyp:iii-F-prime} of Theorem~\ref{thm:tikhonov-trace} for $(\mathcal{T}, \gamma_F, h) = (\Phi, \gamma_t, h)$, through either of two admissible sufficient routes. The operator-level route asks that $u_h^{(\eps)} = \Phi^{*}(\gamma_t + \eps I)^{-1}h$ lie in $\mathbb{D}^{1,2}(\HW)$ for every $\eps \in (0,1]$, as in~\eqref{eq:uheps-D12}. The basis-wise route, which is the one used below, asks instead that $G^{(\eps)} = (\gamma_t + \eps I)^{-1}h$ lie in $\mathbb{D}^{1,2}(H)$ for every $\eps > 0$, together with the H\"older-pair joint integrability~\eqref{eq:linear-sub-Holder} of the pair $(\Phi, G^{(\eps)})$, these being the hypotheses of Theorem~\ref{thm:linear-sub-D12}. Both routes are available automatically at the symmetric edge $p^* \geq 6$, equivalently $2 < q \leq 3$ on the edge choice $p^* = 2q/(q-2)$; sufficient conditions for $q > 3$ are catalogued in Remark~\ref{rem:covering-iiiF-SPDE}.

For each $\eps > 0$, the Tikhonov regularisation
\[
  G^{(\eps)} := (\gamma_t + \eps I)^{-1}\,h
\]
belongs to $\mathbb{D}^{1,2}(H)$ with the exact Leibniz formula
\begin{equation}\label{eq:D-Geps}
  \D_r G^{(\eps)} = -(\gamma_t + \eps I)^{-1}\,(\D_r\gamma_t)\,(\gamma_t + \eps I)^{-1}\,h,
\end{equation}
and the pointwise spectral bound $\norm{G^{(\eps)}}_H \leq \norm{\gamma_t^{\dagger}h}_H$ together with $G^{(\eps)} \to \gamma_t^{\dagger}h$ in $L^q(\Omega; H)$ as $\eps \downarrow 0$. Define the associated covering fields intrinsically via the Hilbert--Schmidt kernel $\Phi_r := D_r X(t) \in L_2(U, H)$:
\[
  v_h^{(\eps)}(r) := \Phi_r^{*}\,G^{(\eps)}, \qquad v_h(r) := \Phi_r^{*}\,\gamma_t^{\dagger}h.
\]

Then the covering fields satisfy the explicit sharp spectral $\HW$-estimate
\begin{equation}\label{eq:vh-L2-estimate}
  \E\!\left[\int_0^t \norm{v_h^{(\eps)}(r) - v_h(r)}_U^2\dr\right] \;\leq\; \frac{\eps}{4}\,\E\!\left[\norm{\gamma_t^{\dagger}h}_H^2\right],
\end{equation}
so that $v_h^{(\eps)} \to v_h$ in $L^2(\Omega;\HW)$ at the rate $O(\eps^{1/2})$. Each $v_h^{(\eps)}$ belongs to $\mathbb{D}^{1,2}(\HW) \subset \Dom(\delta_U)$, the limit $v_h$ belongs to $\Dom(\delta_U)$, and
\begin{equation}\label{eq:delta-Uconv}
  \delta_U\bigl(v_h^{(\eps)}\bigr) \;\longrightarrow\; \delta_U(v_h) \qquad \text{in } L^2(\Omega).
\end{equation}
Consequently, the scalar correction integral, defined intrinsically with the positive sign convention so that $\delta_U(v_h^{(\eps)}) = M^{(\eps)} + C^{(\eps)}$ (i.e., absorbing the Nualart minus sign of~\eqref{eq:nualart-substitution} via the identity $\D G^{(\eps)} = -(\gamma_t+\eps I)^{-1}(\D\gamma_t)(\gamma_t+\eps I)^{-1}h$),
\begin{equation}\label{eq:C-eps-positive-def}
  C^{(\eps)} \;:=\; -\int_0^t \Tr_U\bigl[\Phi_r^{*}\,\D_r G^{(\eps)}\bigr]\dr \;=\; \int_0^t \Tr_U\bigl[\Phi_r^{*}\,(\gamma_t+\eps I)^{-1}\,[\D_r\gamma_t]\,(\gamma_t+\eps I)^{-1}\,h\bigr]\dr,
\end{equation}
converges in $L^2(\Omega)$ as $\eps \downarrow 0$; the limit is independent of the Moore--Penrose freedom of $\D_r(\gamma_t^{\dagger}h)$ on $\ker(\gamma_t)$ and is denoted symbolically by
\begin{equation}\label{eq:D-gamma-inv-h}
  C \;:=\; \lim_{\eps \downarrow 0} C^{(\eps)} \;=\; \int_0^t \Tr_U\bigl[\Phi_r^{*}\,\gamma_t^{\dagger}\,(\D_r\gamma_t)\,\gamma_t^{\dagger}\,h\bigr]\dr.
\end{equation}
\end{lemma}

\begin{remark}[The meaning of $\D_r(\gamma_t^{\dagger}h)$]\label{rem:Dgamma-dagger-symbolic}
The constant $1/4$ in~\eqref{eq:vh-L2-estimate} is sharp, being attained spectrally at the eigendirection on which $\lambda_k = \eps$ in the diagonalisation of $\gamma_t$; the rate $O(\eps^{1/2})$ is therefore not improvable. No claim is made that $\gamma_t^{\dagger}h$ belongs to $\mathbb{D}^{1,2}_{\mathrm{loc}}(H)$, or to any other local Malliavin--Sobolev space \cite{nualart2006malliavin}. The notation
\[
  \D_r(\gamma_t^{\dagger}h) \;=\; -\gamma_t^{\dagger}\,(\D_r\gamma_t)\,\gamma_t^{\dagger}\,h \qquad \text{(modulo $\ker(\gamma_t)$)}
\]
is shorthand for the Tikhonov trace limit~\eqref{eq:D-gamma-inv-h}. What the lemma asserts is the $L^2(\Omega)$-convergence of the regularised correction $C^{(\eps)}$ of~\eqref{eq:C-eps-positive-def}, and that requires only the differentiability of the bounded resolvent $(\gamma_t + \eps I)^{-1}$ applied to the deterministic vector $h$ (the content of Lemma~\ref{lemma:D-inverse}), not Malliavin differentiability of the Moore--Penrose pseudoinverse $\gamma_t^{\dagger}$ as an operator-valued or pseudoinverse-valued random variable. The minus sign is the standard one for resolvent derivatives, and is absorbed by~\eqref{eq:C-eps-positive-def} into the positive-sign convention for $C$. The $\ker(\gamma_t)$-ambiguity does not contribute to~\eqref{eq:D-gamma-inv-h} because the covering-field factor $\Phi_r^{*}$ annihilates $\ker(\gamma_t)$ for a.e.\ $r$ (Remark~\ref{rem:kernel-annih}). The underlying measure-theoretic picture (Fomin differentiability of $\mu_t$) is discussed in the differentiable-measures framework of~\cite{bogachev2010differentiable}.

\end{remark}

\begin{remark}[SPDE-level verification of \ref{hyp:iii-F-prime} at the symmetric edge and beyond]\label{rem:lemma-D-gamma-iiiF-status}
Hypothesis~\ref{hyp:iii-F-prime} holds automatically at $p^* \geq 6$, equivalently $2 < q \leq 3$ on the edge choice $p^* = 2q/(q-2)$, via Theorem~\ref{thm:linear-sub-D12} and Remark~\ref{rem:covering-iiiF-SPDE}.

For $q > 3$, the H\"older-pair joint integrability~\eqref{eq:linear-sub-Holder} (or equivalently the operator-level \ref{hyp:iii-F-prime}) is an additional structural hypothesis on top of Assumption~\ref{ass:SC-raised}, and it is delivered in concrete examples by any of several routes, namely by raising the target to some $q^{**} > q^*$ in Assumption~\ref{ass:SC-raised}, chosen so that the moment chain of Remark~\ref{rem:Phi-all-moments} delivers $\D\gamma_t \in L^r$ at the matched exponent $r = 2p^*/(p^*-2)$; by the sharper $H$-slot bound on the diffusion forcing available when $\calB''_{uu}$ extends to $H \times H$ (Remark~\ref{rem:SC-raised-two-regimes}), which improves the moments of $\calW$ and hence of $\D\gamma_t$; by the bounded-Malliavin-matrix input $\mathcal{T} \in L^\infty(\Omega; L_2(\HW,H))$ of Theorem~\ref{thm:tikhonov-trace}, clause~\ref{hyp:iii-F-prime-c}; or by a structural verification specific to the equation at hand.
The lemma is proved at this maximum level of generality, with the basis-wise hypothesis~\eqref{eq:Geps-D12}+\eqref{eq:linear-sub-Holder} treated as an explicit hypothesis at the SPDE level matching Theorem~\ref{thm:tikhonov-trace}.
\end{remark}

\begin{proof}
Since $\gamma_t \geq 0$ is self-adjoint, the operator $\gamma_t + \eps I$ is pointwise invertible on $H$ with the deterministic bound $\norm{(\gamma_t + \eps I)^{-1}}_{L(H)} \leq \eps^{-1}$. By Proposition~\ref{prop:gamma-properties}, $\gamma_t$ is trace-class self-adjoint non-negative, and $\gamma_t \in \mathbb{D}^{1,2}(L_2(H,H))$ by~\eqref{eq:CM-gamma-D12} of clause~\textup{(iii)} of Assumption~\ref{ass:nondeg} (equivalently, by Proposition~\ref{prop:D-gamma} whenever $\Phi \in \mathbb{D}^{1,p_0}$ with $p_0 \geq 4$, as it is at $p^* \geq 6$). Lemma~\ref{lemma:D-inverse} therefore applies directly with $K = \gamma_t$ and the deterministic direction $h$, yielding $G^{(\eps)} = (\gamma_t + \eps I)^{-1}h \in \mathbb{D}^{1,2}(H)$, the exact Leibniz formula~\eqref{eq:D-Geps}, and the pathwise bound $\norm{\D_r G^{(\eps)}}_{L_2(U, H)} \leq \eps^{-2}\norm{h}_H\norm{\D_r\gamma_t}_{L_2(U, L_2(H,H))}$. (If $\gamma_t \in \mathbb{D}^{1,p}(L_2(H,H))$ for some $p > 2$ under stronger structural input, $G^{(\eps)} \in \mathbb{D}^{1,p}(H)$ likewise; this stronger membership is not needed for the trace-limit argument below, which uses only $L^2(\Omega)$-convergence.) No localisation on $\Omega$ and no operator-level Malliavin calculus for $(\gamma_t + \eps I)^{-1} \in L(H)$ is invoked at this stage; only the vector-level statement for $h$ deterministic is used.

Let $\{(\lambda_k, e_k)\}_{k \geq 1}$ be the (random) spectral resolution of the compact positive self-adjoint operator $\gamma_t$, and write $h = \sum_k h_k\,e_k$ with $h_k = \ip{h}{e_k}_H$. By Assumption \ref{ass:nondeg}, $h \in \Ran(\gamma_t)$ a.s., so in particular $h_k = 0$ whenever $\lambda_k = 0$, and $\sum_{k:\lambda_k > 0}(h_k/\lambda_k)^2 = \norm{\gamma_t^\dagger h}_H^2 < \infty$ a.s. The functional calculus gives
\[
  G^{(\eps)} \;=\; \sum_{k:\,\lambda_k > 0}\frac{h_k}{\lambda_k + \eps}\,e_k, \qquad \gamma_t^{\dagger}h \;=\; \sum_{k:\,\lambda_k > 0}\frac{h_k}{\lambda_k}\,e_k.
\]
Term-by-term inspection yields the pointwise bound
\begin{equation}\label{eq:Geps-spectral-bound}
  \norm{G^{(\eps)}(\omega)}_H^2 \;=\; \sum_{k:\,\lambda_k > 0}\frac{h_k^2}{(\lambda_k + \eps)^2} \;\leq\; \sum_{k:\,\lambda_k > 0}\frac{h_k^2}{\lambda_k^2} \;=\; \norm{\gamma_t^{\dagger}h(\omega)}_H^2
\end{equation}
and the pointwise convergence $G^{(\eps)}(\omega) \to \gamma_t^{\dagger}h(\omega)$ in $H$ for every $\omega$ where $h \in \Ran(\gamma_t)$. Dominated convergence on $\Omega$, using the moment bound of Assumption \ref{ass:nondeg}, upgrades this to $G^{(\eps)} \to \gamma_t^{\dagger}h$ in $L^q(\Omega; H)$.

Setting $w_\eps := G^{(\eps)} - \gamma_t^{\dagger}h$ and using $\int_0^t \norm{\Phi_r^* z}_U^2\dr = \ip{z}{\gamma_t z}_H$, we compute spectrally. With $z_k := h_k/\lambda_k$ for $\lambda_k > 0$,
\[
  G^{(\eps)} - \gamma_t^{\dagger}h \;=\; -\eps\sum_{k:\,\lambda_k > 0}\tfrac{z_k}{\lambda_k+\eps}\,e_k,
\]
so
\begin{equation}\label{eq:gamma-w-eps-spectral}
\begin{split}
  \int_0^t \norm{v_h^{(\eps)}(r) - v_h(r)}_U^2\dr
  &\;=\; \ip{w_\eps}{\gamma_t\,w_\eps}_H
  \;=\; \sum_{k:\,\lambda_k > 0}\frac{\lambda_k\,\eps^2}{(\lambda_k+\eps)^2}\,z_k^2 \\
  &\;\leq\; \frac{\eps}{4}\sum_{k:\,\lambda_k > 0} z_k^2
  \;=\; \frac{\eps}{4}\,\norm{\gamma_t^{\dagger}h}_H^2,
\end{split}
\end{equation}
where the spectral inequality $\lambda\eps^2/(\lambda+\eps)^2 \leq \eps/4$ for all $\lambda \geq 0$ is the elementary maximum at $\lambda = \eps$ of the function $\lambda \mapsto \lambda\eps^2/(\lambda+\eps)^2$ (derivative $\eps^2(\eps-\lambda)/(\lambda+\eps)^3$, zero at $\lambda=\eps$, value $\eps\cdot\eps^2/(2\eps)^2 = \eps/4$). This bound is sharp, attained at $\lambda_k = \eps$ on the eigendirection where $z_k$ concentrates. Taking expectations and invoking Assumption~\ref{ass:nondeg}(ii) (which ensures $\E[\norm{\gamma_t^{\dagger}h}_H^q] < \infty$ for some $q \geq 2$, hence in particular $\E[\norm{\gamma_t^{\dagger}h}_H^2] < \infty$),
\begin{equation}\label{eq:vh-L2-estimate-sharp}
  \E\!\left[\int_0^t \norm{v_h^{(\eps)}(r) - v_h(r)}_U^2\dr\right] \;\leq\; \frac{\eps}{4}\,\E\bigl[\norm{\gamma_t^{\dagger}h}_H^2\bigr],
\end{equation}
so that $v_h^{(\eps)} \to v_h$ in $L^2(\Omega; \HW)$ at rate $O(\eps^{1/2})$.

Under Assumption~\ref{ass:diff2}, $\Phi \in \mathbb{D}^{1, p^*}$ (Remark~\ref{rem:Phi-all-moments}) and $G^{(\eps)} \in \mathbb{D}^{1,2}(H)$ via~\eqref{eq:D-Geps}. The product rule gives
\begin{equation}\label{eq:Dvh-eps-intrinsic}
  \D_s v_h^{(\eps)}(r) \;=\; (\D_s \Phi_r)^{*}\,G^{(\eps)} \;+\; \Phi_r^{*}\,\D_s G^{(\eps)}.
\end{equation}
The first term satisfies, via the HS-norm bound $\|(\D_s\Phi_r)^* z\|_U \leq \|\D_s\Phi_r\|_{L_2(U, H)}\,\|z\|_H$ at $z = G^{(\eps)}$ together with the deterministic bound $\|G^{(\eps)}\|_H \leq \|\gamma_t^{\dagger}h\|_H$,
\[
  \E\!\left[\int_0^t\!\!\int_0^t \|(\D_s \Phi_r)^{*} G^{(\eps)}\|_U^2\,\mathrm{d}r\,\mathrm{d}s\right] \;\leq\; \E\bigl[\|\D \Phi\|_{L_2(\HW; L^2([0,t]; L_2(U, H)))}^2\,\|\gamma_t^{\dagger}h\|_H^2\bigr],
\]
which is finite by H\"older's inequality applied at the conjugate-exponent pair $(p^*/2, p^*/(p^*-2))$ on $(\|\D\Phi\|^2, \|\gamma_t^{\dagger}h\|^2)$, using $\|\D\Phi\| \in L^{p^*}$ from $\Phi \in \mathbb{D}^{1, p^*}$ and $\|\gamma_t^{\dagger}h\| \in L^q = L^{2p^*/(p^*-2)}$ from Assumption~\ref{ass:nondeg}(ii) at the H\"older-triangle pair. The second term satisfies, via the HS-composition bound $\|\Phi_r^{*}\,\D_s G^{(\eps)}\|_{L_2(U, U)} \leq \|\Phi_r\|_{L_2(U, H)}\,\|\D_s G^{(\eps)}\|_{L_2(U, H)}$ together with~\eqref{eq:D-Geps},
\[
  \|\D_s G^{(\eps)}\|_{L_2(U, H)} \;\leq\; \eps^{-2}\,\|h\|_H\,\|\D_s\gamma_t\|_{L_2(U, L_2(H,H))}.
\]
Squaring, integrating over $s \in [0,t]$, and taking $\E[\cdot]$:
\[
  \begin{aligned}
  &\E\!\left[\int_0^t\!\!\int_0^t \|\Phi_r^{*}\,\D_s G^{(\eps)}\|_{L_2(U,U)}^2\,\mathrm{d}s\,\mathrm{d}r\right] \\
  &\qquad \;\leq\; \eps^{-4}\|h\|^2 \,\E\!\left[\int_0^t \|\Phi_r\|_{L_2(U, H)}^2\,\mathrm{d}r\,\cdot\,\|\D\gamma_t\|^2_{L_2(\HW; L_2(H,H))}\right],
  \end{aligned}
\]
which is finite at fixed $\eps > 0$ by H\"older at the matched pair on $(\|\Phi\|^2, \|\D\gamma_t\|^2)$, under hypothesis~\ref{hyp:iii-F-prime} (automatic at $p^* \geq 6$ by Remark~\ref{rem:covering-iiiF-SPDE}). This yields $v_h^{(\eps)} \in \mathbb{D}^{1,2}(\HW) \subset \Dom(\delta_U)$ at fixed $\eps > 0$.

Theorem~\ref{thm:linear-sub-D12} applied to $v_h^{(\eps)} = \Phi^{*}G^{(\eps)}$ gives the exact $\eps$-level identity
\begin{equation}\label{eq:nualart-eps}
  \delta_U\bigl(v_h^{(\eps)}\bigr) \;=\; \underbrace{\bigl[\delta_U(w(\cdot)(z))\bigr]_{z = G^{(\eps)}}}_{M^{(\eps)}} \;-\; \int_0^t \Tr_U\bigl[\Phi_r^{*}\,\D_r G^{(\eps)}\bigr]\,\mathrm{d}r,
\end{equation}
where $M^{(\eps)} = \langle G^{(\eps)}, \mathcal{K}\rangle_H$ via the HS kernel (Lemma~\ref{lemma:HS-kernel}) and
\[
  C^{(\eps)} := -\int_0^t \Tr_U[\Phi_r^{*}\D_r G^{(\eps)}]\,\mathrm{d}r
\]
(positive-sign convention~\eqref{eq:C-eps-positive-def}).

Under Assumption~\ref{ass:SC-raised} at target $q^* = p^*/2$, Remark~\ref{rem:Phi-all-moments} gives $\Phi \in \mathbb{D}^{1, p^*}$, so Lemma~\ref{lemma:HS-kernel} provides a kernel $\mathcal{K} \in L^{p^*}(\Omega; H)$ with $\delta_U(\Phi^* z) = \ip{z}{\mathcal{K}}_H$ for deterministic $z$. At random $Z = G^{(\eps)}$,
\begin{equation}\label{eq:subst-hs}
  M^{(\eps)}(\omega) \;=\; \ip{G^{(\eps)}(\omega)}{\mathcal{K}(\omega)}_H \qquad \text{for }\mathbb{P}\text{-a.e.\ } \omega.
\end{equation}
Subtracting at $Z = \gamma_t^\dagger h$ and applying Cauchy--Schwarz in $H$ then H\"older on $\Omega$ at $(q/2, q/(q-2))$:
\begin{equation}\label{eq:M-eps-conv}
  \norm{M^{(\eps)} - M}_{L^2(\Omega)} \;\leq\; \norm{G^{(\eps)} - \gamma_t^{\dagger}h}_{L^q(\Omega; H)}\,\norm{\mathcal{K}}_{L^{p^*}(\Omega; H)} \;\longrightarrow\; 0 \quad \text{as } \eps \downarrow 0,
\end{equation}
where the first factor tends to zero by~\eqref{eq:Geps-spectral-bound} and $L^q$-dominated convergence, and the second is finite by Lemma~\ref{lemma:HS-kernel}(iii) and Remark~\ref{rem:Phi-all-moments}. Hence $M^{(\eps)} \to M := \ip{\gamma_t^{\dagger}h}{\mathcal{K}}_H$ in $L^2(\Omega)$.

In any measurable spectral resolution $\{(\lambda_k, e_k)\}$ of $\gamma_t$ and deterministic ONB $\{f_j\}$ of $U$,
\begin{equation}\label{eq:Ceps-spectral}
  \Tr_U\!\bigl[\Phi_r^*(\gamma_t + \eps I)^{-1}(\D_r\gamma_t)(\gamma_t + \eps I)^{-1}h\bigr]
  \;=\; \!\!\!\sum_{\substack{j;\,k,l\\ \lambda_k,\lambda_l > 0}}\!\!\frac{h_k\,\langle(\D_r\gamma_t)(f_j)\,e_k,\,e_l\rangle_H\,\langle e_l,\,\Phi_r f_j\rangle_H}{(\lambda_l + \eps)(\lambda_k + \eps)},
\end{equation}
where the $\lambda_k = 0$ and $\lambda_l = 0$ terms vanish by the range condition and the kernel-annihilation property~\eqref{eq:kernel-annih}. By clause~(iii) of Assumption~\ref{ass:nondeg}, the scalar integrands converge, $\mathcal{I}_r^{(\eps)} \to \mathcal{I}_r^{(0)}$ a.s.\ for a.e.\ $r$, and are dominated, $|\mathcal{I}_r^{(\eps)}| \leq \Psi_r$ with $\E[(\int_0^t\Psi_r\,\mathrm{d}r)^2] < \infty$; dominated convergence therefore gives
\begin{equation}\label{eq:Ceps-limit}
  C^{(\eps)} \;=\; \int_0^t \mathcal{I}_r^{(\eps)}\dr \;\longrightarrow\; C \;:=\; \int_0^t \mathcal{I}_r^{(0)}\dr \qquad \text{in } L^2(\Omega),
\end{equation}
the limit being denoted symbolically by $\int_0^t\Tr_U[\Phi_r^{*}\gamma_t^{\dagger}(\D_r\gamma_t)\gamma_t^{\dagger}h]\dr$. Under the spectral sufficient condition of Remark~\ref{rem:CM-compat-spectral}, $\mathcal{I}_r^{(0)}$ is represented by the absolutely convergent series obtained from~\eqref{eq:Ceps-spectral} by replacing $(\lambda_k+\eps)^{-1}(\lambda_l+\eps)^{-1}$ by $(\lambda_k\lambda_l)^{-1}$; without it, \eqref{eq:Ceps-limit} is the definition of the correction and no interchange of an infinite spectral sum with the Tikhonov limit is performed. Clause~(iii) is essential here, and both of its parts are used, since kernel annihilation controls only the $\lambda_k = 0$ terms, neither the convergence nor the small-positive-eigenvalue contributions.

Combining (M) and (C), $\delta_U(v_h^{(\eps)}) = M^{(\eps)} + C^{(\eps)} \to M + C$ in $L^2(\Omega)$. Together with $v_h^{(\eps)} \to v_h$ in $L^2(\Omega; \HW)$, closedness of $\delta_U$ \cite{nualart2006malliavin} yields $v_h \in \Dom(\delta_U)$ with $\delta_U(v_h) = M + C$.

The symbolic expression $\gamma_t^{\dagger}(\D_r\gamma_t)\gamma_t^{\dagger}h$ in~\eqref{eq:D-gamma-inv-h} is shorthand for the Tikhonov limit $C = \lim_{\eps \downarrow 0} C^{(\eps)}$, not a pointwise operator composition (since $(\D_r\gamma_t)(\gamma_t^{\dagger}h)$ need not lie in $\Ran(\gamma_t)$). The limit is well-defined under clause~\textup{(iii)} by the scalar convergence-and-domination argument. Under Remark~\ref{rem:CM-compat-spectral}, spectral dominated convergence is one sufficient verification of that hypothesis. The kernel-annihilation property ensures that any $\ker(\gamma_t)$-ambiguity is killed by $\Phi_r^{*}$ under the trace.
\end{proof}

\begin{remark}[Kernel-annihilation property]\label{rem:kernel-annih}
The Hilbert--Schmidt Malliavin derivative $\Phi_r^* = [D_r X(t)]^*$ satisfies
\begin{equation}\label{eq:kernel-annih}
  \Phi_r^* \circ P_{\ker(\gamma_t)} \;=\; 0 \qquad \text{$\mathbb{P}$-a.s., Leb-a.e.\ } r \in [0,t],
\end{equation}
since for $v \in \ker(\gamma_t)$,
\begin{equation}\label{eq:kernel-annih-energy}
  \int_0^t \norm{\Phi_r^*\,v}_U^2\,\mathrm{d}r \;=\; \langle \gamma_t\,v,\,v\rangle_H \;=\; 0.
\end{equation}
The projection $P_{\ker(\gamma_t)}$ is strongly $\F_t$-measurable via the Tikhonov representation $P_{\ker(\gamma_t)} = I_H - \mathrm{s\text{-}lim}_{\eps \downarrow 0}\,\gamma_t(\gamma_t + \eps I)^{-1}$. Property~\eqref{eq:kernel-annih} is the structural reason the pseudoinverse formula makes sense without closed-range assumptions on $\gamma_t$.
\end{remark}

\begin{remark}[Tikhonov convergence mechanism]\label{rem:trace-conv-tech}
The Tikhonov family $\{v_h^{(\eps)}\}$ is not uniformly bounded in $\mathbb{D}^{1,2}(\HW)$ ($\D_s G^{(\eps)}$ scales as $\eps^{-1}$). Instead, convergence of $\delta_U(v_h^{(\eps)})$ is established directly. The main term converges because $M^{(\eps)} = \ip{G^{(\eps)}}{\mathcal{K}}_H \to \ip{\gamma_t^{\dagger}h}{\mathcal{K}}_H = M$, by H\"older's inequality and the Hilbert--Schmidt kernel representation. The correction converges because clause~(iii) of Assumption~\ref{ass:nondeg} supplies both the a.s.\ convergence $\mathcal{I}_r^{(\eps)} \to \mathcal{I}_r^{(0)}$ and an $L^2$-integrable dominator, so that ordinary dominated convergence identifies the limit; the spectral condition~\eqref{eq:CM-compat-spectral} is one sufficient route to both, through the bound $(\lambda+\eps)^{-1} \leq \lambda^{-1}$. Closedness of $\delta_U$ then places $v_h$ in $\Dom(\delta_U)$ with $\delta_U(v_h) = M + C$.
\end{remark}

\medskip\noindent The correction term $C_h$ involves $\D_r\gamma_t$, which the second Malliavin fibre expresses in one line.

\begin{proposition}[Malliavin derivative of $\gamma_t$]\label{prop:D-gamma}
Under Assumptions~\ref{ass:LR}, \ref{ass:diff2}, and~\ref{ass:SC}, suppose in addition that $\Phi \in \mathbb{D}^{1,p_0}(L^2([0,t]; L_2(U,H)))$ for some $p_0 \geq 2$. Then the Malliavin covariance operator $\gamma_t$ is trace class a.s.\ and, viewed as an $L_2(H,H)$-valued random variable, belongs to $\mathbb{D}^{1,p_0/2}(L_2(H,H))$, with
\begin{equation}\label{eq:D-gamma-formula}
  \D_r\gamma_t \;=\; \int_0^t \Bigl[\calW(t;s,r)\,\Phi_s^{*} \;+\; \Phi_s\,\calW(t;s,r)^{*}\Bigr]\ds
\end{equation}
as an identity in $L_2(U, L_2(H,H))$, where $\Phi_s := \Phi_{t,s} = \D_s X(t)$ and $\calW(t;s,r) = \D_r\Phi_{t,s}$ is the second Malliavin fibre of Definition~\ref{def:second-fibre}, read in the Malliavin slot $r$ and composed with $\Phi_s^{*}$ in the remaining slot.
\end{proposition}

The exponent halves because~\eqref{eq:D-gamma-formula} is quadratic in the pair $(\Phi, \D\Phi)$, so this route delivers $\gamma_t \in \mathbb{D}^{1,2}(L_2(H,H))$ exactly when $p_0 \geq 4$; at smaller $p_0$ that membership is supplied instead by~\eqref{eq:CM-gamma-D12} of clause~\textup{(iii)} of Assumption~\ref{ass:nondeg}.

\begin{proof}
Recall $\gamma_t = \int_0^t \Phi_s\,\Phi_s^{*}\,\mathrm{d}s$. The bilinear map $(A,B) \mapsto AB^{*}$ sends $L_2(U,H) \times L_2(U,H)$ continuously into $L_1(H) \subset L_2(H,H)$, with
\[
  \norm{AB^{*}}_{L_2(H,H)} \;\leq\; \norm{AB^{*}}_{L_1(H)} \;\leq\; \norm{A}_{L_2(U,H)}\,\norm{B}_{L_2(U,H)};
\]
in particular $\gamma_t$ is trace class a.s. The Hilbert-valued Malliavin product rule, obtained by applying the scalar rule basis-wise in deterministic orthonormal bases of $H$ and $U$ and passing to the limit, therefore gives
\begin{equation}\label{eq:D-gamma-interchange}
  \D_r\gamma_t \;=\; \int_0^t \D_r\bigl(\Phi_s\,\Phi_s^{*}\bigr)\ds \;=\; \int_0^t \Bigl[(\D_r\Phi_s)\,\Phi_s^{*} + \Phi_s\,(\D_r\Phi_s)^{*}\Bigr]\ds,
\end{equation}
the interchange with the Bochner integral in $L_2(H,H)$ being justified by
\[
  \norm{\D_r(\Phi_s\,\Phi_s^{*})}_{L_2(U, L_2(H,H))} \;\leq\; 2\,\norm{\D_r\Phi_s}_{L_2(U, L_2(U,H))}\,\norm{\Phi_s}_{L_2(U,H)},
\]
whose right-hand side, after Cauchy--Schwarz in $s$ and H\"older in $\omega$ at the conjugate pair $(2,2)$, lies in $L^{p_0/2}(\Omega)$ because both $\Phi$ and $\D\Phi$ lie in $L^{p_0}$; this is the source of the halved exponent. Lemma~\ref{lemma:DPhi} identifies $\D_r\Phi_s = \D_r\Phi_{t,s}$ with $\calW(t;s,r)$ for almost every $(s,r)$, which turns~\eqref{eq:D-gamma-interchange} into~\eqref{eq:D-gamma-formula}, the $s$-integral being insensitive to a null set.
\end{proof}

Formula~\eqref{eq:D-gamma-formula} is a single symmetric expression valid on the whole square $\{0 \leq s, r \leq t\}$, with no splitting at the diagonal. The two mechanisms that a splitting would exhibit separately are both carried by $\calW$, and they are visible in its own data rather than in the formula for $\D_r\gamma_t$. The second derivatives $\calA''_{uu}$ and $\calB''_{uu}$ enter through the forcings~\eqref{eq:bold-forcings} of~\eqref{eq:W-equation}, and the response of the diffusion coefficient to a perturbation of the state enters through the initial datum~\eqref{eq:W-initial}, whose two branches are the two orders of $r$ and $s$. When $\calA''_{uu} = \calB''_{uu} = 0$ and $\calB'_u = 0$, both the forcing and the datum vanish, so $\calW \equiv 0$ and $\D_r\gamma_t \equiv 0$; this is the linear equation with state-independent diffusion, where the correction term of Theorem~\ref{thm:main} is absent. In finite dimensions the four operator-valued integrands obtained by expanding~\eqref{eq:D-gamma-formula} through the data of $\calW$ are the terms $I_1^{p,q}, I_2^{p,q}, I_3^{p,q}, I_4^{p,q}$ of Theorem 2.1 of \cite{mirafzali2025malliavin}, the first two coming from the forcing and the last two from the initial datum.

\subsection{Assembly}\label{subsec:assembly}

Every ingredient is now available, and it remains only to put them in order.

\begin{proof}[Proof of Theorem \ref{thm:main}]
We verify the hypotheses of Theorems~\ref{thm:abstract-bismut-fomin}--\ref{thm:tikhonov-trace} for the SPDE-specialised quadruple $(F, \mathcal{T}, \gamma_F, \Phi)$ with $F = X(t)$, $\mathcal{T} = DX(t)$, $\gamma_F = \gamma_t$, and $\Phi_r = D_r X(t) = Y(t,r)\calB(r, X(r))$ (the last identification is Proposition~\ref{prop:D-Xt}). The Bismut formula (Theorem \ref{thm:bismut-formula}, which is the SPDE realisation of Theorem~\ref{thm:abstract-bismut-fomin}), the linear-field substitution (Theorem \ref{thm:linear-sub-D12}), and the explicit Malliavin derivative computations of the preceding sections combine as follows.

By Theorem \ref{thm:bismut-formula}, the logarithmic derivative is
\begin{equation}\label{eq:proof-bismut}
  \beta_h(X(t)) = -\E[\delta_U(v_h) \mid X(t)],
\end{equation}
which is precisely the Bismut--Fomin identity \eqref{eq:abstract-bismut-fomin} of Theorem~\ref{thm:abstract-bismut-fomin} specialised at $F = X(t)$, $u_h = v_h = \Phi^{*}\gamma_t^{\dagger}h$ (the canonical pseudoinverse covering field of Theorem~\ref{thm:pseudoinverse-covering}). It remains to obtain an explicit representation of $\delta_U(v_h)$ in terms of variation processes; this is the content of Theorem~\ref{thm:tikhonov-trace} executed in the SPDE-specific notation.

Recall $v_h(r) = \Phi_r^{*}\,\tilde{h}$ with $\tilde{h} = \gamma_t^{\dagger}h$. By Lemma~\ref{lemma:linearity}, $\nabla_z w_r = \Phi_r^{*}$, independent of $z$, satisfying the kernel-annihilation property~\eqref{eq:kernel-annih} by Remark~\ref{rem:kernel-annih}.

For $\eps > 0$, define
\[
  G^{(\eps)} := (\gamma_t + \eps I)^{-1}\,h \in \mathbb{D}^{1,2}(H), \qquad v_h^{(\eps)}(r) := w_r\bigl(G^{(\eps)}\bigr) = \Phi_r^{*}\,G^{(\eps)}.
\]
At fixed $\eps > 0$, clause~\textup{(iii)} of Assumption~\ref{ass:nondeg} and Lemma~\ref{lemma:D-inverse} give
\[
  G^{(\eps)} \;=\; (\gamma_t + \eps I)^{-1}h \;\in\; \mathbb{D}^{1,2}(H),
\]
with
\[
  \D_r G^{(\eps)} \;=\; -(\gamma_t + \eps I)^{-1}(\D_r\gamma_t)(\gamma_t + \eps I)^{-1}h .
\]
No stronger $\mathbb{D}^{1,p}$-membership of $G^{(\eps)}$ is needed in the proof. The operative substitution argument is Theorem~\ref{thm:linear-sub-D12}; the required joint product integrability is precisely the fixed-Tikhonov substitution regularity~\ref{hyp:iii-F-prime}.

The membership $v_h^{(\eps)} = \Phi^{*}G^{(\eps)} \in \mathbb{D}^{1,2}(\HW) \subset \Dom(\delta_U)$ at each fixed $\eps > 0$ is not inferred from the $\mathbb{D}^{1,2}$-regularity of $G^{(\eps)}$ alone. The composition of $\Phi$ with $G^{(\eps)}$ requires joint product integrability of the pair, which is exactly the content of the fixed-Tikhonov substitution regularity~\ref{hyp:iii-F-prime} hypothesised in Part~II. We use that hypothesis directly, through either of its two admissible sufficient routes (operator-level, or the H\"older-pair form~\eqref{eq:linear-sub-Holder} feeding the basis-wise route of Theorem~\ref{thm:linear-sub-D12}).
The fixed-$\eps$ substitution identity is then supplied by Theorem~\ref{thm:tikhonov-trace}, equivalently, in its basis-wise formulation, by Theorem~\ref{thm:linear-sub-D12}, which needs only $G^{(\eps)} \in \mathbb{D}^{1,2}(H)$ together with the H\"older-pair integrability~\eqref{eq:linear-sub-Holder}, in the form
\begin{equation}\label{eq:sub-eps}
  \delta_U\bigl(v_h^{(\eps)}\bigr) \;=\; \bigl[\delta_U(w(\cdot)(z))\bigr]_{z = G^{(\eps)}} \;-\; \int_0^t \Tr_U\bigl[\nabla_z w_r \cdot \D_r G^{(\eps)}\bigr]\dr.
\end{equation}
No finite-dimensional projection and no $\mathbb{D}^{1,4}_{\mathrm{loc}}$-regularity are used, since the basis-wise route collapses the chain rule to a single-derivative operation per coordinate, which is exactly why it operates at every point of the H\"older triangle~\eqref{eq:nondeg-Holder-triangle}.

We pass to the limit in \eqref{eq:sub-eps} term-by-term, invoking Lemma \ref{lemma:D-gamma-inv-h} for each of the three components.

By Lemma \ref{lemma:D-gamma-inv-h}, the covering field converges in $L^2(\Omega; \HW)$ with the quantitative estimate \eqref{eq:vh-L2-estimate}, and the Skorokhod integrals converge in $L^2(\Omega)$:
\[
  \delta_U\bigl(v_h^{(\eps)}\bigr) \;\longrightarrow\; \delta_U(v_h) \qquad \text{in } L^2(\Omega), \qquad v_h \in \Dom(\delta_U).
\]

By the main-term argument of Lemma~\ref{lemma:D-gamma-inv-h}, the map $z \mapsto \delta_U(\Phi^* z)$ is Hilbert--Schmidt from $H$ to $L^2(\Omega; \R)$, represented by a unique kernel $\mathcal{K} \in L^2(\Omega; H)$ via $\delta_U(\Phi^* z) = \ip{z}{\mathcal{K}}_H$ (deterministic $z$), whose extension to random $Z$ is the pointwise inner product $\omega \mapsto \ip{Z(\omega)}{\mathcal{K}(\omega)}_H$, and the $L^q$-convergence $G^{(\eps)} \to \gamma_t^{\dagger}h$ combined with H\"older's inequality and the higher-moment bound on $\mathcal{K}$ under Assumption~\ref{ass:diff2} yields
\[
  \bigl[\delta_U(w(\cdot)(z))\bigr]_{z = G^{(\eps)}} \;=\; \ip{G^{(\eps)}}{\mathcal{K}}_H \;\longrightarrow\; \ip{\gamma_t^{\dagger}h}{\mathcal{K}}_H \;=\; \bigl[\delta_U(w(\cdot)(z))\bigr]_{z = \gamma_t^{\dagger}h} \quad \text{in } L^2(\Omega).
\]

From the substitution identity \eqref{eq:sub-eps} and the two convergences just established, the positive-sign correction term~\eqref{eq:C-eps-positive-def}
\[
  C^{(\eps)} \;=\; -\int_0^t \Tr_U\bigl[\Phi_r^{*}\,\D_r G^{(\eps)}\bigr]\dr \;=\; \delta_U\bigl(v_h^{(\eps)}\bigr) \;-\; M^{(\eps)}
\]
converges in $L^2(\Omega)$ as $\eps \downarrow 0$. The limit equals (positive-sign symbolic form, intrinsic in $\Phi^{*}$)
\[
  C \;=\; \int_0^t \Tr_U\bigl[\Phi_r^{*}\,\gamma_t^{\dagger}\,(\D_r\gamma_t)\,\gamma_t^{\dagger}\,h\bigr]\dr
\]
(interpreted modulo the kernel-annihilation of Remark~\ref{rem:kernel-annih}); the positive sign comes from absorbing the resolvent-derivative minus sign in $\D_r G^{(\eps)} = -(\gamma_t+\eps I)^{-1}(\D_r\gamma_t)(\gamma_t+\eps I)^{-1}h$ via the explicit minus in the definition~\eqref{eq:C-eps-positive-def}.

Passing to the limit $\eps \downarrow 0$ in \eqref{eq:sub-eps} using the three convergences just established,
\begin{equation}\label{eq:proof-sub}
  \delta_U(v_h) \;=\; \bigl[\delta_U(w(\cdot)(z))\bigr]_{z = \tilde{h}} \;+\; \int_0^t \Tr_U\bigl[\Phi_r^{*}\;\gamma_t^{\dagger}\,(\D_r\gamma_t)\,\gamma_t^{\dagger}\,h\bigr]\dr,
\end{equation}
which is precisely \eqref{eq:skorokhod-decomp}, $\delta_U(v_h) = M_h + C_h$.

The correction integral is unfolded by Proposition \ref{prop:D-gamma}, which writes $\D_r\gamma_t$ as the single symmetric expression \eqref{eq:D-gamma-formula} in the second Malliavin fibre. The two mechanisms it contains are read off the data of $\calW$ rather than off a splitting of the $s$-integral. The second derivatives $\calA''_{uu}$ and $\calB''_{uu}$ enter through the forcings \eqref{eq:bold-forcings} of equation \eqref{eq:W-equation}, capturing the response of the fibre to the noise perturbation at time $r$, and the first derivative $\calB'_u$ enters through the initial datum \eqref{eq:W-initial}, capturing the direct response of the noise coefficient to the solution perturbation.
For equations with $\calA''_{uu} = 0$ and $\calB'_u = 0$ (linear drift, state-independent diffusion), both types of correction vanish, and $\delta_U(v_h)$ reduces to the main term, an It\^o integral, recovering the Gaussian result of \cite{mirafzali2025infinite}.

We verify that no abstract Malliavin derivatives remain. The first is $\D_r X(t) = \Phi_{t,r} = Y(t,r)\calB(r,X(r))$ (Proposition \ref{prop:D-Xt}); the second is $\D_r\Phi_{t,s} = \calW(t;s,r)$, which solves \eqref{eq:W-equation} with the datum \eqref{eq:W-initial}, both built from $\calA'_u$, $\calA''_{uu}$, $\calB'_u$, $\calB''_{uu}$ and the fibres alone (Lemma \ref{lemma:DPhi}); $\D_r\gamma_t$ is the symmetric quadratic expression \eqref{eq:D-gamma-formula} in the two fibres (Proposition \ref{prop:D-gamma}); and the effective derivative of $\gamma_t^{\dagger}h$ is $-\gamma_t^{\dagger}(\D_r\gamma_t)\gamma_t^{\dagger}h$ modulo $\ker(\gamma_t)$ (Lemma \ref{lemma:D-gamma-inv-h}), which is the form in which this object enters the trace-class Bismut formula. Every quantity is computable from the solution, the two fibres, the operator derivatives, and the Moore--Penrose pseudoinverse of the Malliavin covariance.

This completes the proof of Theorem \ref{thm:main}.
\end{proof}

\section{The theorem in use}
\label{sec:applications}

A theorem of this kind is to be judged by the cases it settles. We first verify that Theorem~\ref{thm:main} recovers the known results in the limiting cases (finite-dimensional SDEs, Gaussian measures, linear-additive SPDEs), then verify the coefficient-side and analytic structural hypotheses that are available explicitly for the two principal examples, the stochastic $p$-Laplacian and the two-dimensional Navier--Stokes equation, and finally extend the Bismut formula to a class of singular SPDEs via a scalar reduction.

\subsection{The classical formulae recovered}\label{subsec:reduction}

In each of three classical cases the Bismut formula collapses to a known answer. The simplest is the finite-dimensional reduction.

\subsubsection*{Finite-dimensional nonlinear SDEs}

\begin{corollary}\label{cor:finite-dim}
When $V = H = \R^m$ and $U = \R^d$, the Bismut formula \eqref{eq:score-main}--\eqref{eq:skorokhod-decomp} reduces to Theorem 2.1 of \cite{mirafzali2025malliavin}.
\end{corollary}

\begin{proof}
In finite dimensions the Gelfand triple collapses to $V = H = V^* = \R^m$, and $U = \R^d$. We set $\calA(t, u) = -b(t, u) \in \R^m$ (with the sign convention matching \cite{mirafzali2025malliavin}) and $\calB(t, u) = \sigma(t, u) \in \R^{m \times d}$.

The operator $Y(t,r) \in L(\R^m) = \R^{m \times m}$ satisfies
\[
  \mathrm{d}Y(t,r) = \partial_x b(t, X_t)\,Y(t,r)\dt + \sum_{l=1}^d \partial_x \sigma^l(t, X_t)\,Y(t,r)\,\mathrm{d}B_t^l, \qquad Y(r,r) = I_m,
\]
which is equation (2.1) of \cite{mirafzali2025malliavin} (noting $\calA'_u = -\partial_x b$). When $Y(r,0)$ is invertible, the flow factorisation $Y(t,r) = Y(t,0)\,Y(r,0)^{-1}$ gives $Y(t,r) = Y_T Y_r^{-1}$ in the notation of \cite{mirafzali2025malliavin}.

The operator $\gamma_t \in \R^{m \times m}$ becomes (writing $T = t$ for the terminal time, as in \cite{mirafzali2025malliavin})
\[
  \gamma_t = \int_0^t Y(t,r)\,\sigma(r,X_r)\,\sigma(r,X_r)^\top\,Y(t,r)^\top\dr = \int_0^t Y_T Y_r^{-1}\sigma_r\sigma_r^\top(Y_T Y_r^{-1})^\top\dr = \gamma_{X_T},
\]
matching definition (2.3) of \cite{mirafzali2025malliavin}.

For the canonical direction $h = e_k$ ($k$-th standard basis vector of $\R^m$), the covering field is
\[
  v_{e_k}(r) = \sigma(r,X_r)^\top Y(t,r)^\top\gamma_t^{\dagger}e_k \in \R^d,
\]
and the $j$-th component is $[v_{e_k}(r)]_j = \sum_{i=1}^m \sigma^{ij}(r,X_r)[Y(t,r)^\top\gamma_t^{\dagger}e_k]_i$, matching the covering vector field $v_r^{(k)}$ of \cite{mirafzali2025malliavin}.

The infinite-dimensional $\Tr_U$ reduces to a sum over the $d$ standard basis vectors $\{f_1, \ldots, f_d\}$ of $U = \R^d$. In the correction term of \eqref{eq:skorokhod-decomp}, with $M \in L(U, H) = \R^{m \times d}$ symbolically standing for $\gamma_t^{\dagger}(\D_r\gamma_t)\gamma_t^{\dagger}h$ in the effective sense of Lemma~\ref{lemma:D-gamma-inv-h}, we have
\[
  \begin{aligned}
  \Tr_U\bigl[\calB(r,X(r))^*\,Y(t,r)^* \cdot M\bigr]
  &\;=\; \sum_{i=1}^d \bigl\langle \sigma(r,X_r)^\top Y(t,r)^\top (Mf_i),\,f_i\bigr\rangle_{\R^d} \\
  &\;=\; \sum_{i=1}^d \bigl[\sigma_r^\top Y(t,r)^\top M\bigr]_{ii},
  \end{aligned}
\]
i.e., the trace of the $d \times d$ matrix $\sigma_r^\top Y(t,r)^\top M$. By the cyclic invariance of the trace, this equals $\sum_{j=1}^m [M\sigma_r^\top Y(t,r)^\top]_{jj}$, the trace of the corresponding $m \times m$ matrix; both expressions are the same scalar and reproduce the trace appearing in \cite{mirafzali2025malliavin} after the appropriate identifications.

In the notation of Theorem 2.1 of \cite{mirafzali2025malliavin}:
\begin{enumerate}[(a)]
  \item The ``main term'' $[\delta_U(w(\cdot)(z))]_{z = \gamma_t^{\dagger}h}$ corresponds to the Skorokhod integral in equation (2.7) of \cite{mirafzali2025malliavin}, evaluated at the random point $z = \gamma_{X_T}^{-1}e_k$.
  \item The correction term in \eqref{eq:skorokhod-decomp}, when expanded via $\D_r(\gamma_t^{\dagger}h) = -\gamma_t^{\dagger}(\D_r\gamma_t)\gamma_t^{\dagger}h$ and the symmetric formula \eqref{eq:D-gamma-formula} for $\D_r\gamma_t$ (Proposition~\ref{prop:D-gamma}), yields four distinct integrand types, according to which of the two data of $\calW$ each comes from. In the notation of \cite{mirafzali2025malliavin}, the forcing of \eqref{eq:W-equation} carries $\partial_x b$ and $\partial_x^2\sigma$ and produces, after symmetrisation, the integrands $I_1^{p,q}$ and $I_2^{p,q}$ integrated over the full range $s \in [0,t]$; the initial datum \eqref{eq:W-initial} carries $\partial_x\sigma$ and produces, after symmetrisation, the cross-integrands $I_3^{p,q}$ and $I_4^{p,q}$. In finite dimensions the fibre is $\Phi_{t,s} = Y_TY_s^{-1}\sigma_s$ and $\calW(t;s,r) = D_r\Phi_{t,s}$ is the derivative appearing there.
\end{enumerate}
Thus the formula \eqref{eq:skorokhod-decomp} with the score \eqref{eq:score-main} reduces precisely to Theorem 2.1 of \cite{mirafzali2025malliavin}.
\end{proof}

\subsubsection*{Linear equations with additive noise}

\begin{corollary}\label{cor:linear-infinite}
When $\calA(t, u) = -Au$ is linear with $A$ generating an analytic semigroup $S(t)$, and $\calB(t, u) = Q^{1/2}$ is state-independent, then for every deterministic direction $h \in \Ran(\gamma_t)$ (the covariance $\gamma_t$ being deterministic here, so that clause~\textup{(ii)} of Assumption~\ref{ass:nondeg} is automatic) the Bismut formula reduces to
\[
  \beta_h(u) = -\ip{u - S(t)x}{\gamma_t^{\dagger}h}_H,
\]
which is Theorem 4 of \cite{mirafzali2025infinite}.
\end{corollary}

\begin{proof}
With $\calA$ linear and $\calB$ constant, we have $\calA'_u = -A$ (constant in $u$), $\calA''_{uu} = 0$, $\calB'_u = 0$, and $\calB''_{uu} = 0$. The first variation equation \eqref{eq:first-var} becomes $\mathrm{d}Y - AY\dt = 0$ with $Y(r,r) = I$, giving $Y(t,r) = S(t-r)$ (deterministic). The second Malliavin fibre vanishes, $\calW \equiv 0$, since equation \eqref{eq:W-equation} then has zero forcing and zero initial datum. Consequently $\D_r\Phi_{t,s} = 0$, $\D_r\gamma_t = 0$, and all correction terms vanish.

The Malliavin covariance is deterministic:
\[
  \gamma_t = \int_0^t S(t-r)\,Q^{1/2}(Q^{1/2})^*\,S(t-r)^*\dr = \int_0^t S(s)\,Q\,S(s)^*\ds,
\]
which is the Malliavin covariance operator of \cite{mirafzali2025infinite}, equation (12).

The covering field $v_h(r) = (Q^{1/2})^*S(t-r)^*\gamma_t^{\dagger}h$ is deterministic, so $\delta_U(v_h)$ reduces to the It\^o integral:
\[
  \delta_U(v_h) = \int_0^t \ip{(Q^{1/2})^*S(t-r)^*\gamma_t^{\dagger}h}{\mathrm{d}W(r)}_U.
\]

From the mild solution $X(t) = S(t)x + \int_0^t S(t-s)Q^{1/2}\dW(s)$, the stochastic integral is
\[
  \int_0^t S(t-s)Q^{1/2}\dW(s) = X(t) - S(t)x =: z.
\]

Therefore
\[
  \delta_U(v_h) = \ip{z}{\gamma_t^{\dagger}h}_H = \ip{X(t) - S(t)x}{\gamma_t^{\dagger}h}_H,
\]
and since this is already $\sigma(X(t))$-measurable, the conditional expectation is trivial:
\[
  \beta_h(X(t)) = -\ip{X(t) - S(t)x}{\gamma_t^{\dagger}h}_H. \qedhere
\]
\end{proof}

\subsubsection*{Semilinear equations with additive noise}

\begin{corollary}\label{cor:semilinear-additive}
Consider semilinear SPDEs with additive noise,
\[
  \mathrm{d}X(t) = [AX(t) + F(X(t))]\dt + Q^{1/2}\dW(t),
\]
where $A$ generates an analytic semigroup on $H$, $Q^{1/2} \in L_2(U,H)$, and $F: V \to V^*$ is twice continuously Fr\'echet differentiable with $F'_u, F''_{uu}$ satisfying the growth and integrability conditions of Assumptions \textup{(D1)--(D3)} (so that $\calA(t,u) := -[Au + F(u)]$ fits the abstract framework of Assumptions \ref{ass:LR}--\ref{ass:diff2}). Assume that the hypotheses of Theorem~\ref{thm:main}, Part~II, hold for the given deterministic direction $h \in H$, namely the structural Assumptions~\ref{ass:SC} and~\ref{ass:SC-raised}, all three clauses of Assumption~\ref{ass:nondeg} (the range and moment conditions \textup{(i)}--\textup{(ii)} at the chosen H\"older-triangle pair $(q, p^*)$, and the Cameron--Martin compatibility \textup{(iii)}), and the fixed-Tikhonov substitution regularity~\ref{hyp:iii-F-prime}.
Then $\calB'_u = 0$, $\calB''_{uu} = 0$, the second Malliavin fibre \eqref{eq:W-equation} has no stochastic integral term, the correction terms involve only $F''_{uu}$ (not $\calB''_{uu}$), and the formula reads
\[
\begin{aligned}
  \delta_U(v_h) &= \bigl[\delta_U\bigl(r \mapsto \Phi_r^{*}z\bigr)\bigr]_{z = \gamma_t^{\dagger}h}\\
  &\quad + \lim_{\eps \downarrow 0}\int_0^t \Tr_U\bigl[\Phi_r^{*}\,(\gamma_t+\eps I)^{-1}\,[\D_r\gamma_t]\,(\gamma_t+\eps I)^{-1}h\bigr]\dr,
\end{aligned}
\]
with $\Phi_r = \D_r X(t)$, the first line being the frozen-$z$ divergence evaluated after the substitution $z = \gamma_t^{\dagger}h$ and not the composed divergence $\delta_U(r \mapsto \Phi_r^{*}\gamma_t^{\dagger}h)$, which already absorbs the trace correction (Lemma~\ref{lemma:HS-kernel}(iv)).
\end{corollary}

Two features of this specialisation are worth separating from the statement. The frozen-$z$ divergence is anticipating unless $F = 0$, in which case it reduces to an It\^o integral; and if the random evolution family generated by $A + F'(X(\cdot))$ admits an almost surely bounded realisation $Y(t,r) \in L(H)$, which is a structural input rather than a consequence of analyticity of $S(t)$, then $\Phi_r^{*} = (Q^{1/2})^{*}Y(t,r)^{*}$ and the frozen-$z$ field may be written in that factored form, the factorisation being otherwise a shorthand for the Hilbert--Schmidt kernel $\Phi_r$.

The Cameron--Martin hypothesis carries conditional content here. Clause~\textup{(iii)} of Assumption~\ref{ass:nondeg} is a joint convergence-and-integrability condition on the quadruple $(\gamma_t, \D_r\gamma_t, \Phi_r, h)$, and it is not an automatic consequence of $A$ being analytic. Even for state-independent $\calB = Q^{1/2}$, the first variation $Y(t,r)$ is random whenever the drift is nonlinear, so $\D_r\gamma_t \neq 0$ in general and the trace integrand $\mathcal{I}^{(\eps)}_r$ does not vanish identically. It can be verified through parabolic maximal-regularity bounds on the linearised evolution family generated by $A + F'(X(\cdot))$ at the operating point $(q, p^*)$ on the H\"older triangle, together with input on the spectral structure of $\gamma_t$ relative to $h$. The fully linear case $F = 0$ with $\calB$ deterministic is the one subcase in which it is automatic, with $\mathcal{I}^{(0)} \equiv 0$ and $\Psi \equiv 0$ (cf.\ Remark~\ref{rem:CM-compat-role}).

\begin{proof}
With $\calB(t,u) = Q^{1/2}$ state-independent, we have $\calB'_u = 0$ and $\calB''_{uu} = 0$. The first variation equation \eqref{eq:first-var} becomes
\begin{equation}\label{eq:Y-semilinear}
  \mathrm{d}Y(t,r) - [A + F'(X(t))]Y(t,r)\dt = 0, \qquad Y(r,r) = I,
\end{equation}
which is a random linear ODE (random coefficients $A + F'(X(t))$ depending on the SPDE solution, but with no $\dW$ term). Since the coefficients depend on $X(s)$ for $s \in [r,t]$, the propagator $Y(t,r)$ is $\F_t$-measurable but not $\F_r$-adapted; the covering field $v_h(r) = \Phi_r^{*}\gamma_t^{\dagger}h$ is therefore anticipating, and the linear-field substitution formula (Theorem \ref{thm:linear-sub-D12}) is applied as in the general case.

Since $\calB'_u = 0$ and $\calB''_{uu} = 0$, the second Malliavin fibre \eqref{eq:W-equation} reduces to
\begin{equation}\label{eq:W-semilinear}
  \mathrm{d}\calW^{a,u}(\tau;s,r) - [A + F'(X(\tau))]\,\calW^{a,u}(\tau;s,r)\,\mathrm{d}\tau = F''(X(\tau))\bigl(\Phi_{\tau,r}a,\,\Phi_{\tau,s}u\bigr)\,\mathrm{d}\tau,
\end{equation}
which is also a random ODE, with random coefficients and a forcing that is bilinear in the two variational fibres $\Phi_{\tau,r} = Y(\tau,r)Q^{1/2}$ and $\Phi_{\tau,s} = Y(\tau,s)Q^{1/2}$. Both fibres carry the smoothing of $Q$, and by Theorem~\ref{thm:Y-wellposed}\textup{(III)} their fibre energies are finite, so the forcing is square-integrable in time without any bound on $\norm{Y(\tau,r)}_{L(H,V)}$.

For the Bismut formula, by Theorem \ref{thm:main}, the substitution gives:
\begin{align*}
  \delta_U(v_h) &= \bigl[\delta_U(w(\cdot)(z))\bigr]_{z = \gamma_t^{\dagger}h} - \int_0^t \Tr_U[\Phi_r^{*} \cdot \D_r(\gamma_t^{\dagger}h)]\dr.
\end{align*}
The first term is the frozen-$z$ divergence of $w_r(z) = \Phi_r^{*}z$ evaluated at the random point $z = \gamma_t^{\dagger}h$, giving the first line of the claimed formula. (Since $Y(t,r)$ is not $\F_r$-adapted, this is a genuine Skorokhod integral; it reduces to an It\^o integral only in the linear case $F = 0$, where $Y(t,r) = S(t-r)$ is deterministic.)

The second term (the correction) is the Tikhonov limit~\eqref{eq:correction-tikhonov}, written symbolically through $\D_r(\gamma_t^{\dagger}h) = -\gamma_t^{\dagger}(\D_r\gamma_t)\gamma_t^{\dagger}h$. Since $\calB'_u = 0$, the initial datum \eqref{eq:W-initial} of the second Malliavin fibre vanishes and $\calW$ is driven by the $F''$-forcing alone. Substituting \eqref{eq:D-gamma-formula} yields the second line.
\end{proof}

\subsubsection*{Semilinear equations with state-independent diffusion}

\begin{corollary}\label{cor:state-indep-diffusion}
For SPDEs with state-independent diffusion $\calB(t, u) = \calB(t)$ (independent of $u$) but nonlinear drift $\calA(t,u)$, we have $\calB'_u = 0$ and $\calB''_{uu} = 0$. The second Malliavin fibre \eqref{eq:W-equation} simplifies to
\begin{equation}\label{eq:W-state-indep}
\begin{aligned}
  \mathrm{d}\calW^{a,u}(\tau;s,r) + \calA'_u(\tau,X(\tau))\,\calW^{a,u}(\tau;s,r)\,\mathrm{d}\tau
  &= -\calA''_{uu}(\tau,X(\tau))\bigl(\Phi_{\tau,r}a,\,\Phi_{\tau,s}u\bigr)\,\mathrm{d}\tau, \\
  \calW^{a,u}(\max(s,r);s,r) &= 0,
\end{aligned}
\end{equation}
a random PDE with no stochastic integral term but with random coefficients depending on the SPDE solution, and with a forcing bilinear in the two variational fibres. The Bismut formula \eqref{eq:skorokhod-decomp} then involves only the fibres $\Phi$ and $\calW$ and the operator derivatives $\calA'_u$ and $\calA''_{uu}$.
\end{corollary}

\begin{proof}
Since $\calB$ is independent of $u$, we have $\calB'_u = 0$ and $\calB''_{uu} = 0$. The first variation equation \eqref{eq:first-var} retains its stochastic integral term only if $\calB'_u \neq 0$; since $\calB'_u = 0$, it reduces to the random ODE $\mathrm{d}Y(t,r) + \calA'_u(t,X(t))Y(t,r)\dt = 0$, $Y(r,r) = I$. The second Malliavin fibre \eqref{eq:W-equation} similarly loses both its $\dW$ terms and, by \eqref{eq:W-initial}, its initial datum, yielding \eqref{eq:W-state-indep}.

The covering field is $v_h(r) = \Phi_r^{*}\gamma_t^{\dagger}h$, with the factored form $\calB(r)^*Y(t,r)^*\gamma_t^{\dagger}h$ available whenever an a.s.\ bounded realisation of $Y(t,r)$ is separately assumed. The first variation $Y(t,r)$ here solves a random ODE (no $\dW$ term, but with random coefficients $\calA'_u(\tau, X(\tau))$ depending on the SPDE solution), so $Y(t,r)$ is not $\F_r$-adapted; combined with the $\F_t$-measurability of $\gamma_t^{\dagger}h$, the covering field is anticipating, and the Hilbert-space linear substitution identity of Theorem~\ref{thm:linear-sub-D12} applies exactly as in the general nonlinear case.

Since $\calB'_u = 0$, the initial datum \eqref{eq:W-initial} vanishes for every ordering of $r$ and $s$, so $\calW$ is driven by the $\calA''_{uu}$-forcing alone and \eqref{eq:D-gamma-formula} carries a single mechanism. The correction terms in \eqref{eq:skorokhod-decomp} simplify correspondingly, and the formula coincides with the state-independent-diffusion case of \cite{mirafzali2025malliavin} (the specialisation of the nonlinear-drift result with $\calB = \calB(t)$) when restricted to finite dimensions.
\end{proof}

On integrability and ellipticity, Assumption~\ref{ass:nondeg} is strictly weaker than ``$\gamma_t$ injective a.s.'' It requires only that $h \in \Ran(\gamma_t)$ with $\E[\norm{\gamma_t^\dagger h}^q] < \infty$. For additive trace-class noise with $\ker(Q) = \{0\}$, injectivity of $Q^{1/2}$ gives $\ker(\gamma_t) = \{0\}$ a.s., but the pseudoinverse moment bound and Cameron--Martin compatibility require the intrinsic observability (Proposition~\ref{prop:nondeg-additive}) or the Hairer--Mattingly H\"ormander framework (Proposition~\ref{prop:nondeg-NS}).

\subsection{Two equations of the variational class}\label{subsec:verification}

The classical formulae are recovered because the hypotheses hold there for trivial reasons. The test of the theorem is whether they can be verified for equations that were not already understood, and we take the two that motivated the work, a quasilinear diffusion whose linearisation degenerates and a hydrodynamic equation whose nonlinearity is not monotone. The two illustrate the two alternatives of~(SC1) respectively.

\subsubsection*{The stochastic \texorpdfstring{$p$}{p}-Laplacian}

Consider the stochastic $p$-Laplacian equation on a bounded domain $\mathcal{O} \subset \R^d$ with smooth boundary:
\begin{equation}\label{eq:p-Lap-full}
\begin{aligned}
  \mathrm{d}X(t) - \mathrm{div}\bigl(\abs{\nabla X(t)}^{p-2}\nabla X(t)\bigr)\dt
  &= \calB(t, X(t))\dW(t), \\
  X(0) = x, \quad X|_{\partial\mathcal{O}} &= 0,
\end{aligned}
\end{equation}
with $p \geq 2$. Set $V = W_0^{1,p}(\mathcal{O})$, $H = L^2(\mathcal{O})$, $V^* = W^{-1,p'}(\mathcal{O})$ where $p' = p/(p-1)$.

\begin{proposition}[Fully local monotonicity for the $p$-Laplacian]\label{prop:p-Lap-LR}
On the Gelfand triple $W_0^{1,p}(\mathcal{O}) \hookrightarrow L^2(\mathcal{O}) \hookrightarrow W^{-1,p'}(\mathcal{O})$, the drift $\calA(u) = -\mathrm{div}(\abs{\nabla u}^{p-2}\nabla u)$ satisfies the drift-side conditions of Assumption~\ref{ass:LR}. If, in addition, $\calB$ satisfies the $H$-Lipschitz/continuity and linear-growth conditions used in the proof below, in particular for state-independent Hilbert--Schmidt noise, or for $\calB(t,u) = \sigma(u)Q^{1/2}$ with $\sigma$ globally Lipschitz, then the pair $(\calA, \calB)$ satisfies Assumption~\ref{ass:LR}.
\end{proposition}

\begin{proof}

The embedding $V = W_0^{1,p}(\mathcal{O}) \hookrightarrow L^2(\mathcal{O}) = H$ is compact by the Rellich--Kondrachov theorem (for every dimension $d$, since $p \geq 2 > 2d/(d+2)$).

For $u, v, w \in V$, the map $\lambda \mapsto \dual{\calA(u + \lambda v)}{w} = \int_{\mathcal{O}} \abs{\nabla(u + \lambda v)}^{p-2}\nabla(u+\lambda v) \cdot \nabla w\,\mathrm{d}x$ is continuous in $\lambda$ by the dominated convergence theorem (using $\abs{\nabla(u+\lambda v)}^{p-2} \leq C(\abs{\nabla u}^{p-2} + \abs{\lambda}^{p-2}\abs{\nabla v}^{p-2})$).

The classical vector inequality (cf.\ \cite{lindqvist2006plaplace}, Section~10, inequality~(I)): for all $\xi, \eta \in \R^d$ and $p \geq 2$,
\begin{equation}\label{eq:p-Lap-mono}
  (\abs{\xi}^{p-2}\xi - \abs{\eta}^{p-2}\eta) \cdot (\xi - \eta) \geq c_p\,\abs{\xi - \eta}^p,
\end{equation}
where $c_p > 0$ depends only on $p$ (one may take $c_p = 2^{2-p}$). Hence, by~\eqref{eq:p-Lap-mono},
\begin{align*}
  -2\dual{\calA(u) - \calA(v)}{u-v}
  &= -2\int_{\mathcal{O}}(\abs{\nabla u}^{p-2}\nabla u - \abs{\nabla v}^{p-2}\nabla v)\cdot\nabla(u-v)\,\mathrm{d}x \\
  &\leq -2c_p\norm{u-v}_V^p \;\leq\; 0.
\end{align*}
Combining this with the Lipschitz bound $\norm{\calB(t,u) - \calB(t,v)}_{L_2}^2 \leq C\norm{u-v}_H^2$ yields the joint local monotonicity~(H2):
\[
  -2\dual{\calA(u)-\calA(v)}{u-v} + \norm{\calB(t,u)-\calB(t,v)}_{L_2}^2 \;\leq\; C\norm{u-v}_H^2.
\]

Directly,
\[
  \dual{\calA(u)}{u} \;=\; \int_\mathcal{O} \abs{\nabla u}^{p-2}\,\nabla u \cdot \nabla u\,\mathrm{d}x
  \;=\; \int_\mathcal{O} \abs{\nabla u}^p\,\mathrm{d}x \;=\; \norm{u}_V^p,
\]
so $-2\dual{\calA(u)}{u} = -2\norm{u}_V^p \leq -\alpha\norm{u}_V^p$ with $\alpha = 2$.

By H\"older's inequality:
\[
  \abs{\dual{\calA(u)}{v}} = \Bigl\lvert\int_\mathcal{O} \abs{\nabla u}^{p-2}\nabla u \cdot \nabla v\,\mathrm{d}x\Bigr\rvert \leq \norm{\abs{\nabla u}^{p-1}}_{L^{p'}}\,\norm{\nabla v}_{L^p} = \norm{u}_V^{p-1}\norm{v}_V,
\]
giving $\norm{\calA(u)}_{V^*} \leq \norm{u}_V^{p-1} \leq C(1 + \norm{u}_V^{p-1})$.

The $H$-norm continuity~\eqref{eq:B-H-cont} and the growth bound~\eqref{eq:B-H-growth} are conditions on $\calB$, not on $\calA$. For any state-independent noise $\calB(t,u) = \calB(t)$, both hold trivially with $g(t) = \norm{\calB(t)}_{L_2(U,H)}^2$. For noise of the form $\calB(t,u) = \sigma(u)\,Q^{1/2}$ with $\sigma: H \to \R$ globally Lipschitz, $\norm{\calB(t,u) - \calB(t,v)}_{L_2} \leq C\norm{u-v}_H$ gives~\eqref{eq:B-H-cont}, and $\norm{\calB(t,u)}_{L_2}^2 \leq C(1+\norm{u}_H^2)$ gives~\eqref{eq:B-H-growth}.
\end{proof}

\begin{proposition}[The $p$-Laplacian satisfies \textup{(D1)} and \textup{(D3)}]\label{prop:p-Lap-diff}
For $p > 3$ strictly, the operator $\calA(u) = -\mathrm{div}(\abs{\nabla u}^{p-2}\nabla u)$ satisfies conditions (D1) and (D3), the drift side of Assumptions~\ref{ass:diff1} and~\ref{ass:diff2}. The boundary case $p = 3$ falls into the regularised regime alongside $2 < p < 3$ (Proposition~\ref{prop:singular-p-Lap}); see Remark~\ref{rem:p3-borderline} below for the precise reason. The range $3 < p < 4$, where (D1) and (D3) do hold but the domination estimate of Proposition~\ref{prop:p-Lap-Lambda} does not close, is routed through the same regularisation. Conditions (D2)/(D2$'$) and (D4) on $\calB$ are independent of the choice of drift and must be verified for the specific noise coefficient under consideration; they hold trivially for any state-independent $\calB$ uniformly Hilbert--Schmidt-bounded (with $\calB'_u = \calB''_{uu} = 0$, so the unified growth bound~\eqref{eq:B-growth-m} holds with $\kappa_{m_0} = 0$ at any finite $m_0 \geq 2$, and equally trivially for linear multiplicative noise $\calB(t,u) = Lu$ with $\kappa_{m_0} = m_0$).
\end{proposition}

\begin{proposition}[The $p$-Laplacian satisfies the domination clause of \textup{(SC1)(a)} for $p \geq 4$]\label{prop:p-Lap-Lambda}
Let $\calA$ be the $p$-Laplacian $\calA(u) = -\mathrm{div}(\abs{\nabla u}^{p-2}\nabla u)$ with exponent $p \geq 4$, and suppose that the diffusion coefficient does not depend on the state. Then the dissipativity and the second-derivative domination required by alternative \textup{(a)} of~\textup{(SC1)} hold with slack $\delta = 2$, the remaining non-autonomous form-domain and well-posedness datum of~\textup{(SC1)(a)} being a separate structural input, the second derivative of the drift being dominated by
\[
  \Lambda_\sigma(\varphi, \psi) \;=\; C_p\Bigl(\int \abs{\nabla X(\sigma)}^{p-4}\,\abs{\nabla\varphi}^2\,\abs{\nabla\psi}^2\,\mathrm{d}x\Bigr)^{1/2} \;\leq\; C_p\,\norm{X(\sigma)}_V^{(p-4)/2}\,\norm{\varphi}_V\,\norm{\psi}_V.
\]
\end{proposition}

\begin{proof}
The dissipation is
\begin{multline*}
  D_\sigma(w) \;=\; \int\bigl[(p-2)\abs{\nabla X}^{p-4}(\nabla X\cdot\nabla w)^2 + \abs{\nabla X}^{p-2}\abs{\nabla w}^2\bigr]\,\mathrm{d}x \\
  \;\geq\; \int \abs{\nabla X}^{p-2}\abs{\nabla w}^2\,\mathrm{d}x \;\geq\; 0,
\end{multline*}
and the slack $\delta = 2$ is automatic since $\calB'_u \equiv 0$. Every term of $\dual{\calA''_{uu}(X)(\varphi,\psi)}{w}$ carries the weight $\abs{\nabla X}^{p-3}$ against one gradient of each argument, so
\[
  \abs{\dual{\calA''_{uu}(X)(\varphi,\psi)}{w}} \;\leq\; C_p\int\abs{\nabla X}^{p-3}\abs{\nabla\varphi}\,\abs{\nabla\psi}\,\abs{\nabla w}\,\mathrm{d}x ;
\] the weighted Cauchy--Schwarz split $\abs{\nabla X}^{p-3} = \abs{\nabla X}^{(p-2)/2}\cdot\abs{\nabla X}^{(p-4)/2}$ bounds this by
\begin{align*}
  C_p\Bigl(\int\abs{\nabla X}^{p-2}\abs{\nabla w}^2\Bigr)^{1/2}\Bigl(\int\abs{\nabla X}^{p-4}\abs{\nabla\varphi}^2\abs{\nabla\psi}^2\Bigr)^{1/2}
  \;\leq\; \Lambda_\sigma(\varphi, \psi)\,D_\sigma(w)^{1/2}.
\end{align*}
The displayed bound on $\Lambda$ is H\"older with the exponent triple $(\tfrac{p}{p-4}, \tfrac{p}{2}, \tfrac{p}{2})$ on $L^p$-gradients, whose reciprocals sum to one. The joint square-summability~\eqref{eq:Lambda-square-summable} is immediate from the displayed form, since $\Lambda_\sigma^2$ is an integral of the weight against $\abs{\nabla\varphi}^2\abs{\nabla\psi}^2$ and summing over two orthonormal systems replaces those two factors by the gradient densities $\sum_j\abs{\nabla\Theta f_j}^2$ and $\sum_k\abs{\nabla\Xi f_k}^2$, both basis-free.
\end{proof}

\begin{remark}[The forcing clause for the $p$-Laplacian]\label{rem:p-Lap-forcing}
With $\Lambda$ as in Proposition~\ref{prop:p-Lap-Lambda}, the drift part of the total forcing~\eqref{eq:total-forcing} is, by the first line of~\eqref{eq:G1-def} and the pointwise form of the double square sum~\eqref{eq:Lambda-square-summable},
\begin{equation}\label{eq:p-Lap-G1}
  G_1(\tau; r,s) \;=\; C_p^2\int_{\mathcal{O}} \abs{\nabla X(\tau)}^{p-4}\,\bigl\lvert\nabla\Phi_{\tau,r}\bigr\rvert^2_{\mathrm{HS}}\,\bigl\lvert\nabla\Phi_{\tau,s}\bigr\rvert^2_{\mathrm{HS}}\,\mathrm{d}x, \qquad \bigl\lvert\nabla\Theta\bigr\rvert^2_{\mathrm{HS}} := \sum\nolimits_j\abs{\nabla \Theta f_j}^2,
\end{equation}
a single weighted integral of the two fibre gradient densities, and $\mathfrak{J}_{r,s}$ and $\mathbf{F}_2$ vanish for state-independent noise. The clause~(SC1$'$) is the statement that~\eqref{eq:p-Lap-G1}, integrated in time, carries the moment~\eqref{eq:forcing-integrability}, and this is the datum of the degenerate mechanism, in the same sense as the form domain of~(SC1)(a). It is a weighted quantity of the same kind as the tangent energy $\int D_\tau(\Phi_{\tau,r}f_j)\,\mathrm{d}\tau$ that Theorem~\ref{thm:Y-wellposed}\textup{(III)} produces, but at the weight $\abs{\nabla X}^{p-4}$ rather than $\abs{\nabla X}^{p-2}$, and the two weights differ by exactly the degeneracy of the linearisation.

Wherever unweighted fibre regularity is available, the clause reduces to moments already in the assumption set. Bounding the two densities in $L^{p/2}$ and the weight in $L^{p/(p-4)}$ by H\"older at the exponent triple $(p/(p-4), p/2, p/2)$, whose reciprocals sum to one for $p \geq 4$, and then using Minkowski in $L^{p/2}$ to pass from the density to the fibre norm, $\bigl\lVert\,\lvert\nabla\Theta\rvert^2_{\mathrm{HS}}\bigr\rVert_{L^{p/2}} \leq \norm{\Theta}^2_{L_2(U,V)}$, gives
\[
  G_1(\tau;r,s) \;\leq\; C_p^2\,\norm{X(\tau)}_V^{p-4}\,\norm{\Phi_{\tau,r}}^2_{L_2(U,V)}\,\norm{\Phi_{\tau,s}}^2_{L_2(U,V)} ,
\]
and a H\"older triple in $(\tau,\omega)$ at the conjugate exponents $(3,3,3)$ reduces~\eqref{eq:forcing-integrability} to
\[
  \E\!\int_0^T\!\norm{X}_V^{3(p-4)(1+\eps_0)}\,\mathrm{d}\tau,
  \qquad
  \E\Bigl[\Bigl(\int_0^T\!\norm{\Phi_{\tau,r}}^{2}_{L_2(U,V)}\,\mathrm{d}\tau\Bigr)^{3(1+\eps_0)}\Bigr],
\]
\[
  \E\Bigl[\Bigl(\int_0^T\!\norm{\Phi_{\tau,s}}^{2}_{L_2(U,V)}\,\mathrm{d}\tau\Bigr)^{3(1+\eps_0)}\Bigr],
\]
the first supplied by (SC4.1)$_{q^*}$ at the matching exponent. The last two are the unweighted counterpart of the fibre energy, and they follow from~\eqref{eq:Phi-energy} on any event on which $\abs{\nabla X}$ is bounded below, since there the two weights are comparable; this is the regime in which the linearisation is uniformly parabolic, and it is the same regime in which the regularised formulation of Proposition~\ref{prop:singular-p-Lap} converges. Separate moments of the two fibre factors do not suffice in any case, since a product of three factors each controlled at the exponent of the product is not controlled, which is why~\eqref{eq:forcing-integrability} is carried as the single clause~(SC1$'$) and verified at raised exponents rather than assembled from base-exponent bounds.
\end{remark}

\begin{proof}[Proof of Proposition~\ref{prop:p-Lap-diff}]

The first Fr\'echet derivative (in the strong-operator sense; the formal expression is rigorous on a dense set and extends by continuity to $V$ for $p > 3$) is (cf.\ Example \ref{ex:p-Laplace-second-order}):
\[
  \calA'_u(u)v = -\mathrm{div}\bigl((p-2)\abs{\nabla u}^{p-4}(\nabla u \cdot \nabla v)\nabla u + \abs{\nabla u}^{p-2}\nabla v\bigr).
\]
For the growth bound:
\begin{align*}
  \abs{\dual{\calA'_u(u)v}{w}} &\leq (p-1)\int_\mathcal{O} \abs{\nabla u}^{p-2}\abs{\nabla v}\abs{\nabla w}\,\mathrm{d}x\\
  &\leq (p-1)\norm{\abs{\nabla u}^{p-2}}_{L^{p/(p-2)}}\norm{\nabla v}_{L^p}\norm{\nabla w}_{L^p}\\
  &= (p-1)\norm{u}_V^{p-2}\norm{v}_V\norm{w}_V,
\end{align*}
giving $\norm{\calA'_u(u)}_{L(V,V^*)} \leq (p-1)\norm{u}_V^{p-2} \leq C(1 + \norm{u}_V^{p-2})$.

We must show that for any sequence $u_n \to u_0$ in $V = W_0^{1,p}$ and any fixed test direction $v \in V$,
\begin{equation}\label{eq:strong-op-claim}
  \norm{\calA'_u(u_n)v - \calA'_u(u_0)v}_{V^{*}} \;\longrightarrow\; 0 \qquad (n \to \infty),
\end{equation}
uniformly over $w$ in the unit ball of $V$ (this is the content of strong-operator continuity, strictly stronger than scalar-paired/weak-operator continuity, which would establish only $\dual{\calA'_u(u_n)v}{w} \to \dual{\calA'_u(u_0)v}{w}$ for each individual $w$). The duality $V \hookrightarrow V^{*}$ via the divergence form identifies $\calA'_u(u)v \in V^{*} = W^{-1, p'}(\mathcal{O})$, $p' = p/(p-1)$, with the (negative) divergence of the vector field
\begin{equation}\label{eq:p-Lap-kernel}
  K_n^{v}(x) \;:=\; (p-2)\abs{\nabla u_n(x)}^{p-4}(\nabla u_n(x) \cdot \nabla v(x))\,\nabla u_n(x) \;+\; \abs{\nabla u_n(x)}^{p-2}\,\nabla v(x) \;\in\; \R^d
\end{equation}
(at the operating-point sequence $u_n$); concretely, $\dual{\calA'_u(u_n)v}{w} = \int_{\mathcal{O}} K_n^{v}(x) \cdot \nabla w(x)\,\mathrm{d}x$. By the standard duality $W^{-1,p'} = (W_0^{1,p})^{*}$, with the divergence playing the role of the $V^{*}$-isometry,
\begin{equation}\label{eq:V-star-norm-via-Lp-prime}
  \norm{\calA'_u(u_n)v - \calA'_u(u_0)v}_{V^{*}} \;=\; \sup_{\norm{w}_V \leq 1} \abs{\,\textstyle\int_{\mathcal{O}}\bigl(K_n^{v} - K_0^{v}\bigr)\cdot\nabla w\,\mathrm{d}x\,} \;\leq\; \norm{K_n^{v} - K_0^{v}}_{L^{p'}(\mathcal{O};\R^d)},
\end{equation}
with equality up to the Poincar\'e constant if $W^{-1,p'}$ is normed via $V^{*}$-duality and the $V$-norm taken as $\norm{\nabla\cdot}_{L^p}$. Therefore strong-operator continuity~\eqref{eq:strong-op-claim} reduces to convergence of the kernel sequence in $L^{p'}(\mathcal{O};\R^d)$:
\begin{equation}\label{eq:Kn-Lpprime-conv}
  K_n^{v} \;\longrightarrow\; K_0^{v} \qquad \text{in } L^{p'}(\mathcal{O}; \R^d).
\end{equation}
We verify~\eqref{eq:Kn-Lpprime-conv} by Vitali's convergence theorem in $L^{p'}(\mathcal{O}; \R^d)$. Since $u_n \to u_0$ in $V$, we have $\nabla u_n \to \nabla u_0$ in $L^p(\mathcal{O}; \R^d)$, hence $\nabla u_n(x) \to \nabla u_0(x)$ a.e.\ along a subsequence. Continuity of $\xi \mapsto (p-2)\abs{\xi}^{p-4}(\xi \cdot \nabla v)\,\xi + \abs{\xi}^{p-2}\,\nabla v$ at $p > 3$ (the exponent $p - 4 > -1$ ensures integrability at $\xi = 0$) gives $K_n^{v}(x) \to K_0^{v}(x)$ a.e. For uniform integrability, the pointwise bound
\begin{equation}\label{eq:Kn-pointwise-bound}
  \abs{K_n^{v}(x)} \;\leq\; (p-1)\,\abs{\nabla u_n(x)}^{p-2}\,\abs{\nabla v(x)} \qquad \text{a.e.\ in } x,
\end{equation}
(via $\abs{(p-2)\abs{\xi}^{p-4}(\xi\cdot\eta)\xi + \abs{\xi}^{p-2}\eta} \leq (p-2)\abs{\xi}^{p-2}\abs{\eta} + \abs{\xi}^{p-2}\abs{\eta} = (p-1)\abs{\xi}^{p-2}\abs{\eta}$), gives, for any measurable $E \subset \mathcal{O}$,
\begin{equation}\label{eq:Kn-UI-bound}
  \int_E \abs{K_n^{v}(x)}^{p'}\,\mathrm{d}x \;\leq\; (p-1)^{p'}\int_E \abs{\nabla u_n}^{(p-2)p'}\,\abs{\nabla v}^{p'}\,\mathrm{d}x.
\end{equation}
Apply H\"older's inequality at the conjugate pair $(\alpha, \beta)$ with $\alpha := p/((p-2)p')$ and $\beta := p/p'$ (whose reciprocals satisfy $1/\alpha + 1/\beta = (p-2)p'/p + p'/p = p'(p-1)/p = (p/(p-1))(p-1)/p = 1$, so this is a valid conjugate pair):
\begin{align*}
  \int_E \abs{\nabla u_n}^{(p-2)p'}\,\abs{\nabla v}^{p'}\,\mathrm{d}x
  \;&\leq\; \Bigl(\textstyle\int_E\abs{\nabla u_n}^p\,\mathrm{d}x\Bigr)^{(p-2)p'/p}\,\Bigl(\textstyle\int_E\abs{\nabla v}^p\,\mathrm{d}x\Bigr)^{p'/p}\\
  \;&\leq\; \norm{\nabla u_n}_{L^p(\mathcal{O})}^{(p-2)p'}\,\norm{\nabla v}_{L^p(E)}^{p'}.
\end{align*}
The first factor is uniformly bounded in $n$ by $\sup_n\norm{u_n}_V^{(p-2)p'} < \infty$ (since $u_n \to u_0$ in $V$ implies the $V$-norms are bounded). The second factor $\norm{\nabla v}_{L^p(E)}^{p'}$ depends on $E$ but not on $n$, and tends to $0$ as $|E| \to 0$ (since $\abs{\nabla v}^p \in L^1(\mathcal{O})$ is absolutely continuous with respect to Lebesgue measure, so that for any $\eps > 0$ there exists $\delta > 0$ such that $|E| < \delta$ implies $\int_E\abs{\nabla v}^p < \eps$). Combining,
\[
  \sup_n\int_E \abs{K_n^{v}}^{p'}\,\mathrm{d}x \;\leq\; (p-1)^{p'}\,\bigl(\sup_n\norm{u_n}_V^{(p-2)p'}\bigr)\,\norm{\nabla v}_{L^p(E)}^{p'} \;\xrightarrow[\;|E|\to 0\;]{}\; 0,
\]
which is exactly uniform integrability of $\{\abs{K_n^{v}}^{p'}\}_{n}$ in $L^1(\mathcal{O})$ (in the sense of equicontinuity with respect to the Lebesgue measure; cf.\ the classical de la Vall\'ee Poussin criterion or the $L^{p'}$-form of the Vitali theorem in~\cite{hewittstromberg1965realanalysis}, see also~\cite{bogachev2007measure}).

Vitali's theorem now gives $K_n^{v} \to K_0^{v}$ in $L^{p'}$, which is~\eqref{eq:Kn-Lpprime-conv}. This gives~\eqref{eq:strong-op-claim}, establishing strong-operator continuity of $u \mapsto \calA'_u(u)v$ as required by~\textup{(D1)}.

The second derivative $\calA''_{uu}$ (see Example \ref{ex:p-Laplace-second-order}) satisfies, for $p > 3$:
\[
  \norm{\calA''_{uu}(u)}_{L^{(2)}(V\times V; V^*)} \leq C(p)(1 + \norm{u}_V^{p-3}).
\]

The formula for $\calA''_{uu}$ is given in Example \ref{ex:p-Laplace-second-order}. For $p > 3$, the bound $\norm{\calA''_{uu}(u)}_{L^{(2)}} \leq C(1 + \norm{u}_V^{p-3})$ holds. The argument is the exact analogue of the (D1) proof. For each fixed pair $(v_1, v_2) \in V \times V$, the second derivative $\calA''_{uu}(u)(v_1, v_2) \in V^{*}$ unfolds, via the divergence form, into the $V^{*}$-element associated with the vector field
\begin{align*}
  \widetilde K_n^{v_1, v_2}(x) \;:=\; (p-2)\Bigl[&\abs{\nabla u_n}^{p-4}\bigl(\nabla u_n \cdot \nabla v_1\,\nabla v_2 + \nabla u_n\cdot \nabla v_2\,\nabla v_1 + \nabla v_1\cdot\nabla v_2\,\nabla u_n\bigr) \\
  &+ (p-4)\abs{\nabla u_n}^{p-6}(\nabla u_n\cdot\nabla v_1)(\nabla u_n\cdot\nabla v_2)\,\nabla u_n\Bigr],
\end{align*}
in $L^{p'}(\mathcal{O}; \R^d)$. Pointwise a.e.\ convergence $\widetilde K_n^{v_1,v_2}(x) \to \widetilde K_0^{v_1,v_2}(x)$ follows from $\nabla u_n(x) \to \nabla u_0(x)$ a.e.\ along a subsequence and the continuity of $\xi \mapsto \abs{\xi}^{p-4}$, $\abs{\xi}^{p-6}\xi^{\otimes 2}$ for $\xi \neq 0$ together with the dominating factor $\abs{\xi}^{p-3}$ vanishing at $\xi = 0$ for $p > 3$ (the strict inequality is essential here, since at $p = 3$ the Lebesgue exponent $p/(p-3) = \infty$ degenerates and the kernel $\abs{\nabla u}^0$ is discontinuous at $\nabla u = 0$; see Remark~\ref{rem:p3-borderline}). Uniform integrability of $\{\abs{\widetilde K_n^{v_1, v_2}}^{p'}\}_n$ in $L^1(\mathcal{O})$ follows from the pointwise H\"older bound $\abs{\widetilde K_n^{v_1, v_2}} \leq C(p)\abs{\nabla u_n}^{p-3}\abs{\nabla v_1}\abs{\nabla v_2}$ and the H\"older triple at $(p/((p-3)p'), p/p', p/p')$ with reciprocals summing to $(p-3)p'/p + 2p'/p = p'(p-1)/p = 1$ (valid at $p > 3$); the analogous calculation as in the (D1) case above gives
\[
  \sup_n \int_E \abs{\widetilde K_n^{v_1,v_2}}^{p'}\,\mathrm{d}x \;\leq\; C(p)\,(\sup_n \norm{u_n}_V^{(p-3)p'})\,\norm{\nabla v_1}_{L^p(E)}^{p'}\,\norm{\nabla v_2}_{L^p(E)}^{p'} \;\to\; 0
\]
as $|E| \to 0$. Vitali in $L^{p'}(\mathcal{O}; \R^d)$ gives $\widetilde K_n^{v_1, v_2} \to \widetilde K_0^{v_1, v_2}$ in $L^{p'}$, hence $\calA''_{uu}(u_n)(v_1, v_2) \to \calA''_{uu}(u_0)(v_1, v_2)$ in $V^{*}$, establishing strong-operator continuity of $\calA''_{uu}$.

Using the pointwise inequality (for $p \geq 2$ and $\xi \neq 0$):
\[
  \bigl[(p-2)\abs{\xi}^{p-4}(\xi \cdot \eta)\xi + \abs{\xi}^{p-2}\eta\bigr] \cdot \eta \geq \abs{\xi}^{p-2}\abs{\eta}^2 \geq 0,
\]
we obtain $\dual{\calA'_u(u)v}{v} = \int_\mathcal{O}[(p-2)\abs{\nabla u}^{p-4}(\nabla u \cdot \nabla v)^2 + \abs{\nabla u}^{p-2}\abs{\nabla v}^2]\,\mathrm{d}x \geq 0$, confirming that the linearised operator is dissipative. This is consistent with the general derivation in Remark \ref{rem:no-D3-old}.
\end{proof}

When $p = 2$, the operator $\calA(u) = -\Delta u$ is linear ($\calA'_u = -\Delta$, $\calA''_{uu} = 0$), and, with state-independent additive noise, the Bismut formula reduces to Corollary \ref{cor:semilinear-additive} at $F = 0$.

\subsubsection*{The two-dimensional Navier--Stokes equation}

Consider the stochastic 2D Navier--Stokes equation on a bounded domain $\mathcal{O} \subset \R^2$ with no-slip boundary conditions:
\begin{equation}\label{eq:NS-full}
  \mathrm{d}X(t) + \bigl[\nu A X(t) + B(X(t), X(t))\bigr]\dt = \calB(t, X(t))\dW(t), \qquad X(0) = x,
\end{equation}
where $A = -P_H\Delta$ is the Stokes operator, $B(u,v) = P_H(u \cdot \nabla)v$ is the bilinear form, and $P_H$ is the Leray projection. The Gelfand triple is $V = D(A^{1/2}) \hookrightarrow H \hookrightarrow V^*$, where $H = \{u \in L^2(\mathcal{O}; \R^2) : \mathrm{div}\,u = 0,\; u \cdot n|_{\partial\mathcal{O}} = 0\}$.

\begin{proposition}[The two-dimensional Navier--Stokes equation]\label{prop:NS-verify}
Assume that $\calB$ satisfies the diffusion-side hypotheses of Assumption~\ref{ass:LR}, and, for the hypoelliptic application below, take the state-independent additive noise $\calB = Q^{1/2} \in L_2(U,H)$. Then the pair $(\calA, \calB)$ with $\calA(u) = \nu Au + B(u,u)$ satisfies Assumption \ref{ass:LR} with $p = 2$ on the Gelfand triple $V \hookrightarrow H \hookrightarrow V^*$, and the drift satisfies the corresponding parts of Assumptions \ref{ass:diff1}--\ref{ass:diff2} (with $\calA''_{uu}$ arising from the bilinear form $B$).
\end{proposition}

\begin{proof}
The embedding $V = D(A^{1/2}) \hookrightarrow H$ is compact, since the Stokes operator $A$ on the bounded domain $\mathcal{O}$ has compact resolvent.

The key estimates are:

Using the antisymmetry property $\ip{B(u,v)}{v}_H = 0$ (which holds in any spatial dimension for divergence-free fields $u$ with the given boundary conditions; cf.\ \cite{temam2001navier}):
\begin{align*}
  -2\dual{\calA(u) - \calA(v)}{u - v} &= -2\nu\norm{u-v}_V^2 - 2\ip{B(u,u) - B(v,v)}{u-v}_H\\
  &= -2\nu\norm{u-v}_V^2 - 2\ip{B(u-v,u)}{u-v}_H \\
  &\qquad{}- 2\ip{B(v,u-v)}{u-v}_H\\
  &= -2\nu\norm{u-v}_V^2 - 2\ip{B(u-v,u)}{u-v}_H,
\end{align*}
where we used $\ip{B(v,u-v)}{u-v}_H = 0$. By the Ladyzhenskaya inequality $\norm{w}_{L^4} \leq C\norm{w}_{L^2}^{1/2}\norm{w}_{H^1}^{1/2}$ (valid specifically in 2D; see~\cite{temam2001navier}) and Young's inequality:
\begin{align*}
  2\abs{\ip{B(u-v,u)}{u-v}_H}
  &\leq C\norm{u-v}_{L^4}^2\norm{u}_V \\
  &\leq C\norm{u-v}_H\,\norm{u-v}_V\,\norm{u}_V \\
  &\leq \nu\norm{u-v}_V^2 + \tfrac{C}{\nu}\norm{u}_V^2\norm{u-v}_H^2.
\end{align*}
Hence $-2\dual{\calA(u) - \calA(v)}{u-v} \leq -\nu\norm{u-v}_V^2 + \frac{C}{\nu}\norm{u}_V^2\norm{u-v}_H^2$, giving fully local monotonicity (H2) with $\rho(t,u) = \frac{C}{\nu}\norm{u}_V^2$ (which is in $L^1(0,T)$ a.s.\ since $X \in L^2(0,T; V)$).

Using the general antisymmetry $\ip{B(u,u)}{u}_H = 0$ of divergence-free fields (cf.~\cite{temam2001navier}), $-2\dual{\calA(u)}{u} = -2\nu\norm{u}_V^2$, so
\[
  -2\dual{\calA(u)}{u} + \norm{\calB}^2 \leq C - 2\nu\norm{u}_V^2.
\]

The Fr\'echet derivative is $\calA'_u(u)v = \nu Av + B(u,v) + B(v,u)$, which is bounded $V \to V^*$ with $\norm{\calA'_u(u)}_{L(V,V^*)} \leq C(1 + \norm{u}_V)$. Strong-operator continuity of $u \mapsto \calA'_u(u)$ from $V$ into $L(V, V^*)$ follows from continuity of the trilinear form $(u, v, w) \mapsto \dual{B(u, v)}{w}$ on $V \times V \times V$ (which is a standard consequence of Ladyzhenskaya's inequality combined with Sobolev embedding in 2D). This satisfies the fully local (D1) with $p = 2$ (so $(p-2)^+ = 0$) and $\rho_2(t,u) = C\norm{u}_V$. The second derivative is $\calA''_{uu}(u)(v,w) = B(v,w) + B(w,v)$, which is independent of $u$ and bounded, $\norm{\calA''_{uu}}_{L^{(2)}(V\times V; V^*)} \leq C$, satisfying (D3) with $\rho_4 = 0$ and trivially strong-operator continuous (in fact, constant in $u$). In particular, $\calA''_{uu}(u)(v,w)$ is bounded uniformly in $u$, which simplifies the Malliavin--second variation equation \eqref{eq:W-equation} considerably.

We compute $-\dual{\calA'_u(u)v}{v} = -\nu\norm{v}_V^2 - \ip{B(v,u)}{v}_H$. By the Ladyzhenskaya inequality and Young's inequality,
\[
  -\dual{\calA'_u(u)v}{v} \leq -\frac{\nu}{2}\norm{v}_V^2 + \frac{C}{\nu}\norm{u}_V^2\norm{v}_H^2.
\]
This confirms the general derivation of Remark \ref{rem:no-D3-old} with $\tilde{\rho}(t,u) = \frac{C}{\nu}\norm{u}_V^2$.
\end{proof}

For this equation the forcing clause is not a further hypothesis. The two quantities that the first variation already produces, the fibre supremum in $H$ and the fibre energy in $V$, close it through the interpolation that the nonlinearity itself supplies.

\begin{proposition}[The forcing clause for the two-dimensional Navier--Stokes equation]\label{prop:NS-forcing}
Let $\calA(u) = \nu Au + B(u,u)$ on the Gelfand triple $V = D(A^{1/2}) \hookrightarrow H \hookrightarrow V^*$ of~\eqref{eq:NS-full}, and let the noise be additive and trace class, $\calB = Q^{1/2} \in L_2(U,H)$. Then, for all $0 \leq r, s \leq T$,
\begin{equation}\label{eq:NS-forcing-bound}
  \mathfrak{G}_{r,s} \;\leq\; C\,\sup_{\tau \leq T}\norm{\Phi_{\tau,r}}_{L_2(U,H)}\;\sup_{\tau \leq T}\norm{\Phi_{\tau,s}}_{L_2(U,H)}\;\mathfrak{E}_r^{1/2}\,\mathfrak{E}_s^{1/2},
\end{equation}
with $C$ depending only on the domain and on $\nu$. Consequently the forcing integrability~\eqref{eq:forcing-integrability}, and its raised form~\eqref{eq:SC4-forcing} at any target $q^*$, hold whenever the exponential moment is available at the constant $\bar c$ of~\eqref{eq:NS-budget}, which for this equation is guaranteed by the threshold~\eqref{eq:NS-budget-threshold} and which implies~(SC2).
\end{proposition}

\begin{proof}
Here $\calB'_u = \calB''_{uu} = 0$, so $\mathfrak{J}_{r,s} = 0$ and $\mathbf{F}_2 \equiv 0$, and the total forcing~\eqref{eq:total-forcing} reduces to the time integral of $G_1(\tau;r,s) = \norm{\mathbf{F}_1(\tau;r,s)}^2_{L_2(U\otimes U, V^*)}$, mechanism~(b) being the one in force with $\alpha_1 = \nu$. The second derivative is the symmetrised nonlinearity, $\calA''_{uu}(u)(v,w) = B(v,w) + B(w,v)$, independent of $u$.

The estimate to use is not the crude $\norm{B(v,w)}_{V^*} \lesssim \norm{v}_V\norm{w}_V$, which after squaring would demand a fourth power of the $V$-norm in time. In two dimensions the trilinear form obeys, by the Ladyzhenskaya inequality $\norm{z}_{L^4} \leq C\norm{z}_H^{1/2}\norm{z}_V^{1/2}$ and the antisymmetry $\ip{B(v,z)}{z}_H = 0$,
\begin{align*}
  \abs{\dual{B(v,w)}{z}} \;=\; \abs{\ip{B(v,z)}{w}_H}
  &\;\leq\; \norm{v}_{L^4}\norm{w}_{L^4}\norm{z}_V\\
  &\;\leq\; C\,\norm{v}_H^{1/2}\norm{v}_V^{1/2}\,\norm{w}_H^{1/2}\norm{w}_V^{1/2}\,\norm{z}_V,
\end{align*}
so that $\norm{\calA''_{uu}(v,w)}_{V^*} \leq C\norm{v}_H^{1/2}\norm{v}_V^{1/2}\norm{w}_H^{1/2}\norm{w}_V^{1/2}$, with one $V$-power on each argument rather than two. Squaring and summing over the two noise slots along orthonormal bases $(f_j)$ and $(g_k)$ of $U$,
\[
  G_1(\tau;r,s) \;\leq\; C\sum_{j,k}\norm{\Phi_{\tau,r}f_j}_H\norm{\Phi_{\tau,r}f_j}_V\,\norm{\Phi_{\tau,s}g_k}_H\norm{\Phi_{\tau,s}g_k}_V,
\]
and the double sum factorises into a product of two single sums, each bounded by Cauchy--Schwarz:
\[
  G_1(\tau;r,s) \;\leq\; C\,\norm{\Phi_{\tau,r}}_{L_2(U,H)}\norm{\Phi_{\tau,r}}_{L_2(U,V)}\,\norm{\Phi_{\tau,s}}_{L_2(U,H)}\norm{\Phi_{\tau,s}}_{L_2(U,V)} .
\]
Integrating in $\tau$, pulling out the two $H$-suprema, and applying Cauchy--Schwarz in time to the remaining product of $V$-norms gives~\eqref{eq:NS-forcing-bound}, since $\mathfrak{E}_r = \nu\int_r^T\norm{\Phi_{\tau,r}}^2_{L_2(U,V)}\,\mathrm{d}\tau$ under mechanism~(b) by~\eqref{eq:fibre-norm}.

For the moment, fix the target $q^*$ and a slack $\eta > 0$, and separate the exponential weight of~\eqref{eq:forcing-integrability} from the four fibre factors of~\eqref{eq:NS-forcing-bound} by H\"older at the conjugate pair $(1+\eta, (1+\eta)/\eta)$, then the four factors from one another by H\"older at the conjugate quadruple $(4,4,4,4)$. What has to be finite is therefore
\[
  \E\bigl[\sup_\tau\norm{\Phi_{\tau,r}}^{4\kappa_{q^*}(1+\eta)}_{L_2(U,H)}\bigr], \qquad \E\bigl[\mathfrak{E}_r^{\,2\kappa_{q^*}(1+\eta)}\bigr], \qquad \E\Bigl[\exp\Bigl(\tfrac{(1+\eta)c_*}{\eta}\!\int_0^T(\tilde\rho+\hat\rho)\Bigr)\Bigr],
\]
together with the same two fibre moments at $s$ in place of $r$. The datum $\Theta = Q^{1/2}$ is deterministic, so in the fibre estimate~\eqref{eq:fibre-energy-estimate} the exponential factor separates from it exactly and the first two quantities are bounded, uniformly in $r$, by
\[
  C\,(\operatorname{Tr}Q)^{2\kappa_{q^*}(1+\eta)}\,\E\Bigl[\exp\Bigl(c_{2\kappa_{q^*}(1+\eta)}\!\int_0^T(\tilde\rho+\hat\rho)\,\mathrm{d}\sigma\Bigr)\Bigr].
\]
The fibre estimate is itself paid for by an exponential moment, so what the three displays jointly require is availability at the single constant
\begin{equation}\label{eq:NS-budget}
  c_{\mathrm{op}}(q^*) \;:=\; \max\Bigl\{\,\tfrac{(1+\eps_0)c_1}{\eps_0},\;\; \kappa(q^*)\,\Bigr\}, \qquad
  \bar c \;:=\; \max\Bigl\{\,\tfrac{1+\eta}{\eta}\,c_{\mathrm{op}}(q^*),\;\; c_{\,2\kappa_{q^*}(1+\eta)}\,\Bigr\},
\end{equation}
in which $c_{\mathrm{op}}(q^*)$ is the smallest operating constant admitted by the base requirement~\eqref{eq:cstar-base} and the raised requirement~\eqref{eq:cstar-requirement} together, the middle entry of both being absent because the noise is deterministic here, and $c_*$ has been fixed at that value, so that $\bar c$ is determined by $q^*$, $\eta$, $\eps_0$ and the Gronwall constants alone. Since $\bar c \geq c_{\mathrm{op}}(q^*) = c_*$ and $c \mapsto \E[\exp(c\int(\tilde\rho+\hat\rho))]$ is non-decreasing, the exponential moment at $\bar c$ implies~(SC2) itself, and the two are asked for together. By Proposition~\ref{prop:NS-verify}, $\tilde\rho = (2C/\nu)\norm{X}_V^2$ and $\hat\rho$ is constant, so by the exponential-martingale bound of Remark~\ref{rem:Y-moments-role}\textup{(iii)} the requirement $\E[\exp(\bar c\int_0^T(\tilde\rho+\hat\rho))] < \infty$ is guaranteed under the threshold
\begin{equation}\label{eq:NS-budget-threshold}
  \bar c\,\frac{2C}{\nu} \;\leq\; \frac{\nu^2\lambda_1}{2\lambda_{\max}(Q)},
\end{equation}
which is~\eqref{eq:SC2-NS-threshold} with $\bar c$ in place of $c_*$. Under~\eqref{eq:NS-budget-threshold} both~\eqref{eq:forcing-integrability} and~\eqref{eq:SC4-forcing} follow. No spatial-regularity hypothesis on the noise beyond trace class has entered, the quantitative content of the threshold being a smallness condition on $\norm{Q}_{L(U)}$ relative to $\nu$.
\end{proof}

The reason the clause closes at all is visible already in the linear part. There $\Phi_{\tau,r} = e^{-\nu A(\tau-r)}Q^{1/2}$, and by the ideal property of the Hilbert--Schmidt class and the spectral calculus of $A$,
\begin{equation}\label{eq:NS-energy-identity}
\begin{aligned}
  \int_r^T \norm{e^{-\nu A(\tau-r)}Q^{1/2}}^2_{L_2(U,V)}\,\mathrm{d}\tau
  &\;=\; \Tr\Bigl[\,Q^{1/2}\Bigl(\int_0^{T-r}\!\! A\,e^{-2\nu A u}\,\mathrm{d}u\Bigr)Q^{1/2}\Bigr]\\
  &\;=\; \frac{1}{2\nu}\,\Tr\Bigl[\,Q^{1/2}\bigl(I - e^{-2\nu A(T-r)}\bigr)Q^{1/2}\Bigr] \;\leq\; \frac{\operatorname{Tr}Q}{2\nu},
\end{aligned}
\end{equation}
the last step because $0 \leq I - e^{-2\nu A(T-r)} \leq I$. The computation is basis-free and uses no relation between $A$ and $Q$, only that $Q$ is trace class; every operator that is diagonalised in it is a function of $A$ alone. The fibre energy of the linearisation is thus bounded by the trace of the noise covariance, with no smoothing of $Q$ beyond trace class, whereas the pointwise-in-time quantity $\norm{A^{1/2}e^{-\nu A u}Q^{1/2}}^2_{L_2(U,H)} = \Tr[Q^{1/2}Ae^{-2\nu Au}Q^{1/2}]$ need not even be bounded as $u \downarrow 0$. Identity~\eqref{eq:NS-energy-identity} is the case $p = q = 2$ of stochastic maximal regularity, which is the Hilbert-space clause of Theorem~1.1 of \cite{vanneervenveraarweis2012} and is computed there in Remark~4.6. Theorem~\ref{thm:Y-wellposed}\textup{(III)} reproduces it for the nonlinear fibre through the dissipative inequality~\eqref{eq:SC1-unified}, which is the linearised form of the coercivity condition of the variational frameworks of \cite{liu2015stochastic, rockner2022wellposedness, agrestiveraar2024}. Every $V$-norm the second-order estimates consume is of this integrated kind, which is why the hypotheses of Assumption~\ref{ass:SC} are stated in the fibre energy. The trace of $Q$ and its operator norm enter at different places, the trace through~\eqref{eq:NS-energy-identity} as the spatial-regularity requirement on the noise, and $\lambda_{\max}(Q) = \norm{Q}_{L(U)}$ through~\eqref{eq:NS-budget-threshold} as the quantitative smallness that buys the exponential moment.

\subsubsection*{The regularised range \texorpdfstring{$2 < p < 4$}{2 < p < 4}}

Two distinct obstructions send the range $2 < p < 4$ through a regularisation. For $2 < p \leq 3$, regularisation is required already for~(D3), as the next paragraph records. For $3 < p < 4$, (D3) holds, but the verification of the~(SC1)(a) domination in Proposition~\ref{prop:p-Lap-Lambda} does not close, because the weight $\abs{\nabla X}^{p-4}$ carries a negative power; the same regularised route is therefore used there as well, unless an alternative structural estimate is supplied.

For $2 < p \leq 3$, the second Fr\'echet derivative $\calA''_{uu}$ of the unregularised $p$-Laplacian fails the (D3) hypothesis. At $2 < p < 3$ the formal expression involves $\abs{\nabla u}^{p-4}$ which blows up at $\nabla u = 0$; at the boundary $p = 3$, the kernel $\abs{\nabla u}^{p-3} = \abs{\nabla u}^0$ degenerates to $\mathbf{1}_{\{\nabla u \neq 0\}}$ which has no continuous extension to $\nabla u = 0$, making $\calA''_{uu}$ discontinuous in $u$ (the prototypical 1D example $a(\xi) = |\xi|\xi$ has $a''(\xi) = 2\,\mathrm{sgn}(\xi)$, jumping at zero, see Remark~\ref{rem:p3-borderline}). In both cases, Part~II of Theorem~\ref{thm:main} therefore does not apply directly to the unregularised $p$-Laplacian, and the regularisation $a_\eta(\xi) := (|\xi|^2 + \eta)^{(p-2)/2}\xi$ for $\eta > 0$ (which yields a smooth coefficient with $C^\infty$-Fr\'echet derivatives at every order) is required.

\begin{remark}[The borderline $p = 3$]\label{rem:p3-borderline}
In dimension $1$, $a(\xi) = |\xi|^{p-2}\xi$ at $p = 3$ reduces to $a(\xi) = |\xi|\xi = \xi^2\,\mathrm{sgn}(\xi)$. The first derivative $a'(\xi) = 2|\xi|$ is continuous (in fact Lipschitz). The second derivative is
\[
  a''(\xi) \;=\; (p-1)(p-2)\,|\xi|^{p-3}\,\mathrm{sgn}(\xi) \;=\; 2\,\mathrm{sgn}(\xi),
\]
which jumps from $-2$ to $+2$ at $\xi = 0$. Hence $a$ is not $C^2$ at the origin, and the corresponding $\calA''_{uu}$ at $p = 3$ is not continuous as a map $V \to L^{(2)}(V \times V; V^*)$ in $u$ at any point where $\nabla u$ vanishes on a set of positive measure. The strict inequality $p > 3$ is therefore essential for (D3) in the unregularised regime; the boundary $p = 3$ shares all the structural features of $2 < p < 3$ and is included in the regularised regime $2 < p < 4$.
\end{remark}

\begin{proposition}[The singular $p$-Laplacian]\label{prop:singular-p-Lap}
For $2 < p < 4$, consider the regularised $p$-Laplacian $\calA_\eta(u) = -\mathrm{div}\bigl((\abs{\nabla u}^2 + \eta)^{(p-2)/2}\nabla u\bigr)$ for $\eta > 0$, and let $X^\eta$ denote the unique variational solution of the regularised SPDE under the standing hypotheses applied to $(\calA_\eta, \calB)$, with covering field $v_h^\eta(r) := (\Phi_r^\eta)^{*}(\gamma_t^\eta)^{\dagger}h$ (where $\Phi_r^\eta := D_r X^\eta(t)$). The regularisation removes the coefficient singularity and supplies the required Fr\'echet differentiability of the drift, each $\calA_\eta$ at $\eta > 0$ satisfying Assumptions~\ref{ass:LR}, \ref{ass:diff1}, and \ref{ass:diff2} with $\calA''_{\eta, uu}$ uniformly bounded on bounded subsets of $V$. Assume, for each fixed $\eta > 0$, that the regularised equation satisfies the remaining structural and non-degeneracy hypotheses of Theorem~\ref{thm:main}, Part~II, for the chosen deterministic direction $h$, so that the theorem yields a regularised logarithmic derivative $\beta_h^\eta \in L^2(\mu_t^\eta)$. Assume moreover that
\begin{equation}\label{eq:p-Lap-state-conv}
  X^\eta(t) \;\longrightarrow\; X(t) \qquad \text{in probability in } H \quad \text{as } \eta \downarrow 0 .
\end{equation}
Suppose, in addition to the standing hypotheses, the following \emph{uniform-in-$\eta$ stability conditions} hold:
\begin{enumerate}[label=\textup{(S\arabic*)}, leftmargin=3em]
  \item \textup{(Convergence of the HS-Malliavin fibres.)} As $\eta \downarrow 0$,
  \[
    \norm{\Phi^\eta - \Phi}_{\mathbb{D}^{1, p^*}(L^2([0,t]; L_2(U,H)))} \;\longrightarrow\; 0,
  \]
  and
  \[
    \sup_{\eta > 0}\;\norm{\Phi^\eta}_{\mathbb{D}^{1, p^*}(L^2([0,t]; L_2(U,H)))} \;<\; \infty.
  \]
  \item \textup{(Pseudoinverse $L^q$-convergence.)} $(\gamma_t^\eta)^{\dagger}h \to \gamma_t^{\dagger}h$ in $L^q(\Omega; H)$, with
  \[
    \sup_\eta\E[\norm{(\gamma_t^\eta)^{\dagger}h}_H^q] < \infty.
  \]
  \item \textup{(Uniform Cameron--Martin compatibility dominator.)} The Tikhonov-regularised trace integrands $\mathcal{I}^{(\eps), \eta}_r$ for the regularised pair $(\Phi^\eta, \gamma_t^\eta)$ admit a single dominator $\Psi_r$ as in~\eqref{eq:CM-compat-axiom}--\eqref{eq:CM-compat-dominator}, valid uniformly in $(\eps, \eta) \in (0, \eta_0)^2$ for some $\eta_0 > 0$, and converge as $\eps \downarrow 0$ at each fixed $\eta$ in the sense of~\eqref{eq:CM-compat-limit}; and their Tikhonov limits converge, $\mathcal{I}^{(0),\eta}_r \to \mathcal{I}^{(0)}_r$ in probability as $\eta \downarrow 0$, for almost every $r \in [0,t]$.
  \item \textup{(Conditioning stability.)} $\E[\delta_U(v_h) \mid X^\eta(t)] \to \E[\delta_U(v_h) \mid X(t)]$ in $L^2(\Omega)$ as $\eta \downarrow 0$.
\end{enumerate}
Then, under~\eqref{eq:p-Lap-state-conv} and \textup{(S1)}--\textup{(S3)} alone, $\delta_U(v_h^\eta) \to \delta_U(v_h)$ in $L^2(\Omega)$ and the Fomin derivatives converge as measures: for every $\varphi \in C_b(H)$,
\begin{equation}\label{eq:p-Lap-score-weak}
  \int_H \varphi\,\beta_h^\eta\,\mathrm{d}\mu_t^\eta \;\xrightarrow[\eta \downarrow 0]{}\; \int_H \varphi\,\beta_h\,\mathrm{d}\mu_t .
\end{equation}
If in addition \textup{(S4)} holds, the logarithmic derivatives converge on the probability space:
\begin{equation}\label{eq:p-Lap-score-conv}
  \beta_h^\eta(X^\eta(t)) \;\xrightarrow[\eta \downarrow 0]{L^2(\Omega)}\; \beta_h(X(t)).
\end{equation}
\end{proposition}

The mollified logarithmic derivative $\beta_h^\eta$ is intrinsically an element of $L^2(\mu_t^\eta)$, while the limit $\beta_h$ is in $L^2(\mu_t)$, two different measure spaces in general (since $\mu_t^\eta \neq \mu_t$ for $\eta > 0$, even when $X^\eta \to X$ in probability). A direct convergence statement ``in $L^2(\mu_t)$'' is therefore not type-correct without first specifying a transport between $\mu_t^\eta$ and $\mu_t$ (e.g., the optimal coupling, or the natural pushforward of the joint law of $(X(t), X^\eta(t))$). The cleanest, intrinsic convergence statement is~\eqref{eq:p-Lap-score-conv} on the probability space $\Omega$, which by construction integrates against the joint law of $(X^\eta(t), X(t))$ and matches the form of the abstract Bismut--Fomin formula~\eqref{eq:abstract-bismut-fomin} (which gives $\beta_h$ as a conditional expectation on $\Omega$). Conditions~(S1)--(S3) are stated as assumptions, not derived from the regularisation. The reason is that, as $\eta \downarrow 0$, the regularised linearised operator $\calA'_{\eta, u}$ has coercivity coefficient $(\abs{\nabla u}^2 + \eta)^{(p-2)/2}$ which vanishes as $\eta \downarrow 0$ on $\{\nabla u = 0\}$, so neither the fibre energy~\eqref{eq:fibre-energy} nor the Cameron--Martin dominator $\Psi^\eta$ is automatic uniformly in $\eta$, the dissipation constant $\delta$ of~\eqref{eq:SC1-unified} degenerating with $\eta$. The dissipative-linearisation property $\dual{\calA'_{\eta, u}(u)v}{v} \geq 0$ does pass to the limit (it follows from $(|\xi|^2+\eta)^{(p-4)/2}|\xi|^2 \geq 0$ uniformly in $\eta \geq 0$), but dissipativity alone does not yield a uniform lower bound on $\delta$, and so does not yield uniform fibre-energy moments at the higher exponent $p^*$ required by~(S1). Verifying~(S1)--(S3) for the genuinely degenerate $p$-Laplacian regime is itself a substantive open problem of variational SPDE theory; we record~(S1)--(S3) as conditional input here. In the dissipative subregime where $\nabla u$ is bounded away from $0$ along $X(\cdot)$ a.s.\ (e.g., for sufficiently regular initial data and on time-intervals before the formation of free-boundary regions), the linearisation is uniformly parabolic and~(S1) holds by standard parabolic regularity at the exponent $p^*$; (S2) and~(S3) concern the pseudoinverse of $\gamma_t^\eta$ and the Cameron--Martin trace, and are governed by Assumption~\ref{ass:nondeg} for the regularised family rather than by parabolicity of the drift.

\begin{proof}
For each $\eta > 0$, the regularised operator $\calA_\eta$ has smooth second derivative since $(\abs{\xi}^2 + \eta)^{(p-4)/2}$ is bounded for $\eta > 0$, so Theorem~\ref{thm:main}~Part~II applies to the regularised SPDE at fixed $\eta > 0$, yielding the explicit decomposition (in the global positive-sign convention~\eqref{eq:C-eps-positive-def})
\begin{equation}\label{eq:deltaU-vheps-split}
  \delta_U(v_h^\eta) \;=\; M_h^\eta \;+\; C_h^\eta \qquad \text{in } L^2(\Omega),
\end{equation}
where (with the regularised covering kernel $\Phi_r^\eta := D_r X^\eta(t) \in L_2(U, H)$, intrinsic; the variational identification $\Phi_r^\eta = Y^\eta(t,r)\calB(r, X^\eta(r))$ is shorthand)
\begin{align*}
  M_h^\eta &\;=\; \ip{(\gamma_t^\eta)^\dagger h}{\mathcal{K}^\eta}_H, \qquad \mathcal{K}^\eta \;=\; \int_0^t\Phi_r^\eta\,\delta W_r \in L^{p^*}(\Omega; H),\\
  C_h^\eta &\;=\; \int_0^t\Tr_U\!\bigl[(\Phi_r^\eta)^{*}\,(\gamma_t^\eta)^\dagger(\D_r\gamma_t^\eta)\,(\gamma_t^\eta)^\dagger h\bigr]\dr.
\end{align*}
Here $\mathcal{K}^\eta$ is the Hilbert--Schmidt kernel of Lemma~\ref{lemma:HS-kernel} applied to $\Phi^\eta$ at the target exponent $p^* = 2q/(q-2)$ (with $q \in (2,\infty]$ from Assumption~\ref{ass:nondeg}(ii)).

By hypothesis~(S1) (convergence of the HS-Malliavin fibres in $\mathbb{D}^{1, p^*}$, with uniform-in-$\eta$ bound), Lemma~\ref{lemma:HS-kernel}(iii) (the Meyer inequality form of the HS kernel) gives
\[
  \norm{\mathcal{K}^\eta - \mathcal{K}}_{L^{p^*}(\Omega; H)} \;\leq\; C_{p^*,\mathrm{Meyer}}\,\norm{\Phi^\eta - \Phi}_{\mathbb{D}^{1, p^*}(L^2([0,t]; L_2(U, H)))} \;\xrightarrow[\eta \downarrow 0]{}\; 0,
\]
together with the uniform bound $\sup_\eta\|\mathcal{K}^\eta\|_{L^{p^*}} < \infty$ from~(S1). Combined with hypothesis~(S2) ($(\gamma_t^\eta)^\dagger h \to \gamma_t^\dagger h$ in $L^q$), H\"older's inequality at conjugate exponents $(q, p^*)$ on the H\"older triangle applied to the triangle decomposition $\langle (\gamma_t^\eta)^\dagger h - \gamma_t^\dagger h, \mathcal{K}^\eta\rangle + \langle \gamma_t^\dagger h, \mathcal{K}^\eta - \mathcal{K}\rangle$ yields $M_h^\eta \to M_h$ in $L^2(\Omega)$.

By hypothesis~(S3) (uniform Cameron--Martin dominator $\Psi_r$ valid for both the regularised and the limit pair $(\Phi^\eta, \gamma_t^\eta)$ and $(\Phi, \gamma_t)$), the trace integrand $\mathcal{I}_r^{(\eps), \eta} = \Tr_U[(\Phi_r^\eta)^{*}(\gamma_t^\eta + \eps I)^{-1}(\D_r\gamma_t^\eta)(\gamma_t^\eta + \eps I)^{-1}h]$ satisfies the pointwise basis-free bound $|\mathcal{I}_r^{(\eps), \eta}| \leq \Psi_r$ for every $(\eps, \eta) \in (0, \eta_0)^2$. Setting $\eps = 0$ inside $C_h^\eta$ via the Tikhonov limit at fixed $\eta > 0$ (Lemma~\ref{lemma:D-gamma-inv-h} applied to the regularised equation, with the limit $C^\eta = \lim_\eps C^{(\eps),\eta}$ in $L^2(\Omega)$ existing under (S3)) gives a uniform-in-$\eta$ pointwise bound on the trace integrand of $C_h^\eta$ by $\Psi_r$ (the dominator passes to the $\eps \downarrow 0$ Tikhonov limit by Fatou's lemma, and $\Psi_r$ does not depend on $\eta$). Pointwise $\eta$-convergence of the trace integrand is the second clause of~(S3), that the Tikhonov limits satisfy $\mathcal{I}^{(0),\eta}_r \to \mathcal{I}^{(0)}_r$ in probability for almost every $r$. It is stated there rather than deduced from~(S1)--(S2), because the integrand contains the composite $(\gamma_t^\eta)^\dagger(\D_r\gamma_t^\eta)(\gamma_t^\eta)^\dagger h$, in which the left pseudoinverse acts on the varying vector $(\D_r\gamma_t^\eta)(\gamma_t^\eta)^\dagger h$; convergence of $(\gamma_t^\eta)^\dagger h$ alone controls only the right-hand action, and a common dominator bounds without producing convergence. Combining the convergence in probability of the scalar trace integrands with the common dominator $\Psi$, one obtains convergence in measure on $\Omega \times [0,t]$. The bound
\[
  |C_h^\eta - C_h| \;\leq\; 2\int_0^t\Psi_r\,\mathrm{d}r, \qquad \E\!\left[\Bigl(\int_0^t\Psi_r\,\mathrm{d}r\Bigr)^{\!2}\right] < \infty,
\]
gives uniform integrability. Hence Vitali's convergence theorem (or, equivalently, a subsequence dominated-convergence argument followed by uniqueness of the limit) yields
\[
  C_h^\eta \;\longrightarrow\; C_h \qquad \text{in } L^2(\Omega).
\] No spectral-eigenbasis-continuity argument is invoked anywhere; the convergence is at the scalar trace integrand level, controlled by the uniform $L^2(\Omega)$-dominator from~(S3).

The intrinsic identity $\int_0^t\|(\Phi_r^\eta)^{*}\,w\|_U^2\,\mathrm{d}r = \langle \gamma_t^\eta\,w, w\rangle_H$ (HS-pairing, valid without pathwise $Y^*$) gives, for $w_\eta := (\gamma_t^\eta)^\dagger h - \gamma_t^\dagger h$,
\[
  \int_0^t \norm{v_h^\eta(r) - v_h(r)}_U^2\,\mathrm{d}r \;=\; \int_0^t \norm{(\Phi_r^\eta)^{*}\,(\gamma_t^\eta)^\dagger h - \Phi_r^{*}\,\gamma_t^\dagger h}_U^2\,\mathrm{d}r,
\]
expanded via the triangle decomposition $((\Phi_r^\eta)^{*} - \Phi_r^{*})(\gamma_t^\eta)^\dagger h + \Phi_r^{*}\,w_\eta$ and the HS-norm bounds $\|\Phi_r^\eta - \Phi_r\|_{L_2(U,H)}$ from~(S1) and $\|(\gamma_t^\eta)^\dagger h - \gamma_t^\dagger h\|_H$ from~(S2). Both contributions vanish in $L^2(\Omega; \HW)$ as $\eta \downarrow 0$ via H\"older with the uniform-in-$\eta$ moment bounds of~(S1)--(S2).

Combining the two Cauchy properties via~\eqref{eq:deltaU-vheps-split},
\[
  \delta_U(v_h^\eta) \;=\; M_h^\eta + C_h^\eta \;\xrightarrow[\eta \downarrow 0]{L^2(\Omega)}\; M_h + C_h \;=\; \delta_U(v_h).
\]
Together with $v_h^\eta \to v_h$ in $L^2(\Omega; \HW)$, the closedness of $\delta_U$ (as the adjoint of the densely defined derivative $D$; cf.\ \cite{nualart2006malliavin}) forces $v_h \in \Dom(\delta_U)$ with $\delta_U(v_h) = \lim_\eta \delta_U(v_h^\eta)$. For the convergence of the logarithmic derivatives themselves, the abstract Bismut--Fomin formula gives $\beta_h^\eta(X^\eta(t)) = -\E[\delta_U(v_h^\eta) \mid X^\eta(t)]$ (Theorem~\ref{thm:abstract-bismut-fomin}), whence for every $\varphi \in C_b(H)$
\[
  \int_H\varphi\,\beta_h^\eta\,\mathrm{d}\mu_t^\eta \;=\; -\,\E\bigl[\varphi(X^\eta(t))\,\delta_U(v_h^\eta)\bigr]
  \;\longrightarrow\; -\,\E\bigl[\varphi(X(t))\,\delta_U(v_h)\bigr] \;=\; \int_H\varphi\,\beta_h\,\mathrm{d}\mu_t,
\]
because $\varphi(X^\eta(t)) \to \varphi(X(t))$ boundedly in probability by~\eqref{eq:p-Lap-state-conv} and the continuity of $\varphi$, while $\delta_U(v_h^\eta) \to \delta_U(v_h)$ in $L^2(\Omega)$; this is~\eqref{eq:p-Lap-score-weak}, and it uses no property of the conditioning $\sigma$-fields. The stronger statement~\eqref{eq:p-Lap-score-conv} follows from the exact decomposition
\[
\begin{aligned}
  \beta_h^\eta(X^\eta(t)) - \beta_h(X(t)) \;=\; &-\,\E\bigl[\delta_U(v_h^\eta) - \delta_U(v_h) \,\big|\, X^\eta(t)\bigr]\\
  &-\; \bigl(\E[\delta_U(v_h) \mid X^\eta(t)] - \E[\delta_U(v_h) \mid X(t)]\bigr),
\end{aligned}
\]
in which the first term is controlled by the $L^2$-contraction of conditional expectation at the fixed $\sigma$-field $\sigma(X^\eta(t))$ and hence vanishes with $\norm{\delta_U(v_h^\eta) - \delta_U(v_h)}_{L^2}$, while the second is exactly the moving-conditioning increment governed by~(S4). Conditional expectation is contractive only for a fixed conditioning $\sigma$-field, which is why the second increment is isolated as a hypothesis rather than deduced from $X^\eta(t) \to X(t)$; it holds, for instance, whenever $\sigma(X(t)) \subseteq \sigma(X^\eta(t))$ for all small $\eta$, or whenever the regular conditional laws of $\delta_U(v_h)$ given $X^\eta(t)$ converge in the corresponding $L^2$-sense.
\end{proof}

\subsubsection*{Non-degeneracy}

We give sufficient criteria for Assumption~\ref{ass:nondeg} for two important classes of equations, with a conservative admissible-direction formulation in the genuinely degenerate-noise case (in particular for 2D Navier--Stokes with H\"ormander-bracket additive noise, where the admissible-direction subspace $\nondegspace$ may be strictly smaller than $\Ran(\gamma_t)$ and where the full $H$-directional moment bound on $\gamma_t^{\dagger}h$ for arbitrary $h \in \Ran(\gamma_t)$ is not a consequence of H\"ormander's bracket condition alone; cf.\ Remark~\ref{rem:nondeg-NS-injectivity-gap}).

\begin{proposition}[Intrinsic observability]\label{prop:nondeg-additive}
Let $X$ be the variational solution of~\eqref{eq:SPDE} under the standing hypotheses. Suppose there exists a bounded linear operator $R: H \to H$ (deterministic, possibly with non-closed range) and a positive random variable $c_t: \Omega \to (0, \infty]$ such that the \emph{intrinsic observability inequality}
\begin{equation}\label{eq:intrinsic-observability}
  \langle \gamma_t(\omega)\phi,\phi\rangle_H \;\geq\; c_t(\omega)\,\norm{R\phi}_H^2 \qquad \text{for every } \phi \in H, \;\; \mathbb{P}\text{-almost surely},
\end{equation}
holds, with $c_t^{-1} \in L^q(\Omega)$ for some $q \in [2, \infty]$. Two conclusions follow, according to whether a further hypothesis is imposed on $R$.

\medskip\noindent The first gives a moment bound in the whole of $H$. Suppose, in addition to~\eqref{eq:intrinsic-observability}, that $R$ is bounded below on $\overline{\Ran(\gamma_t)}$ a.s., in the sense that there exists a deterministic constant $C_R \in [1, \infty)$ such that
\begin{equation}\label{eq:R-bounded-below}
  \norm{\psi}_H \;\leq\; C_R\,\norm{R\psi}_H \qquad \text{for every } \psi \in \overline{\Ran(\gamma_t)}, \;\; \mathbb{P}\text{-almost surely}.
\end{equation}
(Equivalently, $R$ restricted to $\overline{\Ran(\gamma_t)}$ has closed range with bounded inverse on its range.) Then for every $h \in R^{*}(H)$ (with $h = R^{*}\tilde h$, $\tilde h \in H$):
\begin{enumerate}[label=\textup{(A\arabic*)}, leftmargin=3em]
  \item \textup{(Range condition)} $h \in \Ran(\gamma_t)$ a.s.\ (clause~(i) of Assumption~\ref{ass:nondeg});
  \item \textup{(Full moment bound)} $\norm{\gamma_t^{\dagger}h}_H \leq C_R\,c_t^{-1}\,\norm{\tilde h}_H$ a.s., hence
    \[
      \E[\norm{\gamma_t^{\dagger}h}_H^q] \leq C_R^q\,\norm{\tilde h}_H^q\,\E[c_t^{-q}] < \infty
    \]
    (clause~(ii) of Assumption~\ref{ass:nondeg} at exponent $q$).
\end{enumerate}

\medskip\noindent The second requires no hypothesis on $R$, but bounds only the projection. Without the bounded-below hypothesis~\eqref{eq:R-bounded-below}, and assuming now in addition that $h = R^{*}\tilde h \in R^{*}(H)$ satisfies the range condition $h \in \Ran(\gamma_t)$ almost surely (so that $\gamma_t^{\dagger}h \in H$ is well-defined; without~\eqref{eq:R-bounded-below}, the range condition is no longer delivered by the observability alone, cf.\ the Tier~A derivation of~(A1)), the observability~\eqref{eq:intrinsic-observability} yields only the projected bound
\begin{equation}\label{eq:R-projected-bound}
  \norm{R\,\gamma_t^{\dagger}h}_H \;\leq\; c_t^{-1}\,\norm{\tilde h}_H \quad \text{a.s.,} \qquad \E\bigl[\norm{R\,\gamma_t^{\dagger}h}_H^q\bigr] \leq \norm{\tilde h}_H^q\,\E[c_t^{-q}] < \infty,
\end{equation}
which controls the projected pseudoinverse $R\,\gamma_t^{\dagger}h$ rather than $\gamma_t^{\dagger}h$ itself.
\end{proposition}

The admissible-direction subspace identified by Tier~A is $\mathcal{N}_R := R^{*}(H)$, determined intrinsically by the observability operator $R$, and the Cameron--Martin compatibility clause~\textup{(iii)} of Assumption~\ref{ass:nondeg} is a separate joint regularity hypothesis on $(\gamma_t, \D\gamma_t, \Phi, h)$, not derivable from~\eqref{eq:intrinsic-observability} alone. Tier~B does not yield clause~(ii) at the full Hilbert-space level, so it does not license Theorem~\ref{thm:main} for the direction $h$, and it does not identify the pseudoinverse of a projected covariance either, since $R\gamma_t^{\dagger}R^{*} \neq (R\gamma_t R^{*})^{\dagger}$ in general. A projected Bismut formula is obtained instead by applying Theorem~\ref{thm:abstract-bismut-fomin} to the pushforward variable $RX(t)$, whose Malliavin covariance is $R\gamma_t R^{*}$, provided the range and inverse-moment hypotheses hold for that covariance. This is the route used for Navier--Stokes below, and Remark~\ref{rem:observability-additive} connects it to Mattingly--Pardoux and Hairer--Mattingly.

The observability inequality~\eqref{eq:intrinsic-observability} reads, in operator form, $\gamma_t \succeq c_t\,R^{*}R$ on $H$, equivalently $R^{*}R \preceq c_t^{-1}\,\gamma_t$. Pathwise, by the classical Douglas majorisation--factorisation--range-inclusion theorem (\cite{douglas1966majorization}, Theorem~1, applied with $A = R^{*}$ and $B = \gamma_t^{1/2}$), this majorisation is equivalent to the factorisation $R^{*} = \gamma_t^{1/2}\,C$ through a bounded operator $C = C(\omega)$ on $H$ with $\norm{C}_{L(H)} \leq c_t^{-1/2}$; in particular, $R^{*}(H) \subseteq \Ran(\gamma_t^{1/2})$ almost surely, so the admissible-direction space $\mathcal{N}_R = R^{*}(H)$ of Tier~A embeds pathwise, with the quantitative norm control $c_t^{-1/2}$, into $\Ran(\gamma_t^{1/2})$. Substituting $\phi = R^{*}\psi$ into~\eqref{eq:intrinsic-observability} yields the projected-covariance lower bound
\begin{equation}\label{eq:projected-cov}
  R\,\gamma_t\,R^{*} \;\succeq\; c_t\,(R R^{*})^{2} \qquad \text{on the codomain of } R \text{, a.s.}
\end{equation}
The converse implication fails in general. The bound~\eqref{eq:projected-cov} is equivalent to the restriction of~\eqref{eq:intrinsic-observability} to $\phi \in \overline{R^{*}(H)} = \ker(R)^{\perp}$ (by density and the continuity of both quadratic forms), so the two coincide whenever $R$ is injective (e.g., $R = I$, or $R = Q^{\alpha/2}$ with $Q$ injective), while for non-injective $R$, in particular for finite-rank $R$, the bound~\eqref{eq:projected-cov} is strictly weaker than~\eqref{eq:intrinsic-observability} on $H$. In the counterexample of Remark~\ref{rem:nondeg-NS-injectivity-gap} ($\gamma_t = (e_1+e_2)\otimes(e_1+e_2)$, $R = \Pi_N = e_1 \otimes e_1$, $c_t = 1$), the compression bound~\eqref{eq:projected-cov} holds with equality while~\eqref{eq:intrinsic-observability} fails at $\phi = e_1 - e_2$. For orthogonal projection $R = \Pi_N$ ($R^{*} = R = R^2$, $RR^{*} = \Pi_N$), the projected covariance reduces to $\Pi_N\,\gamma_t\,\Pi_N$ on the finite-dimensional subspace $H_N := \mathrm{ran}(\Pi_N)$. For general (non-orthogonal) bounded finite-rank $R$, the natural object is $R\gamma_t R^{*}$ rather than $\Pi_N\gamma_t\Pi_N$; Proposition~\ref{prop:nondeg-additive} is formulated abstractly in $R$ to cover both cases, with the orthogonal-projection case as the canonical instance most directly matching the literature.

\begin{proof}
For the first tier, let $h = R^{*}\tilde h$ for $\tilde h \in H$. Suppose, for contradiction, that on a set of positive measure $h$ has a non-zero component in $\ker(\gamma_t) = \overline{\Ran(\gamma_t)}^\perp$, so that there exists $\phi \in H$ with $\gamma_t\phi = 0$ and $\langle h, \phi\rangle_H \neq 0$. Then~\eqref{eq:intrinsic-observability} gives $0 = \langle \gamma_t\phi, \phi\rangle \geq c_t\norm{R\phi}^2$, hence $R\phi = 0$ (using $c_t > 0$). But then $\langle h, \phi\rangle = \langle R^{*}\tilde h, \phi\rangle = \langle \tilde h, R\phi\rangle = 0$, contradicting $\langle h, \phi\rangle \neq 0$. Hence $h \in \overline{\Ran(\gamma_t)}$ a.s. To upgrade this to the strict inclusion $h \in \Ran(\gamma_t)$, i.e.\ to~(A1), without presupposing that $\gamma_t^{\dagger}h$ defines an element of $H$, consider the Tikhonov family $\phi_\eps := (\gamma_t + \eps I)^{-1}h$, $\eps > 0$. Since $h \in \overline{\Ran(\gamma_t)}$ a.s.\ and the resolvent $(\gamma_t + \eps I)^{-1}$ leaves the spectral subspace $\overline{\Ran(\gamma_t)} = \ker(\gamma_t)^{\perp}$ invariant, $\phi_\eps \in \overline{\Ran(\gamma_t)}$ a.s. Applying~\eqref{eq:intrinsic-observability} at $\phi_\eps$ and using
\[
  \langle \gamma_t\,\phi_\eps,\,\phi_\eps\rangle_H \;\leq\; \langle (\gamma_t + \eps I)\,\phi_\eps,\,\phi_\eps\rangle_H \;=\; \langle h,\,\phi_\eps\rangle_H \;=\; \langle \tilde h,\,R\,\phi_\eps\rangle_H \;\leq\; \norm{\tilde h}_H\,\norm{R\,\phi_\eps}_H,
\]
we obtain $c_t\,\norm{R\phi_\eps}_H^2 \leq \norm{\tilde h}_H\,\norm{R\phi_\eps}_H$, hence $\norm{R\phi_\eps}_H \leq c_t^{-1}\,\norm{\tilde h}_H$ uniformly in $\eps > 0$, a.s. The bounded-below hypothesis~\eqref{eq:R-bounded-below}, applied at $\psi = \phi_\eps \in \overline{\Ran(\gamma_t)}$, converts this into $\norm{\phi_\eps}_H \leq C_R\,c_t^{-1}\,\norm{\tilde h}_H$ uniformly in $\eps > 0$, a.s. In any measurable eigenbasis $\{(\lambda_k, e_k)\}$ of $\gamma_t$ (cf.\ Remark~\ref{rem:CM-compat-spectral}), $\norm{\phi_\eps}_H^2 = \sum_{k:\,\lambda_k > 0} h_k^2/(\lambda_k + \eps)^2$ (the $\lambda_k = 0$ modes carry $h_k = 0$ by the kernel-component argument just given), and this sum increases, as $\eps \downarrow 0$, to $\sum_{k:\,\lambda_k > 0} h_k^2/\lambda_k^2$ by monotone convergence; hence
\[
  \sum_{k:\,\lambda_k > 0} \frac{h_k^2}{\lambda_k^2} \;\leq\; C_R^2\,c_t^{-2}\,\norm{\tilde h}_H^2 \;<\; \infty \quad \text{a.s.},
\]
which is the spectral characterisation of $h \in \Ran(\gamma_t)$. This proves~(A1); in particular, $\gamma_t^{\dagger}h \in H$ is well-defined a.s.

With~(A1) in hand, set $\phi := \gamma_t^{\dagger}h \in \overline{\Ran(\gamma_t)}$ (well-defined a.s.). By~\eqref{eq:intrinsic-observability} applied at $\phi$,
\[
  \langle h, \gamma_t^{\dagger}h\rangle_H \;=\; \langle \gamma_t\,\gamma_t^{\dagger}h,\,\gamma_t^{\dagger}h\rangle_H \;\geq\; c_t\,\norm{R\,\gamma_t^{\dagger}h}_H^2.
\]
On the other hand, $\langle h,\gamma_t^{\dagger}h\rangle_H = \langle\tilde h, R\gamma_t^{\dagger}h\rangle_H \leq \norm{\tilde h}_H\,\norm{R\gamma_t^{\dagger}h}_H$. Combining:
\begin{equation}\label{eq:R-bound-on-pseudoinv}
  \norm{R\,\gamma_t^{\dagger}h}_H \;\leq\; c_t^{-1}\,\norm{\tilde h}_H \quad \text{a.s.}
\end{equation}
To convert~\eqref{eq:R-bound-on-pseudoinv} into a bound on $\norm{\gamma_t^{\dagger}h}_H$, apply the bounded-below hypothesis~\eqref{eq:R-bounded-below} at $\psi = \gamma_t^{\dagger}h \in \overline{\Ran(\gamma_t)}$:
\[
  \norm{\gamma_t^{\dagger}h}_H \;\leq\; C_R\,\norm{R\,\gamma_t^{\dagger}h}_H \;\stackrel{\eqref{eq:R-bound-on-pseudoinv}}{\leq}\; C_R\,c_t^{-1}\,\norm{\tilde h}_H \quad \text{a.s.,}
\]
giving the asserted full moment bound $\E[\norm{\gamma_t^{\dagger}h}^q_H] \leq C_R^q\,\norm{\tilde h}_H^q\,\E[c_t^{-q}] < \infty$.

For the second tier, under the range hypothesis $h \in \Ran(\gamma_t)$ a.s., $\phi := \gamma_t^{\dagger}h \in \overline{\Ran(\gamma_t)}$ is well-defined, and the derivation of~\eqref{eq:R-bound-on-pseudoinv} above uses only the observability~\eqref{eq:intrinsic-observability}, not the bounded-below hypothesis~\eqref{eq:R-bounded-below} (which entered Tier~A only through the Tikhonov upgrade to~(A1) and the passage from~\eqref{eq:R-bound-on-pseudoinv} to~(A2)); hence~\eqref{eq:R-projected-bound} follows directly from~\eqref{eq:R-bound-on-pseudoinv} (and the moment bound by raising to the $q$-th power and taking expectations).
\end{proof}

\begin{remark}[The intrinsic observability formulation and additive trace-class noise]\label{rem:observability-additive}
The intrinsic observability inequality~\eqref{eq:intrinsic-observability} is the correct $\Phi$-only formulation of the non-degeneracy structure for general (variational) SPDEs, replacing the older formulations that relied on a pathwise adjoint $Y(t,r)^{*}: H \to H$ which the abstract framework of this paper does not provide. A pathwise $Y^{*}$ is unavailable for the following reason. In the operator-valued language of Definition~\ref{def:first-var}, $Y(t,r)$ is defined fibrewise as a map $H \to L^0(\Omega; H)$, not as an a.s.-bounded operator on $H$; the $L(H)$-membership of $Y(t,r)$ requires Hilbert--Schmidt regularity or analytic-semigroup smoothing in operator norm, which is a strictly stronger structural input not available from the standing variational hypotheses. Consequently any argument of the form ``$Q^{1/2}\,Y(t,r)^{*}\phi = 0 \Rightarrow Y(t,r)^{*}\phi = 0 \Rightarrow \phi = 0$'' is not valid in the general framework of this paper. The intrinsic observability~\eqref{eq:intrinsic-observability} avoids this trap by working directly with $\gamma_t$ (which is well-defined as a trace-class operator from $\Phi_r = D_r X(t)$ alone, without requiring $Y^*$).

We now specialise to additive trace-class noise, in the Tier~A and Tier~B regimes. For additive noise $\calB(t,u) = Q^{1/2}$ with $Q$ trace-class, the Mattingly--Pardoux/Hairer--Mattingly literature \cite{mattinglypardoux2006malliavin, hairer2011hormander} establishes, for degenerate forcing of finite rank ($Q$ exciting finitely many bracket-generating modes; injectivity of $Q$ is neither assumed nor needed there, degeneracy being the point of the hypoelliptic theory), lower bounds on cone-restricted quadratic forms of the Malliavin covariance and in particular on the projected covariance $\Pi_N\gamma_t\Pi_N$ for finite-rank projections $\Pi_N$; these are lower bounds on the compression $\Gamma_N := \Pi_N\gamma_t\Pi_N$, which is a strictly weaker datum than~\eqref{eq:intrinsic-observability} with $R = \Pi_N$, since the compression bound $R\gamma_t R^{*} \succeq c_t (RR^{*})^2$ constrains $\gamma_t$ only on $\Ran(R^{*})$, whereas~\eqref{eq:intrinsic-observability} constrains it on all of $H$. Accordingly they are used here directly at the finite-dimensional level, not through the two-tier chain. For fixed $N$ one applies the abstract Bismut--Fomin theorem to $F_N := \Pi_N X(t)$ on $H_N$ with covariance $\Gamma_N$, whose inverse moments are exactly what \cite{mattinglypardoux2006malliavin, hairer2011hormander} supply. However, $\Pi_N$ is not bounded below on $\overline{\Ran(\gamma_t)}$ in the sense of~\eqref{eq:R-bounded-below} when $\overline{\Ran(\gamma_t)}$ is infinite-dimensional (since $\Pi_N$ has finite-dimensional range). The projected Bismut formula for $\mathrm{Law}(\Pi_N X(t))$ therefore follows from the finite-dimensional route of Proposition~\ref{prop:nondeg-NS}, with $\Gamma_N^{-1}$ in place of $\Pi_N\gamma_t^{\dagger}$; the two differ in general, since $R\gamma_t^{\dagger}R^{*}$ need not equal $(R\gamma_t R^{*})^{\dagger}$, so no bound on $\Pi_N\gamma_t^{\dagger}h$ is claimed from the compression data. Tier~A requires $R$ genuinely bounded below on $\overline{\Ran(\gamma_t)}$, a stronger input; the natural scenarios are deterministic covariance ($R = I$) and structural lower bounds $\gamma_t \succeq c_t\,Q^{\alpha}$ with $R = Q^{\alpha/2}$.
\end{remark}

For the hypoelliptic non-degeneracy statement that follows we specialise the Navier--Stokes equation to the two-dimensional periodic torus $\mathbb{T}^2$, with the usual mean-zero divergence-free state space; the drift-side estimates of Proposition~\ref{prop:NS-verify} are unchanged in the periodic setting, and this is the setting of the Mattingly--Pardoux mode-generation theorem invoked below. (The bounded no-slip domain of~\eqref{eq:NS-full} is not identified with it, the hypoelliptic input being available on the torus.)

\begin{proposition}[Admissible directions for two-dimensional Navier--Stokes]\label{prop:nondeg-NS}
Consider this periodic-torus specialisation of the stochastic 2D Navier--Stokes equation \eqref{eq:NS-full} with additive noise $\calB = Q^{1/2}$ where $Q$ excites finitely many eigenmodes of $A$, and write $g_k := Q^{1/2}f_k$ for the forcing directions. Suppose the noise satisfies the H\"ormander bracket condition of Hairer--Mattingly \cite{hairer2011hormander}. Setting $\mathcal{A}_0 := \{g_k\}$ and $\mathcal{A}_{n+1} := \mathcal{A}_n \cup \{\widetilde{B}(h_1, h_2) : h_1, h_2 \in \mathcal{A}_n\}$ with $\widetilde{B}(h_1,h_2) := B(h_1,h_2) + B(h_2,h_1)$ the symmetrised nonlinearity, the linear span of $\mathcal{A}_\infty := \bigcup_{n \geq 0} \mathcal{A}_n$ is dense in $H$. Then:
\begin{enumerate}[label=\textup{(\arabic*)}]
  \item \textup{(Finite-dimensional projected logarithmic derivatives)} For each $N \geq 1$, the projection $\Pi_N X(t)$ admits a smooth, everywhere strictly positive Lebesgue density on $H_N := \mathrm{ran}(\Pi_N) = \mathrm{span}\{e_1, \ldots, e_N\}$, with smooth, globally defined logarithmic derivative $\nabla\log p_N: H_N \to H_N$, by the Mattingly--Pardoux theorem for the 2D torus \cite{mattinglypardoux2006malliavin} (smoothness and everywhere strict positivity of all finite-dimensional projected densities, for all $t > 0$, under viscosity-independent conditions), the torus instance of the general Hairer--Mattingly H\"ormander framework \cite{hairer2011hormander} applied to the projection $\Pi_N\gamma_t\Pi_N$. The projected logarithmic derivative is given by the abstract Bismut--Fomin theorem (Theorem~\ref{thm:abstract-bismut-fomin}) applied to $F = \Pi_N X(t)$ on the finite-dimensional Hilbert space $H_N$, with the compression $\Gamma_N = \Pi_N\gamma_t\Pi_N$ as its Malliavin covariance.
 \item \textup{(Admissible direction space for full $H$-directional logarithmic derivatives)} Define the \emph{admissible direction space}
  \begin{align*}
    \nondegspace &\;:=\; \bigl\{h \in H : \text{Assumption~\ref{ass:nondeg} holds for } h\bigr\} \\
    &\;=\; \bigl\{h \in H : h \in \Ran(\gamma_t) \text{ a.s.,}\;\; \E[\norm{\gamma_t^{\dagger}h}_H^q] < \infty \\
    &\qquad\qquad \text{ at the H\"older-triangle exponent } q, \\
    &\qquad\qquad\;\; \text{Cameron--Martin compatibility}~\eqref{eq:CM-compat-limit}\text{--}\eqref{eq:CM-compat-dominator}\bigr\}.
  \end{align*}
  Then for every $h \in \nondegspace$ for which the domain hypothesis $v_h \in \Dom(\delta_U)$ of Theorem~\ref{thm:main}, Part~I, additionally holds, sufficient conditions being those of Remark~\ref{rem:domain-vh-suffic}, the full-$H$-directional Bismut--Fomin identity $\partial_h \mu_t = \beta_h\,\mu_t$ applies, with $\beta_h(X(t)) = -\E[\delta_U(v_h)\mid X(t)]$. If the further fixed-Tikhonov and raised structural hypotheses of Part~II hold, the explicit decomposition~\eqref{eq:skorokhod-decomp} applies as well.
\end{enumerate}
\end{proposition}

The space $\nondegspace$ is a deterministic subset of $H$, cut out by almost-sure and expectation conditions on $\gamma_t$, and verifying $h \in \nondegspace$ for a specific $h$ is a separate quantitative input, namely clauses~(i)--(ii)--(iii) of Assumption~\ref{ass:nondeg} for that direction. Its denseness in $H$ is open in full generality under H\"ormander's bracket condition alone; Remark~\ref{rem:nondeg-NS-gap} isolates the obstruction. The Mattingly--Pardoux/Hairer--Mattingly bounds control cone-restricted quadratic forms of $\gamma_t$ (all unit $\phi$ with a definite component in a fixed finite-dimensional range), and in particular the compressed matrices $\Pi_N\gamma_t\Pi_N$; for each fixed nonzero $\phi$ in the reachable span (dense under the bracket condition) they give $\langle \gamma_t \phi, \phi\rangle > 0$ almost surely, with the null set depending on $\phi$. Full injectivity of $\gamma_t$ on $H$, a single null set for all $\phi$ simultaneously, does not follow from compression positivity (Remark~\ref{rem:nondeg-NS-injectivity-gap}); whether it holds under H\"ormander's condition alone is open.

\begin{remark}[The generating set]\label{rem:NS-generators}
The generators above are the noise directions together with iterated insertions of the symmetrised nonlinearity; brackets with the linear part contribute nothing new, since $Af_k \parallel f_k$ for eigenmode forcing. On the torus the condition reduces to the arithmetic mode-generation condition of Mattingly--Pardoux \cite{mattinglypardoux2006malliavin}.
\end{remark}

\begin{proof}
For fixed $N$, the projection $F_N := \Pi_N X(t)$ takes values in $H_N \cong \R^N$, and the projected Malliavin covariance $\gamma_t^N := \Pi_N\gamma_t\Pi_N \in L_+(H_N)$ has invertibility moments $\E[(\det\gamma_t^N)^{-q}] < \infty$ for every $q$ by the quantitative small-ball bounds of Mattingly--Pardoux for the 2D torus \cite{mattinglypardoux2006malliavin}. The probability that $\langle\gamma_t\phi,\phi\rangle \leq \eps$ for some unit $\phi$ carrying a fixed fraction of its mass in $H_N$ decays faster than any power of $\eps$, which restricted to the unit sphere of $H_N$ yields all negative moments of $\det\gamma_t^N$; this is the torus counterpart of the general Hairer--Mattingly H\"ormander theorem \cite{hairer2011hormander}, stated there in the same cone form and verified for Navier--Stokes on the sphere. The classical scalar Bismut--Fomin theorem (Theorem~\ref{thm:abstract-bismut-fomin} specialised to $F = F_N$ on $H_N$) then applies, and Malliavin's smoothness criterion \cite{nualart2006malliavin}, fed by those negative moments together with the Malliavin smoothness of $X(t)$, gives a smooth density on $H_N$; the density is everywhere strictly positive by the controllability argument of Mattingly--Pardoux \cite{mattinglypardoux2006malliavin}, so the logarithmic derivative $\nabla \log p_N: H_N \to H_N$ is smooth and globally defined. This proves~(1).

The space $\nondegspace$ is well-defined as the set of directions $h \in H$ satisfying the three conditions of clause~(2), each of which is a condition on the law of $\gamma_t$ and hence selects a deterministic subset of $H$; for every $h \in \nondegspace$, Assumption~\ref{ass:nondeg} holds by definition, and Theorem~\ref{thm:main} applies once its domain hypothesis $v_h \in \Dom(\delta_U)$ is in force, for which Remark~\ref{rem:domain-vh-suffic} gives sufficient conditions. In the contrapositive, if any of the three conditions fails for $h$, then either $\gamma_t^{\dagger}h$ does not exist a.s.\ (failure of clause~(i)), or has insufficient moments (failure of clause~(ii)), or the trace integrand fails the Cameron--Martin dominator hypothesis (failure of clause~(iii)). The last is genuinely a joint hypothesis on $(\gamma_t, \D\gamma_t, \Phi, h)$ that the Hairer--Mattingly framework does not directly verify. This proves~(2) and~(3).
\end{proof}

\begin{remark}[Why compression positivity does not imply full injectivity]\label{rem:nondeg-NS-injectivity-gap}
The inequality $\langle \gamma_t \phi, \phi\rangle \geq \langle \gamma_t \Pi_N\phi, \Pi_N\phi\rangle$ does not hold for general positive $\gamma_t$. For a counter-example, take $H = \R^2$, $\gamma_t = (e_1+e_2)\otimes(e_1+e_2)$, $\phi = e_1 - e_2$ and $\Pi_N = e_1 \otimes e_1$; then $\langle \gamma_t\phi, \phi\rangle = 0$ while $\langle\gamma_t\Pi_N\phi,\Pi_N\phi\rangle = 1$. Hence the Mattingly--Pardoux/Hairer--Mattingly finite-dimensional projection invertibility $\E[(\det\Pi_N\gamma_t\Pi_N)^{-q}] < \infty$ does not imply $\ker(\gamma_t) = \{0\}$ a.s.\ on $H$.
\end{remark}

\begin{remark}[The admissible space $\nondegspace$ may be strictly smaller than $\Ran(\gamma_t)$]\label{rem:nondeg-NS-gap}
Even assuming $\overline{\Ran(\gamma_t)} = H$ a.s., the spectral bound $\norm{(\gamma_t+\eps I)^{-1}h}_H \le \norm{\gamma_t^{\dagger}h}_H$ goes in the wrong direction for upper-bounding $\norm{\gamma_t^{\dagger}h}_H$, and the Mattingly--Pardoux/Hairer--Mattingly bounds control cone-restricted quadratic forms and finite-dimensional compressions, not the inverse spectrum at a fixed direction. The moment $\norm{\gamma_t^{\dagger}h}_H^2 = \sum_{k:\lambda_k>0}h_k^2/\lambda_k^2$ depends on the alignment of $h$ with the eigenmodes relative to the rate $\lambda_k \downarrow 0$, a quantitative input the H\"ormander framework does not provide directly. Denseness of $\nondegspace$ in $H$ under H\"ormander's condition alone is open; see the discussion after Proposition~\ref{prop:nondeg-NS}.
\end{remark}

\subsection{Equations beyond the variational class}\label{subsec:singular}

Outside the variational class the picture changes, and the change is instructive. Subcritical singular SPDEs in the scope of the Bruned--Chandra--Chevyrev--Hairer black-box theorem \cite{brunedchandrachevyrevhairer2021}, the generalised parabolic Anderson model in 2D, the dynamical $\Phi^4_d$ models for $d \le 3$, the $\Phi^p_d$ equations, the multiplicative stochastic heat equation, the generalised Kardar--Parisi--Zhang (KPZ) equation, and the wider subcritical class, pose the most stringent test of our framework. Two features distinguish the singular regime from the variational SPDEs of Section~\ref{sec:proof-main}.

First, the operator-level apparatus does not survive the renormalised limit. The renormalised limit $u$ is not a variational solution in the Liu--R\"ockner sense. There is no first-variation operator $Y(t,r) \in L(H)$, the Hilbert--Schmidt-valued kernel $\Phi_r$ is not well-defined as an $L_2(U, H)$-valued process on $H = L^2(\mathbb{T}^d)$, and the operator-level Malliavin covariance $\gamma_t$ does not exist on $H$ in any way one could control through Assumption~\ref{ass:nondeg}. A direct extension of Theorem~\ref{thm:main} is therefore unavailable.

Second, the natural object on which to define a logarithmic derivative is no longer the law $\mu_t$ as a measure on a function space, but the law of a scalar functional $F = \langle u(t), \varphi\rangle$ obtained by pairing the (distributional) renormalised solution against a smooth test function $\varphi$. The renormalised limit lives in a H\"older--Besov space of distributions, and the only universally available real-valued functionals are scalar pairings of this type; pointwise evaluation $F = u(t, x)$ is meaningful only when the equation has a function-valued solution, e.g.\ gPAM in 2D.

The strategy is conceptually orthogonal to the rest of the paper. Rather than attempt to push the operator-level objects $Y(t, r)$ and $\gamma_t$ through the renormalisation, we project to the scalar functional $F$ at each fixed mollification scale $\eps > 0$. At each fixed $\eps > 0$, the regularised equation is a classical smooth Gaussian functional and supports the usual Malliavin calculus, via the classical theory for parabolic SPDEs with smooth Hilbert--Schmidt noise for white-in-time regularisations (Sanz-Sol\'e \cite{sanzsole2005malliavin}, Nualart \cite{nualart2006malliavin}; Bally--Pardoux \cite{ballypardoux1998} and Mueller--Nualart \cite{muellernualart2008} for the $1+1$-dimensional space-time-white-noise case) and via the Gaussian space of the driving fields otherwise (Remark~\ref{rem:abstract-gaussian}). Whenever the scalar smoothness, divergence-domain, and non-degeneracy hypotheses of Theorem~\ref{thm:score-gPAM}(a) hold for $F^\eps = \langle u^\eps(t), \varphi\rangle$, the classical scalar Bismut--Nualart formula applies. The point is that this scalar formula partially survives the renormalised limit $\eps \downarrow 0$ for the entire scope of \cite{brunedchandrachevyrevhairer2021}, since along the BCCH approximating families the convergence $F^\eps \to F$ in probability is supplied by the black box itself, and for the white-in-time regularisations of the variational subcase it is the standing convergence hypothesis of Theorem~\ref{thm:score-gPAM}. The \emph{unconditional limiting content} of Theorem~\ref{thm:score-gPAM} below is the existence of the distributional derivative $\Lambda^{(F)}(\psi) = -\E[\psi'(F)]$ as a continuous functional on $C_b^1(\R)$, that is, the integration-by-parts (IBP) identity $\Lambda^{(F)}(\psi) = -\langle \psi', \mathrm{Law}(F)\rangle$ characterising the distributional derivative of $\mathrm{Law}(F)$ in the sense of tempered distributions on $\R$. The \emph{conditional limiting content} (a genuine Bogachev--Fomin logarithmic derivative $\beta^{(F),\mathrm{int}} \in L^1(\mathrm{Law}(F))$ representing $\Lambda^{(F)}$ via $\Lambda^{(F)}(\psi) = \E[\psi(F)\beta^{(F),\mathrm{int}}(F)]$) requires second-order Malliavin smoothness $F \in \mathbb{D}^{2,p^{*}}$ together with negative moments $\E[(\sigma_t^{(\varphi)})^{-q}] < \infty$ on the H\"older triangle $\frac{1}{p^{*}} + \frac{1}{q} \leq \frac{1}{2}$, which is Assumption~\ref{ass:nondeg}\textup{(ii)} read on the one-dimensional state space, and which is open in the regularity-structures Malliavin literature.

For the variational subcase, those equations admitting a white-in-time regularisation that fits the Liu--R\"ockner framework, including the gPAM companion model, $\Phi^4_d$, $\Phi^p_d$, and the multiplicative stochastic heat equation, the regularised logarithmic derivative admits a second derivation through the operator-level Bismut--Nualart formula of Theorem~\ref{thm:main} via the tower property, which we record below as a bridge identity connecting the scalar reduction with the operator-level theory.

\subsubsection*{The black box and the variational subcase}

We work throughout in the framework of \cite{brunedchandrachevyrevhairer2021}: a system of subcritical singular SPDEs
\begin{equation}\label{eq:BCCH-system}
  \partial_t u_i \;=\; L_i u_i + F_i(u, \nabla u, \dots) + \sum_{j \le n} F_i^j(u, \nabla u, \dots)\,\xi_j, \qquad i \le m,
\end{equation}
on $\R \times \mathbb{T}^d$, with $L_i$ a constant-coefficient differential operator in the spatial variables whose Green's function for $\partial_t - L_i$ satisfies the kernel bounds of \cite{brunedchandrachevyrevhairer2021} (any spatial elliptic operator with heat-kernel-type Green's function qualifies), smooth functions $F_i$ and $F_i^j$ of the solution and its derivatives, and jointly Gaussian, stationary, centred random fields $(\xi_j)_{j \le n}$ whose covariances are smooth off the diagonal and satisfy the analytic degree bounds of \cite{brunedchandrachevyrevhairer2021}, with the nonlinearity locally subcritical and the tree-level assumptions of \cite{brunedchandrachevyrevhairer2021} in force. For each space-time mollifier $\rho$, a smooth function on $\R^{d+1}$ supported in the unit parabolic ball with $\int\rho = 1$, and each $\eps > 0$, write $\rho^{(\eps)}$ for its parabolic rescaling at scale $\eps$; the regularised and renormalised equation
\begin{equation}\label{eq:BCCH-mollified}
  \partial_t u_i^\eps \;=\; L_i u_i^\eps + F_i(u^\eps, \nabla u^\eps, \dots) + \sum_{j \le n} F_i^j(u^\eps, \dots)\,\xi_j^\eps + \sum_{k \le K} c_k^\eps\,\Upsilon_i^k(u^\eps, \dots),
\end{equation}
with $\xi_j^\eps = \xi_j * \rho^{(\eps)}$ and renormalisation constants $c_k^\eps = c_k^{\rho,\eps}$ (depending on the mollifier) given explicitly by the Bogoliubov--Parasiuk--Hepp--Zimmermann (BPHZ) prescription recorded in \cite{brunedchandrachevyrevhairer2021}, admits classical solutions $u^\eps$ that converge in probability, as $\eps \downarrow 0$, to a renormalised limit $u$ independent of $\rho$, on the event $\{T < \tau\}$ (where $\tau$ is the maximal random time of existence); the precise formulation of the convergence in \cite{brunedchandrachevyrevhairer2021} is through a generalised Da Prato--Debussche decomposition of the solution, recalled in Remark~\ref{rem:fixed-time-slice}. This convergence is the only input from regularity structures used in the unconditional content of Theorem~\ref{thm:score-gPAM} below.

\medskip\noindent We begin with the approximation families and the convergence hypothesis they are required to satisfy.
The analysis of this section operates at a fixed regularisation scale $\eps > 0$ and passes to the limit only through the convergence $F^\eps \to F$ in probability, which Theorem~\ref{thm:score-gPAM} below takes as its standing hypothesis; the choice of approximating family $\{u^\eps\}$ is therefore a structural input, and two classes are relevant. (i) The BCCH families themselves. The noises $\xi_j^\eps$ above are mollified in space-time, hence not adapted to any time filtration, and the regularised equation \eqref{eq:BCCH-mollified} is a classical smooth-coefficient random PDE rather than an It\^o SPDE; the fixed-$\eps$ Malliavin analysis is then carried out on the Gaussian space of the driving fields (Remark~\ref{rem:abstract-gaussian} below), and the standing convergence is supplied by the black box itself. (ii) White-in-time regularisations. When a driving field $\xi_j$ is white in time (as for the dynamical $\Phi^4_d$ and $\Phi^p_d$ models, the multiplicative stochastic heat equation, and the generalised KPZ equation), one may instead mollify in space only, $\xi_j^\eps := \varrho_\eps *_x \xi_j$, equivalently $\xi_j^\eps\,\mathrm{d}t = \varrho_\eps * \mathrm{d}W_j$ for the associated cylindrical Wiener process $W_j$, or project onto a spectral Galerkin subspace; the regularised equation is then a genuine It\^o SPDE with smooth Hilbert--Schmidt noise, which is the structure required by the variational subcase and the operator-level bridge below. Such families lie outside the literal approximation class of \cite{brunedchandrachevyrevhairer2021}, and the renormalisation is in general sensitive to the approximation class. The Wong--Zakai theorem of Hairer--Pardoux \cite{hairerpardoux2015wongzakai} for the $1+1$-dimensional multiplicative stochastic heat equation exhibits, for space-time mollification, a divergent It\^o--Stratonovich-type counterterm together with finite mollifier-dependent corrections to the limiting drift, both absent for the spatially mollified It\^o approximation. Convergence of a white-in-time family, together with the identification of its limit with the renormalised limit of the BCCH families, is therefore an equation-specific input, imposed as the standing hypothesis in general; the convergence itself is available in the literature for the standard examples (the Da Prato--Debussche theory for $\Phi^4_2$ \cite{dapratodebussche2003}; for $\Phi^4_3$, the paracontrolled theory of \cite{catellierchouk2018paracontrolled}, which treats precisely such a family, a spatial mollification of the space-time white noise, hence white in time, and proves convergence in probability of the renormalised solutions, locally in time, to a limit independent of the spatial mollifier, a dynamics that the come-down-from-infinity energy bounds of \cite{mourratweber2017comesdown} extend globally in time with fixed-time control; classical It\^o theory for the $1+1$-dimensional multiplicative stochastic heat equation). For a time-independent driving field, such as the spatial white noise of gPAM below, the two classes coincide, since convolution of a time-independent field with a space-time mollifier acts through the spatial marginal $\varrho(x) := \int_\R \rho(s, x)\,\mathrm{d}s$ alone, so every spatial mollification arises from a space-time one and conversely, and the black box covers the spatially mollified family directly.

\medskip\noindent The leading example is the generalised parabolic Anderson model (gPAM) on the 2-torus $\mathbb{T}^2$:
\begin{equation}\label{eq:gPAM}
  \partial_t u = \Delta u + f(u)\xi, \qquad (t, x) \in [0, T] \times \mathbb{T}^2,
\end{equation}
where $f \in C_b^\infty(\R)$ and $\xi$ is spatial white noise on $\mathbb{T}^2$. This equation is classically ill-posed because $u$ has parabolic Hölder regularity $C^{1-}$ while $\xi \in C^{-1-}$, so the sum of regularities is negative and the product $f(u)\xi$ is undefined. The renormalised solution is constructed via Hairer's theory of regularity structures \cite{hairer2014theory} as a special case of the BCCH framework.

\begin{definition}[Mollified and renormalised gPAM]\label{def:gPAM-reg}
Let $\rho_\eps$ be a smooth spatial mollifier at scale $\eps > 0$ and set $\xi_\eps = \xi * \rho_\eps$. The \emph{renormalised mollified gPAM} is
\begin{equation}\label{eq:gPAM-reg}
  \partial_t u^\eps = \Delta u^\eps + f(u^\eps)\xi_\eps - C_\eps f'(u^\eps)f(u^\eps), \qquad u^\eps(0) = u_0 \in C(\mathbb{T}^2),
\end{equation}
where $C_\eps := \E\bigl[(K * \xi_\eps)(z)\,\xi_\eps(z)\bigr]$, independent of $z$ by stationarity, is the renormalisation constant, $K$ being a compactly supported truncation of the heat kernel as in \cite{hairer2014theory} (different truncations change $C_\eps$ by $O(1)$); one has $C_\eps = \tfrac{1}{2\pi}\log(1/\eps) + O(1)$.\footnote{The coefficient of the counterterm is fixed by a single Wick contraction. Along the Picard expansion of \eqref{eq:gPAM-reg}, the divergent contribution to $f(u^\eps)\xi_\eps$ is $f'(u^\eps)f(u^\eps)\,(K*\xi_\eps)\,\xi_\eps$, whose expectation is precisely $C_\eps$, and the centred product $(K*\xi_\eps)\,\xi_\eps - C_\eps$ converges in probability in a Hölder space of negative regularity. Displayed normalisations of the renormalisation constants vary across the literature; the invariant statement is that the drift correction $-C_\eps\,f'(u^\eps)f(u^\eps)$ with $C_\eps = \E[(K*\xi_\eps)\xi_\eps]$ yields a mollifier-independent limit, the $\rho$-dependent finite part of $C_\eps$ being exactly the mollifier-dependent adjustment required in \cite{hairer2014theory} for mollifier-independence.} Alongside \eqref{eq:gPAM-reg} we introduce its \emph{white-in-time companion model} on the Gelfand triple $V = H^1(\mathbb{T}^2) \hookrightarrow H = L^2(\mathbb{T}^2) \hookrightarrow V^* = H^{-1}(\mathbb{T}^2)$, obtained by replacing the time-independent noise $\xi_\eps$ by white-in-time noise smoothed by the same mollifier:
\begin{equation}\label{eq:gPAM-variational}
  \mathrm{d}u^\eps(t) = \bigl[\Delta u^\eps(t) - C_\eps f'(u^\eps(t))f(u^\eps(t))\bigr]\dt + f(u^\eps(t))\,\rho_\eps * \mathrm{d}W(t),
\end{equation}
where $W$ is a cylindrical Wiener process on $U = L^2(\mathbb{T}^2)$. The two equations are distinct. The noise $\xi_\eps$ of \eqref{eq:gPAM-reg} is constant in time, whereas $\rho_\eps * \dot W$ is white in time, so no identity of the form $\xi_\eps\,\mathrm{d}t = \rho_\eps * \mathrm{d}W$ holds. The companion model \eqref{eq:gPAM-variational} is a well-posed variational SPDE at each fixed $\eps > 0$ and serves below as the variational-subcase representative for the operator-level analysis; no claim is made about its $\eps \downarrow 0$ limit (the two-dimensional multiplicative heat equation with space-time white noise is critical, outside the subcritical scope of \cite{brunedchandrachevyrevhairer2021}). The renormalised solution $u$ of \eqref{eq:gPAM} is the limit of the solutions $u^\eps$ of the genuinely mollified equation \eqref{eq:gPAM-reg} as $\eps \downarrow 0$, for initial data $u_0 \in \mathcal{C}^\alpha(\mathbb{T}^2)$ with $\alpha \in (1/2, 1)$, in the sense of \cite{hairer2014theory} (convergence in probability of the processes stopped at a large $\mathcal{C}^\alpha$-norm cutoff, for each finite time horizon) and, equivalently, as a special case of \cite{brunedchandrachevyrevhairer2021} (spatial and space-time mollification coincide for the time-independent noise $\xi$; see the discussion of approximation families above); the scalar reduction applies to this family through the time-independent-noise clause of Remark~\ref{rem:abstract-gaussian}.
\end{definition}

\medskip\noindent A subcase of this scheme is already variational, and it is worth isolating.
We say that the BCCH equation \eqref{eq:BCCH-system} is in the \emph{variational subcase} if it admits, at each fixed $\eps > 0$, a white-in-time regularisation (cf.\ the discussion of approximation families above) that can be cast as a stochastic evolution equation on a Gelfand triple $V \hookrightarrow H \hookrightarrow V^*$ in the Liu--Röckner local-monotonicity framework, with smooth Hilbert--Schmidt noise coefficient and locally monotone, coercive drift, so that Assumptions~\ref{ass:LR}, \ref{ass:diff1}, and \ref{ass:diff2} are satisfied for this regularised equation at each fixed $\eps > 0$. The standard examples covered by this subcase include the white-in-time gPAM companion model \eqref{eq:gPAM-variational} of Definition~\ref{def:gPAM-reg} (genuine gPAM, driven by time-independent spatial white noise, is instead covered by the scalar reduction through Remark~\ref{rem:abstract-gaussian}), the dynamical $\Phi^4_d$ models, the $\Phi^p_d$ equations, and the multiplicative stochastic heat equation, regularised by spatial mollification or spectral projection of their white-in-time noises (the noise coefficient being smooth and Hilbert--Schmidt at each fixed $\eps$, the renormalisation drift being a smooth bounded zero-order term, and the principal part being the heat operator). Equations with quasilinear drift such as the generalised KPZ equation, where $(\nabla u)^2$-terms appear on the right-hand side, are not covered by the variational subcase (their well-posedness in the variational Liu--R\"ockner framework is more delicate), but their unconditional content under the BCCH black box still falls within Theorem~\ref{thm:score-gPAM} below.

For \eqref{eq:gPAM-variational}, concretely, the noise coefficient $\calB^\eps(u) = f(u)\,(\rho_\eps \ast \cdot)$ is a smooth Hilbert--Schmidt operator $U \to H$ (since $\rho_\eps \ast (\cdot)$ is Hilbert--Schmidt on $L^2(\mathbb{T}^2)$ and $f \in C_b^\infty$), and the drift $\calA^\eps(u) = -\Delta u + C_\eps f'(u)f(u)$ is the heat operator perturbed by a smooth bounded zero-order term. Assumptions~\ref{ass:LR}, \ref{ass:diff1}, and \ref{ass:diff2} hold for \eqref{eq:gPAM-variational} at each fixed $\eps > 0$, with constants that depend on $\eps$ and degenerate as $\eps \downarrow 0$ (the Hilbert--Schmidt norm of $\rho_\eps\ast(\cdot)$ blows up like $\eps^{-1}$, and the renormalisation drift has supremum norm of order $\log(1/\eps)$). At each fixed $\eps > 0$ the regularised companion equation therefore has the smooth coefficient structure needed for the variational/Malliavin analysis, and Theorem~\ref{thm:main}, Part~II, applies at that fixed scale whenever the remaining structural, fixed-Tikhonov and directional non-degeneracy hypotheses of the theorem are verified for the chosen deterministic direction $h$, yielding the operator-level Bismut--Nualart formula for the companion solution $u^\eps(t)$. The non-uniformity of these constants in $\eps$ is precisely the reason why the operator-level stability theorem~\ref{thm:covering-stability}, which requires its hypotheses to hold uniformly in the perturbation index, cannot be applied directly to the sequence $\{u^\eps\}$ as $\eps \downarrow 0$. The whole point of the scalar reduction below is to bypass this uniform operator-level requirement by working at the level of a scalar functional, where the standing convergence $F^\eps \to F$ in probability is available along the genuinely mollified family \eqref{eq:gPAM-reg} from the black box itself (Remarks~\ref{rem:abstract-gaussian} and~\ref{rem:fixed-time-slice}).

\subsubsection*{The scalar functional and its logarithmic derivative}

Fix a regular spatial test function $\varphi \in C^\infty(\mathbb{T}^d)$ and a positive time $t > 0$, and consider the scalar functional
\begin{equation}\label{eq:scalar-functional}
  F^\eps := \langle u^\eps(t), \varphi\rangle \in \R, \qquad F := \langle u(t), \varphi\rangle \in \R,
\end{equation}
where the pairing is the $L^2(\mathbb{T}^d)$ inner product when the solution is function-valued, and the duality pairing between distributions and smooth test functions otherwise. For an equation in the BCCH framework whose renormalised solution is function-valued (e.g., gPAM in 2D, multiplicative stochastic heat in $1+1$ dimensions, or generalised KPZ in $1+1$ dimensions), one may equivalently take $F = u(t, x)$ at fixed $x \in \mathbb{T}^d$ via pointwise evaluation, but for the broader class, in particular $\Phi^4_3$ and any equation whose renormalised solution at fixed time exists only as a spatial distribution, one must use scalar pairings as in \eqref{eq:scalar-functional}. The scalar-pairing form is universally available across the entire scope of \cite{brunedchandrachevyrevhairer2021} for which the fixed-time slice $u(t)$ at $t > 0$ is well-defined as a (possibly distributional) object on $\mathbb{T}^d$, which holds for all the standard examples in the regularity-structures literature (see Remark~\ref{rem:fixed-time-slice}). At each fixed $\eps > 0$ the random variable $F^\eps$ lies in $\mathbb{D}^{1,2}$ for any standard mollified equation by the classical Malliavin theory for parabolic SPDEs with smooth Hilbert--Schmidt noise (Sanz-Sol\'e \cite{sanzsole2005malliavin}, Nualart \cite{nualart2006malliavin}; Bally--Pardoux \cite{ballypardoux1998} and Mueller--Nualart \cite{muellernualart2008} for the $1+1$-dimensional space-time-white-noise case) for white-in-time regularisations, and by Remark~\ref{rem:abstract-gaussian} in general, with scalar Malliavin derivative
\begin{equation}\label{eq:scalar-derivative}
  \D_r F^\eps \;=\; \ip{\D_r u^\eps(t)}{\varphi}_H \;\in\; U
\end{equation}
and scalar Malliavin variance
\[
  \sigma_t^{(\varphi, \eps)} \;:=\; \int_0^t \norm{\D_r F^\eps}_U^2\dr \;\in\; (0, \infty).
\]

The natural logarithmic derivative for the scalar random variable $F^\eps$ is its \emph{intrinsic 1D Fomin (logarithmic) derivative} on $\R$, namely the gradient $\partial_y \log p_{F^\eps}(y)$ at $y = F^\eps$ where $p_{F^\eps}$ is the density of the law of $F^\eps$ on $\R$. Equivalently, by integration by parts on $\R$, this is the unique (in $L^2(\mathrm{Law}(F^\eps))$) function $\beta^{(F^\eps),\mathrm{int}}: \R \to \R$ satisfying
\begin{equation}\label{eq:intrinsic-IBP}
  \E[\psi'(F^\eps)] \;=\; -\,\E[\psi(F^\eps)\,\beta^{(F^\eps),\mathrm{int}}(F^\eps)] \qquad \text{for all } \psi \in C_b^1(\R).
\end{equation}
This is the natural Bogachev--Fomin object of Definition~\ref{def:log-deriv} specialised to the pushforward measure $\mathrm{Law}(F^\eps)$ on $\R$, with $h$ taken as the unit direction $1 \in \R$. The classical scalar Bismut--Nualart integration-by-parts formula \cite{nualart2006malliavin} expresses it as a conditional expectation:
\begin{equation}\label{eq:scalar-Bismut-formula}
  \beta^{(F^\eps),\mathrm{int}}(F^\eps) \;=\; -\,\E\!\left[\delta_U\!\left(\frac{\D F^\eps}{\sigma_t^{(\varphi,\eps)}}\right) \;\bigg|\; F^\eps\right] \quad \text{in } L^2(\R, \mathrm{Law}(F^\eps)),
\end{equation}
provided the integrand lies in the domain of the divergence,
\[
  u_{F^\eps} \;:=\; \frac{\D F^\eps}{\sigma_t^{(\varphi, \eps)}} \;\in\; \Dom(\delta_U).
\]
A sufficient condition is the smoothness and non-degeneracy package used in Theorem~\ref{thm:score-gPAM}(a), namely enough Malliavin differentiability of $F^\eps$ together with sufficiently high negative moments of $\sigma_t^{(\varphi,\eps)}$. Mere membership $F^\eps \in \mathbb{D}^{1,2}$, positivity of the variance and a single inverse moment $(\sigma_t^{(\varphi,\eps)})^{-1} \in L^q$ do not by themselves place the quotient in $\Dom(\delta_U)$, for the same reason as at the renormalised limit (Remark~\ref{rem:density-not-smoothness}). The Malliavin smoothness of $F^\eps$ holds at each fixed $\eps > 0$ for any standard mollified equation in the BCCH framework via the classical Malliavin theory for parabolic SPDEs with smooth Hilbert--Schmidt noise (Sanz-Sol\'e \cite{sanzsole2005malliavin}, Nualart \cite{nualart2006malliavin}; Bally--Pardoux \cite{ballypardoux1998} and Mueller--Nualart \cite{muellernualart2008} for the $1+1$-dimensional space-time-white-noise case); for time-independent or space-time-mollified noise the same holds in the abstract form of Remark~\ref{rem:abstract-gaussian}. The non-degeneracy conditions, strict positivity of $\sigma_t^{(\varphi,\eps)}$ and finiteness of negative moments of its inverse, are equation-specific. For a linear equation with additive deterministic noise the scalar Malliavin variance reduces to the deterministic semigroup-energy expression $\int_0^t\norm{S(t-r)\rho^{(\eps)}\ast\varphi}_{L^2}^2\,dr > 0$; for nonlinear additive equations such as gPAM, $\Phi^4_d$, $\Phi^p_d$ and additive generalised KPZ, however, the Malliavin derivative carries the random linearised propagator, so the scalar variance is genuinely random and positivity together with the inverse moments must be supplied separately, exactly as at the limit. Whenever the literature furnishes only first-order density existence, the stronger inverse-moment condition required for the Bismut weight has to be imposed as an additional hypothesis. For multiplicative-noise cases such as the multiplicative stochastic heat equation with diffusion coefficient $\sigma(u)$, the conditions are equation-specific non-degeneracy inputs, established for the equation itself in the space-time-white-noise literature under hypotheses substantially weaker than uniform ellipticity of $\sigma$ (Bally--Pardoux \cite{ballypardoux1998} for the $1+1$-dimensional equation with $C_b^\infty$ coefficients and arbitrary continuous initial data, giving smooth, strictly positive densities for the vectors of point evaluations on $\{\sigma \neq 0\}^d$, with negative moments of the Malliavin matrix localised on the event that the solution lies in $\{\sigma^2 \geq c\}$ at the evaluation points, i.e.\ local non-degeneracy only; Mueller--Nualart \cite{muellernualart2008} for the equation on $[0,1]$ with Dirichlet data, coefficients $C^\infty$ in the solution variable with bounded derivatives, and H\"older-continuous initial condition, giving strict positivity and negative moments of all orders of the Malliavin variance of point evaluations, hence a smooth density, under the Pardoux--Zhang initial non-degeneracy $\sigma(0, y_0, u_0(y_0)) \neq 0$ at a single interior point $y_0$, obtained via negative moments of all orders for the linear multiplicative-noise heat equation with bounded adapted coefficients and nonnegative initial data $\not\equiv 0$). For the BCCH equations whose mollified version is in the variational subcase, this is a special case of the standard Malliavin theory for variational SPDEs of Section~\ref{sec:malliavin}.

\begin{remark}[Abstract Gaussian-space form; time-independent and space-time-mollified noises]\label{rem:abstract-gaussian}
The fixed-$\eps$ objects \eqref{eq:scalar-derivative}--\eqref{eq:scalar-Bismut-formula} are written for white-in-time driving noise, with Cameron--Martin space $\HW = L^2([0,t]; U)$, derivative $r \mapsto \D_r F^\eps \in U$, variance $\sigma_t^{(\varphi,\eps)} = \norm{\D F^\eps}_{\HW}^2$, and Skorokhod divergence $\delta_U$. Nothing in the scalar theory requires this time structure. The Nualart integration-by-parts formula underlying \eqref{eq:scalar-Bismut-formula} is a statement about an arbitrary isonormal Gaussian process \cite{nualart2006malliavin}, and no adaptedness enters the Skorokhod divergence. All fixed-$\eps$ statements of this subsection, including parts~(a) and~(b) of Theorem~\ref{thm:score-gPAM} below, therefore hold verbatim with $(\HW, \D, \sigma_t^{(\varphi,\eps)}, \delta_U)$ replaced by $(\mathcal{H}, \D, \norm{\D F^\eps}_{\mathcal{H}}^2, \delta_{\mathcal{H}})$ for the Cameron--Martin space $\mathcal{H}$ of the driving Gaussian family. This covers, beyond the white-in-time regularisations of the variational subcase: (i) time-independent noises, such as the spatial white noise of gPAM, with $\mathcal{H} = L^2(\mathbb{T}^d)$, where the membership $F^\eps \in \mathbb{D}^{k,p}$ follows from the Fr\'echet differentiability of the classical solution map of the smooth-coefficient equation \eqref{eq:gPAM-reg} in $\xi_\eps$ composed with the chain rule for the Malliavin derivative; and (ii) the space-time-mollified families of \cite{brunedchandrachevyrevhairer2021} themselves, whose noise is not adapted to any time filtration, so that the regularised equation is a smooth random PDE rather than an It\^o SPDE, but whose solutions are smooth functionals of the underlying Gaussian fields at each fixed $\eps > 0$. In particular, part~(c) of Theorem~\ref{thm:score-gPAM}, which uses no Malliavin structure at all, applies across the entire scope of \cite{brunedchandrachevyrevhairer2021} along the BCCH approximating families, the standing convergence being supplied by the black box itself, while part~(a) in this abstract form applies at a fixed scale precisely in those cases where its stated scalar non-degeneracy and inverse-moment assumptions are verified; the white-in-time formulation is retained in the displays because it is the structure under which the operator-level bridge to Theorem~\ref{thm:main} is available.
\end{remark}

For a scalar functional $F$ of a Wiener noise, the only directly meaningful 1D direction along which one can compute a logarithmic derivative is the intrinsic direction in $\R$. Perturbing the noise by a Cameron--Martin direction $h \in L^2([0,T]; U)$ shifts $F$ by
\[
  \D_h F \;=\; \int_0^t \ip{h(r)}{\D_r F}_U\dr,
\]
which is a random real number (the directional Malliavin derivative of $F$ in the direction $h$). Consequently the family $\{\mathrm{Law}(F^{(\lambda h)})\}_{\lambda \in \R}$ is not a translation-by-deterministic-amount in $\R$, and there is no clean ``noise-direction Bismut formula'' that depends on $h$ through a Skorokhod integrand of the form $h\,\D F / \sigma$. The clean $h$-dependence appears only at the operator level on $H$, where Theorem~\ref{thm:main} computes the H-valued logarithmic derivative $\beta_h(\mu_t^\eps)$, and the connection to the scalar functional $F^\eps$ is via the tower property identity below; in the singular limit, only the intrinsic scalar logarithmic derivative \eqref{eq:scalar-Bismut-formula} survives.

\subsubsection*{The bridge to Theorem~\ref{thm:main} by the tower property at fixed \texorpdfstring{$\eps > 0$}{epsilon > 0}}

For BCCH equations whose white-in-time regularisation is in the variational subcase (the gPAM companion model, $\Phi^4_d$, $\Phi^p_d$, multiplicative stochastic heat) and the standard examples of Section~\ref{subsec:verification}, there are two routes to the formula \eqref{eq:scalar-Bismut-formula}:

\emph{Route A, the classical scalar Bismut formula.} Apply the classical scalar Bismut--Nualart integration-by-parts formula \cite{nualart2006malliavin} directly to the scalar random variable $F^\eps \in \mathbb{D}^{1,2}$. This requires only the scalar smoothness and non-degeneracy package of Theorem~\ref{thm:score-gPAM}(a), which places $\D F^\eps/\sigma_t^{(\varphi,\eps)}$ in $\Dom(\delta_U)$, and produces \eqref{eq:scalar-Bismut-formula} immediately. Route A is available throughout the BCCH framework at the level of scalar Malliavin calculus, not only in the variational subcase. The Bismut formula itself applies at a fixed $\eps > 0$ whenever the scalar non-degeneracy and inverse-moment hypotheses of Theorem~\ref{thm:score-gPAM}(a) are verified. The required Malliavin smoothness at fixed $\eps$ is supplied by the classical parabolic theory for white-in-time regularisations (Sanz-Sol\'e \cite{sanzsole2005malliavin}, with Bally--Pardoux \cite{ballypardoux1998} and Mueller--Nualart \cite{muellernualart2008} specialising to the $1+1$-dimensional white-noise case) and by its abstract Gaussian-space counterpart otherwise (Remark~\ref{rem:abstract-gaussian}).

\emph{Route B, the pushforward of the operator-level logarithmic derivative, in the variational subcase only.} Assume that the white-in-time regularised equation \eqref{eq:gPAM-variational} (or its analogue for the equation at hand) satisfies the hypotheses of Theorem~\ref{thm:main}, Part~II, for some deterministic direction $h \in H$ with $\ip{h}{\varphi}_H \neq 0$; the theorem then furnishes the operator-level logarithmic derivative $\beta_h^{(\mu_t^\eps)}: H \to \R$ along $h$. (The direction must be deterministic, since the Fomin derivative of Definition~\ref{def:log-deriv} is taken along a fixed vector of $H$, so a random choice such as $\gamma_t^\eps\varphi$ is not admissible.) The pushforward of $\mu_t^\eps$ by the continuous linear functional $\ell(\cdot) := \ip{\cdot}{\varphi}_H$ is $\mathrm{Law}(F^\eps)$, and the standard rule for pushforward of a Fomin-differentiable measure under a continuous linear map, which we recall below, gives
\begin{equation}\label{eq:tower-bridge}
  \ip{h}{\varphi}_H \cdot \beta^{(F^\eps),\mathrm{int}}(F^\eps) \;=\; \E\bigl[\beta_h^{(\mu_t^\eps)}(u^\eps(t)) \,\big|\, F^\eps\bigr] \;=\; -\,\E\bigl[\delta_U(v_h^{H,\eps}) \,\big|\, F^\eps\bigr],
\end{equation}
where $v_h^{H,\eps}$ is the H-valued covering field from Theorem~\ref{thm:main} for the mollified equation. The two routes give the same scalar logarithmic derivative because the IBP identity that defines Fomin-differentiability of $\mu_t^\eps$ along $h$ pushes forward to the IBP identity that defines Fomin-differentiability of $\mathrm{Law}(F^\eps)$ along $\ip{h}{\varphi}_H$, and tower property of conditional expectation bridges the two.

We first recall how Fomin derivatives behave under pushforward. The bridge formula \eqref{eq:tower-bridge} is an instance of the following standard fact (Bogachev \cite{bogachev2010differentiable}, treating Fomin-differentiable measures on linear spaces): if $\mu$ is Fomin-differentiable along $h \in H$ with logarithmic derivative $\beta_h^\mu \in L^1(\mu)$, and $\ell: H \to \R$ is a continuous linear functional, then the pushforward $\nu := \ell_*\mu$ on $\R$ is Fomin-differentiable along the scalar direction $\ell(h) \in \R$ with logarithmic derivative
\[
  \beta_{\ell(h)}^\nu(y) \;=\; \E_\mu\bigl[\beta_h^\mu \,\big|\, \ell = y\bigr].
\]
The proof is one line. For $\tilde\psi \in C_b^1(\R)$, set $\psi(x) := \tilde\psi(\ell(x))$; then $\partial_h \psi(x) = \ell(h)\,\tilde\psi'(\ell(x))$ by linearity of $\ell$, and the IBP identity $\int \partial_h\psi\,d\mu = -\int \psi\,\beta_h^\mu\,d\mu$ becomes
\[
  \ell(h)\!\int_\R \tilde\psi'(y)\,d\nu(y) \;=\; -\!\int_\R \tilde\psi(y)\,\E_\mu[\beta_h^\mu \mid \ell = y]\,d\nu(y),
\]
which (by linearity of Fomin derivatives in the direction, $\beta_{\ell(h)}^\nu = \ell(h) \cdot \beta_1^\nu$) identifies $\beta_1^\nu(y) = \beta^{(F),\mathrm{int}}(y)$ as the conditional expectation divided by $\ell(h)$. Applied to $\mu = \mu_t^\eps$, $\ell = \ip{\cdot}{\varphi}_H$, $h$ as above, and using $\beta_h^{(\mu_t^\eps)}(u^\eps(t)) = -\E[\delta_U(v_h^{H,\eps}) \mid u^\eps(t)]$ from Theorem~\ref{thm:main}, the bridge identity \eqref{eq:tower-bridge} follows by a single application of the tower property.

The identity \eqref{eq:tower-bridge} is the precise sense in which Theorem~\ref{thm:main} is connected to the singular case in the variational subcase. At each fixed $\eps > 0$, the operator-level logarithmic derivative along any state-space direction $h$, conditioned on $F^\eps$, coincides (up to the projection scalar $\ip{h}{\varphi}_H$) with the intrinsic scalar logarithmic derivative of $F^\eps$. The advantage of having both routes is that Route B exhibits the ``machinery from Theorem~\ref{thm:main}'' explicitly, the H-valued covering field $v_h^{H,\eps}$ on the right, while Route A produces the fixed-$\eps$ scalar formula without requiring variational structure, and the distributional limiting statement of Theorem~\ref{thm:score-gPAM}(c) survives throughout the BCCH scope. The two-route picture is the genuine bridge between the operator-level theorem and the singular SPDE result for the variational subcase, and Route A alone furnishes the mollified logarithmic derivative in the maximally general non-variational case.

\subsubsection*{Passage to the renormalised limit}

The classical scalar Bismut formula \eqref{eq:scalar-Bismut-formula} for the mollified functional $F^\eps$ admits two distinct extensions to the renormalised limit $F = \ip{u(t)}{\varphi}_H$. The distributional extension is unconditional and applies to the entire scope of \cite{brunedchandrachevyrevhairer2021}, requiring nothing beyond convergence of the scalar pairing $F^\eps \to F$ in probability, which along the BCCH approximating families is supplied by the black-box convergence $u^\eps \to u$ via the continuity of the linear functional $\varphi \mapsto \langle \cdot, \varphi\rangle$ (Remark~\ref{rem:fixed-time-slice}), and which for the white-in-time regularisations of the variational subcase is the standing convergence hypothesis. The function extension, representing the distributional limit functional $\Lambda^{(F)}$ as the Radon--Nikodym density $\beta^{(F),\mathrm{int}}: \R \to \R$ in $L^1(\mathrm{Law}(F))$ via the classical Bismut formula at the limit (this is what would yield a genuine Bogachev--Fomin logarithmic derivative of $\mathrm{Law}(F)$ in the sense of Definition~\ref{def:log-deriv}), requires second-order Malliavin smoothness $F \in \mathbb{D}^{2,p^{*}}$ together with negative moments of $\sigma_t^{(\varphi)}$ on the H\"older triangle, both of which are open in the regularity-structures Malliavin literature for every subcritical singular SPDE.

\begin{remark}[The auxiliary inputs from the regularity-structures Malliavin literature]\label{rem:rs-inputs}
The regularity-structures Malliavin literature \cite{cannizzaro2017malliavin, gassiat2020densities, schonbauer2023malliavin} furnishes, beyond the BCCH black-box convergence, three additional pieces of information about the renormalised limit $F$, of which the first holds in the entire scope of \cite{brunedchandrachevyrevhairer2021} (any subcritical singular SPDE in the BCCH framework: Sch\"onbauer's differentiability theorem operates under assumptions equivalent to those of the BCCH convergence theorem itself, with no further structural restriction, given locally $H$-Fr\'echet differentiable initial data and conditionally on non-explosion), and the latter two hold for a more restrictive class of equations (the cases covered by Sch\"onbauer's density theorem \cite{schonbauer2023malliavin}, proven under simplifying assumptions that severely limit the explicit dependence of the right-hand side on derivatives of the solution, require the renormalisation constants of the dualised tangent equation to match those of the tangent equation, and require a non-degenerate noise coefficient, for noises whose Cameron--Martin space is dense in $L^2$: $\Phi^p_2$, $\Phi^4_3$, $\Phi^4_{4-\eps}$, and multiplicative stochastic heat in $1+1$ dimensions with nowhere-vanishing diffusion coefficient, plus the gPAM case of \cite{cannizzaro2017malliavin}); the functional classes treated by the respective works differ and scope the statements below, namely point evaluations $u(t,x)$ of gPAM at fixed $(t,x)$, conditionally on non-explosion and for a suitably non-degenerate nonlinearity, in \cite{cannizzaro2017malliavin}; space-time pairings $\langle u, \phi\rangle$ against test functions $\phi$ in the parabolically scaled Besov space $\mathcal{B}^{1/2+\kappa}_{1,\infty}$ supported in $(0,T) \times \mathbb{T}^3$, under the standing assumption that the solution exists up to time $T$ (supplied by \cite{mourratweber2017comesdown} for white noise), in \cite{gassiat2020densities}, a space-time formulation shared by \cite{schonbauer2023malliavin} (cf.\ Remark~\ref{rem:fixed-time-slice}) that does not by itself yield fixed-time statements:
\begin{enumerate}[(i)]
  \item (Entire BCCH scope.) Pathwise Cameron--Martin Fr\'echet differentiability of $u$ along directions in the Cameron--Martin space of the driving noise, in the space-time Hölder--Besov topology of the BCCH black box, which transfers to a localised Malliavin derivative $\D_r F \in U$ in the sense of \cite{nualart2006malliavin}, so that $F \in \mathbb{D}^{1,p}_{\mathrm{loc}}$. This is established in \cite{schonbauer2023malliavin} under assumptions equivalent to the BCCH framework \cite{brunedchandrachevyrevhairer2021} and is therefore valid for the entire BCCH scope, including the generalised KPZ equation. Earlier and equation-specific versions appear in \cite{cannizzaro2017malliavin} for gPAM and in \cite{gassiat2020densities} for $\Phi^4_3$.
  \item (Restricted scope.) Strong convergence of the scalar Malliavin derivative under mollification, $\int_0^t \norm{\D_r F^\eps - \D_r F}_U^2\dr \to 0$ in probability, and consequently $\sigma_t^{(\varphi, \eps)} \to \sigma_t^{(\varphi)} := \int_0^t \norm{\D_r F}_U^2\dr$ in probability. (This refines the BCCH convergence $F^\eps \to F$ by additionally controlling the first-order Malliavin object. It is not stated in this integrated form in \cite{cannizzaro2017malliavin, gassiat2020densities}, but is contained in them: \cite{cannizzaro2017malliavin} prove an explicit convergence theorem for the solutions of the renormalised tangent equation for gPAM, and \cite{gassiat2020densities} identify the tangent under the renormalised mollified model with the classical Malliavin tangent of the mollified renormalised equation and prove the tangent map locally Lipschitz jointly in the model and the direction, and linear in the direction, so that convergence, uniformly over Cameron--Martin directions in the unit ball and hence in the integrated form above, follows from the convergence of the renormalised models; cf.\ also \cite{schonbauer2023malliavin}.)
  \item (Restricted scope.) Strict positivity of the scalar Malliavin variance at the limit, $\sigma_t^{(\varphi)} > 0$ almost surely, which by the Bouleau--Hirsch criterion \cite{nualart2006malliavin} implies that $\mathrm{Law}(F)$ is absolutely continuous with respect to Lebesgue measure on $\R$. This is established, for the respective functional classes above, in \cite{schonbauer2023malliavin}, in \cite{cannizzaro2017malliavin} via the strong maximum principle for gPAM, and in \cite{gassiat2020densities} for $\Phi^4_3$.
\end{enumerate}
None of (i)--(iii) is required for the unconditional content of Theorem~\ref{thm:score-gPAM} below, which is established in the entire BCCH scope using only the BCCH convergence $F^\eps \to F$ in probability. Statements (i)--(iii) are nevertheless illuminating context. When (iii) holds, the distributional derivative $\Lambda^{(F)}$ of Theorem~\ref{thm:score-gPAM}(c) is the IBP-distributional derivative of an absolutely continuous probability measure on $\R$, and the question of representing $\Lambda^{(F)}$ as a function (which would yield a Bogachev--Fomin logarithmic derivative in the sense of Definition~\ref{def:log-deriv}) reduces to the question of upgrading first-order to second-order Malliavin smoothness (a distinct open problem). Statements (i)--(iii) are also not enough to obtain the function representation, covering the first-order Malliavin object only, and Sch\"onbauer's results establish density existence but not smoothness of densities for the renormalised solution tested against test functions \cite{schonbauer2023malliavin}.
\end{remark}

\begin{remark}[Fixed-time slices in the BCCH framework]\label{rem:fixed-time-slice}
The convergence result of \cite{brunedchandrachevyrevhairer2021} is formulated through a generalised Da Prato--Debussche decomposition $u = \mathcal{S}^- + \mathcal{S}^+$. The \emph{stationary part} $\mathcal{S}^-$ is a finite sum of explicit renormalised multilinear functionals of the noises, and the mollified stationary parts converge in probability as space-time distributions in Hölder--Besov spaces whose exponents may be negative (e.g., $-\tfrac{1}{2} - \kappa$ for the linear part of $\Phi^4_3$, with corresponding time regularity $-\tfrac{1}{4} - \tfrac{\kappa}{2}$, so that fixed-time slices of $\mathcal{S}^-$ are not directly provided by this topology); the remainder $\mathcal{S}^+$ converges in probability in a topology of maximal solutions whose elements are continuous functions of time valued in a spatial Hölder--Besov space on $\mathbb{T}^d$, up to the random existence time (we refer to the corrected construction of this topology in the arXiv version~4 of \cite{brunedchandrachevyrevhairer2021}, which fixes an error in the published construction). Fixed-time slices of the remainder are therefore supplied by the theorem itself, and defining $u(t)$ reduces to the routine construction of fixed-time versions of the finitely many explicit stationary Gaussian-chaos objects. For equations with positive solution regularity no decomposition is needed. The renormalised solution at fixed time $t > 0$ is automatically a function on $\mathbb{T}^d$ in some Hölder space, the convergence $u^\eps(t) \to u(t)$ in probability follows from the space-time convergence by the embedding of parabolic Hölder spaces into continuous functions of time, and the standing convergence hypothesis $F^\eps \to F$ in probability of Theorem~\ref{thm:score-gPAM} holds with $F = \langle u(t), \varphi\rangle$. This covers gPAM in 2D, multiplicative stochastic heat in $1+1$ dimensions, generalised KPZ in $1+1$ dimensions, the $\Phi^p_d$ equations for $p$ small enough, and many other standard cases. For negative-regularity cases, a solution theory with fixed-time slices can also be constructed directly, without passing through this decomposition; for $\Phi^4_3$, the local-in-time solution theory of \cite{catellierchouk2018paracontrolled} (built over the spatially mollified, white-in-time family; cf.\ the approximation-family discussion above) and the come-down-from-infinity energy bounds of \cite{mourratweber2017comesdown} yield a global dynamics with fixed-time slices in $\mathcal{C}^{-1/2-\kappa}(\mathbb{T}^3)$ for $t > 0$, and analogous fixed-time-slice constructions are known case-by-case for several other negative-regularity standard examples (mass-deformed $\Phi^4_d$, dynamical sine-Gordon, etc.). For all these standard cases, the standing convergence hypothesis of Theorem~\ref{thm:score-gPAM} is verified, and $F = \langle u(t), \varphi\rangle$ has the natural interpretation as the duality pairing between $u(t) \in \mathcal{C}^{\mathrm{reg}(t)}(\mathbb{T}^d)$ and the smooth spatial test function $\varphi$. For BCCH equations where fixed-time slices have not been independently established (no standard example of interest is currently known to lack them), an analogous version of Theorem~\ref{thm:score-gPAM} holds with space-time pairings $F^\eps = \int_{(0,T) \times \mathbb{T}^d} u^\eps(s,x)\,\phi(s,x)\,ds\,dx$ against smooth space-time test functions $\phi \in C_c^\infty((0,T) \times \mathbb{T}^d)$, exactly as in Sch\"onbauer's setup \cite{schonbauer2023malliavin}; the BCCH black-box convergence in space-time topology then directly gives $F^\eps \to F$ in probability without needing fixed-time slices, and the rest of the argument is identical. The fixed-time formulation is preferred here only because it matches the rest of the paper's $X(t)$ framework and makes the bridge identity \eqref{eq:tower-bridge-restated} to Theorem~\ref{thm:main} clean.
\end{remark}

\begin{remark}[Densities, but not necessarily smoothness]\label{rem:density-not-smoothness}
The first-order Malliavin information of Remark~\ref{rem:rs-inputs} establishes existence of densities via Bouleau--Hirsch for the functionals and cases of \cite{cannizzaro2017malliavin, gassiat2020densities, schonbauer2023malliavin} (cf.\ the functional-class scoping in Remark~\ref{rem:rs-inputs}). To pass from existence of the density to existence of a score function for the renormalised limit, two additional ingredients are needed: (a) negative moments
\[
  \E\bigl[(\sigma_t^{(\varphi)})^{-q}\bigr] \;<\; \infty \qquad \text{for some } q \geq 2,
\]
which is the standard ingredient for smoothness of the density (cf.\ Nualart \cite{nualart2006malliavin}), and (b) second-order Malliavin smoothness $F \in \mathbb{D}^{2,p^{*}}$, which is the standard ingredient for the integrand $\D F / \sigma_t^{(\varphi)}$ in the classical scalar Bismut formula to lie in $\mathrm{Dom}(\delta_U)$, so that $\delta_U(\D F / \sigma_t^{(\varphi)})$ is well-defined as an element of $L^2(\Omega)$. The two ingredients enter through the same H\"older triangle as in the operator-level theorem. The Malliavin derivative of the integrand is $\D^2 F / \sigma_t^{(\varphi)} - \D F \otimes \D \sigma_t^{(\varphi)}/(\sigma_t^{(\varphi)})^2$, and by the identity $\norm{\D F}^2_{\HW} = \sigma_t^{(\varphi)}$ its second term is at most twice its first, so that both are governed by the single ratio $\norm{\D^2 F}_{HS}/\sigma_t^{(\varphi)}$ and the requirement is $1/p^{*} + 1/q \leq 1/2$, the relation~\eqref{eq:nondeg-Holder-triangle} at $H = \R$. Neither ingredient is in the regularity-structures Malliavin literature, which establishes only first-order Cameron--Martin Fr\'echet differentiability of the renormalised solution. We therefore treat both as explicit hypotheses in Theorem~\ref{thm:score-gPAM}(b), and obtain the classical Bismut formula for the renormalised limit conditionally on these hypotheses being verified for the equation at hand. The existence of the distributional derivative of the limit law as a continuous functional on $C_b^1(\R)$ does not require either hypothesis and is unconditional, valid in the entire BCCH scope, and is precisely the content of Theorem~\ref{thm:score-gPAM}(c).
\end{remark}
\subsubsection*{Distributional derivative and conditional logarithmic derivative for renormalised scalar functionals}

\begin{theorem}[The scalar law at the renormalised limit]\label{thm:score-gPAM}
Let $u$ be the renormalised solution of any subcritical singular SPDE in the framework of \cite{brunedchandrachevyrevhairer2021}, with regularised solutions $u^\eps$, along an approximating family as in the discussion above, converging to $u$ in probability on $\{T < \tau\}$. Fix $t \in (0, T)$, let $\varphi \in C^\infty(\mathbb{T}^d)$, and set $F^\eps := \langle u^\eps(t), \varphi\rangle_{L^2(\mathbb{T}^d)}$.

Assume throughout that $F^\eps$ converges in probability to a real-valued random variable $F$ as $\eps \downarrow 0$, as it does in the standard cases of the framework (Remark~\ref{rem:fixed-time-slice}). Part~\textup{(c)} uses nothing beyond this; parts~\textup{(a)} and~\textup{(b)} carry the further Malliavin hypotheses stated with them.

\begin{enumerate}[leftmargin=2.4em, label=\textup{(\alph*)}, ref=\textup{(\alph*)}]
  \item[\textup{(c)}] (Distributional derivative at the renormalised limit; unconditional.) Under the standing hypothesis alone, the linear functional
  \[
    \Lambda^{(F)}: C_b^1(\R) \to \R, \qquad \Lambda^{(F)}(\psi) \;:=\; -\,\E[\psi'(F)]
  \]
  is the distributional derivative of $\mathrm{Law}(F)$ on $\R$ in the integration-by-parts sense; the corresponding distributional derivatives of the mollified laws converge weak-$*$ in $(C_b^1(\R))^*$:
  \begin{equation}\label{eq:weak-conv-singular}
  \begin{aligned}
    \Lambda^{(F^\eps)}(\psi) \;=\; -\,\E[\psi'(F^\eps)] &\;\longrightarrow\; -\,\E[\psi'(F)] \;=\; \Lambda^{(F)}(\psi) \\
    &\qquad \text{for every } \psi \in C_b^1(\R), \quad \text{as } \eps \downarrow 0.
  \end{aligned}
  \end{equation}

  \medskip\noindent
  \item[\textup{(a)}] (Mollified score at fixed $\eps$; classical Bismut formula.) Suppose, in addition, that at each fixed $\eps > 0$ the regularised equation supports the classical Malliavin framework (in the white-in-time or abstract Gaussian-space form of Remark~\ref{rem:abstract-gaussian}), i.e., $F^\eps \in \mathbb{D}^{\infty,p}(\Omega)$ for all $p < \infty$, $\sigma_t^{(\varphi,\eps)} > 0$ a.s., and $(\sigma_t^{(\varphi,\eps)})^{-1} \in L^q(\Omega)$ for all $q < \infty$. Then $F^\eps$ admits the intrinsic scalar logarithmic derivative
  \begin{equation}\label{eq:scalar-Bismut-eps}
    \beta^{(F^\eps),\mathrm{int}}(F^\eps) \;=\; -\E\!\left[\delta_U\!\left(\frac{\D F^\eps}{\sigma_t^{(\varphi,\eps)}}\right) \;\bigg|\; F^\eps\right] \quad \text{in } L^2(\R, \mathrm{Law}(F^\eps)).
  \end{equation}
  In the variational subcase, for any deterministic direction $h \in H$ for which the hypotheses of Theorem~\ref{thm:main}, Part~II, hold for the regularised equation and for which $\ip{h}{\varphi}_H \neq 0$, the pushforward identity
  \begin{equation}\label{eq:tower-bridge-restated}
    \ip{h}{\varphi}_H \cdot \beta^{(F^\eps),\mathrm{int}}(F^\eps) \;=\; -\,\E\bigl[\delta_U(v_h^{H,\eps}) \,\big|\, F^\eps\bigr]
  \end{equation}
  expresses the scalar logarithmic derivative in terms of the operator-level Bismut--Nualart formula of Theorem~\ref{thm:main} applied to the mollified variational equation, with $v_h^{H,\eps}$ the $H$-valued covering field of Section~\ref{sec:proof-main}.

  \item[\textup{(b)}] (Score-as-function at the renormalised limit.) Suppose, in addition, that $F$ satisfies clauses~\textup{(i)}--\textup{(ii)} of Assumption~\ref{ass:nondeg} read on the one-dimensional state space $H = \R$ at the unit direction, that is $\sigma_t^{(\varphi)} > 0$ almost surely and, for a pair of exponents $(q, p^{*}) \in [2,\infty]^2$ on the H\"older triangle~\eqref{eq:nondeg-Holder-triangle},
  \begin{equation}\label{eq:scalar-nondeg}
    \E\bigl[(\sigma_t^{(\varphi)})^{-q}\bigr] \;<\; \infty, \qquad F \in \mathbb{D}^{2,p^{*}}(\Omega), \qquad \frac{1}{q} + \frac{1}{p^{*}} \;\leq\; \frac{1}{2},
  \end{equation}
  or, more generally, the intrinsic integrability condition
  \begin{equation}\label{eq:scalar-nondeg-intrinsic}
    \norm{\D^2 F}_{HS}\,\bigl(\sigma_t^{(\varphi)}\bigr)^{-1} \;\in\; L^2(\Omega).
  \end{equation}
  Then the integrand $u_F := \D F / \sigma_t^{(\varphi)}$ lies in $\mathbb{D}^{1,2}(\HW) \subset \Dom(\delta_U)$, and the intrinsic scalar logarithmic derivative of $\mathrm{Law}(F)$ exists in $L^2(\mathrm{Law}(F))$ as a Bogachev--Fomin logarithmic derivative (Definition~\ref{def:log-deriv}), given by
  \begin{equation}\label{eq:scalar-Bismut-limit}
    \beta^{(F),\mathrm{int}}(F) \;=\; -\E\!\left[\delta_U\!\left(\frac{\D F}{\sigma_t^{(\varphi)}}\right) \;\bigg|\; F\right] \quad \text{in } L^2(\R, \mathrm{Law}(F)).
  \end{equation}
\end{enumerate}
\end{theorem}

\begin{remark}[The scalar reduction on the H\"older triangle]\label{rem:score-gPAM-anatomy}
The functional $\Lambda^{(F)}(\psi) = -\E[\psi'(F)]$ in part~(c) is a continuous linear functional on $C_b^1(\R)$ for any probability law on $\R$, requiring no Malliavin regularity of $F$. We reserve the term ``logarithmic derivative'', following Bogachev~\cite{bogachev2010differentiable}, for the Radon--Nikodym density $\beta_h \in L^1(\mu)$ representing this functional via $\Lambda^{(F)}(\psi) = \E[\psi(F)\beta(F)]$; this exists only when the measure is Fomin-differentiable. The weak-$*$ convergence~\eqref{eq:weak-conv-singular} is a trivial consequence of $\mathrm{Law}(F^\eps) \to \mathrm{Law}(F)$ via bounded convergence and is not a substantive convergence-of-scores claim; score convergence under weak convergence of measures generally fails.

The hypothesis of part~(b) is Assumption~\ref{ass:nondeg} specialised to $H = \R$, and the specialisation is exact. At $H = \R$ the Malliavin matrix is $\mathcal{T} = \D F \in \HW$, the covariance is the positive scalar $\gamma_F = \norm{\D F}^2_{\HW} = \sigma_t^{(\varphi)}$, the pseudoinverse applied to the unit direction is $\gamma_F^{\dagger}h = (\sigma_t^{(\varphi)})^{-1}$, and the requirement $\mathcal{T} \in \mathbb{D}^{1,p^{*}}$ reads $F \in \mathbb{D}^{2,p^{*}}$. Clause~\textup{(i)} is the strict positivity of the variance and clause~\textup{(ii)} is~\eqref{eq:scalar-nondeg}; clause~\textup{(iii)} is not needed, because a positive scalar covariance requires no Tikhonov regularisation and the covering field $u_F = \D F/\sigma_t^{(\varphi)}$ is directly $\mathbb{D}^{1,2}$. The singular half of the paper therefore operates on the same H\"older triangle~\eqref{eq:nondeg-Holder-triangle} as the variational half.

That the triangle is the right condition, rather than something stronger, is a consequence of the identity $\sigma_t^{(\varphi)} = \norm{\D F}^2_{\HW}$ itself. The Malliavin derivative of the integrand is
\begin{equation}\label{eq:Du-quotient}
  \D_s u_F(r) \;=\; \frac{\D^2_{s, r} F}{\sigma_t^{(\varphi)}} \;-\; \frac{(\D_r F)\,\otimes\,(\D_s\sigma_t^{(\varphi)})}{(\sigma_t^{(\varphi)})^2}.
\end{equation}
and the second term looks like the more demanding one, carrying two powers of the variance in the denominator. It is not. By the Cauchy--Schwarz bound $\norm{\D\sigma_t^{(\varphi)}}_{\HW} \leq 2\norm{\D F}_{\HW}\norm{\D^2 F}_{HS}$ of~\eqref{eq:Dsigma-bound} and the identity $\norm{\D F}^2_{\HW} = \sigma_t^{(\varphi)}$, the cross term obeys
\[
  \frac{\norm{\D F}_{\HW}\,\norm{\D\sigma_t^{(\varphi)}}_{\HW}}{(\sigma_t^{(\varphi)})^2} \;\leq\; \frac{2\,\norm{\D F}^2_{\HW}\,\norm{\D^2 F}_{HS}}{(\sigma_t^{(\varphi)})^2} \;=\; \frac{2\,\norm{\D^2 F}_{HS}}{\sigma_t^{(\varphi)}},
\]
which is twice the first term. The two terms therefore carry the same integrability requirement, a single one at the ratio $\norm{\D^2 F}_{HS}/\sigma_t^{(\varphi)}$, and H\"older at the pair $(p^{*}, q)$ places it in $L^2(\Omega)$ under the triangle relation of~\eqref{eq:scalar-nondeg}. At the symmetric point $p^{*} = q = 4$ this is $F \in \mathbb{D}^{2,4}$ together with $(\sigma_t^{(\varphi)})^{-1} \in L^4(\Omega)$. When both~(a) and~(b) hold, the distributional and Radon--Nikodym forms of the derivative coincide, $\Lambda^{(F)}(\psi) = \E[\psi(F)\,\beta^{(F),\mathrm{int}}(F)]$ via~\eqref{eq:scalar-Bismut-limit}. Both conditions in~\eqref{eq:scalar-nondeg} remain open in the regularity-structures Malliavin literature (Remark~\ref{rem:density-not-smoothness}).
\end{remark}

\begin{proof}
At each fixed $\eps > 0$, the regularised equation is a classical equation with smooth noise and smooth coefficients (the smoothness of the renormalisation counterterms $\Upsilon_i^k$ being a consequence of the smoothness of the original $F_i$ and $F_i^j$), an It\^o SPDE with smooth Hilbert--Schmidt noise for the white-in-time regularisations, and a smooth-coefficient random PDE analysed on the Gaussian space of the driving fields (Remark~\ref{rem:abstract-gaussian}) for time-independent or space-time-mollified noise. By the additional hypothesis in part~\textup{(a)} (the classical Malliavin framework applies, in the white-in-time or abstract Gaussian-space form, Nualart \cite{nualart2006malliavin}), the regularised solution $u^\eps(t)$ is in $\mathbb{D}^{k,p}$ for all $k, p$, and the scalar Malliavin variance $\sigma_t^{(\varphi, \eps)}$ is strictly positive a.s.\ with finite negative moments of all orders. The classical scalar Bismut--Nualart integration-by-parts formula \cite{nualart2006malliavin} therefore applies to $F^\eps = \langle u^\eps(t), \varphi\rangle$ and gives \eqref{eq:scalar-Bismut-eps}.

Equivalently, this is the scalar instance of the abstract Theorem~\ref{thm:abstract-bismut-fomin} applied with $F = F^\eps$ on $H = \R$, $\mathcal{T} = DF^\eps \in \HW$, and the canonical covering field $u_h := \mathcal{T}^{*}h/\sigma_t^{(\varphi, \eps)} \in \HW$ (with $h = 1 \in \R$ the unit Fomin direction); the Tikhonov and scalar-trace-convergence mechanism of Theorem~\ref{thm:tikhonov-trace} simplifies in the scalar case ($H = \R$, so $\gamma_F^\eps = \sigma_t^{(\varphi,\eps)}$ is a positive scalar with $\gamma_F^{\eps,\dagger} = (\sigma_t^{(\varphi,\eps)})^{-1}$) to the classical division-by-variance form.

For the bridge identity \eqref{eq:tower-bridge-restated} in the variational subcase, observe that Theorem~\ref{thm:main} Part~II applied to the mollified variational equation furnishes the operator-level logarithmic derivative $\beta_h^{(\mu_t^\eps)}: H \to \R$ for any deterministic $h \in H$ verifying Assumption~\ref{ass:nondeg} for $\gamma_t^\eps$, in the form $\beta_h^{(\mu_t^\eps)}(u^\eps(t)) = -\E[\delta_U(v_h^{H,\eps}) \,|\, u^\eps(t)]$. The natural candidate for the bridge to the scalar pairing $\ell(\cdot) = \ip{\cdot}{\varphi}_H$ is the deterministic direction $h := \varphi$ itself (which is in $H = L^2(\mathbb{T}^d)$ since $\varphi \in C^\infty(\mathbb{T}^d) \subset L^2(\mathbb{T}^d)$ on the compact torus), giving $\ip{h}{\varphi}_H = \norm{\varphi}_H^2 > 0$; for it, Assumption~\ref{ass:nondeg} requires $\varphi \in \Ran(\gamma_t^\eps)$ almost surely together with $\E[\norm{\gamma_t^{\eps,\dagger}\varphi}_H^q] < \infty$ and the Cameron--Martin compatibility, which at fixed $\eps > 0$ is an equation-specific verification and is assumed here rather than proved. The pushforward of $\mu_t^\eps$ under the continuous linear functional $\ell$ is the law of $F^\eps$ on $\R$, and by the standard pushforward rule for Fomin derivatives under continuous linear maps \cite{bogachev2010differentiable} (recalled in the discussion preceding \eqref{eq:tower-bridge}), $\ell_*\mu_t^\eps$ is Fomin-differentiable along $\ell(\varphi) = \norm{\varphi}_H^2 \in \R$ with score
\[
  \beta_{\norm{\varphi}_H^2}^{(F^\eps)}(F^\eps) \;=\; \E\bigl[\beta_\varphi^{(\mu_t^\eps)}(u^\eps(t)) \,\big|\, F^\eps\bigr] \;=\; -\,\E\bigl[\delta_U(v_\varphi^{H,\eps}) \,\big|\, F^\eps\bigr],
\]
where the second equality uses the tower property of conditional expectation. By linearity of $\beta_\eta^{(F^\eps)}(\cdot)$ in $\eta \in \R$ (the intrinsic scalar logarithmic derivative satisfies $\beta_\eta^{(F^\eps)}(\cdot) = \eta\,\beta^{(F^\eps),\mathrm{int}}(\cdot)$ for $\eta \in \R$, since in $\R$ the Fomin derivative along scaled directions is itself scaled), the LHS equals $\norm{\varphi}_H^2 \cdot \beta^{(F^\eps),\mathrm{int}}(F^\eps)$. Combining gives the bridge identity in the form
\begin{equation*}
  \norm{\varphi}_H^2 \cdot \beta^{(F^\eps),\mathrm{int}}(F^\eps) \;=\; -\,\E\bigl[\delta_U(v_\varphi^{H,\eps}) \,\big|\, F^\eps\bigr],
\end{equation*}
which is the special case of \eqref{eq:tower-bridge-restated} with $h = \varphi$. More generally, for any deterministic $h \in H$ with $\ip{h}{\varphi}_H \neq 0$ and $h$ verifying Assumption~\ref{ass:nondeg}, the same argument gives \eqref{eq:tower-bridge-restated} as stated.

Under the hypothesis \eqref{eq:scalar-nondeg} (or its intrinsic form \eqref{eq:scalar-nondeg-intrinsic}), we show that the integrand $u_F := \D F / \sigma_t^{(\varphi)}$ is in $\mathbb{D}^{1,2}(\HW)$ as an $\HW$-valued random variable. The quotient rule is derived in two steps, of which the first is that the scalar Malliavin chain rule~\cite{nualart2006malliavin} applied to the smooth scalar function $\psi(y) = 1/y$ on $(0, \infty)$, localised on the full-measure event $\{\sigma_t^{(\varphi)} > 0\}$ (the range condition of part~(b)) with $\psi$ extended smoothly across $\{y = 0\}$ (the resulting identity is independent of the extension by the negative-moment control on $1/\sigma_t^{(\varphi)}$), gives
\[
  \D_s\!\left(\tfrac{1}{\sigma_t^{(\varphi)}}\right) \;=\; -\,\frac{\D_s\sigma_t^{(\varphi)}}{(\sigma_t^{(\varphi)})^2}, \qquad s \in [0, t];
\]
(ii) the operator-valued Leibniz rule for Malliavin derivatives applied to the product $u_F(r) = (\D_r F) \cdot (\sigma_t^{(\varphi)})^{-1}$ of the $U$-valued $\D_r F \in \mathbb{D}^{1,2}(U)$ and the scalar $(\sigma_t^{(\varphi)})^{-1} \in \mathbb{D}^{1,2}(\R)$ (which is the standard Leibniz rule of \cite{nualart2006malliavin} applied to a Hilbert-valued and scalar factor, the scalar/Hilbert version of the basis-wise expansion explained around Lemma~\ref{lemma:HS-kernel}) yields
\[
  (\D_s u_F)(r) \;=\; \frac{\D^2_{s, r} F}{\sigma_t^{(\varphi)}} \;-\; \frac{(\D_r F)\,\otimes\,(\D_s\sigma_t^{(\varphi)})}{(\sigma_t^{(\varphi)})^2}, \qquad s, r \in [0, t],
\]
with $(\D_s u_F)(r)$ taking values in $U \otimes U$, $\D^2 F \in L^2(\Omega; \HW^{\otimes 2})$ (the scalar second-order Malliavin derivative, with $(\D^2_{s,r} F) \in U\otimes U$ at each $(s,r)$; the membership follows from $F \in \mathbb{D}^{2,p^{*}}$ via the inclusion $\mathbb{D}^{2,p^{*}} \hookrightarrow \mathbb{D}^{2,2}$ for $p^{*} \geq 2$), and $\D\sigma_t^{(\varphi)} \in \HW$ given by $\D_s\sigma_t^{(\varphi)} = 2\int_0^t\langle\D_r F,\,\D^2_{s,r} F\rangle_U\,dr \in U$ (with the $\langle\cdot,\cdot\rangle_U$ contracting $\D_r F \in U$ against the second $U$-slot of $\D^2_{s,r} F \in U \otimes U$, leaving a $U$-valued result at each $s$).

Pointwise in $s \in [0, t]$,
\begin{equation}\label{eq:Dsigma-bound}
  \norm{\D_s\sigma_t^{(\varphi)}}_U \;\leq\; 2\int_0^t\norm{\D_r F}_U\,\norm{\D^2_{s,r} F}_{U\otimes U}\,dr,
\end{equation}
and integrating in $s$ by Cauchy--Schwarz with respect to the $\HW$-norm,
\begin{equation}\label{eq:Dsigma-HW-bound}
  \norm{\D\sigma_t^{(\varphi)}}_{\HW} \;=\; \Bigl(\textstyle\int_0^t\norm{\D_s\sigma_t^{(\varphi)}}_U^2\,ds\Bigr)^{1/2} \;\leq\; 2\,\norm{\D F}_{\HW}\,\norm{\D^2 F}_{HS},
\end{equation}
where $\norm{\D^2 F}_{HS} := (\int_0^t\!\int_0^t \norm{\D^2_{s,r} F}_{U\otimes U}^2\,ds\,dr)^{1/2}$ is the Hilbert--Schmidt norm of $\D^2 F$ on $\HW \otimes \HW$. Substituting~\eqref{eq:Dsigma-HW-bound} into the second term of the quotient rule and using the identity $\norm{\D F}^2_{\HW} = \sigma_t^{(\varphi)}$, which is the definition of the scalar Malliavin variance, the two terms collapse onto one another,
\begin{equation}\label{eq:cross-term-envelope}
\begin{aligned}
  \Bigl\lVert\frac{(\D F)\otimes(\D\sigma_t^{(\varphi)})}{(\sigma_t^{(\varphi)})^2}\Bigr\rVert_{\HW^{\otimes 2}}
  &\;=\; \frac{\norm{\D F}_{\HW}\,\norm{\D\sigma_t^{(\varphi)}}_{\HW}}{(\sigma_t^{(\varphi)})^2}\\
  &\;\leq\; \frac{2\,\norm{\D F}_{\HW}^2\,\norm{\D^2 F}_{HS}}{(\sigma_t^{(\varphi)})^2}
  \;=\; \frac{2\,\norm{\D^2 F}_{HS}}{\sigma_t^{(\varphi)}},
\end{aligned}
\end{equation}
so that the whole derivative is dominated by a single quantity,
\begin{equation}\label{eq:Du-envelope}
  \norm{\D u_F}_{\HW^{\otimes 2}} \;\leq\; 3\,\frac{\norm{\D^2 F}_{HS}}{\sigma_t^{(\varphi)}} .
\end{equation}
This is the intrinsic hypothesis~\eqref{eq:scalar-nondeg-intrinsic}, and it places $\D u_F$ in $L^2(\Omega; \HW^{\otimes 2})$ at once. Under the H\"older form~\eqref{eq:scalar-nondeg}, the pair $(\norm{\D^2 F}_{HS}, (\sigma_t^{(\varphi)})^{-1}) \in L^{p^{*}} \times L^q$ and H\"older's inequality give
\begin{equation}\label{eq:cross-term-Holder}
  \E\Bigl[\Bigl(\frac{\norm{\D^2 F}_{HS}}{\sigma_t^{(\varphi)}}\Bigr)^{\!2}\Bigr]^{1/2} \;\leq\; \bigl(\E[\norm{\D^2 F}_{HS}^{p^{*}}]\bigr)^{1/p^{*}}\bigl(\E[(\sigma_t^{(\varphi)})^{-q}]\bigr)^{1/q} \;<\; \infty
\end{equation}
precisely when the two reciprocals obey
\begin{equation}\label{eq:cross-term-recip}
  \frac{1}{p^{*}} + \frac{1}{q} \;\leq\; \frac{1}{2},
\end{equation}
which is the H\"older triangle~\eqref{eq:nondeg-Holder-triangle} of Assumption~\ref{ass:nondeg}\textup{(ii)} at $H = \R$. The field itself is in $L^2(\Omega;\HW)$ by the same identity, since $\norm{u_F}^2_{\HW} = \norm{\D F}^2_{\HW}/(\sigma_t^{(\varphi)})^2 = (\sigma_t^{(\varphi)})^{-1}$, which is integrable because $q \geq 2 \geq 1$; this is the scalar instance of the bound $\E\norm{u_h}^2_{\HW} \leq \norm{h}_H\E[\norm{\gamma_F^{\dagger}h}_H]$ of Theorem~\ref{thm:pseudoinverse-covering}.

Hence $u_F \in \mathbb{D}^{1,2}(\HW) \subset \Dom(\delta_U)$ by \cite{nualart2006malliavin}, and $\delta_U(u_F) \in L^2(\Omega)$. The classical scalar Bismut formula then yields, for every $\psi \in C_b^1(\R)$,
\[
  \E[\psi'(F)] \;=\; \E\!\left[\psi(F)\,\delta_U\!\left(\frac{\D F}{\sigma_t^{(\varphi)}}\right)\right],
\]
which by Definition~\ref{def:log-deriv} characterises the conditional expectation $-\E[\delta_U(\D F / \sigma_t^{(\varphi)}) \,|\, F]$ as the intrinsic scalar logarithmic derivative $\beta^{(F),\mathrm{int}}(F)$, giving \eqref{eq:scalar-Bismut-limit}.

By hypothesis, $F^\eps \to F$ in probability. For any $\psi \in C_b^1(\R)$, the function $\psi'$ is bounded continuous, so bounded convergence gives
\begin{equation}\label{eq:bdd-conv-step}
  \E[\psi'(F^\eps)] \;\longrightarrow\; \E[\psi'(F)] \qquad \text{as } \eps \downarrow 0.
\end{equation}
Define $\Lambda^{(F)}(\psi) := -\E[\psi'(F)]$, a continuous linear functional on $C_b^1(\R)$ with $\norm{\Lambda^{(F)}}_{(C_b^1)^*} \leq 1$. Then
\begin{equation}\label{eq:Lambda-pointwise-conv}
  \Lambda^{(F^\eps)}(\psi) \;\longrightarrow\; \Lambda^{(F)}(\psi) \qquad \text{for every } \psi \in C_b^1(\R),
\end{equation}
which is~\eqref{eq:weak-conv-singular}. The limit $\Lambda^{(F)}$ is the distributional derivative of $\mathrm{Law}(F)$, established using only bounded convergence and the standing convergence $F^\eps \to F$ — no Malliavin regularity of $F$ is required.

If part~(a) additionally holds at each $\eps > 0$, the IBP identity~\eqref{eq:intrinsic-IBP} via~\eqref{eq:scalar-Bismut-eps} represents
\begin{align*}
  \Lambda^{(F^\eps)}(\psi) &\;=\; \int_\R \psi(y)\,\beta^{(F^\eps),\mathrm{int}}(y)\,d\mathrm{Law}(F^\eps)(y),
\end{align*}
so the score measures $\nu_\beta^{(F^\eps)} := \beta^{(F^\eps),\mathrm{int}}\,\mathrm{Law}(F^\eps)$ converge weak-$*$ to $\Lambda^{(F)}$ in $(C_b^1(\R))^*$. When part~(b) holds, $\Lambda^{(F)}$ is represented by $\beta^{(F),\mathrm{int}} \in L^1(\mathrm{Law}(F))$ via~\eqref{eq:scalar-Bismut-limit}; otherwise it remains distributional.
\end{proof}

The operator-level machinery of Sections~\ref{sec:proof-main}--\ref{sec:variation} does not survive the renormalised limit, since neither $\gamma_t$ as a trace-class operator on $L^2(\mathbb{T}^d)$ nor the identification $\Phi_r = Y(t,r)\calB(r,X(r))$ is available. The scalar reduction replaces operator-level non-degeneracy with positivity of $\sigma_t^{(\varphi)}$. The unconditional content (part~(c)) uses only the standing convergence $u^\eps \to u$, supplied along the BCCH approximating families by \cite{brunedchandrachevyrevhairer2021}; the function-valued representation (part~(b)) requires the additional Malliavin hypotheses, which remain open.

\begin{remark}[Scope]\label{rem:scope}
Part~(c) of Theorem~\ref{thm:score-gPAM} holds throughout the BCCH scope \cite{brunedchandrachevyrevhairer2021}, covering the generalised parabolic Anderson model, the dynamical $\Phi^4_d$ and $\Phi^p_d$ equations, the multiplicative stochastic heat equation, the generalised KPZ equation \cite{schonbauer2023malliavin}, and all other subcritical singular SPDEs satisfying the assumptions of that theorem, with any $\varphi \in C^\infty(\mathbb{T}^d)$. The variational subcase (the gPAM companion model, $\Phi^4_d$, $\Phi^p_d$, multiplicative stochastic heat) additionally admits the operator-level bridge~\eqref{eq:tower-bridge-restated}; the generalised KPZ equation falls outside this subcase but is covered by part~(c) directly.
\end{remark}

\begin{remark}[Pointwise evaluation]\label{rem:pointwise-gPAM}
The bridge identity~\eqref{eq:tower-bridge-restated} requires $\varphi \in H$, so it does not directly cover pointwise evaluation $F = u(t,x)$ (for which $u \mapsto u(x)$ is not continuous on $H$). This case is recovered directly. At each fixed $\eps > 0$, whenever the mollified point evaluation $u^\eps(t,x)$ satisfies the Malliavin smoothness and scalar non-degeneracy hypotheses of Theorem~\ref{thm:score-gPAM}(a), the same scalar Bismut formula~\eqref{eq:scalar-Bismut-eps} applies verbatim, and the renormalised limit proceeds as in Theorem~\ref{thm:score-gPAM}(c) using $u^\eps(t,x) \to u(t,x)$ in probability. Pointwise evaluation in function-valued cases (gPAM, multiplicative stochastic heat in $1+1$ dimensions) is therefore not an exception but simply bypasses the operator-level bridge.
\end{remark}

The regularity-structure framework enters only through the BCCH convergence $u^\eps \to u$ in probability, which is the sole input for part~(c). Parts~(a) and~(b) use classical scalar Malliavin calculus at fixed $\eps$ and additional smoothness hypotheses at the limit, respectively; neither requires the regularity-structures Malliavin literature \cite{cannizzaro2017malliavin, gassiat2020densities, schonbauer2023malliavin}.

\subsection{The formula tested numerically}\label{subsec:numerics}

We now subject the formula of Theorem~\ref{thm:main} to direct numerical test in the additive-noise regime, on spectral discretisations of the equations to which it applies and, for the three singularly motivated rows, on finite-$K$ mollified surrogates, where the drift nonlinearity $\calA'_u \ne -A$ is active and is the sole source of the correction term $C_h$. Two features of the formula make the experiment of independent mathematical interest. First, no closed-form expression for $\beta_h$ exists once $\calA'_u \ne -A$. Outside the Gaussian and linear-additive regimes \cite{mirafzali2025infinite}, the Skorokhod integral $\delta_U(v_h)$ is a path-dependent random variable, and $\beta_h$ is its conditional expectation $-\E[\delta_U(v_h)\mid X(t)]$ across realisations. Second, the marginal $\mu_t$ on $H$ is infinite-dimensional, and kernel-density methods together with their numerical derivatives, which give meaning to ``$\partial_x \log p$'' in $\R^m$, are not available on $H$ except after projection to a scalar pairing. The Bismut formula \eqref{eq:skorokhod-decomp}--\eqref{eq:correction-symbolic} therefore provides what the classical numerical apparatus does not extract from samples of $X(t)$ alone, namely a closed-form per-path expression $\delta_U(v_h)(\omega)$ whose conditional expectation across realisations is $\beta_h$. The experiments of this subsection evaluate the formula in this regime and check it against four independent tests, each isolating a distinct failure mode, together with ablation, direction, and refinement studies that quantify the power of the tests rather than assume it. The complementary state-dependent-diffusion regime $\calB'_u \ne 0$, admitted by the theory through Assumption~\ref{ass:LR}, is not exercised numerically here; we comment on it where the one multiplicative example (gPAM) enters.

The experimental scope is the principal examples of Section~\ref{subsec:verification} (the stochastic $p$-Laplacian in 1D and 2D, and the 2D Navier--Stokes equation in vorticity form), augmented by three semilinear cases of Section~\ref{subsec:reduction} (Allen--Cahn 2D, Cahn--Hilliard 2D, and $\Phi^4_2$) and three singular cases of Section~\ref{subsec:singular} tested at fixed mollification (gPAM 2D, $\Phi^4_3$ in 1D radial form, and KPZ in $1+1$D). The estimator is described in \S\ref{subsubsec:estimator}, the validation framework in \S\ref{subsubsec:validation}, the test cases in \S\ref{subsubsec:suite}, and the results in \S\ref{subsubsec:results}. The full source code, including all parameter choices and seeds, is included with the submission.

\subsubsection{The estimator}\label{subsubsec:estimator}

Every quantity in the formula must first be turned into something a computer can evaluate. In space we work in spectral Galerkin form on the basis $(e_k)_{k}$ of Fourier modes on $\mathbb{T}^d$ (or sine modes on the appropriate domain), retaining the lowest $K^d$ modes in dimension $d$. The cylindrical noise on $U = L^2$ is realised mode-by-mode as $W = \sum_k W_k\,e_k$ with $(W_k)_k$ independent scalar Brownian motions. Throughout the numerical suite the diffusion coefficient is taken state-independent and additive, $\calB(r, u) = q^{1/2}$, diagonal in this basis with eigenvalues $q_k = (1 + m^2 + \lambda_k)^{-s}$ for a mass parameter $m^2 > 0$ and decay exponent $s > d/2$; the shift by $1 + m^2$ keeps $q_k$ finite on the zero mode. The condition $s > d/2$ gives $\sum_k q_k < \infty$, so the covariance $Q = q$ is trace-class and $\calB = q^{1/2} \in L_2(U, H)$ is Hilbert--Schmidt, the setting in which Assumption~\ref{ass:LR} holds throughout the suite. Restricting the experiments to additive noise is a deliberate simplification that isolates the drift nonlinearity $\calA'_u \ne -A$, the source of the correction term $C_h$ and the principal novelty of Theorem~\ref{thm:main}, from the state-dependent-diffusion mechanism $\calB'_u \ne 0$, whose discrete second-fibre implementation would introduce a separate validation burden. The one genuinely multiplicative member of the scope, gPAM, is therefore tested in a frozen-coefficient additive surrogate; this is made precise in \S\ref{subsubsec:suite}. Nonlinear products are de-aliased by the standard $2/3$-rule.

In time, on a uniform grid $\{t_i = i\Delta t\}_{i=0}^{n}$ with $\Delta t = t/n$, the SPDE is integrated by an exponential integrator \cite{lord2014introduction},
\begin{equation}\label{eq:etd-Lawson-discrete}
  c^{(i+1)} \;=\; E\,c^{(i)} \;-\; \varphi_1(-L\Delta t)\,\Delta t\,N(c^{(i)}) \;+\; E\,q^{1/2}\,\Delta W^{(i)},
\end{equation}
with $E := e^{-L\Delta t}$, $\varphi_1(z) := (e^z - 1)/z$, $-Lu := \Delta u$ for the heat-type linear part (replaced by $-\Delta^2 u$ for Cahn--Hilliard and by the Stokes operator for the vorticity Navier--Stokes equation), and $N(u) := \calA(u) - Lu$ the nonlinear residual; the linear semigroup is integrated exactly. For the semilinear, additive-noise members of the suite the scheme attains strong order $\tfrac{1}{2}$ in $\Delta t$ under the trace-class noise and globally Lipschitz hypotheses of \cite{lord2014introduction}. For the quasilinear and superlinearly-growing members ($p$-Laplacian, Navier--Stokes, and the cubic and quadratic reaction terms), the nonlinearity is not globally Lipschitz and this rate is not guaranteed by the semilinear theory; indeed full-discrete exponential and linear-implicit Euler schemes diverge strongly and numerically weakly already for stochastic Allen--Cahn equations \cite{beccari2019divergence}, so we treat $\Delta t$ as a controlled discretisation parameter and confirm the empirical convergence rate directly in test~(V4) rather than invoking an a priori order.

\medskip\noindent With the discretisation fixed, each ingredient of the formula \eqref{eq:skorokhod-decomp}--\eqref{eq:correction-symbolic} is implemented in discrete-exact form along a single realisation $\{c^{(i)}\}_{i=0}^{n}$, that is, as the differential of the discrete map \eqref{eq:etd-Lawson-discrete} rather than as a discretisation of the continuous formula. The first variation $Y(t, r)v$ comes from the discrete linearisation of \eqref{eq:etd-Lawson-discrete}, namely
\begin{equation*}
  \delta c^{(i+1)} \;=\; E\,\delta c^{(i)} \;-\; \varphi_1(-L\Delta t)\,\Delta t\,\bigl(\calA'_u(c^{(i)}) - L\bigr)\,\delta c^{(i)};
\end{equation*}
the backward adjoint $w_r := Y(t, r)^*\,\tilde h$ comes from the discrete-adjoint version of the same recursion. The Malliavin covariance applied to a vector, $h \mapsto \gamma_t h = \int_0^t \Phi_r\Phi_r^* h\dr$, is taken in Duhamel form and used as the matvec inside a conjugate-gradient solve $\gamma_t \tilde h = h$ with relative residual tolerance $\le 10^{-4}$ and Tikhonov regularisation $\ge 10^{-9}$ (with an equation-specific floor of $10^{-6}$ for $p$-Laplacian 1D, where $\gamma_t$ is most ill-conditioned). Since $w_r = Y(t,r)^{*}\tilde h$ anticipates the driving path whenever the drift is nonlinear, the forward sum $\sum_{i}\langle q^{1/2}w_{t_i}, \Delta W^{(i)}\rangle_U$ is not by itself a discrete divergence; the estimator used throughout is the full discrete Skorokhod divergence
\begin{equation}\label{eq:discrete-skorokhod}
  M_h^{\Delta} \;=\; \sum_{i=0}^{n-1}\bigl\langle q^{1/2}\,w_{t_i},\,\Delta W^{(i)}\bigr\rangle_U \;-\; \Delta t\sum_{i=0}^{n-1}\Tr_U\bigl(q^{1/2}\,\mathrm{D}_{\Delta W^{(i)}}\,w_{t_i}\bigr),
\end{equation}
in which the integrand being differentiated is the full covering field $q^{1/2}w_{t_i}$, not $w_{t_i}$ alone; for deterministic $q$ the factor pulls out of the derivative, $\mathrm{D}_{\Delta W^{(i)}}(q^{1/2}w_{t_i}) = q^{1/2}\mathrm{D}_{\Delta W^{(i)}}w_{t_i}$, and for a general state-dependent coefficient $q^{1/2}$ is replaced by $\calB(t_i, X_{t_i})^{*}$ throughout. The anticipating-diagonal term is assembled from the same discrete backward-adjoint pass that produces $w$, by automatic differentiation of $w_{t_i}$ in the increment $\Delta W^{(i)}$, and traced by the same probe estimator used for $C_h$. The diagonal term belongs to the main frozen-$z$ divergence and is distinct from $C_h$, which corrects substitution of the random value $z = \gamma_t^{\dagger}h$. Its measured magnitude at every operating point of Table~\ref{tab:numerics-main} lies below one paired standard error of the reported value (at the per-mille level for gPAM), so it is resolved but not dominant in the reported comparisons. Finally, the discrete realisation of $[\D_r\gamma_t](\xi)$ comes from the second-fibre construction of Proposition~\ref{prop:D-gamma}, assembled on a subgrid $\{r_l\}_{l=1}^{n_{\mathrm{sub}}}\subset\{t_i\}$ via the discrete-adjoint variant of the first-variation recursion together with a discrete forward Duhamel pass with source $-\calA''_{uu}(c)\cdot(\cdot,\cdot)$; the time integral over $r$ in $C_h$ is taken with composite Simpson weights on the uniform-stride subgrid, and a paired-path quadrature audit released with the code verifies its convergence at every operating point of the suite.
The trace $\Tr_U[\Phi_r^* (\gamma_t + \eps I)^{-1}[\D_r \gamma_t](\cdot)(\gamma_t + \eps I)^{-1}h]$ in $C_h$ is estimated stochastically. The default is the Hutchinson trace estimator
$\Tr(A) \approx n_{\mathrm{H}}^{-1}\sum_{j=1}^{n_{\mathrm{H}}}\langle \xi_j, A\xi_j\rangle$,
used here in its Gaussian-probe form, in which the $(\xi_j)$ are i.i.d.\ standard Gaussian on the retained $K^d$-mode spectral truncation of $U$, the estimator is unbiased for any square $A$, and for $A$ with symmetric part $A_{\mathrm{s}}$ its variance is $2\norm{A_{\mathrm{s}}}_{\mathrm{HS}}^2/n_{\mathrm{H}}$, the Gaussian-probe variance recorded in \cite{hutchinson1990stochastic} and credited there to Girard. The estimator of \cite{hutchinson1990stochastic} proper draws the probes instead from the Rademacher distribution ($\pm 1$ with equal probability), the unique minimum-variance unbiased choice among i.i.d.\ zero-mean probes, whose variance $2\sum_{i \neq j}(A_{\mathrm{s}})_{ij}^2/n_{\mathrm{H}}$ (entries taken in the probe basis) lies below the Gaussian value by exactly the diagonal contribution $2\sum_{i}(A_{\mathrm{s}})_{ii}^2/n_{\mathrm{H}}$; we retain the Gaussian probes, matching the Gaussian-probe form of Hutch++ below, and, as is now customary, use ``Hutchinson estimator'' for either probe law. For two equations whose assembled trace operator $A_r := \Phi_r^*(\gamma_t + \eps I)^{-1}[\D_r\gamma_t](\cdot)(\gamma_t + \eps I)^{-1}$ is observed to have rapidly decaying singular values, we replace Hutchinson by Hutch++ \cite{meyer2021hutchpp}, splitting the trace as $\Tr(A) = \Tr(Q^\top A Q) + \Tr((I - QQ^\top)A)$, where $Q$ is the $Q$-factor of a randomised sketch of dimension $k$ for the dominant subspace of $A$; the first term is computed exactly, the second by a Hutchinson estimator of dimension $m$ restricted to the orthogonal complement of $Q$. The Gaussian-probe form of Hutch++ is unbiased for any square $A$ and satisfies $\operatorname{Var} \leq \tfrac{8}{m-2}\norm{A - A_k}_{\mathrm{HS}}^2 \leq \tfrac{16}{(m-2)^2}\norm{A}_{*}^2$ at the coupling $k = \tfrac{m-2}{8} - 1$ of \cite{meyer2021hutchpp}, where $A_k$ is the best rank-$k$ approximation of $A$ and $\norm{\cdot}_{*}$ is the nuclear norm (equal to the trace for positive-semidefinite $A$); in particular no symmetry or definiteness is needed for unbiasedness or for this variance form. Our $A_r$ is neither symmetric nor positive-semidefinite, and our sketch and tail dimensions $(k, m)$ below sit outside the regime of these explicit constants, so we invoke Hutch++ as an empirical variance-reduction device rather than through the quantitative guarantee, the gain being realised whenever $\norm{A_r - A_{r,k}}_{\mathrm{HS}} \ll \norm{A_r}_{\mathrm{HS}}$, precisely the mechanism of the general-matrix variance bound of \cite{meyer2021hutchpp}, which we observe to hold for the assembled operator of gPAM and $\Phi^4_3$ at fixed $r$ (Remark~\ref{rem:numerics-gPAM} below), and we report raw and variance-reduced residuals side by side so that the conclusions do not rest on the estimator choice.

\medskip\noindent To reduce the variance of the resulting estimator, a linear control variate is constructed from each realised Brownian path $W(\omega)$ by simultaneously integrating the linearised SPDE $\mathrm{d}X^{\mathrm{lin}} = -L X^{\mathrm{lin}}\dt + q^{1/2}\dW$ on the same noise; the discrete-exact Bismut Skorokhod $\delta_U^{\mathrm{lin}}(v_h^{\mathrm{lin}})$ of the linear equation has closed form by Theorem~4 of \cite{mirafzali2025infinite} (specialised to the matched discrete scheme), and the pair $(\delta_U(v_h), \delta_U^{\mathrm{lin}}(v_h^{\mathrm{lin}}))$ is strongly correlated path-by-path. The classical regression-coefficient construction \cite{glasserman2003monte} replaces $\delta_U(v_h)(\omega)$ by $\delta_U(v_h)(\omega) - \hat c^{*}\,\bigl(\delta_U^{\mathrm{lin}}(v_h^{\mathrm{lin}})(\omega) - \E[\delta_U^{\mathrm{lin}}(v_h^{\mathrm{lin}})]\bigr)$ with $\hat c^{*}$ the sample regression slope. The construction is unbiased asymptotically (with a residual $O(M^{-1})$ correlation bias from estimating $\hat c^{*}$ on the same sample, dominated by the leading $M^{-1/2}$ Monte-Carlo error throughout our range), and reduces the empirical variance by a factor of $(1 - \rho^2)^{-1}$, where $\rho$ is the sample correlation between $\delta_U(v_h)$ and $\delta_U^{\mathrm{lin}}(v_h^{\mathrm{lin}})$. It does not alter the conclusion of any validation test below, and we report both the raw and the variance-reduced residuals (cf.\ \S\ref{subsubsec:results}).

\subsubsection{What is being tested}\label{subsubsec:validation}

In the absence of a closed-form $\beta_h$, four independent tests probe distinct aspects of the formula. Each is a precise mathematical statement at fixed $M$, and each fails under a different defect, namely implementation error, formula error, density-level disagreement, or finite-$M$ bias.

\medskip\noindent The first test, which we label (V1), concerns the trace identity, and compares the second-fibre route against finite differences. The most algorithmically intricate piece of the formula is the Malliavin derivative $[\D_r\gamma_t](\xi)$ entering $C_h$, delivered by the second-fibre construction of Proposition~\ref{prop:D-gamma}; an undetected error at this layer, whether in the proposition itself or in its discrete implementation, would propagate silently through every downstream test. The purpose of (V1) is to isolate this layer from the rest of the validation chain. For each path and each Hutchinson direction $\xi \in U$, the per-path trace contribution
\begin{equation*}
  \mathcal{I}_r(\xi) \;:=\; \langle \xi,\;\Phi_r^*\,(\gamma_t + \eps I)^{-1}\,[\D_r\gamma_t](\xi)\,(\gamma_t + \eps I)^{-1}\,h\rangle_U
\end{equation*}
is computed along two algorithmically independent routes. The second-fibre pipeline (SV) assembles $\mathcal{I}_r(\xi)$ from the discrete adjoint Jacobian and the second-fibre construction of Proposition~\ref{prop:D-gamma}. The finite-difference pipeline (FD) computes $[\D_r\gamma_t](\xi)$ by perturbing the path at time $r$ along the state-space direction $\calB(r, X(r))\,\xi \in H$ (the natural lift of $\xi \in U$ to a path perturbation, cf.\ \eqref{eq:trace-unfolded}) and reading off the difference quotient. Agreement of SV and FD to within the conjugate-gradient stopping tolerance verifies simultaneously that Proposition~\ref{prop:D-gamma} holds at the discrete level and that the SV pipeline implements it correctly; the subsequent tests (V2)--(V4) then probe Theorem~\ref{thm:main} alone, with no remaining ambiguity between implementation error and formula error.

\medskip\noindent The second test, (V2), is the integration-by-parts identity itself. The defining property of $\beta_h$ from Definition~\ref{def:log-deriv},
\[
  \E[\langle\nabla\varphi(X(t)), h\rangle_H] = -\E[\varphi(X(t))\,\beta_h(X(t))],
\]
combined with the conditional-expectation form of $\beta_h$ in \eqref{eq:score-main} and the tower property, yields
\begin{equation}\label{eq:exp-IBP}
  \E\bigl[\,\langle\nabla\varphi(X(t)), h\rangle_H\,\bigr] \;=\; \E\bigl[\,\varphi(X(t))\,\delta_U(v_h)\,\bigr] \qquad \text{for all } \varphi \in C^1_b(H).
\end{equation}
We choose six test functionals $\{\varphi_k\}_{k=1}^{6}$, of which five depend only on the projected coordinate $z = \langle h, c\rangle$, namely $\varphi_1(c) = \tfrac{1}{2}z^2$, $\varphi_2(c) = \tfrac{1}{6}z^3$, $\varphi_3(c) = \cos z$, $\varphi_5(c) = z\exp(-z^2/2)$, $\varphi_6(c) = \tanh z$, while $\varphi_4(c) = \exp(-\norm{c}_H^2/2)$ depends on the full state through its $H$-norm. The family combines a polynomial pair ($\varphi_1, \varphi_2$) with three bounded-smooth profiles ($\varphi_3, \varphi_6$, and the full-state Gaussian $\varphi_4$) and a Gaussian-localised Schwartz profile ($\varphi_5$), so that the directional gradients $\langle\nabla\varphi_k, h\rangle_H$ span the regimes of linear, quadratic, oscillatory, and exponentially decaying response. The four bounded functionals lie in $C_b^1(H)$ and so satisfy \eqref{eq:exp-IBP} directly; the two polynomials have gradients of at most quadratic growth, and the identity \eqref{eq:exp-IBP} extends to them by truncation under the moment bounds on $\langle h, X(t)\rangle$ following from Assumption~\ref{ass:nondeg} together with $\delta_U(v_h) \in L^2(\Omega)$ (Theorem~\ref{thm:main}). We estimate both sides by their empirical means $\hat L_k, \hat R_k$ over $M$ Monte-Carlo paths and report two standardised residuals, the pooled and the paired,
\begin{equation*}
  z_k^{\mathrm{pool}} \;:=\; \frac{\hat L_k - \hat R_k}{\sqrt{\widehat{\operatorname{Var}}(\hat L_k) + \widehat{\operatorname{Var}}(\hat R_k)}}, \qquad
  z_k^{\mathrm{pair}} \;:=\; \frac{\overline{D_k}}{\widehat{\mathrm{SE}}(\overline{D_k})}, \qquad D_k^{(j)} := L_k^{(j)} - R_k^{(j)},
\end{equation*}
for $k = 1, \dots, 6$,
where $\widehat{\mathrm{SE}}(\overline{D_k}) = \widehat{\operatorname{Var}}(D_k)^{1/2}/\sqrt{M}$. Because both sides are evaluated on the same $M$ paths, the paired form $z_k^{\mathrm{pair}}$ is the correct two-sample statistic. It retains the path-by-path covariance through $\operatorname{Var}(D_k) = \operatorname{Var}(L_k) + \operatorname{Var}(R_k) - 2\operatorname{Cov}(L_k, R_k)$, whereas the pooled form discards it, and for positively correlated $(\hat L_k, \hat R_k)$ it has the smaller denominator, hence is the more stringent statistic. We report both (Table~\ref{tab:numerics-main}) and rest the acceptance decision on $\max_k|z_k|^{\mathrm{pair}}$. Under \eqref{eq:exp-IBP}, the central limit theorem (CLT) gives $\sqrt{M}\,\overline{D_k} \Rightarrow \mathcal{N}(0, \operatorname{Var}(D_k))$, so $z_k^{\mathrm{pair}} \Rightarrow \mathcal{N}(0,1)$ for each $k$; conversely, a violation of the identity for any test functional would force $|z_k^{\mathrm{pair}}| \to \infty$ at rate $\sqrt{M}$ on that functional. Acceptance is the family-wise criterion $\max_k|z_k|^{\mathrm{pair}} < q_{95}$, where $q_{95}$ is the $95\%$ quantile of $\max_k|z_k|$ under the null hypothesis \eqref{eq:exp-IBP}, estimated by a wild bootstrap that multiplies the six per-path residuals $D_k^{(j)}$ by a common Rademacher sign and so preserves their empirical correlation; across the suite $q_{95}$ lies in $2.33$--$2.54$, against $1.96$ for a single test, so the criterion is exact for the correlated family rather than liberal or conservative by an unquantified margin. The per-functional correlations $\rho_k = \operatorname{Corr}(\hat L_k, \hat R_k)$ are recorded with the source. This is the principal distribution-free test of the formula; \eqref{eq:exp-IBP} is the integration-by-parts identity in the form of \cite{bogachev2010differentiable}, and is equivalent to the assertion that $\beta_h$ is the logarithmic derivative of $\mu_t$ along $h$, modulo the resolution provided by the test family $\{\varphi_k\}$. The relative residual $|\hat L_k - \hat R_k|/\max(|\hat L_k|, |\hat R_k|)$ is reported alongside as a scale-aware secondary diagnostic.

\medskip\noindent It is fair to ask what (V2) can actually reject. An identity test is informative only against the alternatives it can reject. We therefore recompute $\max_k|z_k|^{\mathrm{pair}}$ on the same Monte-Carlo paths with the correction term deleted ($\delta_U = M_h$) and with its sign inverted ($\delta_U = M_h - C_h$), and record the resolution $\mathrm{pow} := \max_k |\hat\E[\varphi_k\,C_h]| / \widehat{\mathrm{SE}}{}^{\mathrm{pair}}$, the number of paired standard errors by which the $C_h$ contribution to the right-hand side of \eqref{eq:exp-IBP} exceeds the Monte-Carlo noise floor. Where $\mathrm{pow}$ is small the test accepts the formula without resolving these variants against it; three resolved operating points, at which amplitude, ensemble, and time grid are enlarged until $\mathrm{pow} \ge 4$, carry the affirmative claim (Table~\ref{tab:numerics-power}). As a family-free corroboration we report a kernelised Stein discrepancy of the Nadaraya--Watson-smoothed projected score against the sample \cite{liu2016kernelized}, with a wild-bootstrap $p$-value \cite{chwialkowski2016kernel}; the score there is itself kernel-estimated, so the load-bearing comparison is the ordering of the discrepancy across the full, deleted, and inverted variants on the same sample.

\medskip\noindent The third test, (V3), compares the Bismut score pointwise against a kernel density estimate. The pushforward of $\mu_t$ under the continuous linear functional $\ell(c) = \langle h, c\rangle$ has a marginal density $p_z$ on $\R$, and by the chain rule for Fomin derivatives along bounded linear maps \cite{bogachev2010differentiable},
\begin{equation*}
  |h|_H^2\,\partial_z\,\log p_z(z) \;=\; -\E\bigl[\,\delta_U(v_h) \,\bigm|\, \langle h, X(t)\rangle = z\,\bigr].
\end{equation*}
The right-hand side is estimated from the sample by Nadaraya--Watson regression with Gaussian kernel and Silverman bandwidth $h_{\mathrm{NW}} = 1.06\,\hat\sigma_z\,M^{-1/5}$ \cite{silverman1986density}; the left-hand side is estimated independently by kernel-density estimation (KDE) of $p_z$ followed by numerical logarithmic differentiation, with the same Silverman bandwidth. The bandwidth choice $M^{-1/5}$ is asymptotically optimal for density estimation; for the derivative of the density it is slightly oversmoothed (the optimal rate is $M^{-1/7}$ \cite{wandjones1995kernel}), so the KDE-derivative side carries a residual bias of order $h_{\mathrm{NW}}^2$ that decays with $M$ at the price of variance. The Bismut and KDE-derivative curves are compared on the support of $\hat p_z$ by least-squares regression slope and Pearson correlation. Agreement to slope close to unity and correlation above $0.93$ (Table~\ref{tab:numerics-main}) places the formula within both the Monte-Carlo and the kernel-smoothing tolerance.

\medskip\noindent The fourth test, (V4), examines the Monte-Carlo convergence rate. An unbiased finite-variance estimator converges at the canonical Monte-Carlo rate $M^{-1/2}$ \cite{glasserman2003monte}; a biased estimator plateaus at finite $M$. Define the IBP residual vector $R(M) := (\hat L_k - \hat R_k)_{k=1}^6 \in \R^6$. By the multivariate CLT applied componentwise, $\sqrt{M}\,R(M) \Rightarrow N(0, \Sigma)$ as $M \to \infty$, where $\Sigma_{kl} := \operatorname{Cov}(L_k^{(j)} - R_k^{(j)},\,L_l^{(j)} - R_l^{(j)})$ is the per-sample covariance; consequently $\E\norm{R(M)}_2 \sim C\,M^{-1/2}$ as $M \to \infty$, with $C = \E\norm{N(0, \Sigma)}_2$ a finite constant determined by $\Sigma$. We sweep $M \in \{100, 200, 400, 800, 1600\}$ on Allen--Cahn 2D with three independent seeds per $M$ and fit
\begin{equation}\label{eq:conv-rate-fit}
  \log\Bigl(\,\norm{R(M)}_2 / \norm{\hat L}_2\,\Bigr) \;=\; \log C \;+\; \alpha\,\log M
\end{equation}
by least squares in $(\log M, \log \norm{R}_2/\norm{\hat L}_2)$. We report $\hat\alpha$ with its bootstrap $95\%$ confidence interval ($400$ resamples) and the coefficient of determination of the fit, on four representatives spanning the semilinear, quasilinear, transport, and rough blocks (Allen--Cahn 2D, $p$-Laplacian 2D, Navier--Stokes 2D, KPZ). The acceptance criterion, that the interval contain $-\tfrac12$, applies where the residual is variance-dominated; a flattened exponent at large $M$ is instead the signature of the $M$-independent discretisation floor, and is cross-examined directly by a time-refinement sweep in \S\ref{subsubsec:results} that fits the weak order of the residual bias in $\Delta t$ on the equation carrying the largest correction share.

\medskip\noindent These four tests are independent of one another. (V1) certifies the implementation pathwise, without appeal to the validity of Theorem~\ref{thm:main}; it would pass even if the formula were wrong, provided the algorithm correctly computed its right-hand side. (V2) tests the formula as an integration-by-parts identity at fixed $M$ across six functionals simultaneously; it would fail if any structural piece of the formula were missing or mis-signed, since the resulting bias would drive $|z_k|$ unbounded with $M$ on at least one functional, a sensitivity the ablation of \S\ref{subsubsec:results} verifies empirically rather than assumes. (V3) tests the formula against an external estimator of the same scalar quantity; passing requires the Bismut Nadaraya--Watson and the kernel-density-derivative estimators to coincide on the support of $\hat p_z$. (V4) tests for the absence of an $O(1)$ bias at finite $M$ by confirming the $M^{-1/2}$ convergence rate; it would fail if Theorem~\ref{thm:main} held only modulo a discretisation-induced bias. Together, (V1)--(V4) assemble a validation suite in the tradition of the Malliavin Monte-Carlo literature (Malliavin-weight sensitivity estimators \cite{fournie1999applications}, their absolute-continuity and Malliavin-weight foundations for (jump-)diffusions in finite and infinite dimensions \cite{forsterluetkebohmertteichmann2009}, and Monte-Carlo estimator methodology \cite{glasserman2003monte}), adapted to a setting in which no closed-form benchmark $\beta_h$ is available.

\subsubsection{The equations tested}\label{subsubsec:suite}

Table~\ref{tab:numerics-suite} lists the nine SPDEs and their structural data. We work throughout on the periodic torus $\mathbb{T}^d$ (replaced by an interval with sine modes for the 1D radial $\Phi^4_3$), with the state-independent additive noise $\calB(r, u) = q^{1/2}$ of \S\ref{subsubsec:estimator}, diagonal in the spectral basis with covariance $Q = q$ trace-class for $s > d/2$. The initial datum is fixed at the small spectral perturbation $c_0 = \tfrac{1}{20}\,q^{1/2}$ in $H$; the direction $h$ is the lowest-mode unit vector. On each finite-dimensional discretisation the corresponding regularised covariance solve is well posed at the numerical operating point used in the experiment; no continuum verification of Assumption~\ref{ass:nondeg} is inferred from this numerical fact. The basis size $K$, the time-step count $n$, and the spectral exponent $s$ are recorded in the table.

The first six equations (Allen--Cahn 2D, $p$-Laplacian in 1D and 2D, Cahn--Hilliard 2D, Navier--Stokes 2D in vorticity form, and $\Phi^4_2$) are instances of Theorem~\ref{thm:main} in the following sense. Their fully-local-monotone and differentiability hypotheses (H1)--(H5) and (D1)--(D4) are verified in Section~\ref{subsec:verification} (the principal cases) or Section~\ref{subsec:reduction} (the semilinear cases), as are the structural clauses, namely (SC1)(b) with $\alpha_1 = \nu$ for Navier--Stokes and the domination clause of (SC1)(a) with the dominated second derivative of Proposition~\ref{prop:p-Lap-Lambda} for the $p$-Laplacian at $p = 4$, (SC2) at the operating constant $c_*$ per Remark~\ref{rem:Y-moments-role}(iii), and (SC1$'$) for Navier--Stokes by Proposition~\ref{prop:NS-forcing}, the same clause being identified explicitly for the $p$-Laplacian in Remark~\ref{rem:p-Lap-forcing} and remaining a structural input there unless the unweighted fibre regularity recorded in that remark is available. What is tested numerically is the finite-dimensional spectral discretisation of each equation. At fixed spectral dimension the resulting stochastic system is finite-dimensional and smooth, so the Malliavin covariance, Tikhonov resolvent, second Malliavin fibre, and trace correction entering the implemented Bismut identity are well defined at the numerical operating points used below. The diagnostics \textup{(V1)}--\textup{(V4)} test the resulting finite-dimensional identity and its implementation; they are not intended as a verification of the continuum hypotheses (SC3) and (SC4.1)$_{q^*}$--(SC4.3)$_{q^*}$, the inverse-moment condition, or the Cameron--Martin compatibility condition. The continuum statement of Theorem~\ref{thm:main} applies whenever those hypotheses are verified for the corresponding continuum equation. The remaining three (gPAM 2D, $\Phi^4_3$ in 1D radial form, and KPZ $1+1$D) are singular equations in the BCCH scope of Section~\ref{subsec:singular}; we test them at fixed mollification scale $\eps_{\mathrm{mol}} > 0$ (written with a subscript to distinguish it from the trace ridge $\eps$ of \eqref{eq:correction-symbolic}), with $\eps_{\mathrm{mol}}$ proportional to the smallest resolved spatial scale of the spectral basis. For these three, what is tested is the finite-$K$ Bismut--Tikhonov identity on the mollified, spectrally truncated system at fixed $\eps_{\mathrm{mol}} > 0$, an identity of finite-dimensional Malliavin calculus (Theorem~\ref{thm:abstract-bismut-fomin} on $H_N$, cf.\ Corollary~\ref{cor:finite-dim}); no renormalised limit is tested, that limit falling under Theorem~\ref{thm:score-gPAM}(c) and being, by construction, a distributional rather than function-valued object. See Remark~\ref{rem:singular-numerics} below. For gPAM and $\Phi^4_3$ the mollified equation is additionally a member of the variational subcase (for gPAM, through the companion model \eqref{eq:gPAM-variational}); KPZ in $1+1$D is not, because the quasilinear term $\tfrac{1}{2}\lambda|\nabla u|^2$ takes the mollified KPZ equation outside the Liu--R\"ockner variational subcase (as recorded in Section~\ref{subsec:singular}). In neither case is a verified continuum instance of Theorem~\ref{thm:main} claimed for these three rows.

\begin{table}[pos=ht]
\centering
\small
\setlength{\tabcolsep}{4pt}
\renewcommand{\arraystretch}{1.10}
{\setlength{\tabcolsep}{3pt}
\resizebox{\textwidth}{!}{%
\begin{tabular}{@{}lllccc@{}}
\toprule
Equation & $\calA(u)$ & State space & $K\!/\!n\!/\!s$ & Theorem & Section \\
\midrule
Allen--Cahn 2D                  & $-\Delta u + m^2 u + u^3$                          & $L^2(\mathbb{T}^2)$    & $10/25/1.7$  & \ref{thm:main}                  & \ref{subsec:reduction} \\
$p$-Laplacian 1D                & $-\!\operatorname{div}(|\nabla u|^{p-2}\nabla u) + \eta u$ & $L^2(\mathbb{T}^1)$ & $16/30/1.0$  & \ref{thm:main}                  & \ref{subsec:verification} \\
$p$-Laplacian 2D                & $-\!\operatorname{div}(|\nabla u|^{p-2}\nabla u) + \eta u$ & $L^2(\mathbb{T}^2)$ & $10/25/1.7$  & \ref{thm:main}                  & \ref{subsec:verification} \\
Cahn--Hilliard 2D               & $\Delta(-\eps_{\mathrm{CH}}^2\Delta u + u^3 - u)$  & $L^2(\mathbb{T}^2)$    & $10/25/1.7$  & \ref{thm:main}                  & \ref{subsec:reduction} \\
Navier--Stokes 2D (vorticity)   & $u\!\cdot\!\nabla\omega - \nu\Delta\omega$         & $L^2(\mathbb{T}^2)$    & $10/25/1.7$  & \ref{thm:main}                  & \ref{subsec:verification} \\
$\Phi^4_2$                      & $-\Delta u + m^2 u + u^3$                          & $L^2(\mathbb{T}^2)$    & $10/25/1.0$  & \ref{thm:main}                  & \ref{subsec:reduction} \\
\midrule
gPAM 2D (surrogate)             & $-\Delta u + m^2 u - \alpha\,\Pi(u^2)$              & $L^2(\mathbb{T}^2)$    & $10/100/1.7$ & finite-$K$ Bismut, $\eps_{\mathrm{mol}}$-fix & \ref{subsec:singular} \\
$\Phi^4_3$ (1D rad., mollified) & $-\Delta u + m^2 u + u^3 - 3 C_\eps u$             & $L^2([0,1])$           & $20/25/1.5$  & finite-$K$ Bismut, $\eps_{\mathrm{mol}}$-fix & \ref{subsec:singular} \\
KPZ $1+1$D (mollified)          & $-\nu\Delta u - \tfrac{1}{2}\lambda|\nabla u|^2$   & $L^2(\mathbb{T}^1)$    & $20/20/1.0$  & finite-$K$ Bismut, $\eps_{\mathrm{mol}}$-fix & \ref{subsec:singular} \\
\bottomrule
\end{tabular}%
}}
\caption{The test suite. For the first six equations the tested object is the verified spectral discretisation, the continuum statement of Theorem~\ref{thm:main} applying under the structural hypotheses (SC3) and (SC4.1)$_{q^*}$--(SC4.3)$_{q^*}$ and Assumption~\ref{ass:nondeg}, the clauses (SC1), (SC1$'$) and (SC2) being verified in Sections~\ref{subsec:reduction}--\ref{subsec:verification}. The last three are singular SPDEs in the BCCH scope; we test them at fixed mollification scale $\eps_{\mathrm{mol}} > 0$, the renormalised limit of Theorem~\ref{thm:score-gPAM}(c) being a distributional rather than function-valued object beyond the scope of the numerical test. For the first six rows, the ``Theorem'' column records the continuum variational theorem whose finite-dimensional spectral discretisation is being tested, with the continuum statement remaining conditional on the structural hypotheses stated in the text. For the last three singularly motivated rows, ``finite-$K$ Bismut'' indicates that the actual numerical test is performed on the finite-dimensional mollified/truncated stochastic system. The notation ``$\eps_{\mathrm{mol}}$-fix'' emphasises that no $\eps_{\mathrm{mol}} \downarrow 0$ renormalised-limit claim is tested. The column ``$K/n/s$'' records the spectral basis size $K$ per dimension, time-step count $n$ over $[0, T]$, and the additive-noise decay exponent $s$ in $q_k = (1 + m^2 + \lambda_k)^{-s}$. The exponent $s > d/2$ makes the covariance $Q$ trace-class; the $\Phi^4_2$ row at $s = 1$ in $d = 2$ sits at the borderline $s = d/2$, where the modal sum diverges only logarithmically, and is run on the finite $K^2$-mode truncation where the trace is finite, without claiming the uniform-in-$K$ trace-class bound. The gPAM row is the frozen-coefficient additive surrogate of \S\ref{subsubsec:suite}, not the true multiplicative equation; its drift retains a genuine quadratic Anderson-type nonlinearity $-\alpha\,\Pi(u^2)$, with $\alpha > 0$ a coupling constant and $\Pi$ the mean-zero spectral projection, the Wick mass being absorbed into the linear part.}
\label{tab:numerics-suite}
\end{table}

For each equation in the variational block (Allen--Cahn through $\Phi^4_2$) and for the singular KPZ at fixed $\eps_{\mathrm{mol}} > 0$, we take $M = 400$ Monte-Carlo paths, $n_{\mathrm{sub}} = 6$ subgrid points for the correction term, and one Hutchinson direction per subgrid point. For gPAM we take $M = 300$, $n_{\mathrm{sub}} = 15$, and Hutch++ with sketch and tail dimensions $(k, m) = (2, 2)$; for $\Phi^4_3$ we take $M = 400$, $n_{\mathrm{sub}} = 6$, and the same Hutch++ parameters.

\begin{figure}[pos=t]
\centering
\includegraphics[width=0.95\textwidth]{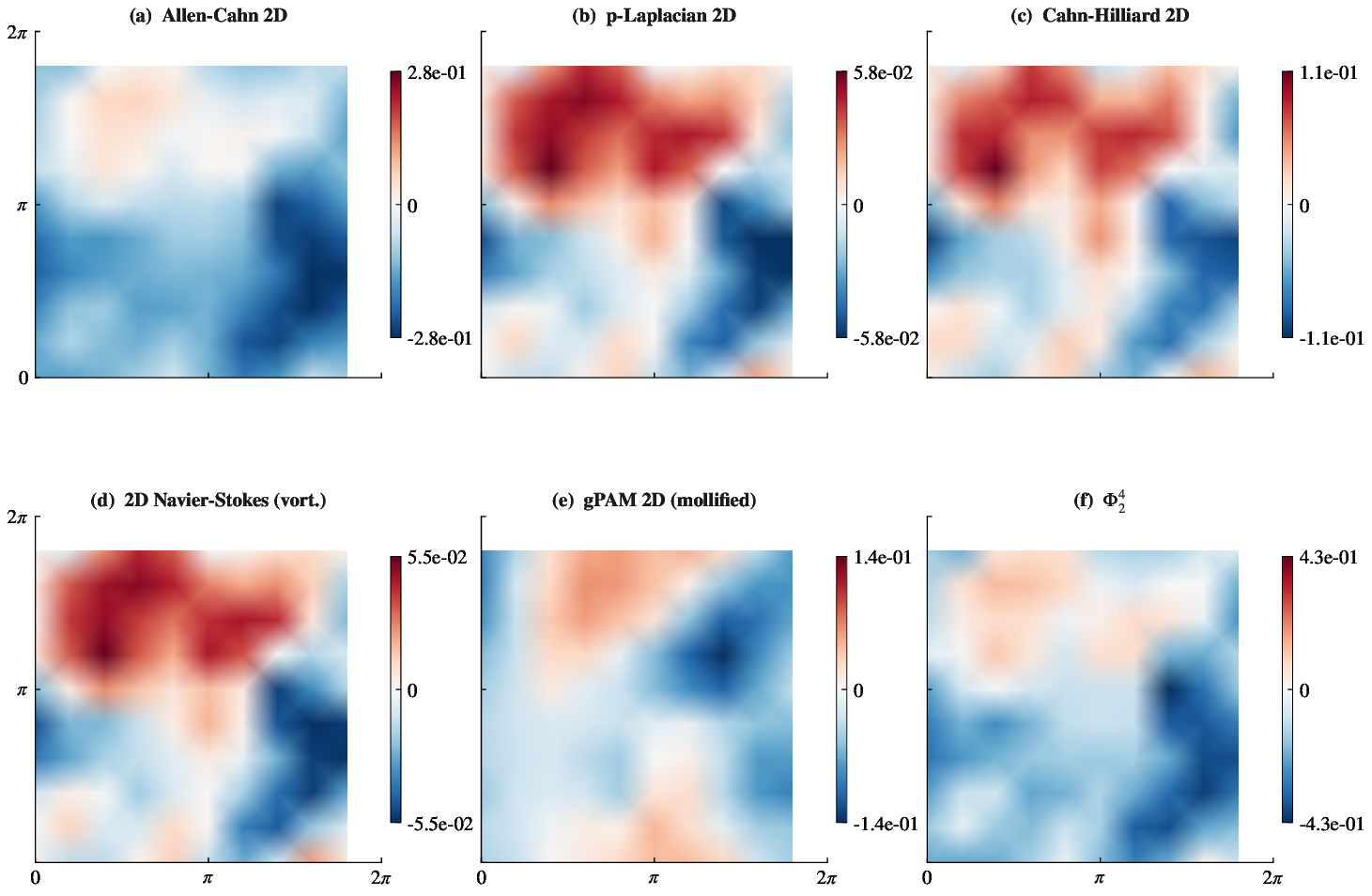}
\caption{Sample terminal states $X(T, \omega)$ at one Monte-Carlo realisation for the six two-dimensional equations of the suite. The panels show (a) Allen--Cahn 2D, (b) $p$-Laplacian 2D, (c) Cahn--Hilliard 2D, (d) 2D Navier--Stokes in vorticity form, (e) gPAM 2D (surrogate), and (f) $\Phi^4_2$. The structural diversity of the suite is visible in the regularising character of Cahn--Hilliard and Navier--Stokes, the rougher profile of the $p$-Laplacian, and the mollified-noise texture of the gPAM surrogate. The estimator of \S\ref{subsubsec:estimator} treats all six on the same footing, as spectral discretisations for the first five and as the finite-$K$ mollified surrogate for gPAM.}
\label{fig:numerics-2d-states}
\end{figure}

A clarification on gPAM is in order, since it is the one member of the scope with genuinely multiplicative noise. The mollified equation $\mathrm{d}X = (\Delta - m^2)X\dt + X\,(\rho_{\eps_{\mathrm{mol}}} Q^{1/2}\,\dW)$ has state-dependent diffusion $\calB^{\eps_{\mathrm{mol}}}(u) = u\cdot\rho_{\eps_{\mathrm{mol}}} Q^{1/2}$, so $\calB'_u \ne 0$ and the full diffusion contribution to the second Malliavin fibre would be required to test it as stated. We do not do this. Consistent with the additive scope of the present experiments, the gPAM row tests a \emph{frozen-coefficient additive surrogate}, in which the multiplicative channel $u\,\rho_{\eps_{\mathrm{mol}}}Q^{1/2}\,\dW$ is replaced by the additive channel $\rho_{\eps_{\mathrm{mol}}}Q^{1/2}\,\dW$ (the diffusion coefficient frozen at a deterministic background), while a genuine quadratic Anderson-type drift nonlinearity $N(u) = -\alpha\,\Pi(u^2)$, with $\Pi$ the mean-zero projection and the Wick mass absorbed into $L$, is retained. The surrogate is therefore not the heat equation; it carries a real local nonlinearity that activates the drift second variation $\calA''_{uu}$ exactly as the cubic does for $\Phi^4$, and it is this drift nonlinearity, not the multiplicative noise, that the gPAM row exercises. What it does not test is the state-dependent-diffusion mechanism $\calB'_u \ne 0$ of the true equation; a genuine multiplicative test is deferred to the state-dependent-diffusion implementation noted in \S\ref{subsubsec:estimator}. The choice of Hutch++ for gPAM and $\Phi^4_3$ follows from an empirical observation about the assembled trace operator $A_r = \Phi_r^*(\gamma_t + \eps I)^{-1}[\D_r\gamma_t](\cdot)(\gamma_t + \eps I)^{-1}$. On the smooth mollifier of scale $\eps_{\mathrm{mol}} > 0$ its singular spectrum at a typical $r$ is concentrated in the first few modes, so $\norm{A_r - A_{r,k}}_{\mathrm{HS}} \ll \norm{A_r}_{\mathrm{HS}}$ for small $k$, the regime in which Hutch++ reduces variance relative to Hutchinson \cite{meyer2021hutchpp}. KPZ falls outside this regime, since the quasilinear nonlinearity $|\nabla u|^2$ couples all spatial frequencies, the assembled operator is full-rank, and vanilla Hutchinson is used. Sample terminal states for the six two-dimensional equations of the suite are shown in Figure~\ref{fig:numerics-2d-states}.

The equation-specific physical parameters are collected in Table~\ref{tab:numerics-params}; together with the $K/n/s$ data of Table~\ref{tab:numerics-suite} and the per-test Monte-Carlo settings above, they fix the experiment. The shared numerical settings are as follows. The trace regularisation is the ridge $\eps$ in $(\gamma_t + \eps I)^{-1}$ of \eqref{eq:correction-symbolic}, the same $\eps$ as in the theory; to avoid collision with it, the singular-SPDE mollification scale is written $\eps_{\mathrm{mol}}$ throughout this subsection and in Table~\ref{tab:numerics-params}. We take $\eps = 10^{-8}$ for the main-term solve and $\eps = 10^{-9}$ for the correction-term solve, raised to an equation-specific floor $\eps = 10^{-6}$ for $p$-Laplacian 1D, with conjugate-gradient relative tolerance $\le 10^{-4}$. The (V1) finite-difference baseline uses a two-sided central difference with step $5\times10^{-5}$ in the path-perturbation amplitude. The stochastic convolution is integrated by the exact-in-law increment $E\,q^{1/2}\Delta W^{(i)}$ of the Lawson scheme \eqref{eq:etd-Lawson-discrete} rather than an Euler approximation. The (V3) kernel-density estimates use the Silverman bandwidth on the empirical support of $\hat p_z$ with the outer $2.5\%$ of mass trimmed at each tail to suppress boundary bias. All runs use a fixed seed recorded with the source; the seed and the complete parameter file accompany the submission.

\begin{table}[pos=t]
\centering
\small
\setlength{\tabcolsep}{6pt}
\renewcommand{\arraystretch}{1.12}
\begin{tabular}{lll}
\toprule
Equation & Terminal $T$ & Physical parameters \\
\midrule
Allen--Cahn 2D              & $0.30$ & $m^2 = 1.5$ \\
$p$-Laplacian 1D            & $0.10$ & $p = 4$, $\eta = 2.0$ \\
$p$-Laplacian 2D            & $0.20$ & $p = 4$, $\eta = 0.5$ \\
Cahn--Hilliard 2D           & $0.20$ & $\eps_{\mathrm{CH}} = 0.30$ \\
Navier--Stokes 2D (vort.)   & $0.25$ & $\nu = 0.50$ \\
$\Phi^4_2$                  & $0.30$ & $m^2 = 1.0$ \\
gPAM 2D (surrogate)         & $0.20$ & $m^2 = 0.5$, $\eps_{\mathrm{mol}} = 0.10$ \\
$\Phi^4_3$ (1D radial)      & $0.30$ & $m^2 = 1.0$ \\
KPZ $1+1$D                  & $0.10$ & $\nu = 1.0$, $\lambda = 0.10$ \\
\bottomrule
\end{tabular}
\caption{Equation-specific parameters. $T$ is the terminal time; remaining columns give the drift parameters ($m^2$ mass, $p$ the $p$-Laplacian exponent and $\eta$ the coefficient of the linear zeroth-order term $\eta u$ in its drift, $\eps_{\mathrm{CH}}$ the Cahn--Hilliard interface width, $\nu$ viscosity, $\lambda$ the KPZ coupling) and, for the singular cases, the mollification scale $\eps_{\mathrm{mol}}$. The additive-noise decay exponent $s$ and the discretisation $K, n$ are in Table~\ref{tab:numerics-suite}.}
\label{tab:numerics-params}
\end{table}

\subsubsection{Outcomes}\label{subsubsec:results}

We report the four tests in turn. Turning first to (V1), the trace identity holds at numerical precision. On a grid of four subgrid times $r$ and two Hutchinson directions $\xi$ per equation, the SV and FD routes to $\mathcal{I}_r(\xi)$ agree to relative discrepancy between $2.8 \times 10^{-10}$ (gPAM) and $1.6 \times 10^{-7}$ ($p$-Laplacian 2D), the residual set by the conjugate-gradient tolerance for $(\gamma_t + \eps I)^{-1}$ rather than by any SV/FD disagreement; Proposition~\ref{prop:D-gamma} and its implementation are in agreement at numerical precision across the suite. Repeating the grid at doubled field amplitude leaves eight of nine equations at the same precision; the exception is $p$-Laplacian 1D, whose degenerate-gradient regime at large amplitude degrades the agreement to $4.5 \times 10^{-3}$, a recorded regime boundary of the second-fibre construction there, immaterial to every statistic below since the $C_h$ share of that equation is $0.1$ per cent. As a further guard, three algorithmically independent assemblies of $C_h$ itself are compared on identical paths for every equation, namely the SV pipeline, the restricted FD route, and a full-product FD route that additionally differentiates the transfer factor $Y(t, r)$ under the path perturbation. The first two agree within five per cent of scale wherever $C_h$ is material, and the third measures the anticipating-diagonal term of the discrete Skorokhod divergence~\eqref{eq:discrete-skorokhod}, the second sum of that display, carried inside the reported main term, at no more than one third of $C_h$ at doubled amplitude, and at the per-mille level for gPAM, below one paired standard error at every operating point of Table~\ref{tab:numerics-main}.

\medskip\noindent For (V2), the integration-by-parts identity is tested across the whole suite. Table~\ref{tab:numerics-main} collects, for each equation, the minimum scale-aware relative residual across the six functionals, the paired statistic $\max_k|z_k|^{\mathrm{pair}}$ against its wild-bootstrap quantile $q_{95}$, the (V3) regression diagnostics, and the share of $\delta_U(v_h)$ carried by the correction term. All nine equations satisfy the family-wise criterion, with $\max_k|z_k|^{\mathrm{pair}}$ between $0.76$ and $1.86$ against thresholds between $2.33$ and $2.43$; the verdict is unchanged under the pooled standard error (maxima $0.59$--$1.85$, recorded with the source).

\begin{table}[pos=t]
\centering
\small
\setlength{\tabcolsep}{4pt}
\renewcommand{\arraystretch}{1.08}
\begin{tabular}{lrrrrrrr}
\toprule
Equation & $M$ & rel.\ res. & $\max_k|z_k|^{\mathrm{pair}}$ & $q_{95}$ & slope & corr. & $C_h$ share \\
\midrule
Allen--Cahn 2D                & $400$ & $4.6 \times 10^{-3}$ & $0.99$ & $2.34$ & $0.915$ & $0.987$ & $3.4\%$ \\
$p$-Laplacian 1D              & $400$ & $8.2 \times 10^{-2}$ & $1.85$ & $2.42$ & $0.859$ & $0.993$ & $0.1\%$ \\
$p$-Laplacian 2D              & $400$ & $6.9 \times 10^{-2}$ & $1.79$ & $2.42$ & $0.950$ & $0.955$ & $0.6\%$ \\
Cahn--Hilliard 2D             & $400$ & $6.4 \times 10^{-2}$ & $1.66$ & $2.43$ & $0.959$ & $0.960$ & $2.3\%$ \\
Navier--Stokes 2D (vort.)     & $400$ & $6.9 \times 10^{-2}$ & $1.86$ & $2.43$ & $0.949$ & $0.953$ & $0.1\%$ \\
$\Phi^4_2$                    & $400$ & $6.3 \times 10^{-3}$ & $0.96$ & $2.33$ & $0.906$ & $0.984$ & $5.4\%$ \\
\midrule
gPAM 2D (surrogate)           & $300$ & $4.6 \times 10^{-2}$ & $1.17$ & $2.35$ & $0.918$ & $0.986$ & $13.4\%$ \\
$\Phi^4_3$ (1D radial)        & $400$ & $1.6 \times 10^{-2}$ & $0.88$ & $2.38$ & $0.926$ & $0.989$ & $7.2\%$ \\
KPZ $1+1$D                    & $400$ & $2.6 \times 10^{-3}$ & $0.76$ & $2.37$ & $0.981$ & $0.939$ & $1.1\%$ \\
\bottomrule
\end{tabular}
\caption{Tests (V2) and (V3) on the nine SPDEs of Table~\ref{tab:numerics-suite}. The columns report the minimum scale-aware relative residual of \eqref{eq:exp-IBP} across the six functionals; the paired maximum $\max_k|z_k|^{\mathrm{pair}}$ against its multiplicity-correct wild-bootstrap null quantile $q_{95}$ (the acceptance criterion); the regression slope and Pearson correlation of the Bismut Nadaraya--Watson estimator against the kernel-density derivative along $z = \langle h, X(t)\rangle$; and the share $\overline{|C_h|}/(\overline{|M_h|} + \overline{|C_h|})$ of the Skorokhod integral carried by the correction. All nine equations pass.}
\label{tab:numerics-main}
\end{table}

\medskip\noindent The power of (V2) against the deleted and sign-inverted corrections is measured next. At the operating points of Table~\ref{tab:numerics-main} the resolution of $C_h$ is at most $2.1$ paired standard errors (gPAM), so that table accepts the formula without resolving the deleted and sign-inverted variants against it. Even at those points the structure is visible on the equation with the largest correction share. Inverting the sign of $C_h$ on gPAM's own $300$ paths fails the criterion at $3.06$ against $q_{95} = 2.35$, and at doubled noise amplitude the two semilinear representatives pass in full form ($0.91$, $1.47$) while their sign-inverted variants fail ($2.49$, $3.54$). The designed points of Table~\ref{tab:numerics-power} complete the demonstration. At resolutions of $4.2$--$5.3$ standard errors, the identical Monte-Carlo paths that accept $\delta_U = M_h + C_h$ reject the deleted correction at $3.0$--$3.6$ and the sign-inverted one at $7.2$--$7.8$. The kernelised Stein discrepancy orders the three variants identically (for gPAM, $0.068$, $0.299$, and $1.06$ for the full, deleted, and sign-inverted forms). The $\Phi^4_2$ point is placed at $1.5\times$ rather than $2\times$ amplitude so that the measured It\^o-sum remainder of (V1), which grows with amplitude, sits a full standard error inside the acceptance margin while the resolution exceeds five; the choice is recorded with the source.

\begin{table}[pos=t]
\centering
\small
\setlength{\tabcolsep}{3pt}
\renewcommand{\arraystretch}{1.10}
\begin{tabular}{lccrrrrrc}
\toprule
Equation & amplitude & $M$ & $n$ & full & $q_{95}$ & deleted & inverted & resolution \\
\midrule
Allen--Cahn 2D      & $2\times$   & $3600$ & $200$ & $\mathbf{1.57}$ & $2.48$ & $3.30$ & $7.27$ & $5.1\,\sigma$ \\
$\Phi^4_2$          & $1.5\times$ & $4400$ & $200$ & $\mathbf{2.04}$ & $2.54$ & $2.96$ & $7.16$ & $5.3\,\sigma$ \\
gPAM 2D (surrogate) & $1\times$   & $1200$ & $200$ & $\mathbf{1.38}$ & $2.41$ & $3.58$ & $7.78$ & $4.2\,\sigma$ \\
\bottomrule
\end{tabular}
\caption{Ablation of the correction term at the three resolved operating points (refined time grid, enlarged ensembles; noise amplitude relative to Table~\ref{tab:numerics-params}). All three variants are evaluated on identical Monte-Carlo paths, with ``full'' denoting $\max_k|z_k|^{\mathrm{pair}}$ for the estimator $\delta_U = M_h + C_h$ of Theorem~\ref{thm:main}, ``deleted'' the same statistic for $\delta_U = M_h$, and ``inverted'' that for $\delta_U = M_h - C_h$; ``resolution'' is the number of paired standard errors by which the $C_h$ contribution to \eqref{eq:exp-IBP} exceeds the Monte-Carlo noise floor. The full formula passes at every point at which either variant fails.}
\label{tab:numerics-power}
\end{table}

\begin{figure}[pos=t]
\centering
\includegraphics[width=0.95\textwidth]{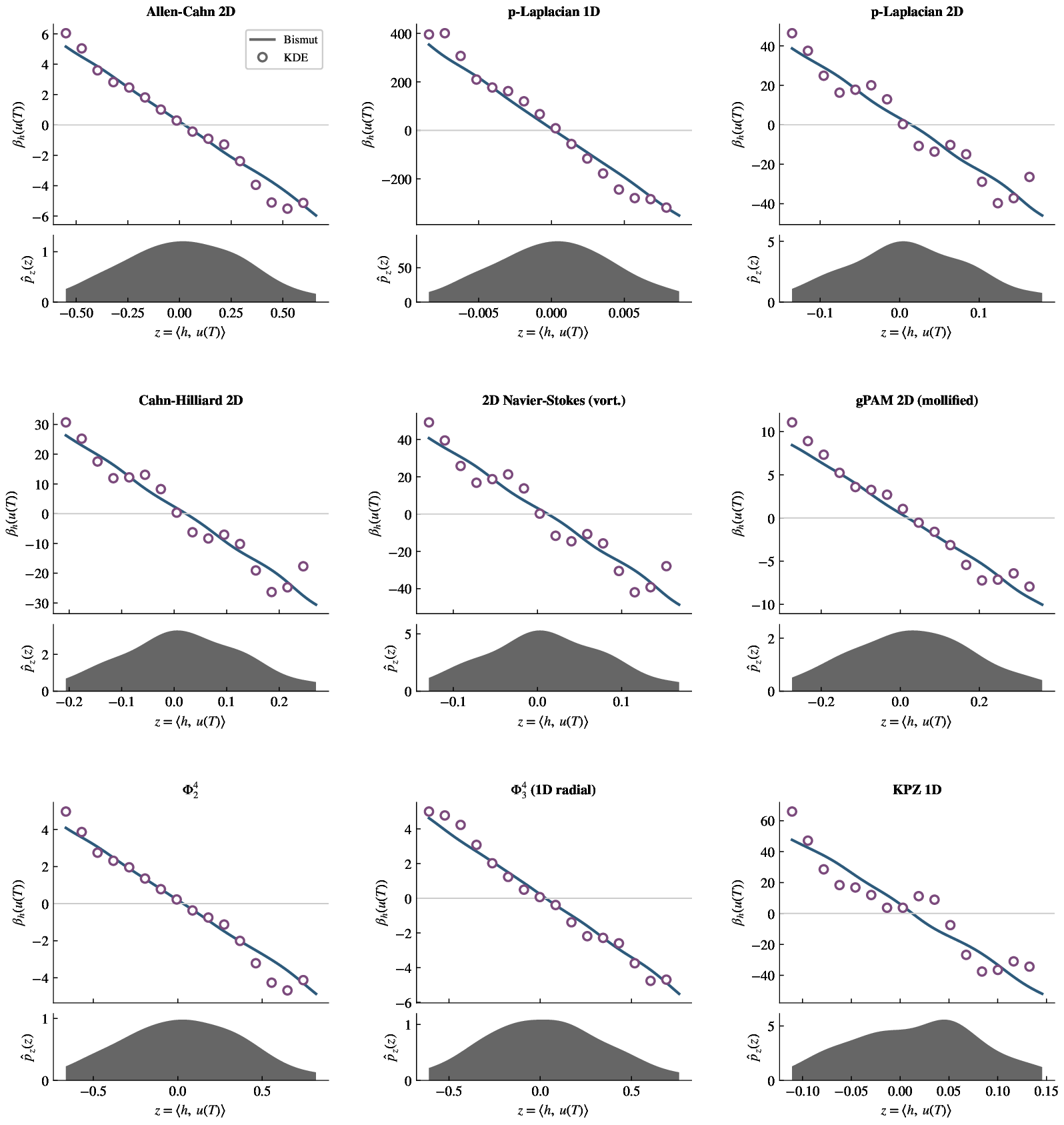}
\caption{Test (V3). For each equation, the upper panel overlays the Bismut Nadaraya--Watson estimator (solid line) on the kernel-density logarithmic derivative (open circles), as functions of the projected pairing $z = \langle h, X(T)\rangle$. The lower panel shows the kernel-density estimate $\hat p_z(z)$; the Bismut estimator is reliable where $\hat p_z$ is supported. The quantitative slope and correlation are in Table~\ref{tab:numerics-main}.}
\label{fig:numerics-pointwise}
\end{figure}

\medskip\noindent For (V3), the pointwise comparison, the attenuation can be attributed. The regression slope of the Bismut Nadaraya--Watson estimator against the kernel-density derivative ranges over $0.859$ ($p$-Laplacian 1D) to $0.981$ (KPZ), with Pearson correlation $0.939$--$0.993$ (Table~\ref{tab:numerics-main}; overlay in Figure~\ref{fig:numerics-pointwise}). We attribute the sub-unit slopes by direct measurement. Sweeping the trace ridge $\eps$ over $10^{-4}$--$10^{-8}$ on the worst case ($p$-Laplacian 1D), at a fixed reference sample of $M = 300$, drives the slope from $0.139$ to $0.797$ with the correlation pinned at $0.992$ throughout and $\norm{\tilde h}$ growing accordingly, identifying Tikhonov shrinkage of $(\gamma_t + \eps I)^{-1}h$ as the dominant component; holding the ridge at its floor and enlarging the sample gives slopes $0.797$, $0.900$, and $0.917$ at $M = 300$, $400$, and $1200$, identifying the remainder as the kernel-bandwidth attenuation of the $M^{-1/5}$ Silverman rule on the KDE-derivative side. The same machinery evaluated at the production operating point ($\eps = 10^{-6}$, $M = 400$) reproduces the Table~\ref{tab:numerics-main} row for this equation exactly, so the three experiments are mutually consistent, each moving one component at a time. Both components vanish in the joint limit, and neither is a property of the formula.

\medskip\noindent The outcome does not depend on the chosen direction. Replacing the suite direction by the third spectral mode and by a fixed unit mixture of the top four modes, on four representative equations with matched noise paths across directions, gives fifteen of eighteen rows inside their per-row $q_{95}$. The two third-mode exceedances (Allen--Cahn at $3.13$, KPZ at $3.31$) halve into the acceptance band when the time grid is doubled ($1.60$ and $1.10$), attributing them to the $O(\Delta t)$ scheme bias through the stiffer mode couplings. The KPZ mixture row persists ($2.72$, and $3.45$ at the doubled grid) and coincides with the largest conjugate-gradient iteration counts of the study, consistent with the conditioning of $\gamma_t^{-1}h$ for non-eigenvector directions under the quasilinear nonlinearity; with a $C_h$ share of one per cent it is not a correction-term effect, and we record it as a solver-conditioning boundary at this operating point rather than as evidence about the representation.

\begin{figure}[pos=t]
\centering
\includegraphics[width=0.85\textwidth]{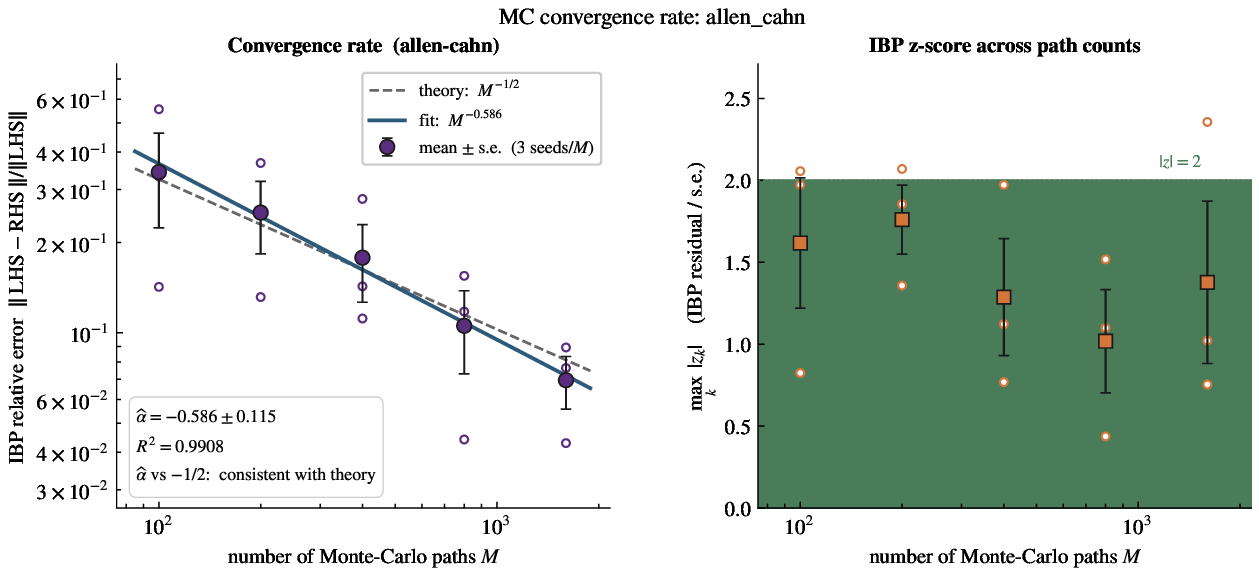}
\caption{Test (V4) on Allen--Cahn 2D. The left panel shows the log--log plot of the vector relative IBP residual $\norm{R(M)}_2/\norm{\hat L}_2$ against the path count $M$; the dashed grey line is the theoretical Monte-Carlo rate $M^{-1/2}$, and the solid blue line the least-squares fit \eqref{eq:conv-rate-fit} with $\hat\alpha = -0.586 \pm 0.115$, $R^2 = 0.991$, and $95\%$ bootstrap confidence interval $[-0.794, -0.353]$. The right panel shows the corresponding maximum paired $z$-statistic across the six test functionals at each $M$, with the multiplicity-correct acceptance quantile shaded.}
\label{fig:numerics-conv}
\end{figure}

\medskip\noindent For (V4), the convergence rates in $M$ and $\Delta t$ come out as predicted. The fit \eqref{eq:conv-rate-fit} on Allen--Cahn 2D yields
\begin{equation}\label{eq:conv-rate-result}
  \hat\alpha \;=\; -0.586 \pm 0.115, \qquad \text{95\% bootstrap interval } [-0.794,\, -0.353], \qquad R^2 = 0.991,
\end{equation}
bracketing the unbiased prediction $\alpha = -\tfrac{1}{2}$ with no visible plateau (Figure~\ref{fig:numerics-conv}); KPZ gives $\hat\alpha = -0.305 \pm 0.211$ with interval $[-0.670, +0.157]$, compatible though noise-limited at this cost. For Navier--Stokes and $p$-Laplacian 2D the fitted exponents flatten to $\approx -0.20$ with intervals excluding $-\tfrac{1}{2}$; at those operating points the sweep reaches the $M$-independent discretisation floor, the exponent interpolates between $-\tfrac{1}{2}$ and $0$ accordingly, and the floor itself is then measured directly. On the equation carrying the largest correction share, a four-rung refinement $n \in \{25, 50, 100, 200\}$ at $M = 200$ with three seeds per rung drives the control-variate-resolved relative bias from $2.31 \times 10^{-2}$ to $2.44 \times 10^{-3}$, a factor of $9.5$ against the factor of $8$ for an exact first-order scheme, with fitted weak order $p = 1.12 \pm 0.51$ (standard error; $R^2 = 0.976$); at the two finest rungs the bias falls below the seed-scatter floor of the probe. The residual bias of the estimator is thus the scheme's weak order one, not the formula's.

\begin{remark}[The gPAM correction term]\label{rem:numerics-gPAM}
The assembled trace operator $A_r$ for the mollified gPAM surrogate has rapidly decaying singular values; of the leading $100$ singular values at a typical $r$, the first three account for approximately $97$ per cent of the Hilbert--Schmidt norm, the regime in which Hutch++ with $(k, m) = (2, 2)$ essentially eliminates the stochastic trace error, and the reason gPAM and $\Phi^4_3$ use it. With that variance suppressed, gPAM is the suite member on which the correction term is largest ($13.4$ per cent of $\delta_U$) and the structure of the formula most exposed. It passes at $1.17$ against $q_{95} = 2.35$; its sign-inverted variant fails on the same paths already at $M = 300$; and at $(M, n) = (1200, 200)$ the full, deleted, and inverted forms stand at $1.38$, $3.58$, and $7.78$ (Table~\ref{tab:numerics-power}). The time-refinement sweep of (V4) then locates the residual entirely in the discretisation; the bias is weak order one in $\Delta t$ and $2.4 \times 10^{-3}$ at $n = 200$, a property of the discrete scheme \eqref{eq:etd-Lawson-discrete} and not of the formula \eqref{eq:skorokhod-decomp}.
\end{remark}

\begin{remark}[Three singularly motivated cases at fixed $\eps_{\mathrm{mol}}$]\label{rem:singular-numerics}
Three rows of Table~\ref{tab:numerics-suite}, gPAM, $\Phi^4_3$, and KPZ, are motivated by singular SPDEs in the scope of \cite{brunedchandrachevyrevhairer2021}. The experiments are carried out at a fixed spatial mollification $\eps_{\mathrm{mol}} > 0$ and a fixed finite spectral truncation.

At this finite-dimensional level the implemented systems are smooth stochastic differential equations and the tested integration-by-parts and Tikhonov-correction identities are finite-dimensional Malliavin identities. The gPAM row is the frozen-coefficient additive surrogate described in the numerical setup. No claim is made that the KPZ row is a continuum instance of Theorem~\ref{thm:main}, since KPZ lies outside the variational Liu--R\"ockner subcase considered there.

No numerical claim is made about the $\eps_{\mathrm{mol}} \downarrow 0$ renormalised limit. The unconditional limiting statement available in the present paper is the distributional derivative of Theorem~\ref{thm:score-gPAM}(c); a function-valued logarithmic derivative at the limit additionally requires the second-order Malliavin and inverse-moment hypotheses of Theorem~\ref{thm:score-gPAM}(b).
\end{remark}

\medskip\noindent Finally, we turn to variance reduction and the bias probe. The linear control variate of \S\ref{subsubsec:estimator} reduces the empirical variance of the IBP residual by factors between $1.2 \times 10^{3}$ ($\Phi^4_3$) and $2.1 \times 10^{6}$ ($p$-Laplacian 1D, where the small-noise regime keeps the solution close to its linearisation and the control variate correspondingly tightly correlated with the full estimator). Its uncentred residual has exactly zero mean under the discrete measure, so beyond variance reduction it furnishes a bias probe whose noise floor sits roughly the square root of the reduction factor below the raw Monte-Carlo error; it is this probe that resolves the discretisation floor in (V4) and supplies the weak-order fit above. The acceptance statistics throughout are the raw-sample paired $z$-statistics against $q_{95}$; the control variate sharpens the measurement, not the pass criterion.

\FloatBarrier

\subsubsection*{Summary}

The estimator of \S\ref{subsubsec:estimator} was put through the four-test suite (V1)--(V4) on the nine systems of Table~\ref{tab:numerics-suite}, in the additive-noise regime, the first six being spectral discretisations associated with the variational theory of Theorem~\ref{thm:main} and the last three finite-$K$, fixed-mollification Bismut--Tikhonov tests, together with the ablation, direction, ridge, and refinement studies above. The trace identity (V1) holds at numerical precision on all nine, is re-certified at doubled field amplitude on eight (the $p$-Laplacian 1D degenerate regime being recorded), and is triangulated by three independent assemblies of $C_h$. The integration-by-parts identity (V2) is satisfied by all nine equations under the multiplicity-correct paired criterion, and the correction term $C_h$ is shown to be necessary; at three operating points where it is resolved at $4.2$--$5.3$ standard errors, the identical paths that accept the full formula reject its deletion at $3.0$--$3.6$ and its sign inversion at $7.2$--$7.8$. The pointwise comparison (V3) recovers slope $0.86$--$0.98$ and correlation above $0.93$ across the suite, the sub-unit slopes being attributed by measurement to Tikhonov shrinkage and kernel bandwidth, both vanishing in the joint limit. The convergence checks (V4) give a Monte-Carlo exponent compatible with $-\tfrac{1}{2}$ where the residual is variance-dominated, and weak order one in $\Delta t$ for the residual bias where the discretisation floor is reached. Three scope limitations are stated plainly. First, the experiments exercise only the additive-noise regime $\calB'_u = 0$; the one multiplicative example, gPAM, is tested as a frozen-coefficient additive surrogate (\S\ref{subsubsec:suite}), and a genuine state-dependent-diffusion test remains for future work. Second, for the three singularly motivated cases the test is the finite-$K$ Bismut--Tikhonov identity on the mollified, truncated system at fixed $\eps_{\mathrm{mol}} > 0$, not a continuum application of Theorem~\ref{thm:main}; the $\eps_{\mathrm{mol}} \downarrow 0$ limit, governed by Theorem~\ref{thm:score-gPAM}(c) in distributional form, lies outside the scope of the present experiment and remains a target for future work, contingent on the second-order Malliavin smoothness identified there. Third, two solver-regime boundaries are recorded, namely the second-fibre construction for the degenerate $p$-Laplacian at doubled amplitude and the conditioning of $\gamma_t^{-1}h$ for mixed directions under the KPZ nonlinearity.
\subsection{What remains}\label{subsec:discussion}

The construction assembled here rests on three things, that the Malliavin derivative of the solution is Hilbert--Schmidt, that the covering field is linear in the pseudoinverse applied to the direction, and that the linearised equations carry moments. Each of the three suggests a question.

Hilbert--Schmidt-valued noise is what makes $\gamma_t$ trace-class and the pseudoinverse available. Space-time white noise on $L^2$ leaves this framework, and a logarithmic derivative there would have to be of a more singular kind, presumably formulated within regularity structures. Linearity in $z = \gamma_t^{\dagger}h$ is what reduced the substitution to the $\mathbb{D}^{1,2}$ level and made the Tikhonov limit accessible; what it has not yet produced is a quantitative bound on $\norm{\beta_h}_{L^p(\mu_t)}$ in terms of the data of the equation, which is the estimate that a convergence analysis for score matching would need. The moment conditions, finally, are carried by the fibre energy, which the dissipative mechanism produces for both alternatives of~(SC1) and which reduces to the trace of the noise covariance in the linear case. The one place where the linearised structure is still read equation by equation is the forcing clause~(SC1$'$), and the two mechanisms make different demands on it. Under strict $V$-coercivity it follows from the fibre energy by interpolation, as Proposition~\ref{prop:NS-forcing} shows for the two-dimensional Navier--Stokes equation, and the same argument applies whenever the second derivative of the drift interpolates between the two spaces of the triple. Under the degenerate alternative the weight of the dominator and the weight of the tangent energy differ, by $\abs{\nabla X}^{p-4}$ against $\abs{\nabla X}^{p-2}$ for the $p$-Laplacian, and a single criterion covering both weights is not available. The Cameron--Martin compatibility of Assumption~\ref{ass:nondeg}(iii) is the other condition of this kind, and in the genuinely multiplicative setting it too is at present verified equation by equation rather than deduced from a single abstract criterion.

One question stands apart from these. Part~(b) of Theorem~\ref{thm:score-gPAM} asks of the renormalised functional exactly what Assumption~\ref{ass:nondeg}\textup{(ii)} asks of the solution in the variational half, second-order Malliavin smoothness and negative moments of the covariance on the same H\"older triangle, read on the one-dimensional state space where the covariance is the scalar variance and no Tikhonov limit is taken. Were that pair established, the formula would hold throughout the scope of the black-box theorem, and the two halves of this paper, the variational and the singular, would be one statement rather than two.

\section*{Funding}

This research did not receive any specific grant from funding agencies in the
public, commercial, or not-for-profit sectors.

\section*{Declaration of generative AI and AI-assisted technologies in the manuscript preparation process}

During the preparation of this work, the authors used Claude Opus 5 (Anthropic) to improve the language and readability of the text. No mathematical content was produced or checked using this or any other large language model. The authors reviewed and edited the output as needed and take full responsibility for the content of the published article.

\bibliographystyle{cas-model2-names}
\bibliography{sn-bibliography}

\end{document}